\documentclass[12pt]{article}

\usepackage{amsmath}

\usepackage{stmaryrd}
\usepackage{enumerate}

\usepackage[pagebackref]{hyperref}\renewcommand*\backref[1]{}\renewcommand*\backrefalt[4]{  \ifcase #1  No citations  \or  (Cited on page #2)   \else  (Cited on pages #2)   \fi}

\usepackage{graphicx}
\usepackage{amsmath}
\usepackage{amssymb}
\usepackage{latexsym}
\usepackage{url}
\usepackage{epsfig}
\usepackage{mathrsfs}
\usepackage{amsfonts}
\usepackage{amsthm}

\def\choose#1#2{\left (\!\!\begin{array}{c}#1\\#2\end{array}\!\!\right )}

\usepackage{color}

\newcommand{\MAGENTA}[1]{\textcolor{magenta}{#1}}

\long\def\delete#1{}

\newcommand{\sm}[1]{\MAGENTA{\bf [Sanming: #1]}}

\usepackage{xcolor}
\usepackage[normalem]{ulem}

\newcommand{\pmat}[1]{\begin{pmatrix}#1\end{pmatrix}}

\newcommand{\be}{\begin{equation}}
\newcommand{\ee}{\end{equation}}
\newcommand{\bea}{\begin{eqnarray}}
\newcommand{\eea}{\end{eqnarray}}
\newcommand{\bean}{\begin{eqnarray*}}
\newcommand{\eean}{\end{eqnarray*}}

\def\CCC{\mathbb{C}}
\def\FFF{\mathbb{F}}
\def\ZZZ{\mathbb{Z}}
\def\NNN{\mathbb{N}}
\def\PPP{\mathbb{P}}

\def\QQQ{\mathbb{Q}}
\def\RRR{\mathbb{R}}

\def\AA{{\cal A}}

\def\EE{{\cal E}}
\def\FF{{\cal F}}
\def\GG{{\cal G}}
\def\II{{\cal I}}

\def\NN{{\cal N}}
\def\OO{{\cal O}}
\def\PP{{\cal P}}

\def\SS{{\cal S}}
\def\TT{{\cal T}}
\def\XX{{\cal X}}

\def\b0{{\bf 0}}

\def\bfe{{\bf e}}

\def\bx{{\bf x}}

\def\bv{{\bf v}}

\def\diam{{\rm diam}}

\def\diag{{\rm diag}}

\def\codim{{\rm codim}}
\def\Cay{{\rm Cay}}

\def\gcd{{\rm gcd}}
\def\Spec{{\rm Spec}}
\def\mod{{\rm mod}}
\def\ICG{{\rm ICG}}

\def\Sud{{\rm Sud}}
\def\SudP{{\rm SudP}}
\def\PLSG{{\rm PLSG}}

\def\Cy{{\rm Cy}}

\def\Tr{{\rm Tr}}

\def\Cay{{\rm Cay}}
\def\Aut{{\rm Aut}}

\def\wr{{\rm\,wr\,}}

\def\PSL{{\rm PSL}}

\def\GL{{\rm GL}}
\def\AGL{{\rm AGL}}

\def\GO{{\rm GO}}

\def\Dic{{\rm Dic}}
\def\Dih{{\rm Dih}}
\def\Gal{{\rm Gal}}

\def\E{\mathcal {E}}

\def\la{\langle}
\def\ra{\rangle}

\def\lc{\lceil}
\def\rc{\rceil}
\def\lf{\lfloor}
\def\rf{\rfloor}

\def\De{\Delta}
\def\Ga{\Gamma}

\def\Si{\Sigma}
\def\Om{\Omega}

\def\a{\alpha}
\def\b{\beta}
\def\d{\delta}
\def\g{\gamma}
\def\l{\lambda}
\def\om{\omega}
\def\s{\sigma}
\def\t{\tau}
\def\ve{\varepsilon}

\theoremstyle{plain}
\newtheorem{thm}{Theorem}
\newtheorem{cor}[thm]{Corollary}

\newtheorem{assump}[thm]{Assumption}
\newtheorem{conj}[thm]{Conjecture}

\numberwithin{equation}{section}

\newcommand{\NEPS}{\mathrm{NEPS}}
\newcommand*{\legendre}[2]{\left(\frac{#1}{#2}\right)}

\theoremstyle{remark}
\newtheorem{prob}{{Problem}}

\date{}

\title{\textbf{Eigenvalues of Cayley graphs}}

\author{Xiaogang Liu\thanks{Supported by the National Natural Science Foundation of China (Nos. 11601431 and 11871398), the Natural Science Foundation of Shaanxi Province (No. 2020JM-099), and the Natural Science Foundation of Qinghai Province  (No. 2020-ZJ-920).}\\
{\small  School of Mathematics and Statistics}\\[-0.8ex]
{\small  Northwestern Polytechnical University}\\[-0.8ex]
{\small  Xi'an, Shaanxi 710072, P. R. China}\\
\emph{{\small \tt xiaogliu@nwpu.edu.cn}}
\and
Sanming Zhou\thanks{Supported by a Future Fellowship (FT110100629) of the Australian Research Council.}\\
{\small School of Mathematics and Statistics}\\[-0.8ex]
{\small  The University of Melbourne}\\[-0.8ex]
{\small  Parkville, VIC 3010, Australia}\\
\emph{{\small \tt sanming@unimelb.edu.au}}}

\begin{document}

\maketitle

\begin{abstract}
We survey some of the known results on eigenvalues of Cayley graphs and their applications, together with related results on eigenvalues of Cayley digraphs and generalizations of Cayley graphs.


\bigskip

\end{abstract}


\tableofcontents



\section{Introduction}
\label{sec:int}

The study of eigenvalues of graphs is an important part of modern graph theory. In particular, eigenvalues of Cayley graphs have attracted increasing attention due to their prominent roles in algebraic graph theory and applications in many areas such as expanders \cite{Hoory06,Lubotzky95}, chemical graph theory \cite{Trinajstic92} and quantum computing \cite{Bernasconi08,Stevanovic11}. A large number of results on spectra of Cayley graphs have been produced over the last more than four decades. This paper is a survey of the literature on eigenvalues of Cayley graphs and their applications.

All definitions below are standard and can be found in, for example, \cite{BondyM08, Cvetkovic10, GodsilG2001}. A finite undirected \emph{graph} consists of a finite set whose elements are called the \emph{vertices} and a collection of unordered pairs of (not necessarily distinct) vertices each called an \emph{edge}. As usual, for a graph $G$, we use $V(G)$ and $E(G)$ to denote its vertex and edge sets, respectively, and we call the size of $V(G)$ the \emph{order} of $G$. An edge $\{u, v\}$ of $G$ is usually denoted by $uv$ or $vu$. If $uv$ is an edge of $G$, we say that $u$ and $v$ are joined by this edge, $u$ and $v$ are \emph{adjacent} in $G$, both $u$ and $v$ are \emph{incident} to the edge $uv$, and $u$ and $v$ are the \emph{end-vertices} of the edge $uv$. An edge joining a vertex to itself is called a \emph{loop}, and two or more edges joining the same pair of distinct vertices are called \emph{parallel edges}. The \emph{degree} of $u$ in $G$, denoted by $d_{G}(u)$ or simply $d(u)$ if there is no risk of confusion, is the number of edges of $G$ incident to $u$, each loop counting twice. A graph without loops is called a \emph{multigraph}, and a graph is \emph{simple} if it has no loops or parallel edges.

Let $G$ be a finite undirected graph of order $n$. The {\em adjacency matrix} of $G$, denoted by $A(G)$, is the $n \times n$ matrix with rows and columns indexed by the vertices of $G$ such that the $(u, v)$-entry is equal to the number of edges joining $u$ and $v$, with each loop counting as two edges. The eigenvalues and eigenvectors (eigenfunctions) of $A(G)$ are called the \emph{eigenvalues} and \emph{eigenvectors} (\emph{eigenfunctions}) of $G$, respectively, and the collection of eigenvalues of $G$ with multiplicities is called the \emph{spectrum} of $G$. Since $G$ is undirected, $A(G)$ is a real symmetric matrix and hence the eigenvalues of $G$ are all real numbers. If $\lambda_1,\lambda_2,\ldots,\lambda_r$ are the distinct eigenvalues of $G$ and $m_1,m_2,\ldots,m_r$ the corresponding multiplicities, then the spectrum of $G$ is denoted by
\begin{eqnarray*}
\Spec (G) = \pmat{\lambda_1 & \lambda_2 & \ldots  & \lambda_r \\
m_1  & m_2  & \ldots & m_r}
\end{eqnarray*}
or
$$
\Spec (G) = (\lambda_1^{m_1}, \lambda_2^{m_2}, \ldots, \lambda_r^{m_r}),
$$
where in the latter notation we usually omit $m_i$ if $m_i = 1$ for some $i$. An eigenvalue with multiplicity $1$ is called a \emph{simple} eigenvalue. The \emph{spectral radius} of a graph is the maximum modulus of its eigenvalues. Two graphs are said to be \emph{cospectral} if they have the same spectrum. Cospectral graphs need not be isomorphic, but isomorphic graphs are always cospectral.

Let $D(G)$ be the $n \times n$ diagonal matrix whose $(u, u)$-entry is equal to the degree $d(u)$ of $u$ in $G$, for each $u \in V(G)$. The matrices $D(G) - A(G)$ and $D(G) + A(G)$ are called the \emph{Laplacian matrix} and \emph{signless Laplacian matrix} of $G$, respectively, and their eigenvalues are called the \emph{Laplacian eigenvalues} and \emph{signless Laplacian eigenvalues} of $G$, respectively. It is readily seen that the smallest Laplacian eigenvalue of $G$ is $0$. The second smallest Laplacian eigenvalue of $G$, with multiple eigenvalues counted separately, is called the \emph{algebraic connectivity} of $G$. Since the multiplicity of $0$ is equal to the number of connected components of $G$, the algebraic connectivity of $G$ is positive if and only if $G$ is connected. In the case when $G$ is \emph{$k$-regular} for some integer $k \ge 0$, that is, each vertex has degree $k$, the largest eigenvalue of $G$ is $k$, and moreover a real number $x$ is an eigenvalue of $G$ if and only if $k - x$ is a Laplacian eigenvalue of $G$. Since all graphs considered in this paper are regular, this implies that all results to be reviewed in the paper can be presented in terms of Laplacian eigenvalues. Henceforth we will mostly talk about eigenvalues of graphs.

Similar to undirected graphs, a finite \emph{directed graph} (or \emph{digraph}) can be defined by specifying a finite set of vertices and a collection of ordered pairs of (not necessarily distinct) vertices; each ordered pair $(u, v)$ in the collection is called an \emph{arc} with \emph{tail} $u$ and \emph{head} $v$ (from $u$ to $v$). We also say that $u$ and $v$ are the end-vertices of the arc $(u, v)$. The \emph{in-degree} of a vertex $u$ is the number of arcs with $u$ as their head, and the \emph{out-degree} of $u$ is the number of arcs with $u$ as their tail. The \emph{underlying graph} of a digraph is the undirected graph obtained by replacing each arc by an edge with the same end-vertices. The \emph{adjacency matrix} of a finite digraph $G$ is the matrix whose $(u, v)$-entry is equal to the number of arcs from $u$ to $v$, and the eigenvalues of this matrix are called the \emph{eigenvalues} of $G$. Note that, unlike the undirected case, a digraph may have complex eigenvalues as its adjacency matrix is not necessarily symmetric.

A digraph without loops is called a \emph{multidigraph}. A digraph is \emph{symmetric} if whenever $(u, v)$ is an arc, $(v, u)$ is also an arc. We often identify a loopless symmetric digraph with its underlying simple graph, which is obtained by replacing each pair of arcs $(u, v), (v, u)$ by a single edge joining $u$ and $v$. Note that a loopless symmetric digraph and its underlying simple graph have the same adjacency matrix and hence the same spectrum.

Let $\Gamma$ be a finite group with identity element $1$, and let $S$ be a subset of $\Gamma$. The \emph{Cayley (di)graph} on $\Gamma$ with \emph{connection set} $S$, denoted by $\Cay(\Gamma,S)$, is defined to be the digraph with vertex set $\Gamma$ and arcs $(x, y)$ for all pairs $x, y \in \Gamma$ such that $xy^{-1} \in S$. Clearly, this digraph has order $|\Ga|$ and degree $k$ in the sense that each vertex has both in-degree and out-degree $k$, where $k = |S|$. Moreover, $\Cay(\Gamma,S)$ is connected if and only if $\la S \ra = \Gamma$, where $\la S \ra$ is the subgroup of $\Gamma$ generated by $S$. In the case when $1 \not \in S$ and $S$ is \emph{inverse-closed} (that is, $S = S^{-1} := \left\{s^{-1}: s \in S\right\}$), the digraph $\Cay(\Gamma,S)$ is symmetric with no loops and hence may be identified with its underlying simple graph. In other words, if $S \subseteq \Ga \setminus \{1\}$ is inverse-closed, then $\Cay(\Gamma,S)$ is understood as a $k$-regular simple graph.

Throughout this paper, unless stated otherwise whenever we say a graph we mean a finite undirected simple graph. In particular, unless stated otherwise Cayley graphs considered in this paper are finite, undirected and simple.

\subsection{Outline of the paper}
\label{subsec:outline}

It is well known that eigenvalues of Cayley (di)graphs can be expressed in terms of the irreducible characters of the underlying groups. In Section \ref{sec:cha} we will review these fundamental results together with their counterparts for vertex-transitive graphs.

An important problem in spectral graph theory is to understand when all eigenvalues of a graph are integers. A graph with this property is called integral, or more specifically, integral over the field $\mathbb{Q}$. Similarly, one can consider graphs which are integral over some other algebraic number fields such as the Gaussian field $\mathbb{Q}(i)$. A number of results on integral Cayley graphs have been produced in the past more than ten years. In Section \ref{sec:ICG} we will review these results, along with a few other results on Cayley graphs which are integral over $\mathbb{Q}(i)$ or a general algebraic number field. In Section \ref{sec:CoCayG} we will survey results on cospectral Cayley graphs and Cayley graphs which are determined by their spectra. In Section \ref{sec:ComRing} we will discuss several families of Cayley graphs on finite commutative rings, including unitary Cayley graphs and quadratic unitary Cayley graphs of finite commutative rings, some Cayley graphs on finite chain rings, and generalized Paley graphs on finite fields.

The energy of a graph is a concept that arises from chemical graph theory. This notion has been studied extensively over the past more than four decades. In particular, a number of results on energies of Cayley graphs have been produced in recent years. We will give an account of such results in Section \ref{sec:EnCay}.

It is well known \cite{Hoory06} that the expansion of a regular graph is determined by its second largest eigenvalue and that Cayley graphs play an important role in constructing expander graphs. Roughly, regular graphs achieving best possible expansion (in terms of the well-known Alon-Boppana bound) are called Ramanujan graphs. In Sections \ref{sec:Ram} and \ref{sec:2ndCay} we will review some results on Ramanujan Cayley graphs and the second largest eigenvalue of Cayley graphs, respectively.

Perfect state transfer in graphs is an important concept that arises from quantum computing. In general, it is challenging to construct graphs admitting perfect state transfer or prove the existence of perfect state transfer in a graph. In Section \ref{sec:pst} we will survey recent results on perfect state transfer in several families of Cayley graphs along with related results on the periodicity of such graphs. We will also review recent results on pretty good state transfer in a few families of Cayley graphs.

A number of distance-regular graphs (and in particular strongly regular graphs) are known to be Cayley graphs. In Section \ref{sec:drCay} we will review some results on distance-regular Cayley graphs. This discussion will be continued in Section \ref{subsec:srnCay}, where we will review several results on strongly regular $n$-Cayley graphs, where $n \ge 2$. In Section \ref{sec:genCay} we will discuss a few notions of generalized Cayley graphs, including $n$-Cayley graphs, Cayley sum graphs and group-subgroup pair graphs, with a focus on their eigenvalues.

In Section \ref{sec:diCay} we will give a brief account of results on eigenvalues of directed Cayley graphs. In Section \ref{sec:miscell} we will mention a number of miscellaneous results, including results on eigenvalues of random Cayley graphs, distance eigenvalues of Cayley graphs, eigenvalues of mixed Cayley graphs, etc. We will conclude this paper in Section \ref{sec:open} with a collection of open problems and research topics.

Despite its length this survey is far from being a comprehensive treatise on eigenvalues of Cayley graphs. Some important omissions are as follows.
\begin{itemize}
\item Only finite Cayley graphs are considered in our paper. Some results on eigenvalues of infinite Cayley graphs can be found in \cite[Section 5]{MW89}.
\item We do not survey results on families of expanders with a fixed degree but increasing orders (see, for example, \cite{BourgainG08, Davidoff03, Lubotzky88, Lubotzky12, Margulis73}). We believe that this important area involving deep group theory and other mathematical tools deserves a separate treatment by experts in the area. The reader is referred to \cite{Hoory06} for a survey on expanders and their applications and \cite{Lubotzky95} for a survey on eigenvalues of Cayley graphs in the context of expanders and random walks. See also \cite{Li1996, Murty03} for two surveys on Ramanujan graphs and \cite{Tao2015} for a book on expansion in Cayley graphs and finite simple groups of Lie type. A survey on analytic, geometric and combinatorial properties of finitely generated groups that are related to expansion properties of Cayley graphs as measured by the isoperimetric number (Cheeger constant) can be found in \cite{Zuk16}. See also \cite{Kowalski19, Krebs11} for two introductory books on expanders and \cite{Davidoff03} for a monograph focusing on the Ramanujan graphs discovered by Lubotzky, Phillips and Sarnak \cite{Lubotzky88} and Margulis \cite{Margulis73}.
\item Strongly regular Cayley graphs are essentially partial difference sets in groups. A large number of results on partial difference sets exist in the literature (see \cite{Ma94} for a not-so-recent survey on this topic), but we believe that a survey of them should be a subject of its own. So this area is mostly omitted in our paper, with the exception of a few results whose proofs use methods which seem to be nontypical in the area of partial difference sets. The reader is referred to \cite{BrouwerCN89} for a classic treatment and \cite{vanDamKT16} for a recent survey on distance-regular graphs. Most results on distance-regular Cayley graphs reviewed in our paper are not covered in \cite{vanDamKT16}.
\item By the well-known matrix-tree theorem due to Kirchhoff (see \cite{Biggs93,BrouwerH12,GodsilG2001}), the number of spanning trees in a regular graph is determined by its eigenvalues. Neverthless, results on the number of spanning trees in Cayley graphs will not be surveyed in our paper.
\end{itemize}

We conclude this section by mentioning a few references which complement the present paper. In \cite{ChanG97}, Chan and Godsil surveyed results on automorphisms of a graph in relation to its eigenvalues and eigenvectors. A survey on Laplacian eigenvalues of graphs can be found in \cite{Mohar92}, and a survey on spectra of digraphs is given in \cite{Brualdi10}. Basic results on eigenvalues of Cayley graphs can be found in two widely used textbooks \cite{Biggs93, GodsilG2001} on algebraic graph theory and a recent book \cite{GodsilM16} on Erd\H{o}s-{K}o-{R}ado theorems.

\subsection{Terminology and notation}
\label{subsec:notation}

The reader is referred to \cite{BondyM08}, \cite{Isaacs08}, \cite{Hardy08} and \cite{Atiyah69, BiniF2002} for terminology and notation on graph theory, group theory, number theory and the theory of finite commutative rings, respectively.

The \emph{complement} $\overline{G}$ of a graph $G$ is the graph with vertex set $V(G)$ in which two vertices are adjacent if and only if they are not adjacent in $G$. The \emph{line graph} $L(G)$ of $G$ is the graph with vertices the edges of $G$ such that two vertices are adjacent if and only if the corresponding edges have a common end-vertex. The \emph{distance} between two vertices $u$ and $v$ in a graph $G$, denoted by $d_{G}(u, v)$ or simply $d(u, v)$, is the length of a shortest path between them in $G$; if there is no path between $u$ and $v$ in $G$, then we set $d(u, v) = \infty$. The \emph{diameter} of $G$, denoted by $\diam(G)$, is the maximum distance between two vertices of $G$.

The \emph{Cartesian product} $G \boxempty H$ of two graphs $G$ and $H$ is the graph with vertex set $V(G)\times V(H)$ in which $(u, v)$ and $(x, y)$ are adjacent if and only if either $u = x$ and $vy \in E(H)$, or $v = y$ and $ux \in E(G)$. The \emph{tensor product} (also known as direct product, Kronecker product or categorical product in the literature) of $G$ and $H$, denoted by $G\otimes H$, is the graph with vertex set $V(G)\times V(H)$ in which $(u,v)$ is adjacent to $(x,y)$ if and only if $u$ is adjacent to $x$ in $G$ and $v$ is adjacent to $y$ in $H$. These two operations are associative and so the Cartesian (tensor) product of any finite number of graphs is well defined.

A group $\Ga$ is said to {\em act} on a set $\Om$ if each pair $(\a, g) \in \Om \times \Ga$ corresponds to some $\a^g \in \Om$ such that $\a^1 = \a$ and $(\a^g)^{h} = \a^{gh}$ for $g, h \in \Ga$, where $1$ is the identity element of $\Ga$. The {\em stabilizer} of $\a \in \Om$ in $\Ga$ is the subgroup $\Ga_{\a} := \{g \in \Ga: \a^g = \a\}$ of $\Ga$. The {\em $\Ga$-orbit} containing $\a$ is defined to be $\a^{\Ga} := \{\a^g: g \in \Ga\}$. The \emph{rank} of $\Ga$ is the number of $\Ga$-orbits on $\Om \times \Om$ under the induced action defined by $(\a, \b)^g = (\a^g, \b^g)$ for $\a, \b \in \Om$ and $g \in \Ga$. If $\Ga_{\a} = \{1\}$ for all $\a \in \Om$, then $\Ga$ is called {\em semiregular} on $\Om$. If $\a^{\Ga} = \Om$ for some (and hence all) $\a \in \Om$, then $\Ga$ is {\em transitive} on $\Om$. If $\Ga$ is both transitive and semiregular on $\Om$, then it is said to be {\em regular} on $\Om$.

An \emph{automorphism} of a graph $G$ is a permutation of $V(G)$ which maps adjacent vertices to adjacent vertices and nonadjacent vertices to nonadjacent vertices. The \emph{automorphism group} of $G$, denoted by $\Aut(G)$, is the group of automorphisms of $G$ under the usual composition of permutations. Of course $\Aut(G)$ acts on $V(G)$ in a natural way, and this action induces natural actions on the set of edges, the set of arcs and the set of $2$-arcs of $G$, where an \emph{arc} is an ordered pair of adjacent vertices, and a \emph{$2$-arc} is an ordered triple $(u, v, w)$ of distinct vertices such that $v$ is adjacent to both $u$ and $w$. If a subgroup $\Ga$ of $\Aut(G)$ is transitive on $V(G)$, then $G$ is said to be \emph{$\Ga$-vertex-transitive}. Similarly, if $\Ga$ is transitive on the set of edges of $G$, then $G$ is \emph{$\Ga$-edge-transitive}; if $\Ga$ is transitive on the set of arcs of $G$, then $G$ is \emph{$\Ga$-arc-transitive}; and if $\Ga$ is transitive on the set of $2$-arcs of $G$, then $G$ is \emph{$(\Ga, 2)$-arc-transitive}. A graph $G$ is called \emph{vertex-transitive}, \emph{edge-transitive}, \emph{arc-transitive} or \emph{$2$-arc-transitive} if it is $\Aut(G)$-vertex-transitive, $\Aut(G)$-edge-transitive, $\Aut(G)$-arc-transitive or $(\Aut(G), 2)$-arc-transitive, respectively. Clearly, any vertex-transitive graph must be regular and any arc-transitive graph without isolated vertices must be vertex-transitive. It is known that any edge-transitive but not vertex-transitive graph must be bipartite (see, for example, \cite[Proposition 15.1]{Biggs93}). It is well known that all Cayley graphs are vertex-transitive, but the converse is not true, the Petersen graph being a counterexample. It is also well known that a graph is isomorphic to a Cayley graph if and only if its automorphism group contains a subgroup which is regular on the vertex set (see, for example, \cite[Lemma 16.3]{Biggs93}). A \emph{circulant graph}, or simply a \emph{circulant}, is a Cayley graph on a cyclic group.

Throughout the paper we use the following notation:

$\NNN,\ \ZZZ,\ \mathbb{P}$: Sets of positive integers, integers, primes, respectively

$\QQQ,\ \RRR,\ \CCC$: Fields of rational numbers, real numbers, complex numbers, respectively

$|V|$: Cardinality of a set $V$ (in particular, if $\Ga$ is a group, then $|\Gamma|$ is the order of $\Ga$)

$n, d, r$: Positive integers

$\omega_n = \exp(2\pi i/n)$: An $n$th primitive root of unity, where $i^2=-1$

$\phi(n)$: Euler's totient function, which gives the number of positive integers up to $n$ that are coprime to $n$

$e_p(n)$: Exponent of prime $p$ in $n$

$\Cay(\Ga, S)$: Cayley (di)graph of a group $\Ga$ with respect to connection set $S \subseteq \Ga$

$\ZZZ_n$: Ring of integers modulo $n$

$\ZZZ_n^{\times}$: Set of units of ring $\mathbb{Z}_n$

$C_n,\ K_n$: Cycle and complete graph of order $n$, respectively

$K_{m,n}$: Complete bipartite graph with $m$ and $n$ vertices in the bi-parts of the bipartition, respectively

$H(d,q)$: Hamming graph, namely, the Cartesian product $K_q \Box \cdots \Box K_q$ ($d$ factors)

$H(d,2)$: Hypercube of dimension $d$

$\NEPS(G_1,\ldots,G_d;\mathcal{B})$: NEPS of graphs $G_1,\ldots,G_d$ with respect to basis $\mathcal{B}\subseteq \mathbb{Z}_2^{d}\setminus \{\mathbf{0}\}$

$L(G)$: Line graph of a graph $G$

$\overline{G}$: Complement of a graph $G$

$\Spec(G)$: Spectrum of a graph $G$

$\EE(G)$: Energy of a graph $G$

$D(n)$: Set of positive divisors of $n$

$\ICG(n, D)$: gcd graph of the cyclic group $\ZZZ_n$ with respect to $D \subseteq D(n) \setminus \{n\}$

$\ICG(R/(c), D)$: gcd graph of a UFD $R$ with respect to a set $D$ of proper divisors of $c \in R \setminus \{0, 1\}$

$G_{R} = \Cay(R, R^{\times})$: Unitary Cayley graph of a finite commutative ring $R$, where $R^{\times}$ is the set of units of $R$

$G_{\ZZZ_n} = \Cay(\ZZZ_n, \ZZZ_n^{\times}) = \ICG(n, \{1\})$: Unitary Cayley graph of ring $\ZZZ_n$

$\mathcal{G}_{R}$: Quadratic unitary Cayley graph of a finite commutative ring $R$

$\FFF_q$: Finite field of $q$ elements, $q$ being a prime power

$w(\mathbf{a})$: Hamming weight, namely the number of nonzero coordinates of $\mathbf{a}$

$S_n$: Symmetric group of degree $n$

$A_n$: Alternating group of degree $n$

$D_{2n}$: Dihedral group $\langle a,b \mid a^n=b^2=1, b^{-1}ab=a^{-1}\rangle$ of order $2n$

$\Dih(\Si)$: Generalized dihedral group of an abelian group $\Si$, namely, $\la \Si, b \mid b^2 = 1, b^{-1} a b = a^{-1}, a \in \Si \ra$

$SD_{8n}$: Semi-dihedral group $\la a, b \mid a^{4n} = b^2 = 1, bab = a^{2n-1}\ra$ of order $8n$

$\Dic_{n}$: Dicyclic group $\langle a, b \mid a^{2n}=1, b^2=a^n, b^{-1}ab=a^{-1}\rangle$ of order $4n$

$\langle S\rangle$: Subgroup of a group $\Ga$ generated by a subset $S \subseteq \Ga$ (in particular, $\langle x\rangle = \langle \{x\}\rangle$, $\langle x, y\rangle = \langle \{x, y\}\rangle$, and so on)

$XY$: The subset $\{xy: x \in X, y \in Y\}$ of a group $\Ga$, for given $X, Y \subseteq \Ga$

$\Gamma \rtimes \Sigma$: Semidirect product of a group $\Gamma$ by a group $\Sigma$

$\Gamma \wr \Sigma$: Wreath product of a group $\Gamma$ by a group $\Sigma$

$\Aut(G)$: Automorphism group of a graph $G$

$\Aut(\Ga)$: Automorphism group of a group $\Ga$

$\mathbb{C}\Ga$: Group algebra of a group $\Ga$ over $\mathbb{C}$

$\ZZZ \Ga$: Group ring of a group $\Ga$ over $\ZZZ$


\section{Eigenvalues of Cayley graphs}
\label{sec:cha}

\subsection{Characters}
\label{subsec:cha}

The reader is referred to \cite{Isaacs76,James01,Sagan01,Serre77,Steinberg12} for representation theory of finite groups and properties of characters. Here we recall only a few basic definitions for ordinary representations of finite groups. Let $V$ be an $n$-dimensional vector space over $\CCC$, where $n \ge 1$ is an integer. A \emph{representation} of a finite group $\Gamma$ on $V$ is a group homomorphism $\pi: \Gamma \rightarrow \GL(V)$, where $\GL(V)$ is the \emph{general linear group} of $V$ defined as the group of invertible linear transformations of $V$. The dimension $n$ of $V$ is called the \emph{degree} or \emph{dimension} of $\pi$. Two representations $\pi_i: \Gamma \rightarrow \GL(V_i)$, $i=1, 2$ are said to be \emph{equivalent}, written $\pi_{1} \sim \pi_{2}$, if there exists a vector space isomorphism $T: V_1 \rightarrow V_2$ such that, for all $g \in \Ga$, $\pi_{1}(g) T = T \pi_{2}(g)$. As usual, by identifying each element of $\GL(V)$ with its matrix with respect to a chosen basis for $V$, we may identify $\GL(V)$ with the group $\GL(n,\mathbb{C})$ of all invertible $n \times n$ matrices over $\CCC$ with operation the product of matrices, and we may view a representation of $\Ga$ on $V$ as a group homomorphism $\pi: \Gamma \rightarrow \GL(n,\mathbb{C})$. The \emph{character} of $\pi$ is the mapping $\chi_{\pi}: \Ga \rightarrow \CCC$ defined by
$$
\chi_{\pi}(g) = \Tr(\pi(g)),\; g \in \Gamma,
$$
where $\Tr$ denotes the trace of a matrix. It is readily seen that $\chi_{\pi}(1) = \Tr(\pi(1)) = \Tr(I_{n}) = n$, which is also called the \emph{degree} of $\chi_{\pi}$, where $1$ is the identity element of $\Gamma$ and $I_n$ the identity matrix of size $n \times n$. The representation $\pi: \Gamma \rightarrow \GL(1,\mathbb{C})$ defined by $\pi(g) = (1),\, g \in \Ga$ is called the \emph{trivial representation} of $\Ga$ and the corresponding character is called the \emph{trivial character} of $\Ga$.

Let $X$ and $Y$ be matrices. Define $X \oplus Y = \pmat{X & 0 \\ 0 & Y}$, where the two zero-blocks are all-$0$ matrices of appropriate sizes. Let $\pi_1$ and $\pi_2$ be representations of $\Gamma$. The \emph{direct sum} of $\pi_1$ and $\pi_2$, denoted by $\pi_1\oplus\pi_2$, is defined to be the representation of $\Gamma$ given
by $(\pi_1\oplus \pi_2)(g) = \pi_1(g) \oplus \pi_2(g)$ for $g \in \Ga$. A representation $\pi$ of $\Gamma$ is \emph{reducible} if there exist representations $\pi_1, \pi_2$ of $\Gamma$ such that $\pi \sim \pi_1\oplus\pi_2$, and \emph{irreducible} otherwise. In the latter case we say that $\chi_{\pi}$ is an \emph{irreducible character}.

It is well known \cite{Isaacs76,James01} that every irreducible character of any finite abelian group is one-dimensional. In other words, the characters of any finite abelian group $\Ga$ are precisely homomorphisms $\chi: \Ga \rightarrow \CCC^*$ (that is, $\chi(gh) = \chi(g) \chi(h)$ for $g, h \in \Ga$). Note that $|\chi(g)| = 1$ for all $g \in \Ga$. It is also well known that any finite abelian group $\Ga$ has $|\Ga|$ distinct characters and the dual group $\hat{\Ga}$ of $\Ga$ is isomorphic to $\Ga$, where $\hat{\Ga}$ is the group of all characters of $\Ga$ with operation defined by $(\chi \psi)(g) = \chi (g) \psi(g)$ for $\chi, \psi \in \hat{\Ga}, g \in \Ga$. In particular, the characters of the cyclic group $\ZZZ_n$ are given by $\chi_{k}(g) = \om_n^{kg},\; 0 \le k \le n-1$, where $\omega_n = \exp(2\pi i/n),\, i^2=-1$ is an $n$th primitive root of unity.


\subsection{Eigenvalues of Cayley graphs}
\label{subsec:spec-Cay}

In this section we survey a few basic results on eigenvalues of Cayley (di)graphs and Cayley colour (di)graphs. Many results on eigenvalues of Cayley graphs to be surveyed in later sections are based on these fundamental results.

Let $\Gamma$ be a finite group of order $n$ and $\alpha: \Gamma\rightarrow \mathbb{C}$ a function. The \emph{Cayley colour digraph} of $\Gamma$ with \emph{connection function} $\alpha$ (see, for example,  \cite{Babai79}), denoted by $\Cay(\Gamma,\alpha)$,  is defined to be the directed graph with vertex set $\Ga$ and arc set $\{(x, y): x, y \in\Gamma\}$ such that each arc $(x, y)$ is coloured by $\alpha(yx^{-1})$. (The colour of $(x, y)$ was defined as $\alpha(xy^{-1})$ in \cite{Babai79}, but we use $\alpha(yx^{-1})$ here as eigenvectors are treated as column vectors in this survey.) The \emph{adjacency matrix} of $\Cay(\Gamma,\alpha)$ is defined as the matrix with rows and columns indexed by the elements of $\Ga$ such that the $(x, y)$-entry is equal to $\alpha(yx^{-1})$. The \emph{eigenvalues} of $\Cay(\Gamma,\alpha)$ are simply the eigenvalues of its adjacency matrix. In the special case when $\alpha: \Gamma \rightarrow \{0,1\}$ and the set $S = \{g:\alpha(g)=1\}$ satisfies $1\notin S$ and $S = S^{-1}$, the Cayley colour digraph $\Cay(\Gamma,\alpha)$ can be identified with the Cayley graph $\Cay(\Gamma, S)$, and the adjacency matrix of $\Cay(\Gamma,\alpha)$ agrees with that of $\Cay(\Gamma, S)$.

In 1975, Lov\'{a}sz \cite{Lovasz75} proved the following result.

\begin{thm}
\emph{(\cite{Lovasz75}; also \cite[Corollary 3.2]{Babai79})}
\label{BabaiS1}
Let $\Gamma$ be an abelian group of order $n$ with irreducible characters $\chi_1,\chi_2,\ldots,\chi_n$. Let $\alpha: \Gamma\rightarrow \mathbb{C}$ be a function. Then the eigenvalues of the Cayley colour digraph $\Cay(\Gamma,\alpha)$ of $\Ga$ are given by
$$
\lambda_i=\sum_{g\in\Gamma}\alpha(g)\chi_i(g),\; i=1,2,\ldots,n.
$$
\end{thm}

Moreover, the eigenvector corresponding to $\lambda_i$ is $(\chi_i(x): x \in \Ga)^T$ (that is, the column vector whose $x$-entry is $\chi_i(x)$) because for any $x \in \Ga$ we have $\sum_{y \in \Ga} \alpha(yx^{-1}) \chi_i(y) = \sum_{g \in \Ga} \alpha(g) \chi_i(gx) = (\sum_{g \in \Ga} \alpha(g) \chi_i(g)) \chi_i(x) = \lambda_i \chi_i(x)$.

Now consider the general case when $\Ga$ is not necessarily abelian. We may express the adjacency matrix $A = (\alpha(yx^{-1}))_{x, y \in \Ga}$ of $\Cay(\Gamma,\alpha)$ as $A = \sum_{g \in \Ga} \alpha(g) L_{g}$, where $L_{g}$ is the matrix of the left multiplication by $g$ in the group algebra $\mathbb{C}\Ga$. Set $x = \sum_{g \in \Ga} \alpha(g) g \in \mathbb{C}\Ga$ and let $x^*$ be the left multiplication by $x$. By decomposing $\mathbb{C}\Ga$ into minimal left $\mathbb{C}\Ga$-modules and then considering the spectrum of the restriction of $x^*$ to each summand left $\mathbb{C}\Ga$-module in this decomposition, Babai \cite{Babai79} obtained the following result which gives Theorem \ref{BabaiS1} in the special case when $\Ga$ is abelian.

\begin{thm}
\emph{(\cite[Theorem 3.1]{Babai79})}
\label{BabaiS2}
Let $\Ga$ be a finite group. Let $\chi_1,\chi_2,\ldots,\chi_h$ be the irreducible characters of $\Ga$ and $n_1,n_2,\ldots,n_h$ be their degrees, respectively. Let $\alpha: \Gamma\rightarrow \mathbb{C}$ be a function. Then the eigenvalues of the Cayley colour digraph $\Cay(\Gamma,\alpha)$ of $\Ga$ can be arranged as
$$
\Lambda=\{\lambda_{i j k}: i=1,2,\ldots,h;j,k=1,2,\ldots,n_i\}
$$
such that $\lambda_{i j 1}=\lambda_{i j 2}=\cdots=\lambda_{i j n_i}$ (this common value is denoted by $\lambda_{i,j}$), and for any positive integer $t$,
$$
\lambda_{i,1}^t+\cdots+\lambda_{i,n_i}^t = \sum_{g_1,\ldots,g_t\in\Gamma}\left(\prod_{s=1}^t\alpha(g_s)\right) \chi_i\left(\prod_{s=1}^tg_s\right).
$$
\end{thm}

In the case when $\alpha: \Gamma \rightarrow \{0,1\}$, Theorem \ref{BabaiS1} yields the following corollary.

\begin{cor}
\label{BabaiS1a}
Let $\Gamma$ be an abelian group of order $n$ with irreducible characters $\chi_1,\chi_2,\ldots,\chi_n$. Then the eigenvalues of any Cayley (di)graph $\Cay(\Gamma, S)$ of $\Ga$ are given by
$$
\lambda_i = \sum_{g \in S} \chi_i(g),\; i=1,2,\ldots,n.
$$
\end{cor}

As explained after Theorem \ref{BabaiS1}, the eigenvector corresponding to $\lambda_i$ is $(\chi_i(x): x \in \Ga)^T$.

In the case when $\alpha: \Gamma \rightarrow \{0,1\}$, Theorem \ref{BabaiS2} gives rise to the following corollary.

\begin{cor}
\label{BabaiS2a}
Let $\Ga$ be a finite group. Let $\chi_1,\chi_2,\ldots,\chi_h$ be the irreducible characters of $\Ga$ and $n_1,n_2,\ldots,n_h$ be their degrees, respectively. Then the eigenvalues $\lambda_{i,j}$ ($1 \le i \le h$, $1 \le j \le n_i$) of any Cayley (di)graph $\Cay(\Gamma, S)$ of $\Ga$ satisfy
$$
\lambda_{i,1}^t+\cdots+\lambda_{i,n_i}^t = \sum_{g_1,\ldots,g_t\in S} \chi_i\left(\prod_{s=1}^t g_s \right)
$$
for any positive integer $t$.
\end{cor}

As an example we obtain the eigenvalues of any circulant graph from Corollary \ref{BabaiS1a} immediately (see \cite{Lovasz75}; see also \cite{Biggs93} for a treatment using circulant matrices). Recall that a circulant graph is a Cayley graph $\Cay(\ZZZ_n, S)$ on a cyclic group $\ZZZ_n$, where $n \ge 3$. Since the characters of $\ZZZ_n$ are given by $\chi_{k}(g) = \om_n^{kg},\; 0 \le k \le n-1$, from Corollary \ref{BabaiS1a} it follows that the eigenvalues of $\Cay(\ZZZ_n, S)$ are given by
\be
\label{eq:eigen-circ}
\lambda_k = \sum_{g \in S} \om_n^{k g},\; k = 0, 1, \ldots, n-1.
\ee
In particular, the eigenvalues of the cycle $C_n$ of length $n$ are $2 \cos (2\pi k/n),\; 0 \le k \le n-1$. Eigenvalues of circulants were reviewed in \cite{Moreno2000} as part of a survey of some aspects of quadratic Gauss sums over finite rings and fields with applications. In the case when $S$ is the set of nonzero quadratic residues, formula \eqref{eq:eigen-circ} together with the explicit formula for the Gauss sum gives a complete determination of the spectrum of the corresponding circulant $\Cay(\ZZZ_n, S)$.

As another example let us look at the $d$-dimensional hypercube $H(d, 2)$, $d \ge 1$, which can be defined as the Cayley graph $\Cay(\ZZZ_2^d, S)$, where $S = \{\mathbf{e}_1, \ldots, \mathbf{e}_d\}$ is the standard basis of $\ZZZ_2^d$. The characters of $\ZZZ_2^d$ are given by $\chi_{\mathbf{a}}(x) = (-1)^{\la \mathbf{a}, \mathbf{x} \ra}$, $\mathbf{a} = (a_1, \ldots, a_d) \in \ZZZ_2^d$. Hence $\chi_{\mathbf{a}}(\mathbf{e}_i) = -1$ if $a_i = 1$ and $\chi_{\mathbf{a}}(\mathbf{e}_i) = 1$ if $a_i = 0$. The eigenvalue of $H(d, 2)$ associated with $\mathbf{a}$ is $\sum_{i} \chi_\mathbf{a} (\mathbf{e}_i) = (d-w(\mathbf{a}))-w(\mathbf{a}) = d - 2w(\mathbf{a})$, where $w(\mathbf{a})$ is the Hamming weight of $\mathbf{a}$. Therefore, the eigenvalues of $H(d, 2)$ are $d-2i$ with multiplicity $\choose{d}{i}$, $0 \le i \le d$. In general, using Corollary \ref{BabaiS1a}, Lov\'{a}sz \cite{Lovasz75} determined the eigenvalues of any Cayley graph on $\ZZZ_2^d$ and coined the term ``cubelike graphs" for such graphs.

In general, for $d \ge 1$ and $q \ge 2$, the eigenvalues of the Hamming graph $H(d,q)$ are $q(d-j)-d$ with multiplicities $\choose{d}{j}(q-1)^j$, $0 \le j \le d$ (see, for example, \cite{BrouwerCN89}), where by definition $H(d, q)$ is the Cayley graph $\Cay(\ZZZ_q^d, S)$ with $S$ the set of elements of $\ZZZ_q^d$ having exactly one nonzero coordinate. The eigenvalues of some generalized Hamming graphs can also be computed explicitly as will be seen in Sections \ref{HammingE} and \ref{subsec:Hamming-Johnson}.

A Cayley graph $\Cay(\Gamma,S)$ is called \emph{normal}\footnote{Note that this definition is different from the following notion which is also widely used in the literature: A Cayley graph $G = \Cay(\Gamma, S)$ is called a normal Cayley graph if the right regular representation of $\Ga$ is normal in $\Aut(G)$.} if $S$ is closed under conjugation (that is, $S$ the union of some conjugacy classes of the group $\Gamma$). The following result enables us to compute explicitly the eigenvalues of any normal Cayley graph using character values of the underlying group. This result was proved by Diaconis and Shahshahani in \cite{Diaconis81} and Zieschang in \cite{Zieschang88}, and we adapt its presentation from the latter. It was also implied in the work of Babai \cite{Babai79}.

\begin{thm}
\emph{(\cite[Theorem 1]{Zieschang88})}
\label{Conjeign}
Let $\Ga$ be a finite group and let $\{\chi_1,\chi_2,\ldots,\chi_h\}$ be the set of all irreducible characters of $\Gamma$. Then the eigenvalues of any normal Cayley graph $\Cay(\Gamma,S)$ of $\Ga$ are given by
\be
\label{eq:Conjeign}
\lambda_j=\frac{1}{\chi_j(1)}\sum_{g \in S}\chi_j(g),\; j = 1, 2, \ldots, h.
\ee
Moreover, the multiplicity of $\lambda_j$ is equal to $\sum_{1\le k\le h,\, \lambda_k = \lambda_j}\chi_k(1)^2$.
\end{thm}

A consequence of this result is that the number of distinct eigenvalues of a normal Cayley graph on a finite group $\Ga$ does not exceed the number of irreducible complex representations of $\Ga$. Another application of Theorem \ref{Conjeign} is a formula (see \cite[Theorem 2]{Zieschang88}) for the number of walks of length $n$, $n \in \NNN$, between any two vertices in a normal Cayley graph. In the special case when $\Ga$ contains a cyclic self-normalized subgroup $W$ of order $pq$ for distinct primes $p, q$, and $S = \cup_{g \in \Ga} W_0^g$, with $W_0$ the set of elements of $W$ with order $pq$,
$\Cay(\Gamma,S)$ has $4$ or $5$ distinct eigenvalues and for any $n \in \NNN$ the number of walks of length $n$ between two adjacent vertices in $\Cay(\Gamma,S)$ does not depend on the choice of the two adjacent vertices (see \cite[Theorem 9]{Zieschang88}).

Theorem \ref{Conjeign} can be used to compute explicitly the eigenvalues of specific normal Cayley graphs in the case when the irreducible characters of the underlying groups are known. For example, all groups of order $p^2 q$ together with their character tables are known, where $p$ and $q$ are distinct odd primes. Using this and Theorem \ref{Conjeign}, the eigenvalues of each normal Cayley graph of order $p^2 q$ ($p < q$) were determined in \cite{GhorbaniL16}. Similarly, for distinct odd primes $p, q, r$, the eigenvalues of all normal Cayley graphs on groups of order $kpq$ ($k = 2, 3$) or $pqr$ were determined in \cite{GhorbaniL18} and \cite{GhorbaniN17}, respectively. In addition, the spectra of all connected minimal arc-transitive ``normal'' Cayley graphs of order $p^2 q$ for distinct odd primes $p, q$ were determined in \cite{GhorbaniSN20}, but here the word ``normal'' means that for the Cayley graph $\Cay(\Gamma, S)$ considered the right regular representation of $\Ga$ is normal in the automorphism group of $\Cay(\Gamma, S)$.

\begin{thm}
\emph{(\cite[Theorem 4.3 and Corollary 4.5]{Foster-GreenwoodK16})}
\label{thm:FG16}
Let $\Ga$ be a finite group and let $f: \Ga \rightarrow \CCC$ be a class function. Then the eigenvalues of the matrix $M_f = (f(xy^{-1}))_{x, y \in \Ga}$ are
$$
\theta_{\chi} = \frac{1}{\chi(1)} \sum_{x \in \Ga} f(x) \chi(x),\; \text{ with multiplicity } \chi(1)^2,
$$
where $\chi$ ranges over the irreducible characters of $\Ga$.

Moreover, the spectral radius of $M_f$ is equal to $\sum_{x \in \Ga} |f(x)|$. In particular, if $f: \Ga \rightarrow \RRR^{+}$ takes nonnegative real numbers, then the spectral radius of $M_f$ is the eigenvalue $\theta_1$ corresponding to the trivial character, and in addition $\theta_1$ is the largest eigenvalue of $M_f$ provided that all eigenvalues $\theta_{\chi}$ are real.
\end{thm}


\subsection{Eigenvalues of vertex-transitive graphs}
\label{Sec:VTG}

Vertex-transitive graphs can be considered as generalizations of Cayley graphs. As such we review several results about eigenvalues of vertex-transitive graphs in this section. Since Cayley graphs are vertex-transitive, all these results apply to Cayley graphs.

In \cite{Lovasz75}, Lov\'{a}sz proved that the determination of the spectrum of any vertex-transitive graph can be easily reduced to that of a Cayley graph. In fact, for any $\Ga$-vertex-transitive graph $G$ of order $n$, the lexicographic product $G'$ of $G$ with the empty graph of order $s = |\Ga|/n$ admits $\Ga$ as a vertex-regular group of automorphisms. Thus $G'$ is a Cayley graph of $\Ga$ and so its eigenvalues can be determined using Corollary \ref{BabaiS2a}. On the other hand, one can see that $A(G')$ is the Kronecker product of $A(G)$ and the all-$1$ matrix of size $s \times s$, and hence the eigenvalues of $G'$ are $0$ (with multiplicity $(s-1)n$) and $s \lambda_1, \ldots, s \lambda_n$, where $\lambda_1, \ldots, \lambda_n$ are the eigenvalues of $G$. In other words, the eigenvalues of $G$ can be obtained from that of $G'$ by throwing away $(s-1)n$ $0$'s and dividing each of the rest by $s$. In this way one can obtain, as stated in \cite[Theorem]{Lovasz75}, the eigenvalues of $G$ using the irreducible characters of $\Ga$.

In 1969, Petersdorf and Sachs \cite{PetersdorfS69} proved that for any vertex-transitive graph $G$ with degree $k$ and order $n$, if $n$ is odd then $k$ is the only simple eigenvalue of $G$, and if $n$ is even then any simple eigenvalue of $G$ is contained in the set $\{2i-k: 0 \le i \le k\}$. This result was strengthened by Sachs and Stiebitz \cite{SachsM82} in the follow way.

\begin{thm}
\emph{(\cite[Theorem 13]{SachsM82})}
Let $G$ be a vertex-transitive graph, loops and parallel edges being allowed, with degree $k$ and order $n = 2^e m$, where $m$ is odd. If $e = 0$, then $k$ is the only simple eigenvalue of $G$; if $e=1$, then $G$ has at most one simple eigenvalue other than $k$, and moreover it is of the form $4i-k$ for some $i \in \{0, 1, \ldots, (k-1)/2\}$ if it exists; if $e \ge 2$, then $G$ has at most $2^e$ simple eigenvalues including $k$, and each of them is contained in the set $\{2i-k: 0 \le i \le k\}$.
\end{thm}

Several other interesting results on simple eigenvalues of (directed or undirected) vertex-transitive graphs have also been proved in \cite{SachsM82, SachsM83}. In the following we present a few of them for undirected graphs possibly with loops and parallel edges. Let $\Ga$ be a permutation group of degree $n$. Then each element $g \in \Ga$ corresponds to a permutation matrix $P_g$. Denote by $z(\Ga)$ the number of vectors $\bx = (x_1, \ldots, x_n)^T$ with $x_1 = 1$ such that for every $g \in \Ga$ there is a number $a_{g}(\bx)$ with $P_{g} \bx = a_{g}(\bx) \bx$.

\begin{thm}
\emph{(\cite[Theorems 5--7]{SachsM82})}
Let $G$ be a vertex-transitive graph, loops and parallel edges being allowed. Then the following hold:
\begin{itemize}
\item[\rm (a)] $G$ has at most $z(\Aut(G))$ simple eigenvalues;
\item[\rm (b)] if the order of $G$ is $n = 2^e m$, where $m$ is odd, then $G$ has at most $2^e$ simple eigenvalues;
\item[\rm (c)] for any integer $n = 2^e m$, where $m$ is odd, there exist connected vertex-transitive graphs of order $n$ which have $2^e$ simple eigenvalues.
\end{itemize}
\end{thm}

\begin{thm}
\emph{(\cite[Proposition 4 and Theorem 14]{SachsM82})}
Let $G$ be a vertex-transitive graph of order $n$ and degree $k$ in which loops and parallel edges are allowed. Then the following hold:
\begin{itemize}
\item[\rm (a)] $G$ has at most $k+1$ simple eigenvalues;
\item[\rm (b)] if $G$ has more than $n/2$ simple eigenvalues, then $\Aut(G)$ is abelian.
\end{itemize}
\end{thm}

\begin{thm}
\emph{(\cite[Theorem 16]{SachsM82})}
Let $G$ be a connected vertex-transitive graph of order $n$ and degree $k$ in which loops and parallel edges are allowed. Suppose that $\Aut(G)$ is primitive on $V(G)$. Then the following hold:
\begin{itemize}
\item[\rm (a)] if $n$ is not a prime, then $G$ has exactly one simple eigenvalue, which is $k$;
\item[\rm (b)] if $n=p$ is a prime, then either $G$ has exactly one simple eigenvalue or it has exactly $p$ simple eigenvalues, and the latter occurs if and only if $\Aut(G)$ is a cyclic group of order $p$.
\end{itemize}
\end{thm}

\begin{thm}
\emph{(\cite[Corollary, p.88]{SachsM82})}
Let $G$ be a $\Ga$-vertex-transitive graph, loops and parallel edges being allowed, where $\Ga \le \Aut(G)$. Then $G$ has at most $t+1$ simple eigenvalues, where $t$ is the number of involutions of $\Ga$.
\end{thm}

\begin{thm}
\emph{(\cite[Theorem 19]{SachsM82})}
Let $G$ be a vertex-transitive graph with order $n$ all of whose eigenvalues are simple, loops and parallel edges being allowed. Then the following hold:
\begin{itemize}
\item[\rm (a)] $G$ is connected;
\item[\rm (b)] $n = 2^e$ for some nonnegative integer $e$, $G$ is regular of degree $k \ge 2^e - 1$, and $G$ contains the hypercube $H(e, 2)$ as a spanning subgraph;
\item[\rm (c)] if $-k$ is an eigenvalue of $G$, then the underlying simple graph of $G$ is $H(e, 2)$;
\item[\rm (d)] if $n > 2$, then $G$ contains parallel edges;
\item[\rm (e )] if $n > 1$, then $\Aut(G)$ is an elementary abelian $2$-group.
\end{itemize}
\end{thm}

As mentioned earlier, arc-transitive graphs are vertex-transitive. It is well known (see, for example, \cite[Proposition 16.7]{Biggs93}) that any arc-transitive graph has at most two simple eigenvalues, namely, $\pm k$, with $k$ the degree of the graph.

Recently, Guo and Mohar \cite{GuoM20} studied cubic vertex-transitive graphs with a nontrivial simple eigenvalue. Of course, a cubic vertex-transitive graph $G$ has largest eigenvalue $3$, and this eigenvalue is simple if and only if $G$ is connected. It is known that if $G$ is connected and $-3$ is also an eigenvalue of $G$, then $G$ is bipartite and $-3$ is a simple eigenvalue. Moreover, from the above-mentioned result of Petersdorf and Sachs \cite{PetersdorfS69} one can see that the only possible simple eigenvalues of $G$ are $\pm 3$ and $\pm 1$. In \cite{GuoM20} it was proved among other things that a cubic vertex-transitive graph having both $-1$ and $1$ as simple eigenvalues must be bipartite. Several families of cubic vertex-transitive graphs with $1$ as a simple eigenvalue and a classification of generalized Petersen graphs with $1$ as a simple eigenvalue were also given in \cite{GuoM20}.

Finally, it is easy to see that for any graph $G$ a permutation of $V(G)$ with corresponding permutation matrix $P$ is an automorphism of $G$ if and only if $A(G)P = PA(G)$. Thus, for such a matrix $P$ and any eigenvalue $\l$ of $G$ with eigenvector $\bv$, $P\bv$ is also an eigenvector of $A(G)$ with eigenvalue $\l$. This simple observation foreshadows close connections between the eigenvalues and automorphisms of a graph as one can find in \cite{CriscuoloKM80, Teranishi09}. We only mention the following result as it is most relevant to this survey.

\begin{thm}
\emph{(\cite[Theorem 2]{CriscuoloKM80})}
Let $G$ be a vertex-transitive graph. Then the number of distinct eigenvalues of $G$ is at most the rank of $\Aut(G)$.
\end{thm}


\section{Integral Cayley graphs}
\label{sec:ICG}

A graph is called \emph{integral} if its eigenvalues are all integers. This is equivalent to saying that all eigenvalues are rational numbers, because one can easily show that if all eigenvalues of a graph are rational then they must be integers. Which graphs are integral? This question was first asked by Harary and Schwenk \cite{Harary74} in 1973, with an immediate remark that the general problem appears to be challenging and intractable. It is known that there exist infinitely many integral graphs. However, in general it is nontrivial to construct integral graphs or determine all integral graphs in a given class of graphs. See \cite{Balinska02, Wang05} for two surveys on integral graphs. In this section we give a survey of integral Cayley graphs, and in Section \ref{sec:ComRing} we will focus on three families of integral Cayley graphs on finite commutative rings. Our treatment here has almost no overlap with the survey papers \cite{Balinska02, Wang05}.


\subsection{Characterizations of integral Cayley graphs}

\subsubsection{Integral circulant graphs, unitary Cayley graphs, and gcd graphs}
\label{CirUgcd}

In view of \eqref{eq:eigen-circ} one can see easily that not every circulant graph is integral. So it is natural to ask when a circulant graph is integral. This question was answered completely by So in \cite{So06}. Given an integer $n \ge 1$ and a positive divisor $d$ of $n$, define
\begin{equation}
\label{eq:gnd}
S_n(d) = \{a: 1 \le a \le n,\ \gcd(a,n)=d\}.
\end{equation}

\begin{thm}
\emph{(\cite[Theorem 7.1]{So06})}
\label{SoResult}
Let $n \ge 3$ be an integer. A circulant graph $\Cay(\mathbb{Z}_n, S)$ is integral if and only if $S$ is the union of $S_n(d)$ for some proper divisors $d$ of $n$.
\end{thm}

As a corollary it was noted \cite[Corollary 7.2]{So06} that there are at most $2^{\tau(n)-1}$ integral circulants on $n$ vertices, where $\tau(n)$ is the number of divisors of $n$. It was further conjectured that there are exactly $2^{\tau(n)-1}$ integral circulants on $n$ vertices (see \cite[Conjecture 7.3]{So06}). This conjecture is still open in its general form and partial results about it will be presented in Section \ref{subsec:cosp-isom}.

The sufficiency in Theorem \ref{SoResult} was obtained independently by {Klotz} and {Sander} in \cite{Klotz07}  using a slightly different language. Denote by $\ZZZ_n^{\times}$ the set of units (multiplicatively invertible elements) of ring $\mathbb{Z}_n$. The \emph{unitary Cayley graph} of ring $\mathbb{Z}_n$ is defined as the Cayley graph $\Cay(\mathbb{Z}_n, \ZZZ_n^{\times})$. This graph was first introduced in \cite{BerrizbeitiaG96, BerrizbeitiaG04} in the study of induced cycles of a given length and related chromatic uniqueness problem for the graph. In \cite{Fuchs04, Fuchs05}, longest induced cycles in $\Cay(\mathbb{Z}_n, \ZZZ_n^{\times})$ were further studied. In \cite{Klotz07}, {Klotz} and {Sander} studied several combinatorial properties of unitary Cayley graphs such as connectivity, diameter, chromatic number, clique number, independence number and perfectness. In the same paper they introduced the concept of gcd graphs as a generalization of unitary Cayley graphs. Define
\be
\label{eq:Dn}
D(n) = \{a: 1 \le a \le n,\text{ $a$ is a divisor of } n\}
\ee
to be the set of positive divisors of $n$. So $D(n) \setminus \{n\}$ is the set of positive proper divisors of $n$. Given $D \subseteq D(n) \setminus \{n\}$, where $n \ge 2$, the \emph{gcd graph} of the cyclic group $\mathbb{Z}_n$ with respect to $D$, denoted by $\ICG(n,D)$, is defined by
\be
\label{eq:ICG}
V(\ICG(n,D)) = \mathbb{Z}_n,\;\, E(\ICG(n,D)) = \{\{x, y\}: x,y\in \mathbb{Z}_n, \,\gcd(x-y,n)\in D\}.
\ee
In other words,
$$
\ICG(n,D) = \Cay(\mathbb{Z}_n, S_{n}(D)),
$$
where
\begin{equation}
\label{eq:SnD}
S_{n}(D) = \{a \in \mathbb{Z}_n: \gcd(a,n)\in D\}.
\end{equation}
In particular,
$$
\Cay(\mathbb{Z}_n,\ZZZ_n^{\times})\cong\ICG(n,\{1\}).
$$

It is evident that the graphs in Theorem \ref{SoResult} are precisely the gcd graphs of cyclic groups. Thus the following result obtained by {Klotz} and {Sander} in \cite{Klotz07} gives the sufficiency in Theorem \ref{SoResult}.

\begin{thm}
\emph{(\cite[Theorem 16]{Klotz07})}
\label{thm:intICG}
The gcd graphs of cyclic groups are all integral. In particular, the unitary Cayley graphs of rings $\ZZZ_n$ for $n \ge 3$ are integral circulant graphs.
\end{thm}

Moreover, in \cite{Klotz07} the eigenvalues of any gcd graph $\ICG(n,D)$ were given in terms of Euler's totient function $\varphi(n)$ and the M\"{o}bius function $\mu(n)$, where $\mu(n) = 0$ if $n$ is divisible by the square of a prime and $\mu(n) = (-1)^k$ otherwise with $k$ the number of prime factors of $n$. It is well known that
$$
\mu(n) = \sum_{1\le k \le n,\,\gcd(k,n)=1}\omega_n^{k}.
$$
The well-known \emph{Ramanujan sum} \cite{Hardy08} is defined as
\begin{equation}
\label{eq:RamSum}
c(k,n)=\sum_{1\le j\le n,\,\gcd(j,n)=1}\omega_n^{kj} = \mu\left(\frac{n}{\gcd(k,n)}\right)\cdot\frac{\varphi(n)}{\varphi\left(\frac{n}{\gcd(k,n)}\right)},\; k = 0, 1, \ldots, n-1.
\end{equation}
We know from number theory that $c(k,n)$ takes only integral values. In \cite[Theorem 13]{Klotz07}, {Klotz} and {Sander} proved that the eigenvalues of $\Cay(\mathbb{Z}_n,\ZZZ_n^{\times})$ are $c(k,n)$, $0 \le k \le n-1$. In general, they proved \cite[Theorem 16]{Klotz07} further that the eigenvalues of $\ICG(n,D)$ are
\begin{equation}\label{gcdEig1}
\lambda_l(n,D)=\sum_{d\in D} c(l,n/d),\; l = 0, 1, \ldots, n-1.
\end{equation}
Therefore, $\ICG(n,D)$ is integral as claimed in Theorem \ref{thm:intICG}.

Using the language of gcd graphs, Theorem \ref{SoResult} can be restated as follows.

\begin{cor}
\label{cor:intICG}
A circulant graph is integral if and only if it is a gcd graph of a cyclic group.
\end{cor}

\subsubsection{Integral Cayley graphs on abelian groups}
\label{InUniGCD}

Now that integral circulant graphs have been characterized, it is natural to move on to Cayley graphs on abelian groups. Let $A$ be a set and $\mathcal {F}$ a family of subsets of $A$. The \emph{Boolean algebra} $B(\mathcal {F})$ generated by $\mathcal {F}$ is the lattice of those subsets of $A$ each obtained by taking unions, intersections and complements of members of $\mathcal{F}$ in an arbitrary way but a finite number of times. The following result, due to {Klotz} and {Sander} \cite{Klotz10}, gives a sufficient condition for a Cayley graph  on an abelian group to be integral.

\begin{thm}
\emph{(\cite[Theorem 8]{Klotz10})}
\label{InteBoo}
A Cayley graph $\Cay(\Gamma, S)$ on a finite abelian group $\Gamma$ is integral provided that $S$ belongs to the Boolean algebra generated by the family of subgroups of $\Gamma$.
\end{thm}

With the help of Theorem \ref{SoResult}, {Klotz} and {Sander} also proved that this sufficient condition is necessary in the special case when $\Ga$ is a cyclic group.

\begin{thm}
\emph{(\cite[Theorem 10]{Klotz10})}
A circulant graph $\Cay(\ZZZ_n, S)$ with order $n \ge 2$ is integral if and only if $S$ belongs to the Boolean algebra generated by the family of subgroups of $\ZZZ_n$.
\end{thm}

In the same paper, {Klotz} and {Sander} further conjectured that the same should be true for any finite abelian group. One year later they confirmed their own conjecture in the following special case.

\begin{thm}
\emph{(\cite[Theorems 2 and 4]{Klotz11})}
Let $\Gamma=\ZZZ_{n_1}\oplus \cdots \oplus \ZZZ_{n_r}$ be an abelian group such that $\gcd(n_i,n_j)\le 2$ for $i\ne j$, where each $n_i \ge 2$. Then a Cayley graph $\Cay(\Ga, S)$ on $\Ga$ is integral if and only if $S$ belongs to the Boolean algebra generated by the family of subgroups of $\Ga$.
\end{thm}

The above-mentioned conjecture of Klotz and Sander \cite{Klotz10} was finally confirmed by {Alperin} and {Peterson} in \cite{Alperin12} in its general form. To present their result we need a few definitions first. Let $\Gamma$ be a finite group. Let $\widehat{\Gamma}$ be the set of characters of representations of $\Gamma$ over $\CCC$. A subset $A\subseteq \Gamma$ is called \emph{integral} if $\chi(A) \in \mathbb{Z}$ for every $\chi \in \widehat{\Gamma}$, where $\chi(A)=\sum_{a\in A}\chi(a)$. The group $\Ga$ is called \emph{cyclotomic} \cite{Alperin12} if the Boolean algebra generated by the family $\II(\Ga)$ of integral sets of $\Gamma$ equals the Boolean algebra generated by the family $\FF(\Ga)$ of subgroups of $\Ga$. {Alperin} and {Peterson} proved that any finite abelian group is cyclotomic \cite[Theorem 5.1]{Alperin12}, and that the converse is also true \cite[Theorem 8.3]{Alperin12}. Using the former statement, they obtained the following result, which confirms the above-mentioned conjecture of {Klotz} and {Sander} \cite{Klotz10}.

\begin{thm}
\label{thm:Alperin12}
\emph{(\cite[Corollary 7.2]{Alperin12})}
A Cayley graph $\Cay(\Gamma,S)$ on a finite abelian group $\Gamma$ is integral if and only if $S$ is an integral set. Equivalently, a Cayley graph $\Cay(\Gamma,S)$ on a finite abelian group $\Gamma$ is integral if and only if $S$ belongs to the Boolean algebra generated by the family of subgroups of $\Ga$.
\end{thm}

{Alperin} and {Peterson} also proved that for any group $\Ga$ (not necessarily abelian), all atoms of the Boolean algebra $B(\FF(\Ga))$ are integral and so $B(\FF(\Ga)) \subseteq B(\II(\Ga))$ (see \cite[Section 4]{Alperin12}), and for the dihedral group $D_{2n}$ of order $2n$, $B(\II(D_{2n}))$ equals the power set of $D_{2n}$ (see \cite[Theorem 6.1]{Alperin12}). (An \emph{atom} of a Boolean algebra is an element $x$ such that there exist exactly two elements $y$ satisfying $y \le x$, namely $x$ and $0$, where $\le$ is the partial order on the Boolean algebra such that $a \le b$ precisely when $a = a \wedge b$. In particular, a non-empty subset $X$ of $\Ga$ is an atom of $B(\FF(\Ga))$ if and only if for any $Y \in B(\FF(\Ga))$ such that $Y \subseteq X$ we have $Y = X$ or $Y = \emptyset$.)

Let $G$ be a graph and $\emptyset \neq N \subseteq \NNN \cup \{\infty\}$. The \emph{distance power} of $G$ with respect to $N$, denoted by $G^N$, is the graph with vertex set $V(G)$ such that $u$ and $v$ are adjacent if and only if the distance between them in $G$ belongs to $N$. In particular, for any positive integer $k$, $G^{\{1,2,\ldots,k\}}$ is called the \emph{$k$th power} of $G$ and is usually denoted by $G^k$. The second power $G^2$ is usually called the \emph{square} of $G$.

\begin{thm}\emph{(\cite[Theorem 1]{KlotzS12})}
\label{DisPow1}
Let $G = \Cay(\Gamma, S)$ be a Cayley graph on a finite abelian group $\Gamma$. If $G$ is integral, then for any set $N$ of nonnegative integers (possibly including $\infty$),  the distance power $G^N$ of $G$ is also an integral Cayley graph on $\Gamma$.
\end{thm}

In a recent paper Liu and Li \cite{LiuL16} studied distance powers of the unitary Cayley graph $\Cay(\ZZZ_n, \ZZZ_n^{\times})$. Since $\Cay(\ZZZ_n, \ZZZ_n^{\times})$ is integral, by Theorem \ref{DisPow1}, such distance powers are all integral circulant graphs.

Integral Cayley graphs on abelian groups $\Ga$ have been studied in \cite{BridgesM82a} via the algebra $\CCC \Ga$ over $\CCC$ generated by the permutation matrices representing the elements of $\Ga$ and the algebra over $\QQQ$ of those elements of $\CCC \Ga$ with rational entries and rational eigenvalues. The objects of study in \cite{BridgesM82a} are graphs admitting an abelian group $\Ga$ as an automorphism group regular on the vertex set. Note that such graphs are precisely Cayley graphs on $\Ga$.

\subsubsection{Integral normal Cayley graphs}

Recall that a Cayley graph $\Cay(\Gamma,S)$ is normal if $S$ is the union of some conjugacy classes of the group $\Gamma$. A subset $S$ of $\Ga$ is said to be \emph{power-closed} \cite{Godsil14} if, for every $x\in S$ and $y\in  \langle x\rangle$ with $\langle y\rangle = \langle x\rangle$, we have $y\in S$. A power-closed subset is called an Eulerian subset in \cite{GuoLMR19}. The next result follows from a more general result \cite[Theorem 1.3]{Godsil14} proved by Godsil and Spiga. The ``if'' part of this result was also proved in \cite{GuoLMR19} and \cite{KonstantinovaL20}.

\begin{thm}
\label{thm:GodsilP14}
\emph{(\cite[Theorem 1.1]{Godsil14}; see also \cite[Theorem]{GuoLMR19} and \cite[Theorem 1]{KonstantinovaL20})}
A normal Cayley graph $\Cay(\Gamma,S)$ is integral if and only if $S$ is power-closed.
\end{thm}

As noted in \cite{KonstantinovaL20}, this implies in particular that any normal Cayley graph on any $2$-group with connection set a set of involutions is integral. Another corollary of Theorem \ref{thm:GodsilP14} is that any normal Cayley graph $\Cay(\Gamma,S)$ such that the order of every element of $S$ belongs to $\{2,3,4,6\}$ is integral (see \cite[Corollary 1]{GuoLMR19}). This answers affirmatively Question 19.50(a) in \cite{Kourovka}, which was also answered by A. Abdollahi according to the comments below this question in  \cite{Kourovka}. Using Theorem \ref{thm:GodsilP14}, it was proved in \cite[Corollary 3]{GuoLMR19} that $\Cay(\Ga, S \setminus \Si)$ is integral for any subgroup $\Si$ of a group $\Ga$ and any power set $S$ of $\Ga$ closed under conjugation.

\subsubsection{Integral Cayley multigraphs}\label{ICMAHG}

A \emph{multiset} is a set $S$ together with a multiplicity function $\mu_S:S\rightarrow \mathbb{N}$, where for each $x\in S$ the positive integer $\mu_S(x)$ indicates the number of times that $x$ occurs in the multiset. As a convention we set $\mu_S(x) = 0$ for $x \notin S$. Let $\Gamma$ be a finite group. A multiset $S$ of some elements of $\Gamma$ is called \emph{inverse-closed} if $\mu_S(x) = \mu_S(x^{-1})$ for every $x \in S$. In this case the \emph{Cayley multigraph} on $\Gamma$ with connection set $S$, $\Cay(\Gamma, S)$, is defined to be the multigraph with vertex set $\Gamma$ such that the number of edges joining $x, y \in \Gamma$ is equal to $\mu_S(x y^{-1})$. In other words, $\Cay(\Gamma, S)$ is precisely the Cayley colour graph $\Cay(\Gamma, \mu_S)$. The \emph{adjacency matrix} of $\Cay(\Gamma, S)$ is the matrix with $(x, y)$-entry $\mu_S(x y^{-1})$, and its eigenvalues are called the \emph{eigenvalues} of $\Cay(\Gamma, S)$. A Cayley multigraph is callled \emph{integral} if all its eigenvalues are integers. In the special case when $\mu_S(x) = 1$ for every $x \in S$, $\Cay(\Gamma, S)$ is a Cayley graph in the usual sense. Hence we use the same notation for both Cayley graphs and Cayley multigraphs.

Denote by $\NN(\Ga)$ the set of normal subgroups of $\Gamma$. Define $B(\NN(\Ga))$ to be the Boolean algebra generated by $\NN(\Ga)$. This notion can be extended to multisets by including multiset operations. Formally, we take all atoms of the Boolean algebra $B(\NN(\Ga))$ and take all multisets that can be expressed as nonnegative integer combinations of these atoms. This defines the collection $\mathcal{C}(\Gamma)$ of multisets, which is called the \emph{integral cone} over $B(\NN(\Ga))$.

In 1982, Bridges and Mena \cite{BridgesM82} completely characterized integral Cayley multigraphs on abelian groups.

\begin{thm}\label{MCayInt1}
\emph{(\cite[Theorem 2.4]{BridgesM82}); see also \cite[Theorem 1]{DeVosKMA13}}
Let $\Gamma$ be an abelian group and $S$ an inverse-closed multiset of $\Gamma$. Then the Cayley multigraph $\Cay(\Gamma,S)$ is integral if and only if $S\in \mathcal{C}(\Gamma)$.
\end{thm}

In \cite{DeVosKMA13}, DeVos, Krakovski, Mohar and Ahmady generalized the sufficiency part of this result to Cayley multigraphs on any finite group.

\begin{thm}\label{MCayInt2}
\emph{(\cite[Theorem 2]{DeVosKMA13})}
Let $\Gamma$ be a finite group. Then for any $S\in \mathcal{C}(\Gamma)$ the Cayley multigraph $\Cay(G, S)$ is integral.
\end{thm}

Of course, by Theorem \ref{MCayInt1}, the converse statement in Theorem \ref{MCayInt2} holds for abelian groups. In \cite{DeVosKMA13}, DeVos, Krakovski, Mohar and Ahmady also investigated to what extent the converse would hold for some other groups. They provided necessary and sufficient conditions for the integrality of Cayley multigraphs over Hamiltonian groups (see \cite[Theorem 11]{DeVosKMA13}), where a \emph{Hamiltonian group} is a non-abelian Dedekind group, whilst a \emph{Dedekind group} is a group in which every subgroup is normal. It is known that every finite Hamiltonian group is the direct product $Q_8 \times A$ of the quaternion group $Q_8 = \{\pm 1, \pm i, \pm j, \pm k\}$ and a finite abelian group $A$. Given a Hamiltonian group $\Gamma = Q_8 \times A$ and an inverse-closed multiset $S$ of elements of $\Gamma$, for every $q \in Q_8$, define $B_q = \{a \in A: (q, a) \in S\}$ to be the multiset in which the multiplicity of $a \in B_q$ is equal to the multiplicity of $(q, a)$ in $S$. Then $B_1 = B_1^{-1}$, $B_{-1} = B_{-1}^{-1}$ and $B_{-q} = B_q^{-1}$ for $q \in Q_8 \setminus \{-1, 1\}$. For an irreducible character $\l$ of $A$, define $\hat{\l}(B_q) = \l(B_q) - \l(B_q^{-1}) = \l(B_q) - \l(B_{-q})$ for $q \in Q_8$.

\begin{thm}
\label{MCayInt3}
\emph{(\cite[Theorem 11]{DeVosKMA13})}
Let $\Gamma = Q_8 \times A$ be a Hamiltonian group, where $A$ is an abelian group, and let $S$ be an inverse-closed multiset of elements of $\Gamma$. Then the Cayley multigraph $\Cay(\Gamma, S)$ is integral if and only if the following hold:
\begin{itemize}
\item[\rm (a)] $B_1, B_{-1} \in \mathcal{C}(A)$;
\item[\rm (b)] the multiset union $B_q \cup B_{-q} \in \mathcal{C}(A)$, for every $q \in Q_8 \setminus \{-1, 1\}$;
\item[\rm (c)] $\hat{\l}(B_i)^2 + \hat{\l}(B_j)^2 + \hat{\l}(B_k)^2$ is a negative square of an integer, for every irreducible character $\l$ of $A$.
\end{itemize}
\end{thm}

Using this, DeVos, Krakovski, Mohar and Ahmady also obtained a few results on the integrality of Cayley graphs and Cayley multigraphs on some special families of Hamiltonian groups, including $Q_8 \times \ZZZ_p$ and $Q_8 \times \ZZZ_p^d$, where $p$ is a prime and $d \ge 2$. See \cite[Section 6]{DeVosKMA13} for details.


\subsection{Other families of integral Cayley graphs on abelian groups}

\subsubsection{Unitary finite Euclidean graphs}\label{UnFiEuGr}

In \cite{LiV13}, Li and Vinh defined the \emph{unitary finite Euclidean graphs} $T_n^{(d)}$ by
$$
V(T_n^{(d)})=\mathbb{Z}_n^d,\;\, E(T_n^{(d)})=\left\{(\mathbf{a},\mathbf{b}): \mathbf{a},\mathbf{b}\in \mathbb{Z}_n^d, ~ \sum_{i=1}^d(a_i-b_i)^2\in\mathbb{Z}_n^{\times}\right\},
$$
where $n \ge 2$, $d \ge 1$ and $\mathbb{Z}_n^\times$ is the set of units of $\mathbb{Z}_n$. Clearly,
$$
T_n^{(d)} =\Cay(\mathbb{Z}_n^d, S_n^{(d)}),
$$
where
$$
S_n^{(d)}=\left\{\mathbf{x} \in \mathbb{Z}_n^d: \sum_{i=1}^dx_i^2\in\mathbb{Z}_n^{\times}\right\}.
$$
The following result was proved by Li and Vinh in \cite{LiV13}.

\begin{thm}
\emph{(\cite[Theorem 3.6]{LiV13})}
Let $n \ge 2$ and $d \ge 1$ be integers. Then the unitary finite Euclidean graph $T_n^{(d)}$ is integral when $n$ is odd or $d$ is even.
\end{thm}

In the same paper Li and Vinh also conjectured that $T_n^{(d)}$ is integral for any $n \ge 2$ and $d \ge 1$ (see \cite[Conjecture 3.7]{LiV13}). This conjecture was proved in \cite{DaiHHV13} as a special case of a more general result. Given $U \subseteq \ZZZ_n$, the \emph{distance graph} over $\ZZZ_n^d$ generated by $U$, denoted by $T_n^{(d)}(U)$, was defined \cite{DaiHHV13} in the same way as $T_n^{(d)}$ except that $\mathbb{Z}_n^{\times}$ is replaced by $U$. In particular, $T_n^{(d)} = T_n^{(d)}(\mathbb{Z}_n^{\times})$. Obviously,
$$
T_n^{(d)}(U) = \Cay(\mathbb{Z}_n^d, S_n^{(d)}(U)),
$$
where $S_n^{(d)}(U)$ is defined in the same way as $S_n^{(d)}$ with $\mathbb{Z}_n^{\times}$ replaced by $U$. Since this is a Cayley graph on an abelian group, by Corollary \ref{BabaiS1a} its eigenvalues can be easily found to be
$$
\l_{\mathbf{b}}(U) = \sum_{\mathbf{x} \in S_n^{(d)}(U)} \exp\{2 \pi i (\mathbf{b}^T \mathbf{x})/n\},\, \mathbf{b} \in \ZZZ_n^d,
$$
where $\mathbf{b}^T$ is the transpose of $\mathbf{b}$ and $\mathbf{b}^T \mathbf{x}$ is the dot product of $\mathbf{b}$ and $\mathbf{x}$.

Recall from \eqref{eq:gnd} that $S_{n}(d) = \{1 \le a \le n: \gcd(a, n) = d\}$ for any divisor $d$ of $n$.

\begin{thm}
\emph{(\cite[Theorem 3]{DaiHHV13})}
Let $n \ge 2$ and $d \ge 1$ be integers. Let $D$ be a set of divisors of $n$ and $U = \cup_{d \in D} S_{n}(d)$. Then $T_n^{(d)}(U)$ is integral.
\end{thm}

In the special case when $D = \{1\}$, this implies that the above-mentioned conjecture (\cite[Conjecture 3.7]{LiV13}) is true.

\begin{thm}
\emph{(\cite{DaiHHV13})}
The unitary finite Euclidean graph $T_n^{(d)}$ is integral for any integers $n \ge 2$ and $d \ge 1$.
\end{thm}

\subsubsection{NEPS of complete graphs, gcd graphs of abelian groups, and generalized Hamming graphs}
\label{HammingE}

Let $G_1, \ldots,G_d$ be graphs, and let $\emptyset\neq \mathcal{B}\subseteq \mathbb{Z}_2^{d}\setminus \{\mathbf{0}\}$, where $d \ge 1$ and $\mathbf{0}=(0,\ldots,0) \in \mathbb{Z}_2^{d}$. The \emph{non-complete extended $p$-sum} (NEPS) of $G_1,\ldots,G_d$ with \emph{basis} $\mathcal{B}$ (see \cite[Definition 2.5.1]{Cvetkovic10}), denoted by $\NEPS(G_1,\ldots,G_d;\mathcal{B})$, is the graph with vertex set $V(G_{1})\times \cdots \times V(G_{d})$ in which two vertices $(x_{1}, \ldots, x_{d})$ and $(y_{1}, \ldots, y_{d})$ are adjacent if and only if there exists $\beta=(\beta_{1}, \ldots, \beta_{d})\in \mathcal{B}$ such that $x_{i}=y_{i}$ whenever $\beta_{i}=0$ and $x_{i}$ is adjacent to $y_{i}$ in $G_{i}$ whenever $\beta_{i}=1$. This notion is a generalization of several graph operations such as tensor product, Cartesian product, and strong product. In fact, $\NEPS(G_1,\ldots,G_d; \{(1, \ldots, 1)\})$ is simply the \emph{tensor product} $G_1 \otimes \cdots \otimes G_d$; if $\mathcal{B}_0 = \{(1, 0, \ldots, 0), \ldots, (0, 0, \ldots, 1)\}$ is the standard basis of $\mathbb{Z}_2^{d}$, then $\NEPS(G_1,\ldots,G_d;\mathcal{B}_0)$ is the \emph{Cartesian product} $G_1 \Box \cdots \Box G_d$. Thus $\NEPS(K_{n_1},\ldots,K_{n_d}; \mathcal{B}_0)$ is the Hamming graph $H(n_1, \ldots, n_d) = K_{n_1} \Box \cdots \Box K_{n_d}$, where $n_1, \ldots, n_d \ge 2$. In particular, $\NEPS(K_{2},\ldots,K_{2}; \mathcal{B}_0)$ is the hypercube $H(d,2)$ of dimension $d$. In general, for any $\emptyset \neq \mathcal{B}\subseteq \mathbb{Z}_2^{d} \setminus \{\mathbf{0}\}$, $\NEPS(K_{2},\ldots,K_{2}; \mathcal{B})$ is a cubelike graph. (Recall from Section \ref{subsec:spec-Cay} that a \emph{cubelike graph} \cite{Lovasz75} is precisely a Cayley graph on an elementary abelian $2$-group.) The NEPS construction is a useful tool in graph theory. For example, it was used by Razborov \cite{Razborov92}, Huang and Sudakov \cite{HuangS12} and Cioab\u{a} and Tait \cite{CioabaT11} to construct counterexamples to the Alon-Saks-Seymour and the rank-coloring conjectures.

A fundamental result (see \cite[Theorem 2.5.4]{Cvetkovic10}) on NEPS asserts that, if $\l_{i 1}, \ldots, \l_{i n_i}$ are eigenvalues of $G_i$ for $1 \le i \le d$, where $n_i = |V(G_i)|$, then $\NEPS(G_1,\ldots,G_d;\mathcal{B})$ has eigenvalues
\be
\label{eq:mu}
\mu_{i_1, \ldots, i_d} = \sum_{(\beta_{1}, \ldots, \beta_{d})\in \mathcal{B}} \l_{1 i_1}^{\beta_1} \cdots \l_{d i_d}^{\beta_d},\ i_t = 1, 2, \ldots, n_t,\ t = 1, 2, \ldots, d.
\ee
As an immediate consequence, we get the following result.

\begin{cor}
\label{NEPSint}
Every NEPS of integral graphs is integral.
\end{cor}

This result is relevant since NEPS graphs of complete graphs form an interesting family of integral Cayley graphs on abelian groups, as we now explain. Let $\mathbf{m}=(m_1, \ldots, m_d)$, $\mathbf{n}=(n_1, \ldots, n_d)$ and $\mathbf{k}=(k_1, \ldots, k_d)$ be $d$-tuples of positive integers. If $k_i = \gcd(m_i, n_i)$ for $i = 1, \ldots , d$, then $\mathbf{k}$ is called the \emph{greatest common divisor} of $\mathbf{m}$ and $\mathbf{n}$, written $\mathbf{k} = \gcd(\mathbf{m}, \mathbf{n})$. Define
\begin{equation}
\label{eq:DDn}
\mathbf{D}(\mathbf{n}) = \{(a_1, \ldots, a_d): 1 \le a_i \le n_i \text{ and } a_i \text{ is a divisor of } n_i \text{ for } 1 \le i \le d\}.
\end{equation}

Let $\Gamma=\mathbb{Z}_{n_1}\oplus\cdots\oplus\mathbb{Z}_{n_d}$ be an abelian group, where $d \ge 1$ and each $n_i \ge 2$. Let $\mathbf{n} = (n_1, \ldots, n_d)$ and $\mathbf{D} \subseteq \mathbf{D}(\mathbf{n}) \setminus \{\mathbf{n}\}$. Denote by $S_{\Ga}(\mathbf{D})$ the set of elements $\mathbf{x} = (x_1, \ldots, x_d)$ of $\Gamma$ such that $0 \le x_i \le n_i - 1$ for $1 \le i \le d$ and $\gcd(\mathbf{x},\mathbf{n}) \in \mathbf{D}$. (If $\mathbf{D} = \{\mathbf{k}\}$, then we write $S_{\Ga}(\mathbf{k})$ in place of $S_{\Ga}(\{\mathbf{k}\})$.) The \emph{gcd graph} of the abelian group $\Gamma$ with respect to $\mathbf{D}$, denoted\footnote{As seen in Section \ref{CirUgcd}, the acronym $\ICG$ arose from Integral Circulant Graphs. We use the same acronym for gcd graphs of abelian groups here and gcd graphs of unique factorization domains later (Section \ref{subsec:eng-ring-gcd}) to indicate that all these graphs are defined in the same fashion, though they are not necessarily circulants in general.} by $\ICG(\mathbf{n}, \mathbf{D})$, is defined \cite{Klotz11, KlotzS13} to be the Cayley graph $\Cay(\Gamma, S_{\Ga}(\mathbf{D}))$. This definition is a generalization of the gcd graphs of cyclic groups, because in the special case when $\Gamma=\mathbb{Z}_{n}$ (so $\mathbf{n} = (n)$) and $\mathbf{D} = D \subseteq D(n) \setminus \{n\}$ is a set of positive proper divisors of $n$, $\ICG(\mathbf{n}, \mathbf{D})$ is exactly the gcd graph $\ICG(n, D)$ of the cyclic group $\ZZZ_n$. Note that $\ICG(\mathbf{n}, \mathbf{D})$ relies on the specific expression of the abelian group $\Ga$ as the direct sum $\mathbb{Z}_{n_1}\oplus\cdots\oplus\mathbb{Z}_{n_d}$.

In \cite[Theorems 2.5 and 2.6]{KlotzS13}, it was shown that a graph isomorphic to a gcd graph of an abelian group if and only if it is isomorphic to the NEPS of some complete graphs. In fact, the proof of \cite[Theorem 2.6]{KlotzS13} implies that
\begin{equation}
\label{eq:NEPSICG}
\NEPS(K_{n_1},\ldots, K_{n_d}, \mathcal{B}) \cong \ICG(\mathbf{n}, \mathbf{D}),
\end{equation}
where $\mathbf{D} =\cup_{\beta \in \mathcal{B}} \mathbf{D}(\beta)$, where for each $\beta=(\beta_1,\ldots,\beta_d)\in \mathcal{B}$, $\mathbf{D}(\beta)$ consists of those $(a_1, \ldots, a_d)$ such that $a_i = n_i$ if $\beta_i = 0$ and $a_i \ge 1$ is a proper divisor of $n_i$ if $\beta_i = 1$. Since complete graphs are integral, this together with Corollary \ref{NEPSint} implies the following result.

\begin{cor}
\label{gcdint}
\emph{(\cite[Proposition 3.2]{KlotzS13}; see also \cite{Klotz11})}
NEPS graphs of complete graphs are integral. In other words, every gcd graph of any abelian group is integral.
\end{cor}

In \cite{KlotzS13}, Klotz and Sander proved further that gcd graphs of abelian groups have a particularly simple eigenspace structure, namely, every eigenspace of the adjacency matrix of a gcd graph has a basis with entries $-1,0,1$ only. Corollary \ref{gcdint} implies the following result due to Lov\'{a}sz.

\begin{cor}
\label{cubelike-int}
\emph{(\cite{Lovasz75})}
Every cubelike graph is integral.
\end{cor}

Let $\Gamma=\mathbb{Z}_{n_1}\oplus\cdots\oplus\mathbb{Z}_{n_d}$ be an abelian group, where $d \ge 1$ and $n_i \ge 2$ for $1 \le i \le d$. Let $T=\{t_1,\ldots,t_k\}$ be a set of integers such that $1\le t_i\le d$ for each $i$. The \emph{generalized Hamming graph} $H(n_1,\ldots,n_d;T)$ is defined to have vertex set $\Gamma$ such that $\mathbf{x} = (x_1, \ldots, x_d) \in \Gamma$ and $\mathbf{y} = (y_1, \ldots, y_d) \in \Gamma$ are adjacent if and only if the \emph{Hamming distance} $w(\mathbf{x}-\mathbf{y})$ between them belongs to $T$. As observed in \cite[Example 2.8]{KlotzS13}, $H(n_1,\ldots,n_d;T)$ is an NEPS graph of $K_{n_{1}},\ldots,K_{n_{d}}$. Thus it is integral by Corollary \ref{cubelike-int} and its eigenvalues can be computed using formula \eqref{eq:mu}.

\begin{thm}
\emph{(\cite[Proposition 14]{Klotz10})}
Every generalized Hamming graph is an integral Cayley graph.
\end{thm}

In the special case when $T=\{1\}$, we have $H(n_1,\ldots,n_d;\{1\})\cong K_{n_1}\boxempty\cdots\boxempty  K_{n_d} = H(n_1, \ldots, n_d)$. In particular, $H(q,\ldots,q;\{1\})$ is the classical Hamming graph $H(d, q)$. It is well known \cite{Cvetkovic10} that the spectra of the Cartesian product of graphs can be expressed in terms of that of the factor graphs. Since the eigenvalues of any complete graph are known, this implies \cite{BrouwerCN89} that the eigenvalues of $H(d,q)$ are $q(d-j)-d$ with multiplicities $\choose{d}{j}(q-1)^j$ for $0 \le j \le d$. In \cite{Sander10}, {T. Sander} proved that $H(n_1,\ldots,n_d;\{d\})$ is isomorphic to the complement graph of $H(n_1,\ldots,n_d;\{1,\ldots,d-1\})$, and for distinct primes $p_1,\ldots, p_d$, $H(p_1,\ldots,p_d;\{d\})$ is isomorphic to the unitary Cayley graph $\Cay(\mathbb{Z}_n,\ZZZ_n^{\times})$, where $n = \prod_{i=1}^d p_i$.

In \cite{LiLZZ17}, Li, Liu, Zhang and Zhou characterized isomorphisms between the squares of generalized Hamming graphs and the NEPS of some complete graphs. As an application they determined the eigenvalues of the squares of generalized Hamming graphs. In view of Theorem \ref{DisPow1} and Corollary \ref{NEPSint}, such graphs are all integral Cayley graphs.

\subsubsection{Sudoku graphs and positional Sudoku graphs}

Let $n \ge 2$ be an integer. An \emph{$n$-Sudoku} is an arrangement of $n\times n$ square blocks each consisting of $n\times n$ cells, with each cell filled with a number (colour) from $\{1, 2, \ldots, n^2\}$ such that every block, row or column contains all of the colours $1, 2, \ldots, n^2$. The \emph{Sudoku graph} $\Sud(n)$ is defined \cite{Sander09, Klotz10} to have vertices the $n^4$ cells of an $n$-Sudoku such that distinct vertices (cells) are adjacent if and only if they are in the same block, row or column. The \emph{positional Sudoku graph} $\SudP(n)$ is defined \cite{Klotz10} to have vertices the $n^4$ cells of an $n$-Sudoku such that distinct vertices (cells) are adjacent if and only if they are in the same block, row, column or in the same position of their respective blocks. It can be shown that both $\Sud(n)$ and $\SudP(n)$ are Cayley graphs on $\ZZZ_n^4$; see \cite[Section 4.2]{Klotz10} for an explicit expression of the corresponding connection sets. In \cite[Lemma 3.1]{Sander09}, it was shown that $\Sud(n)$ is the NEPS of four copies of $K_n$ with basis $\{(0,1,0,1), (1,1,0,0), (0,0,1,1), (1,0,0,0), (0,1,0,0), (0,0,1,0), (0,0,0,1)\}$. Using this it was proved in \cite[Theorem 3.2]{Sander09} that $\Sud(n)$ is integral and each of its eigenspaces admits a basis whose entries are all from the set $\{-1, 0, 1\}$.

\begin{thm}
\emph{(\cite{Sander09, Klotz10})}
Let $n \ge 2$ be an integer. Then both $\Sud(n)$ and $\SudP(n)$ are integral Cayley graphs. Moreover, the following hold:
\begin{itemize}
\item[\rm (a)] the spectrum of the Sudoku graph $\Sud(n)$ is
\begin{eqnarray*}
\left(\begin{array}{cccccc}
3n^2-2n-1 & 2n^2-2n-1   &  n^2-n-1   &  n^2-2n-1    & -1                  &   -1-n\\
1                 &   2(n-1)         &  2n(n-1)    &  (n-1)^2       &n^2(n-1)^2  &   2n(n-1)^2
\end{array}
\right);
\end{eqnarray*}
\item[\rm (b)] the spectrum of the positional Sudoku graph $\SudP(n)$ is
\begin{eqnarray*}
\left(\begin{array}{cccccc}
4n(n-1)     & 2n^2-3n       &  n(n-2)        &  0                 & -n             &   -2n  \\
1                 &   4(n-1)         &  4(n-1)^2    &  (n-1)^4      &4(n-1)^3  &   2(n-1)^2
\end{array}
\right).
\end{eqnarray*}
\end{itemize}
\end{thm}

In the expository paper \cite{KlotzSander10}, Klotz and Sander used symmetry properties of Sudoku graphs to determine their eigenvalues. In particular, they recovered the above-mentioned result that all eigenvalues of Sudoku graphs are integral.

\subsubsection{Pandiagonal Latin square graphs}

Let $n \ge 2$ be an integer. The \emph{pandiagonal Latin square graph} $\PLSG(n)$ is defined \cite{Klotz10} to have vertices the $n^2$ positions of an $n \times n$ matrix such that distinct vertices (positions) are adjacent if and only if they are in the
same row, the same column, the same (broken) parallel to the main diagonal, or the same (broken) parallel to the secondary diagonal. This is a Cayley graph on $\ZZZ_n^2$ whose connection set can be found in \cite[Section 4.3]{Klotz10}. In contrast, the usual Latin square graphs (see, for example, \cite[Section 10.4]{GodsilG2001}) are not in general Cayley graphs.

\begin{thm}
\emph{(\cite[Proposition 16]{Klotz10})}
Let $n \ge 2$ be an integer. The pandiagonal Latin square graph $\PLSG(n)$ is an integral Cayley graph. Moreover, the following hold:
\begin{itemize}
\item[\rm (a)] if $n$ is odd, then the spectrum of $\PLSG(n)$ is
\begin{eqnarray*}
\left(\begin{array}{ccc}
4n-4          &   n-4              &  -4  \\
1                 &   4(n-1)         &  n^2-4n+3
\end{array}
\right);
\end{eqnarray*}
\item[\rm (b)] if $n$ is even, then the spectrum of $\PLSG(n)$ is
\begin{eqnarray*}
\left(\begin{array}{cccccc}
4n-5          & 2n-5       &  n-3       &  n-5        & -3             &  -5 \\
1                 &   1           &  n           &  3n-6      &n^2/2-n  &   n^2/2-3n+4
\end{array}
\right).
\end{eqnarray*}
\end{itemize}
\end{thm}


\subsection{Cayley graphs on symmetric and alternating groups}
\label{subsec:sn}

Cayley graphs on symmetric groups arise from various fields of discrete mathematics, including the study of Erd\H{o}s-Ko-Rado theorems \cite{GodsilM16} for permutations. For example, a long-standing conjecture on $k$-intersecting subsets of $S_n$ was proved by Ellis, Friedgut and Pilpel in \cite{EllisFP11} with the help of the smallest eigenvalue of the Cayley graph on $S_n$ with respect to the set of permutations with less than $k$ fixed points. Here we focus on eigenvalues of several interesting families of Cayley graphs on symmetric groups. We will also review a recent result on the integrality of a family of Cayley graphs on alternating groups.

As usual we use $S_n$ and $A_n$, respectively, to denote the symmetric and alternating groups on $[n] = \{1, 2, \ldots, n\}$, where $n \ge 2$. The Cayley graph $\Cay(S_{n}, S)$ on $S_n$ with connection set $S=\{(12),(13),\ldots,(1n)\}$ is known as the \emph{star graph} of degree $n-1$ in theoretical computer science and probability theory. Using GAP \cite{TGap04} and the GRAPE package of L. H. Soicher, the spectra of $\Cay(S_{n}, S)$  for $3 \le n \le 6$ was computed in the proof of \cite[Lemma 2.13]{Abdollahi09}:
$$
\Spec(\Cay(S_{n}, S)) =
\begin{cases}
(-2, 2,  (-1)^2, 1^2), & n = 3 \\
(-3, 3, (-2)^{6}, 2^6, (-1)^{3}, 1^3, 0^4), & n = 4 \\
(-4, 4, (-3)^{12}, 3^{12}, (-2)^{28}, 2^{28}, (-1)^{4}, 1^4, 0^{30}), & n = 5 \\
(-5, 5, (-4)^{20}, 4^{20}, (-3)^{105}, 3^{105}, (-2)^{120}, 2^{120}, (-1)^{30}, 1^{30}, 0^{168}), & n = 6.
\end{cases}
$$
In particular, this implies that $\Cay(S_{n}, S)$ is integral for $3 \le n \le 6$. Suggested by this computational result, {Abdollahi} and {Vatandoost} posed the following conjecture.

\begin{conj}
\emph{(\cite[Conjecture 2.14]{Abdollahi09})}
\label{SnConj1}
Let $S = \{(12), (13), \ldots , (1n)\}$, where $n \ge 4$. Then the star graph $\Cay(S_{n},S)$ is integral. Moreover, $\{0,\pm1, \ldots ,\pm(n - 1)\}$ is the set of all distinct eigenvalues of $\Cay(S_{n},S)$.
\end{conj}

This conjecture was confirmed by {Krakovski} and {Mohar} in \cite{Krakovski12}.

\begin{thm}
\label{thm:star}
\emph{(\cite[Theorem 1]{Krakovski12})}
Let $S = \{(12), (13), \ldots , (1n)\}$, where $n \ge 4$. Then $n-l$ and $-(n-l)$ are eigenvalues of $\Cay(S_{n},S)$ with multiplicity at least $\choose{n-2}{l-1}$, for $1\le l \le n-1$, and $0$ is an eigenvalue of $\Cay(S_{n},S)$ with multiplicity at least $\choose{n-1}{2}$.
\end{thm}

In \cite[Theorem 1]{AvgustinovichKK16}, analytic formulas for the multiplicities of eigenvalues $\pm (n-k)$ of $\Cay(S_n, S)$ for $k = 2, 3, 4, 5$ were obtained. In \cite{KhomyakovaK15}, a Chapuy-Feray combinatorial approach was used to obtain the multiplicities of eigenvalues of $\Cay(S_n, S)$, and the exact values are calculated for $n \le 10$.

Later it was noted in \cite{Chapuy12} (see also \cite{Renteln11} and \cite{Krakovski12}) that the truth of Conjecture \ref{SnConj1} was implied by certain properties of the Jucys-Murphy elements, discovered by Jucys \cite{Jucys74} and independently by Flatto, Odlyzko and Wales \cite{FOW85}. Let $\lambda$ be a partition of $n$. A \emph{standard Young tableau (SYT) of shape $\lambda$} is a filling of the Ferrers diagram of $\lambda$ with elements $\{1, 2, \ldots, n\}$ in such a way that elements increase along rows and columns (in particular all elements are distinct, and each element appears exactly once). Denote by $\TT(\l)$ the set of SYTs of shape $\l$ and $f_{\l} = |\TT(\l)|$. Given $T \in \TT(\l)$ and a box $\Box$ of the Ferrers diagram of $T$, the content of $\Box$ is the difference between its ordinate and abscissa. Define $c_{T}(i)$ to be the content of the box in which label $i$ appears in $T$, and set $I_{\l}(k) = |\{T \in \TT(\l): c_{T}(n)=k\}|$. As noticed in \cite[Corollary 2.1]{Chapuy12}, the above-mentioned result in \cite{FOW85, Jucys74} implies the following result.

\begin{thm}
\emph{(\cite[Corollary 2.1]{Chapuy12})}
\label{thm:star-int}
Let $S = \{(12), (13), \ldots , (1n)\}$, where $n \ge 2$. Then the star graph $\Cay(S_{n},S)$ is integral, and the multiplicity of any integer $k$ as an eigenvalue of $\Cay(S_{n},S)$ is equal to $\sum_{\lambda\in\mathcal {P}(n)}f_{\lambda}I_{\lambda}(k)$, where $\mathcal {P}(n)$ is the set of partitions of $n$.
\end{thm}

For any positive integers $n \ge 3$ and $l$ with $n > 2l$, Goryainov, Kabanov, Konstantinova, Shalaginov and Valyuzhenich presented in \cite{GoryainovKKSV20} a family of eigenfunctions of the star graph of degree $n-1$ corresponding to eigenvalue $n - l - 1$ and linked these functions to the standard basis of a Specht module. In the same paper they also proved that any eigenfunction of the star graph of degree $n-1$ corresponding to the second largest eigenvalue $n-2$ can be reconstructed from its values on the second neighbourhood of a vertex.

For an integer $r$ with $2\le r \le n$, let $\Cy(r)$ be the set of all $r$-cycles in $S_n$ which do not fix $1$. That is,
$$
\Cy(r)=\{\alpha\in S_n: \alpha(1)\ne1 \text{~and $\alpha$ is an $r$-cycle}\}.
$$
Note that $\Cy(2)=\{(12),(13),\ldots,(1n)\}$ and so $\Cay(S_n,\Cy(r))$ can be thought as a generalization of the star graph $\Cay(S_{n},S)$ above. In \cite{Chen14}, Chen, Ghorbani and Wong proved that the graphs $\Cay(S_n,\Cy(r))$ are all integral. In fact, they proved a more general result (Theorem \ref{thm:Chen14} below). For any $T \subseteq [n]$, denote by $S_n(T)$ the set of elements of $S_n$ which fix each point in $[n] \setminus T$. As usual, for $\sigma \in S_n$, let $\mathrm{Cl}(\sigma) = \{\alpha^{-1} \sigma \alpha: \alpha \in S_n\}$ be the conjugacy class of $\sigma$. A subset $S$ of $S_n$ is said to be \emph{nicely separated} \cite{Chen14} if there exists a partition $\{T_1, \ldots, T_l\}$ of $[n]$ with $l \ge 2$ such that $S = (\cup_{\sigma \in S}\mathrm{Cl}(\sigma)) \cap (S_n \setminus (\cup_{i=1}^{l} S_n(T_i)))$. Note that any nicely separated subset $S$ of $S_n$ gives rise to an undirected Cayley graph $\Cay(S_n, S)$.

\begin{thm}
\label{thm:Chen14}
\emph{(\cite[Theorem 3.4]{Chen14})}
Let $S$ be a nicely separated subset of $S_n$, where $n \ge 2$. Then $\Cay(S_n, S)$ is integral.
\end{thm}

Since $\Cy(r)$ is nicely separated with $T_1 = \{1\}$ and $T_2 = \{2, \ldots, n\}$, this implies the following corollary.

\begin{cor}
\emph{(\cite[Theorem 1.4]{Chen14})}
$\Cay(S_n,\Cy(r))$ is integral for $2 \le r \le n$.
\end{cor}

It is known that central characters are algebraic integers \cite[Theorem 3.7]{Isaacs76} and the characters of $S_n$ are integers \cite[Corollary 22.17]{James01}. Combining these with \eqref{eq:Conjeign}, we obtain:
\begin{thm}
\emph{(\cite[Corollary 1.2]{Chen14})}
\label{nomalSn}
All normal Cayley graphs on $S_n$ are integral.
\end{thm}

It is well known that the irreducible representations of $S_n$ are in one-to-one correspondence with the partitions of integer $n$, where a partition of $n$ is an $r$-tuple $\l = (\l_1, \ldots, \l_{r})$ of positive integers, for some $r \ge 1$, such that $\l_1 \ge \cdots \ge \l_{r}$ and $n = \l_1 + \cdots + \l_{r}$. So the irreducible characters of $S_n$ can be indexed by $\l$ and written as $\chi_{\l}$, and in view of \eqref{eq:Conjeign} the eigenvalues of any normal Cayley graph $\Cay(S_n, S)$ on $S_n$ are
$$
\eta_{\l}(\Cay(S_n, S)) = \frac{1}{\chi_{\l}(1)} \sum_{s \in S}\chi_{\l}(s),
$$
with $\l$ running over all partitions of $n$. Moreover, the multiplicity of $\eta_{\l}$ is $\chi(1)^2$. A partition $\l$ of $n$ is said to be associated with the smallest eigenvalue of $\Cay(S_n, S)$ if $\eta_{\l}$ is the smallest eigenvalue of $\Cay(S_n, S)$.

An interesting normal Cayley graph on $S_n$ is the \emph{derangement graph} on $[n]$, which is defined as the Cayley graph $\Cay(S_n, \mathscr {D}_n)$, where $\mathscr  {D}_n =\{\sigma \in S_n: \sigma(i)\ne i \text{ for each } i \in [n]\}$ is the set of derangements of $[n]$. This Cayley graph is normal as $\mathscr {D}_n$ is closed under conjugation. Thus $\Cay(S_n, \mathscr  {D}_n)$ is integral by Theorem \ref{nomalSn}. In \cite{Renteln07}, {Renteln} gave several interesting formulas for the eigenvalues of $\Cay(S_n, \mathscr  {D}_n)$ (see Theorems 3.2, 3.3, 4.2, 6.1 and 6.5 in \cite{Renteln07}), including the formulas in the next two theorems.

\begin{thm}
\emph{(\cite[Theorem 3.2]{Renteln07})}
The eigenvalues of $\Cay(S_n, \mathscr  {D}_n)$ are given by
$$
\eta_{\lambda}=\sum_{k=0}^n(-1)^{n-k}\frac{n!}{(n-k)!}\frac{f_{\lambda/(k)}}{f_{\lambda}},
$$
where $\lambda$ runs over all partitions of  $n$.
\end{thm}

The \emph{complete factorial symmetric function} $\omega_k$ is defined as
$$
\omega_k(z_1,\ldots,z_r)=\sum_{1\le i_1\le \cdots \le i_k\le r}\left(z_{i_1}-i_1+2-1\right)\cdots \left(z_{i_k}-i_k+2-k\right).
$$

\begin{thm}
\emph{(\cite[Theorem 6.1]{Renteln07})}
The eigenvalues of $\Cay(S_n, \mathscr  {D}_n)$ are given by
$$
\eta_{\lambda}=\sum_{k=0}^n(-1)^{n-k}\omega_k(\mu_1,\ldots,\mu_r),
$$
where $\lambda=(\lambda_1,\ldots,\lambda_r)$ runs over all partitions of $n$ and $\mu_i=\lambda_i+r-i$ for each $i$.
\end{thm}

Using these formulas, Renteln \cite{Renteln07} proved part (a) of the following theorem which settled affirmatively a conjecture of Ku and Wong \cite[Conjecture 1]{Ku07}. Part (b) of this theorem was proved by Deng and Zhang \cite{Deng11} with the help of the work in \cite{Renteln07}.

\begin{thm}
\label{thm:RentDeng}
The following hold:
\begin{itemize}
\item[\rm (a)] the smallest eigenvalue of $\Cay(S_n, \mathscr  {D}_n)$ is equal to $-\frac{1}{n-1}|\mathscr  {D}_n|$ \emph{(\cite[Theorem 7.1]{Renteln07})};
\item[\rm (b)] the second largest eigenvalue of $\Cay(S_n, \mathscr  {D}_n)$ is equal to $\frac{n-1}{n-3}|\mathscr {D}_{n-2}|$ \emph{(\cite[Theorem 1.1]{Deng11})}.
\end{itemize}
\end{thm}

Further spectral properties of $\Cay(S_n, \mathscr  {D}_n)$ were obtained in \cite{Ku10, Ku13}. It is well known that the eigenvalues of this graph can be indexed by partitions of $n$. The question about how these eigenvalues are determined by the shape of their corresponding partitions was investigated in \cite{Ku10}. Lower and upper bounds on the absolute values of these eigenvalues were also obtained in \cite{Ku10}. In \cite{Ku13}, Ku and Wong gave a new recurrence formula for the eigenvalues of $\Cay(S_n, \mathscr  {D}_n)$ and proved a conjecture of Ku and Wales \cite[Conjecture 1.1]{Ku10} which give lower and upper bounds for the absolute values of these eigenvalues. In \cite{Deng11}, a lower bound on the connectivity of $\Cay(S_n, \mathscr  {D}_n)$ and lower and upper bounds on the Cheeger constant (see Section \ref{sec:Ram}) of $\Cay(S_n, \mathscr  {D}_n)$ were obtained.

Let $n$ and $k$ be integers with $1 \le k \le n-1$. Let $\mathscr{D}_n^{(k)}$ be the set of permutations in $S_n$ without any $i$-cycle in their cyclic decomposition for every $i$ from $1$ to $k$. In \cite{KuLW16}, the Cayley graph $\Cay(S_n, \mathscr {D}_n^{(k)})$ was studied. Obviously, this is a spanning subgraph of $\Cay(S_n, \mathscr {D}_n)$, and in particular $\Cay(S_n, \mathscr {D}_n^{(1)}) = \Cay(S_n, \mathscr {D}_n)$. The following result is a generalization of part (a) in Theorem \ref{thm:RentDeng}.

\begin{thm}
\label{thm:KuLW16}
\emph{(\cite[Theorem 1.3]{KuLW16})}
Let $n \ge 2$ and $k$ be integers such that $1 \le k \le n^{\delta}$ with $0 < \delta < \frac{2}{3}$. Then for sufficiently large $n$, the smallest eigenvalue of $\Cay(S_n, \mathscr {D}_n^{(k)})$ is equal to $\eta_{(n-1, 1)} = -\frac{1}{n-1}|\mathscr{D}_n^{(k)}|$, and moreover $(n-1, 1)$ is the only partition of $n$ associated with the smallest eigenvalue of $\Cay(S_n, \mathscr {D}_n^{(k)})$.
\end{thm}

In \cite{KuLW16}, it was also proved that under the same condition as in Theorem \ref{thm:KuLW16} the set of maximum independent sets of $\Cay(S_n, \mathscr {D}_n^{(k)})$ coincides the set of maximum independent sets of the derangement graph $\Cay(S_n, \mathscr {D}_n)$, and all such maximum independent sets have size $(n-1)!$.

Given an integer $k$ with $0 \le k \le n-1$, let $S(n, k)$ be the set of elements $\sigma$ of $S_n$ such that $\sigma$ fixes exactly $k$ points in $[n]$. The \emph{$k$-point-fixing graph} is defined \cite{KuLW15} to be the Cayley graph $\Cay(S_n, S(n, k))$, that is, two vertices $\s, \t$ are adjacent if and only if $\s \t^{-1}$ fixes exactly $k$ points. Clearly, the $0$-point-fixing graph is the derangement graph. So $k$-point-fixing graphs can be regarded as a generalization of derangement graphs. Since $S(n, k)$ is closed under conjugation, Theorem \ref{nomalSn} implies that all $k$-point-fixing graphs are integral. In \cite{KuLW15}, Ku, Lau and Wong obtained a recursive formula for the eigenvalues of $\Cay(S_n, S(n, k))$, and using this formula they determined the signs of the eigenvalues of the $1$-point-fixing graph $\Cay(S_n, S(n, 1))$. In \cite{KuLW16a}, they obtained the smallest eigenvalue of $\Cay(S_n, S(n, 1))$.

\begin{thm}
\emph{(\cite[Theorem 1.3]{KuLW16a})}
\label{thm:KuLW16a}
For $n \ge 7$, the smallest eigenvalue of $\Cay(S_n, S(n, 1))$ is given by
$$
\eta_{(n-2,2)} = -\frac{1}{n-3}\left(|\mathscr  {D}_{n-1}| + (-1)^n (n-2)\right).
$$
Moreover, $(n-2,2)$ is the only partition of $n$ associated with the smallest eigenvalue of $\Cay(S_n, S(n, 1))$.
\end{thm}

In general, Ku, Lau and Wong obtained the following result in \cite{KuLW17}, where the existence of $t_{0}$ was proved as well.

\begin{thm}
\emph{(\cite[Theorem 1.4]{KuLW17})}
\label{thm:KuLW17}
Let $n$ and $k$ be integers with $0 \le k < n$. Let $t_{0} = t_{0}(k)$ be the smallest positive integer such that $\sum_{i=0}^{t_{0}}{\choose{k}{i}}\frac{(-1)^{t_{0}-i}}{(t_{0}-i)!} < 0$. Then there exists a positive integer $n_{0} = n_{0}(t_0, k)$ such that for all $n \ge n_{0}$, $(n-t_{0},t_{0})$ is the only partition of $n$ associated with the smallest eigenvalue of $\Cay(S_n, S(n, k))$.
\end{thm}

As an application of this result, an upper bound on the independence number of $\Cay(S_n, S(n, k))$ was given in \cite{KuLW17}.

For $0 \le k \le n-1$ and $2 \le j \le n-1$, let $S^{(j)}(n, k)$ be the set of elements of $S_n$ fixing exactly $k$ points and containing no $i$-cycle in their cyclic decomposition for every $i$ from $2$ to $j$. Then we have a sequence of normal Cayley graphs
$$
\Cay(S_n, S^{(n-1)}(n, k)), \Cay(S_n, S^{(n-2)}(n, k)), \ldots, \Cay(S_n, S^{(2)}(n, k)), \Cay(S_n, S(n, k)),
$$
where every graph is a subgraph of the ones after it. Note that $\Cay(S_n, S^{(j)}(n, 0)) = \Cay(S_n, \mathscr {D}_n^{(j)})$ for $2 \le j \le n-1$. Similar to Theorem \ref{thm:KuLW16}, Ku, Lau and Wong proved the following result in \cite{KuLW18}.

\begin{thm}
\label{thm:KuLW18}
\emph{(\cite[Theorem 1.4]{KuLW18})}
Let $n$ and $j$ be integers such that $2 \le j \le n^{\delta}$ with $0 < \delta < \frac{1}{3}$. Then for sufficiently large $n$, the smallest eigenvalue of $\Cay(S_n, S^{(j)}(n, 1))$ is equal to $\eta_{(n-2, 2)} = -\frac{2}{n-3}|\mathscr{D}_{n-1}^{(j)}|$, and moreover $(n-2, 2)$ is the only partition of $n$ associated with the smallest eigenvalue of $\Cay(S_n, S^{(j)}(n, 1))$.
\end{thm}

A recurrence formula for the eigenvalues of $\Cay(S_n, S^{(j)}(n, k))$ and the exact values of $\eta_{\l}(\Cay(S_n, S^{(j)}(n, 1)))$ for eight special partitions of $n$ were also obtained in \cite{KuLW18}.


The \emph{complete transposition graph} (also known as the transposition network) $T_n$ is the Cayley graph on $S_n$ with connection set consisting of all transpositions in $S_n$. With motivation from computing the bisection width of $T_n$, the following result was proved in \cite{KalpakisYesha97} using Theorem \ref{Conjeign}. The integrality of $T_n$ was recently rediscovered in \cite[Theorem 2]{KonstantinovaL20} and \cite[Corollary 2]{GuoLMR19} using different approaches.

\begin{thm}
\label{thm:KalpakisYesha97}
\emph{(\cite[Lemma 3]{KalpakisYesha97})}
Let $n \ge 2$ be an integer. Then $T_n$ is an integral Cayley graph, and the algebraic and geometric
multiplicity of each eigenvalue of $T_n$ are equal. The largest eigenvalue of $T_n$ is $n(n-1)/2$ with multiplicity 1; the second largest eigenvalue of $T_n$ is $n(n-3)/2$ with multiplicity $(n-1)^2$; and for $1 \le k \le n$, $n(n-2k+1)/2$ is an eigenvalue of $T_n$ with multiplicity at least $n!/(n(n-k)! (k-1)!)$.
\end{thm}

We finish this section by mentioning the following two families of integral Cayley graphs of alternating groups.

\begin{thm}
\emph{(\cite[Theorem 3]{KonstantinovaL20})}
Let $n \ge 4$ and $1 \le k \le n$ be integers. Let $T_k$ be the set of all $3$-cycles of the form $(k, i, j)$. Then $\Cay(A_n, T_k)$ is integral. Moreover, $\{i^2 -n + 1: i = 0, 1, \ldots, n-1\}$ is the set of all distinct eigenvalues of $\Cay(A_n, T_k)$.
\end{thm}

\begin{thm}
\emph{(\cite[Corollary 5]{GuoLMR19})}
Let $S = \{(1, 2, i), (1, 2, i)^{-1}: 3 \le i \le n\}$, where $n \ge 3$. Then $\Cay(A_n, S)$ is integral.
\end{thm}

This answers affirmatively \cite[Question 19.50(b)]{Kourovka} which according to \cite{GuoLMR19} was also answered independently by M. Muzychuk.


\subsection{Integral Cayley graphs on other families of non-abelian groups}

\subsubsection{Dihedral, semi-dihedral and generalized dihedral groups}

In this section $D_{2n} = \langle a,b\left|a^n=b^2=1, b^{-1}ab=a^{-1}\right \rangle$ is the dihedral group of order $2n \ge 4$ and $S$ is a subset of $D_{2n}\setminus \{1\}$ with $S^{-1}=S$. The Cayley graph $\Cay(D_{2n},S)$ is called a \emph{dihedrant} in the literature, and its eigenvalues can be computed using Theorem \ref{BabaiS2} and the character table of $D_{2n}$ (see \cite[Section 7.5]{Krebs11} and \cite[Corollary 2.7]{HuangHL17}). Recently, Gao and Luo \cite{Gao10} gave a simpler way to compute the spectra of $\Cay(D_{2n},S)$ using a more general result (see Theorem \ref{thm:bi-Cay}) on eigenvalues of semi-Cayley (bi-Cayley) graphs on abelian groups. We present their result in the theorem below but defer our discussion on bi-Cayley graphs to Section \ref{subsec:eignCay}.

\begin{thm}\label{Dneigen1}
\emph{(\cite[Theorem 4.3]{Gao10})}
The following hold:
\begin{itemize}
\item[\rm (a)] if $S\cap b \la a \ra =\emptyset$, then the eigenvalues of $\Cay(D_{2n},S)$ are
$$
\sum\nolimits_{a^i\in S}\omega_n^{ir},\, r = 0, 1, \ldots, n-1,
$$
each with multiplicity $2$;
\item[\rm (b)] if $S\cap b \la a \ra \ne\emptyset$, say $a^{i_0}b \in S$, then the eigenvalues of  $\Cay(D_{2n},S)$ are
$$
\sum\nolimits_{a^i \in S}\omega_n^{ir} \pm \left|\sum\nolimits_{a^{i}b\in S}\omega_n^{(i_0-i)r}\right|,\,  r = 0, 1, \ldots, n-1.
$$
\end{itemize}
\end{thm}

Since $\left|\sum_{a^{i}b\in S}\omega_n^{(i_0-i)r}\right|
  = \left|\omega_n^{i_0r}\right|\left|\sum_{a^{i}b\in S}\omega_n^{-ir}\right|
   =\left|\sum_{a^{i}b\in S}\omega_n^{-ir}\right|$, Theorem \ref{Dneigen1} implies the following characterization of integral dihedrants.

\begin{thm}
\emph{(\cite[Theorem 3.5.1]{Ahmdy13})}
\label{DihInte1}
A dihedrant $\Cay(D_{2n},S)$ is integral if and only if $\sum\nolimits_{a^i \in S}\omega_n^{ir}$ and $\left|\sum\nolimits_{a^ib\in S}\omega_n^{-ir}\right|$ are all integers for $0 \le r \le n-1$.
\end{thm}

We have $\sum\nolimits_{j\in S_n(d)}\omega_n^{jr}\in \mathbb{Z}$ for $0 \le r \le n-1$
(see \cite[Theorem 4.4]{So06}).
It can be easily shown that, for $A \subseteq \mathbb{Z}_n$ and $0 \le r \le n-1$, if $\sum\nolimits_{j\in A}\omega_n^{jr}\in\mathbb{Z}$, then $\sum\nolimits_{j\in A}\omega_n^{-jr}=\sum\nolimits_{j\in A}\omega_n^{jr}$. This together with Theorem \ref{DihInte1} implies the following result, where $S_n(d)$ is as in \eqref{eq:gnd}.

\begin{thm}
\emph{(\cite{LiuZhou})}
\label{IntegDH11}
The following hold:
\begin{itemize}
  \item[\rm (a)] $\Cay(D_{2n},S)$ is integral provided that each of $\{i: a^i\in S\}$ and $\{i: a^{i}b \in S\}$ is a union of $S_n(d)$ for some divisors $d$ with $d<n$;
  \item[\rm (b)] $\Cay(D_{2n},S)$ is integral provided that $\{i: a^i\in S\}$ is a union of $S_n(d)$ for some divisors $d$ with $d<n$ and $|\{i: a^{i}b \in S\}|=1$.
\end{itemize}
\end{thm}

Since $S_n(d)=d S_{n/d}(1)$ for any divisor $d$ of $n$, in the special case when $\{i: a^i\in S\} = \{i: a^{i}b \in S\}=S_n(1)\cup S_n(d)$, part (a) of Theorem \ref{IntegDH11} implies the following result.

\begin{thm}
\emph{(\cite[Theorem 1.2]{Abdollahi09})}
If $n \ge 3$ is odd and
$$
S=\{a^i: i\in S_n(1)\}\cup\{a^{di}: i \in S_{n/d}(1)\}\cup\{ba^i: i\in S_n(1)\}\cup\{ba^{di}: i \in S_{n/d}(1)\}
$$
for some proper divisor $d$ of $n$, then the dihedrant $\Cay(D_{2n},S)$ is integral.
\end{thm}

Recently, Lu, Huang and Huang \cite{LuHH2018} obtained several results about integral dihedrants. By the character table of dihedral groups, for each $h$, $0 \le h \le n-1$, the mapping $\chi_{h}: D_{2n} \rightarrow \mathbb{C}$ defined by $\chi_{h}(a^i) = 2 \cos(2hi\pi/n)$ and $\chi_{h}(ba^i) = 0$ for $0 \le i \le n-1$ is a character of $D_{2n}$. Recall from Section \ref{InUniGCD} that $B(\FF(\Ga))$ is the Boolean algebra generated by the family $\FF(\Ga)$ of subgroups of a group $\Ga$.

\begin{thm}\label{DnLHH18}
\emph{(\cite{LuHH2018})}
\label{thm:LuHH2018}
Let $S_1 = S \cap \la a \ra$ and $S_2 = S \cap b \la a \ra$. Then the following hold:
\begin{itemize}
\item[\rm (a)] $\Cay(D_{2n},S)$ is integral if and only if both $\chi_{h}(S_1)$ and $\chi_{h}(S_1^2) + \chi_{h}(S_2^2)$ are integers, and $2(\chi_{h}(S_1^2) + \chi_{h}(S_2^2)) - (\chi_{h}(S_1^2))^2$ is a square, for $1 \le h \le \lfloor (n-1)/2 \rfloor$;
\item[\rm (b)] $\Cay(D_{2n},S)$ is integral if and only if $S_1 \in B(\FF(\la a \ra))$ and $2 \chi_{h}(S_2^2)$ is a square, for $1 \le h \le \lfloor (n-1)/2 \rfloor$;
\item[\rm (c)] if $S_1 \in B(\FF(\la a \ra))$ and $b S_2 \in B(\FF(\la a \ra))$, then $\Cay(D_{2n},S)$ is integral;
\item[\rm (d)] if $n=p$ is an odd prime, then $\Cay(D_{2p},S)$ is integral if and only if $S_1 \in B(\FF(\la a \ra))$ and $S_2 = b \la a \ra$, $\{ba^i\}$ or $b \la a \ra \setminus \{ba^i\}$ for some $0 \le i \le p-1$.
\end{itemize}
\end{thm}

The statement in part (d) of Theorem \ref{DnLHH18} was also proved in \cite[Theorem 3.5.5]{Ahmdy13} independently.

The following characterization of integral normal dihedrants was obtained by Cao and Feng in \cite{CaoF19} during their study of perfect state transfer in dihedrants (see Section \ref{subsec:PSTdih} for details).

\begin{thm}
\emph{(\cite[Theorem 4.1]{CaoF19})}
A normal dihedrant $\Cay(D_{2n},S)$ is integral if and only if $S_1^j = S_1$ for all $j \in \ZZZ_{n}^{\times}$, where $S_1 = S \cap \la a \ra$ and $S_1^j = \{a^{ij}: a^i \in S_1\}$.
\end{thm}

The \emph{semi-dihedral group} of order $8n$, where $n \ge 2$, is the group
\begin{equation}
\label{eq:SD}
SD_{8n} = \la a, b \mid a^{4n} = b^2 = 1, bab = a^{2n-1}\ra = \la a \ra \cup b \la a \ra.
\end{equation}

\begin{thm}
\emph{(\cite[Propositions 4.1 and 4.2]{LuoCWW21})}
Let $SD_{8n}$ be the semi-dihedral group as shown in \eqref{eq:SD}, and let $\Cay(SD_{8n}, S)$ be a connected normal Cayley graph on $SD_{8n}$. Then the following hold:
\begin{itemize}
\item[\rm (a)] if $n$ is odd, then $\Cay(SD_{8n}, S)$ is integral if and only if $\{x^i: x \in S \cap \la a \ra\} = S \cap \la a \ra$ for every $i \in \mathbb{Z}^{\times}_{4n}$;
\item[\rm (b)] if $n$ is even, then $\Cay(SD_{8n}, S)$ is integral if and only if $\{x^i: x \in S\} = S$ for every $i \in \mathbb{Z}^{\times}_{4n}$.
\end{itemize}
\end{thm}

Let $\Si$ be a finite abelian group. The \emph{generalized dihedral group} of $\Si$, written $\Dih(\Si)$, is the semidirect product of $\Si$ by $\ZZZ_2$, with $\ZZZ_2$ acting on $\Si$ by inverting elements. Equivalently, one can define $\Dih(\Si) = \la \Si, b \mid b^2 = 1, b^{-1} a b = a^{-1}, a \in \Si \ra$. In \cite{HuangL21}, Huang and Li obtained a necessary and sufficient condition for a Cayley graph on $\Dih(\Si)$ to be integral. Recall from Section \ref{InUniGCD} that, for a group $\Ga$, $B(\FF(\Ga))$ is the Boolean algebra generated by the family $\FF(\Ga)$ of subgroups of $\Ga$, and $\chi(A) = \sum_{a \in A} \chi(a)$ for $A \subseteq \Ga$ and $\chi \in \widehat{\Ga}$.

\begin{thm}
\emph{(\cite[Theorem 3.1]{HuangL21})}
\label{thm:HuangL21}
Let $\Si$ be a finite abelian group, and let $\Ga = \Dih(\Si)$. Let $S = S_1 \cup bS_2$ be an inverse-closed subset of $\Ga \setminus \{1\}$, where $S_1, S_2$ are subsets of $\Si$. Then $\Cay(\Ga, S)$ is integral if and only if $S_1 \in B(\FF(\Si))$ and $\chi(S_{2} S_{2}^{-1})$ is a square number for all irreducible representations $\chi$ of $\Si$ except those satisfying $\chi(a^2) = 1$ for all $a \in \Si$.
\end{thm}

As consequences the following special cases were noted in \cite[Corollaries 3.4, 3.5 and 3.6]{HuangL21}: (i) if $|S_2| = 1$, then $\Cay(\Ga, S)$ is integral if and only if $S_1 \in B(\FF(\Si))$; (ii) if both $S_1$ and $S_2$ are inverse-closed, then $\Cay(\Ga, S)$ is integral if and only if $S_1, S_2 \in B(\FF(\Si))$; (iii) if $\Si$ is of odd order $n$, $S_1 \in B(\FF(\Si))$, and $S_2$ is an $(n, |S_2|, \lambda)$-difference set of $\Si$ with $|S_2|-\lambda$ a square number, then $\Cay(\Ga, S)$ is integral.

In the special case when $\Si = \ZZZ_n$, $\Dih(\Si)$ is the dihedral group $D_{2n}$ and Theorem \ref{thm:HuangL21} gives rise to part (a) of Theorem \ref{thm:LuHH2018}, a result proved by Lu, Huang and Huang in \cite{LuHH2018}.

\subsubsection{Cayley graphs on dicyclic groups}

In this section $$\Dic_{n}=\langle a, b \mid a^{2n}=1, b^2=a^n, b^{-1}ab=a^{-1}\rangle$$ is the \emph{dicyclic group} of order $4n \ge 8$.

\begin{thm}
\emph{(\cite[Theorem 4.6]{Gao10})}
\label{CayT4n}
Let $S$ be a subset of $\Dic_{n} \setminus \{1\}$ with $S^{-1}=S$.
\begin{itemize}
\item[\rm (a)] If $S$ contains no element of $\{a^ib: 0\le i\le 2n-1\}$, then $\Cay(\Dic_{n},S)$ has eigenvalues $\sum\nolimits_{a^i\in S}\omega_{n}^{ir}$ each with multiplicity $2$, for $0 \le r \le 2n-1$.
\item[\rm (b)] If $S$ contains an element of $\{a^ib: 0\le i\le 2n-1\}$, say $a^{i_0}b \in S$, then $\Cay(\Dic_{n},S)$ has eigenvalues $\sum\nolimits_{a^i \in S}\omega_{n}^{ir} \pm \left|\sum\nolimits_{a^ib \in S}\omega_{n}^{(i_0-i)r}\right|$, for $0 \le r \le 2n-1$.
\end{itemize}
\end{thm}

Since $\left|\sum_{a^{i}b\in S}\omega_n^{(i_0-i)r}\right| =\left|\sum_{a^{i}b\in S}\omega_n^{-ir}\right|$, this implies the following characterization of integral Cayley graphs on $\Dic_{n}$.

\begin{thm}\label{Dicyclic4n}
Let $S$ be a subset of $\Dic_{n} \setminus \{1\}$ with $S^{-1}=S$. Then $\Cay(\Dic_{n},S)$ is integral if and only if $\sum\nolimits_{a^i \in S}\omega_n^{ir}$ and $\left|\sum\nolimits_{a^ib\in S}\omega_n^{-ir}\right|$ are all integers for $0 \le r \le 2n-1$.
\end{thm}

Theorem \ref{CayT4n} also implies the following result.

\begin{thm}
\emph{(\cite[Theorem 1.3]{Abdollahi09})}
If $n \ge 3$ is odd and
$$
S=\{a^{k}: 1\le k \le 2n-1, k\ne n\}\cup\{ab, a^{n+1}b\},
$$
then $\Cay(\Dic_{n},S)$ is integral and its spectrum is
\begin{eqnarray*}
\left(\begin{array}{cccc}
2n              & 2n-4          &  0                  &  -4    \\
1                 &   1              &  3n-1             & n-1
\end{array}
\right).
\end{eqnarray*}
\end{thm}

Theorem \ref{Dicyclic4n} implies the following result.

\begin{thm}
Let $S$ be a subset of $\Dic_{n} \setminus \{1\}$ with $S^{-1}=S$. The following hold:
\begin{itemize}
  \item[\rm (a)] $\Cay(\Dic_{n},S)$ is integral provided that each of $\{i: a^i\in S\}$ and $\{i: a^{i}b \in S\}$ is a union of $S_n(d)$ for some divisors $d$ with $d<n$;
  \item[\rm (b)] $\Cay(\Dic_{n},S)$ is integral provided that $\{i: a^i\in S\}$ is a union of $S_n(d)$ for some divisors $d$ with $d<n$ and $|\{i: a^{i}b \in S\}|=1$.
\end{itemize}
\end{thm}

Recently, integral Cayley graphs on $\Dic_{n}$ were further studied by Cheng, Feng and Huang in \cite{ChengFH19}. All characters of $\Dic_{n}$ are known (see, for example, \cite{James01}), and among them the following irreducible characters $\psi_j$, $1 \le j \le n-1$ are most relevant to the integrality of Cayley graphs on $\Dic_{n}$: $\psi_j(1) = 2$, $\psi_j(a^n) = 2(-1)^j$, $\psi_j(a^r) = \om_{n}^{jr/2} + \om_{n}^{-jr/2}$ ($1 \le r \le n-1$), $\psi_j(b) = 0$, $\psi_j(ab) = 0$. Set $\psi_j(X^2) = \sum_{x_1, x_2 \in X} \psi_j(x_1 x_2)$ for any subset $X$ of $\Dic_{n}$.

\begin{thm}
\label{DicCFH19}
\emph{(\cite[Theorem 3.2]{ChengFH19})}
Let $S = S_1 \cup S_2$ be a subset of $\Dic_{n} \setminus \{1\}$ such that $S^{-1}=S$, where $S_1 \subseteq \la a \ra \setminus \{1\}$ and $S_2 \subseteq b \la a \ra$. Then $\Cay(\Dic_{n},S)$ is integral if and only if the following hold for $1 \le j \le n-1$:
\begin{itemize}
  \item[\rm (a)] $\psi_j(S_1)$ and $\psi_j(S_1^2) + \psi_j(S_2^2)$ are integers;
  \item[\rm (b)] $2\left(\psi_j(S_1^2) + \psi_j(S_2^2)\right) - \psi_j(S_1)^2$ is a square number.
\end{itemize}
\end{thm}

\begin{thm}
\label{DicCFH19a}
\emph{(\cite[Theorem 4.3]{ChengFH19})}
Let $S = S_1 \cup S_2$ be a subset of $\Dic_{n} \setminus \{1\}$ such that $S^{-1}=S$, where $S_1 \subseteq \la a \ra \setminus \{1\}$ and $S_2 \subseteq b \la a \ra$. Then $\Cay(\Dic_{n},S)$ is integral if and only if $S_1$ belongs to the Boolean algebra $B(\FF(\la a \ra))$ and $2\psi_j(S_2^2)$ is a square number for all $1 \le j \le n-1$.
\end{thm}

In particular, by Theorem \ref{DicCFH19a}, $\Cay(\Dic_{n},S)$ is integral if both $S_1$ and $b S_2$ belong to the Boolean algebra $B(\FF(\la a \ra))$ (see \cite[Corollary 4.5]{ChengFH19}). Theorem \ref{DicCFH19a} also implies a necessary condition involving the integral cone over $B(\FF(\la a \ra))$ for a Cayley on $\Dic_{n}$ to be integral (see \cite[Corollary 4.6]{ChengFH19}). With the help of this condition the following result was proved in \cite{ChengFH19}.

\begin{thm}
\label{DicCFH19b}
\emph{(\cite[Theorem 4.3]{ChengFH19})}
Let $p$ be an odd prime and $S = S_1 \cup S_2$ a subset of $\Dic_{p} \setminus \{1\}$ such that $S^{-1}=S$, where $S_1 \subseteq \la a \ra \setminus \{1\}$ and $S_2 \subseteq b \la a \ra$. Then $\Cay(\Dic_{n},S)$ is integral if and only if $S_1$ belongs to $B(\FF(\la a \ra))$ and $S_2$ is equal to $b \la a \ra \setminus \{ba^k, ba^{p+k}\}$, $b \la a \ra$ or $\{ba^k, ba^{p+k}\}$ for some $0 \le k \le p-1$.
\end{thm}

\subsubsection{Cayley graphs on a family of groups with order $6n$}
\label{subsec:order6n}

Let $U_{6n}=\langle a,b~|~ a^{2n}=b^3=1, a^{-1}ba=b^{-1}\rangle$, where $n \ge 1$. Using the character table of this group (see \cite[p.187]{James01}) and Theorem \ref{BabaiS2}, one can compute the spectrum of $\Cay(U_{6n},S)$ for any subset $S$ of $U_{6n}\setminus \{1\}$ with $S^{-1}=S$. This enabled Abdollahi and Vatandoost \cite{Abdollahi09} to prove the following result.

\begin{thm}
\label{thm:U6n}
\emph{(\cite[Theorem 1.4]{Abdollahi09})}
If $n \ge 3$ is odd and
$$
S=\{a^{2k}b:1\le k \le n-1\}\cup\{a^{2k}b^2: 1\le k \le n-1\}\cup\{a^{2k+1}b: 0\le k \le n-1\},
$$
then $\Cay(U_{6n},S)$ is integral and its spectrum is
\begin{eqnarray*}
\left(\begin{array}{ccccc}
3n-2         & n-2          &  1                  &  -2          & 1-2n              \\
1                 &   1           &  4n-2           &  2n-2      &2
\end{array}
\right).
\end{eqnarray*}
\end{thm}


\subsection{Integral Cayley graphs of small degrees}

The following result due to Abdollahi and Vatandoost gives a characterization of connected cubic integral Cayley graphs.

 \begin{thm}
 \emph{(\cite[Theorem 1.1]{Abdollahi09})}
There are exactly seven connected cubic integral Cayley graphs. In particular, a connected cubic Cayley graph $\Cay(\Gamma,S)$ is integral if and only if $\Gamma$ is isomorphic to one the following groups:
$$
\ZZZ_2^2, \ZZZ_4, \ZZZ_6, S_3, \ZZZ_2^3, \ZZZ_2 \times \ZZZ_4, D_8, \ZZZ_2 \times \ZZZ_6, D_{12}, A_4, S_4, D_8\times \ZZZ_3, D_6\times \ZZZ_4, A_4\times \ZZZ_2.
$$
\end{thm}

The following theorem also due to Abdollahi and Vatandoost determines the orders of connected $4$-regular integral Cayley graphs on finite abelian groups.

\begin{thm}
\emph{(\cite[Theorem 1.1]{Abdollahi11})}
The order of any connected $4$-regular integral Cayley graph on a finite abelian group must be one of the following:
$$
5, 6, 8, 9, 10, 12, 16, 18, 20, 24, 25, 32, 36, 40, 48, 50, 60, 64, 72, 80, 96, 100, 120, 144.
$$
\end{thm}

As a side result it was also proved in \cite[Theorem 1.2]{Abdollahi11} that there are precisely 27 connected integral Cayley graphs of order up to $11$.

In 2015, {Minchenko} and {Wanless} \cite{Minchenko15} proved that up to isomorphism there are precisely $32$ connected $4$-regular integral Cayley graphs, $17$ of which are bipartite. They also proved that there are exactly 27 arc-transitive $4$-regular integral graphs of degree $4$, 16 of which are bipartite. See \cite[Tables 3, 5 and 7]{Minchenko15} for the groups and connection sets for these graphs. From their results they noticed that integral Cayley graphs can be cospectral to integral non-Cayley graphs, and integral arc-transitive graphs can be cospectral to integral non-arc-transitive graphs.

Recently, Ghasemi \cite{Ghasemi17} determined all possible orders of $5$-regular integral Cayley graphs on abelian groups other than cyclic groups.

\begin{thm}
\emph{(\cite[Theorem 3.3]{Ghasemi17})}
Let $\Gamma$ be a finite abelian group which is not cyclic, and let $S\subseteq \Gamma \setminus \{1\}$ be such that $|S| = 5$, $S = S^{-1}$ and $\Gamma = \langle S\rangle$.  If $\Cay(\Gamma, S)$ is integral, then the order of $\Gamma$ must be one of the following:
$$
8, 16, 18, 24, 32, 36, 40, 48, 50, 64, 72, 80, 96, 100, 120, 128, 144, 160, 192, 200, 240,
288.
$$
\end{thm}


\subsection{Automorphism groups of integral Cayley graphs}

\subsubsection{Cayley integral groups}
\label{sec:GCInt}

As seen in Corollary \ref{cubelike-int}, every Cayley graph on $\mathbb{Z}^n_2$ is integral. In general, a group $\Gamma$ is called \emph{Cayley integral} \cite{Klotz10} if every Cayley graph on $\Gamma$ is integral, or equivalently any graph admitting $\Ga$ as a regular group of automorphisms is integral. In \cite[Theorem 13]{Klotz10}, Klotz and Sander gave a classification of all finite abelian Cayley integral groups, which can be regarded as a generalization of Corollary \ref{cubelike-int}.

\begin{thm}
\label{thm:abelianCI}
\emph{(\cite[Theorem 13]{Klotz10})}
All nontrivial abelian Cayley integral groups are represented by
$$
\mathbb{Z}^n_2,~\mathbb{Z}^n_3,~\mathbb{Z}^n_4,~\mathbb{Z}^m_2 \oplus \mathbb{Z}^n_3,~\mathbb{Z}^m_2\oplus \mathbb{Z}^n_4,\quad m\ge1,~ n\ge1.
$$
\end{thm}

In \cite[Section 5]{Klotz10}, Klotz and Sander found three non-abelian Cayley integral groups, namely $S_3$, $Q_8$ and $\mathbb{Z}_3 \rtimes \mathbb{Z}_4$, and in \cite[Problem 3]{Klotz10} they posed the problem of determining all non-abelian Cayley integral groups. This problem was settled by Abdollahi and Jazaeri \cite{AbdollahiJ14} and independently by Ahmady, Bell and Mohar \cite{AhmadyBM14}.

\begin{thm}
\label{thm:nonabelianCI}
\emph{(\cite[Theorem 1.1]{AbdollahiJ14} and \cite[Theorem 4.2]{AhmadyBM14})}
A finite non-abelian group is Cayley integral if and only if it is isomorphic to one of the
following groups:
\begin{itemize}
  \item[\rm (a)] the symmetric group $S_3$ of degree $3$;
  \item[\rm (b)] the nontrivial semidirect product $\mathbb{Z}_3 \rtimes \mathbb{Z}_4 = \langle x, y ~| ~x^3 = y^4 = 1, y^{-1}xy = x^{-1}\rangle$;
  \item[\rm (c)] $Q_8 \times  \mathbb{Z}^n_2$ for some integer $n \ge 0$, where $Q_8$ is the quaternion group.
\end{itemize}
\end{thm}

Theorems \ref{thm:abelianCI} and \ref{thm:nonabelianCI} together give a complete classification of Cayley integral groups.

Let $\Ga$ be an abelian group of order at least $3$ having a unique involution $t$. The group $\Dic(\Ga) = \la \Ga, x\ra$, where $x^2 = t$ and $x^{-1}ax = a^{-1}$ for every $a \in \Ga$, is known as a \emph{generalized dicyclic group}. In the case when $\Ga \cong \ZZZ_n$ it is the dicyclic group of order $2n$, and when $\Ga \cong \ZZZ_{2^n}$ it is known as the \emph{generalized quaternion group} $Q_{2^{n+1}}$ of order $2^{n+1}$.

Let $k \ge 1$ be an integer. In \cite{EstelyiK14}, Est\'{e}lyi and Kov\'{a}cs studied the family $\mathcal{G}_{k}$ of finite groups $\Ga$ such that every Cayley graph on $\Ga$ with degree at most $k$ is integral. Clearly, $\mathcal{G}_{1}$ is the family of all finite groups, and $\mathcal{G}_{2}$ consists of those groups whose non-identity elements are of order $2, 3, 4$ or $6$ and contain no subgroup isomorphic to $D_{2n}$ for any $n \ge 4$.

\begin{thm}
\emph{(\cite[Theorem 3]{EstelyiK14})}
\label{thm:EstelyiK14}
If $k \ge 6$, then $\mathcal{G}_{k}$ consists of all Cayley integral groups. Moreover, $\mathcal{G}_{4}$ and $\mathcal{G}_{5}$ are equal and consist of the following groups:
\begin{itemize}
\item[\rm (a)] the Cayley integral groups;
\item[\rm (b)] the generalized dicyclic groups $\Dic(\ZZZ_{3}^n \times \ZZZ_6)$, where $n \ge 1$.
\end{itemize}
\end{thm}

In \cite{EstelyiK14}, Est\'{e}lyi and Kov\'{a}cs also gave a characterization of non-abelian $2$-groups in $\mathcal{G}_{3}$. As usual we use $[a, b] = a^{-1}b^{-1}ab$ to denote the commutator of two elements $a, b$ of a group.

\begin{thm}
\emph{(\cite[Proposition 12]{EstelyiK14})}
Let $\Ga$ be a non-abelian 2-group of exponent 4. Then $\Ga \in \mathcal{G}_{3}$ if and only
if every minimal non-abelian subgroup of $\Ga$ is isomorphic to $Q_8$, the metacyclic group $\la a, b\ | \ a^4 = b^4 = 1,\ b^{-1} a b = a^{-1} \ra$, or the non-metacyclic group $\la a, b, c \ | \ a^4 = b^4 = c^2 = 1,\ [a, b] = c,\ [c, a] = [c, b] = 1 \ra$.
\end{thm}

A complete characterization of $\mathcal{G}_{3}$ was recently obtained by Ma and Wang in \cite{MaW16a}. For an integer $k \ge 2$, let $\mathcal{A}_k$ be the family of finite groups $\Ga$ such that every Cayley graph on $\Ga$ with degree exactly $k$ is integral.

\begin{thm}
\emph{(\cite[Theorem 2.6 and Corollary 2.7]{MaW16a})}
The following hold:
\begin{itemize}
\item[\rm (a)] A finite group $\Ga$ belongs to $\mathcal{A}_{3}$ if and only if $\Ga \cong S_3$, or for any involution $x$ and element $y$ of $\Ga$, $\la x, y \ra$ is isomorphic to one of the following groups:
$$
\ZZZ_2,\ \ZZZ_2^2,\ \ZZZ_4,\ \ZZZ_6,\ \ZZZ_2 \times \ZZZ_4,\ \ZZZ_2 \times \ZZZ_6,\ A_4;
$$
\item[\rm (b)] $\mathcal{G}_{3}$ consists of all finite $3$-groups of exponent $3$ and all groups in $\mathcal{A}_{3}$.
\end{itemize}
\end{thm}

A second proof of Theorem \ref{thm:EstelyiK14} was also given in \cite{MaW16a}.

\subsubsection{Cayley integral simple groups}
\label{sec3.6}

A finite group $\Gamma$ is called \emph{Cayley simple} if the only connected integral Cayley graph on $\Gamma$ is the complete graph of order $|\Gamma|$. This definition was first introduced in \cite{Abdollahi11}, where the following question was posed (\cite[Question 2.21]{Abdollahi11}): Is any finite simple group, Cayley simple? A negative answer to this question is implied in the following result.

\begin{thm}\emph{(\cite[Proposition 2.6]{AbdollahiJ13})}\label{CaySim1}
Let $\Gamma$ be a finite group and $S$ an inverse-closed subset of $\Gamma \setminus \{1\}$. Then
$\Gamma \setminus S$ is a subgroup of $\Gamma$ if and only if $\Cay(\Gamma, S)$ is a complete multipartite graph. Moreover, such a complete multipartite Cayley graph is integral, with $|\Gamma \setminus S|$ vertices in each part of the corresponding multipartition.
\end{thm}

Note that complete graphs form a special class of complete multipartite graphs. So in \cite{AbdollahiJ13} Abdollahi and Jazaeri modified the definition of Cayley simple groups and introduced the notion of  \emph{Cayley integral simple groups} (or \emph{CIS-group} for short) by replacing ``complete graph'' by ``complete multipartite graph''. They proved the next two results.

\begin{thm}\emph{(\cite[Theorem 1.3]{AbdollahiJ13})}
Let $\Gamma$ be a finite non-simple group. Then $\Gamma$ is a CIS-group if and
only if $\Gamma\cong\mathbb{Z}_{p^2}$ for some prime $p$ or $\Gamma\cong\mathbb{Z}_2^2$.
\end{thm}

\begin{thm}\emph{(\cite[Theorem 4.10]{AbdollahiJ13})}\label{CIS-g1}
Let $\Gamma$ be a finite group. Let $H$ and $K$ be proper subgroups of $\Ga$ such
that $HK = \Ga$ and $H \cap K =\{1\}$. If there exists an integral Cayley graph $\Cay(H, S)$ such that $S \cup \{1\}$ is not a subgroup of $H$, then $\Gamma$ is not a CIS-group.
\end{thm}

It is known that any cyclic group of prime order is a CIS-group (see \cite[Corollary 2.19]{Abdollahi11}). Since these are the only simple abelian groups, it is natural to ask the following question (\cite[Question 1.4]{AbdollahiJ13}): Which finite non-abelian simple groups are CIS-groups? Regarding this question, Abdollahi and Jazaeri proved the following result using Theorem \ref{CIS-g1} and the fact that $A_{p} = A_{p-1} \la (1,2,\ldots,p) \ra$.

\begin{thm}\emph{(\cite[Corollary 4.12]{AbdollahiJ13})}
For every prime $p \ge 5$, the alternating group $A_p$ is not a CIS-group.
\end{thm}

Thus the alternating groups $A_p$ for primes $p \ge 5$ form an infinite family of finite non-abelian simple groups which are not CIS-groups. It turns out that this is not a coincidence since the following result shows that there does not exist any non-abelian finite CIS-group, answering the question mentioned above.

\begin{thm}
\emph{(\cite[Theorem 3.3]{AhmadyBM14})}\label{CISnon-NAbl1}
Let $\Ga$ be a CIS group. Then $\Ga$ is abelian and in particular is isomorphic to $\ZZZ_p$ or $\ZZZ_{p^2}$ for some prime $p$ or is isomorphic to $\ZZZ_2^2$.
\end{thm}

In particular, this implies that every finite non-abelian group admits a nontrivial integral Cayley graph.

\subsubsection{Automorphism groups of integral circulant graphs}

Let $r$ be a positive integer. A poset $([r], \preceq)$ on $[r] = \{1, 2, \ldots, r\}$ is called \emph{increasing} if $i \preceq j$ implies $i \le j$. Given $([r], \preceq)$ and positive integers $n_1, \ldots, n_r$, the generalized wreath product $\prod_{([r], \preceq)} S_{n_i}$ introduced in \cite{KlinK12} is a certain permutation group acting on $[n_1] \times \cdots \times [n_r]$. The following result was proved in \cite{KlinK12} using the techniques of Schur rings.

\begin{thm}
\emph{(\cite[Theorem 1.1]{KlinK12})}
Let $\Ga$ be a permutation group acting on the cyclic group $\ZZZ_n$, where $n \ge 2$. The following conditions are equivalent:
\begin{itemize}
\item[\rm (a)] $\Ga = \Aut(\Cay(\ZZZ_n, S))$ for some integral circulant graph $\Cay(\ZZZ_n, S)$;
\item[\rm (b)] $\Ga$ is a permutation group which is permuationally isomorphic to a generalized wreath product $\prod_{([r], \preceq)} S_{n_i}$, where $([r], \preceq)$ is an increasing poset and $n_1, \cdots, n_r$ are integers no less than $2$ such that $n = n_1 \ldots n_r$ and $\gcd(n_i, n_j) = 1$ whenever $i \not \preceq j$.
\end{itemize}
\end{thm}

In \cite{BasicI11}, Ba\v{s}i\'c and Ili\'c determined the automorphism groups of certain classes of integral circulant graphs: the unitary Cayley graph $\Cay(\ZZZ_n, \ZZZ_n^{\times})$ and more generally the gcd graphs $\ICG(n, \{d\})$ for proper divisors $d$ of $n$. They also determined the automorphism groups of those gcd graphs $\ICG(n, D)$ of $\ZZZ_n$ for which $D=\{1, p^k\}$ and $n$ is either square-free or a prime power.

\begin{thm}
\emph{(\cite[Theorem 3.4]{BasicI11})}
Let $n = p_1^{\a_1} p_2^{\a_2} \ldots p_k^{\a_k} \ge 2$ be an integer in canonical factorization, and let $m = p_1 p_2 \ldots p_k$. Then
$$
\Aut(\Cay(\ZZZ_n, \ZZZ_n^{\times})) \cong (S_{p_1} \times \cdots \times S_{p_k}) \wr S_{n/m}.
$$
In particular,
$$
|\Aut(\Cay(\ZZZ_n, \ZZZ_n^{\times}))| = p_{1}! p_{2}! \ldots p_{k}! \left(\left(\frac{n}{m}\right)!\right)^{m}.
$$
\end{thm}


\subsection{Cayley digraphs integral over the Gauss field or other number fields}

A circulant digraph is called \emph{Gaussian integral} if all its eigenvalues are algebraic integers of the Gauss field $\mathbb{Q}(i)$ (that is, all eigenvalues are Gaussian integers $a+bi$, where $a, b \in \ZZZ$). In \cite{XuM11}, circulant digraphs which are Gaussian integral were studied, and those with order $n = k$, $2k$ or $4k$ for odd integers $k$ were characterized. Among others the following results were proved in \cite{XuM11}, where $S_{n}(D)$ is as defined in \eqref{eq:SnD} and $D(n)$ is the set of positive divisors of $n$ as defined in \eqref{eq:Dn}.

\begin{thm}
\emph{(\cite[Theorems 2.1 and 2.2]{XuM11})}
Let $n \ge 3$ be an integer not divisible by $4$. Then a circulant digraph $\Cay(\ZZZ_n, S)$ is Gaussian integral if and only if $S = S_{n}(D)$ for some $D \subseteq D(n) \setminus \{n\}$.
\end{thm}

This together with Corollary \ref{cor:intICG} implies that a circulant graph of order not divisible by $4$ is Gaussian integral if and only if it is integral.

In general, given an algebraic number field $K$ (that is, a finite extension of $\mathbb{Q}$), a circulant digraph is called \emph{integral over $K$} if all its eigenvalues are algebraic integers of $K$. Denote $F = K \cap \mathbb{Q}(\om_n)$. Using Galois theory, Li \cite{LiFei13} obtained a characterization of circulant digraphs which are integral over $K$. He considered the orbits of a certain action of the Galois group $\Gal(\mathbb{Q}(\om_n)/F)$ on $S_{n}(d)$. Setting $m = n/d$, there are exactly $r_d = [F \cap \QQQ(\om_m): \QQQ]$ such orbits $M_i(d)$, for $1 \le i \le r_d$, and each of them has size $[\QQQ(\om_m):F \cap \QQQ(\om_m)]$. Since $\{1, 2, \ldots, n-1\} = \cup_{d \in D(n)} S_{n}(d)$, the set $\{1, 2, \ldots, n-1\}$ is partitioned into $r(n,K) = \sum_{d \in D(n)} r_d$ such orbits $M_i(d)$, with $d$ running over $D(n)$ and $1 \le i \le r_d$.

\begin{thm}
\label{thm:LiFei13}
\emph{(\cite[Theorem 1 and Corollary 1]{LiFei13})}
Let $K$ be an algebraic number field. A circulant digraph $\Cay(\ZZZ_n, S)$, $n \ge 3$ is integral over $K$ if and only if $S$ is the union of some (not necessarily all) of the orbits $M_i(d), d \in D(n), 1 \leq i \le r_d$, possibly from different proper divisors $d$ of $n$. Therefore, there are at most $2^{r(n,K)}$ pairwise non-isomorphic circulant digraphs of order $n$ which are integral over $K$.
\end{thm}

In the special case when $n$ is a multiple of $8$ and $K = \mathbb{Q}(i)$ is the Gauss field, Theorem \ref{thm:LiFei13} yields the following corollary which proves a conjecture in \cite[Conjecture 3.3]{XuM11}.

\begin{cor}
\emph{(\cite[Corollary 2]{LiFei13})}
A circulant digraph $\Cay(\ZZZ_n, S)$ with order $n$ a multiple of $8$ is Gaussian integral if and only if $S$ is the union of some (not necessarily all) of the orbits $M_i(d), d \in D(n), 1 \leq i \le r_d$, possibly from different proper divisors $d$ of $n$.
\end{cor}

Consider an abelian group $\Ga = \ZZZ_{n_1} \oplus \cdots \oplus \ZZZ_{n_d}$ of order $n = n_1 \ldots n_d$, where $n_1, \cdots, n_d$ are positive integers. Let $K$ be an algebraic number field. Then $\Gal(\QQQ(\om_n)/K) \subseteq \Gal(\QQQ(\om_n)/\QQQ) $ $\cong \ZZZ_{n}^*$. Thus, for each $\s \in \Gal(\QQQ(\om_n)/K)$, there exists an element $a \in \ZZZ_{n}^*$ such that $\s(\om_n) = \om_n^a$. This element $a$ gives rise to an action of $H$ on each $\ZZZ_{i}^*$ defined by $\pi_{i}(x) = a x$ ($\mod~n_i$) for $x \in \ZZZ_{i}^*$. Hence $\Gal(\QQQ(\om_n)/K)$ acts on $\Ga$ by $\pi(x_1, \ldots, x_d) = (\pi_{1}(x), \ldots, \pi_{d}(x)) = (ax~\mod~n_1, \ldots, ax~\mod~n_d)$ for $(x_1, \ldots, x_d) \in \Ga$. The following result was proved by Li in \cite{LiFei12}.

\begin{thm}
\label{thm:LiFei12}
\emph{(\cite[Theorem 2.1]{LiFei12})}
Let $\Ga = \ZZZ_{n_1} \oplus \cdots \oplus \ZZZ_{n_d}$ be an abelian group and $K$ an algebraic number field, where each $n_i \ge 2$. A Cayley digraph $\Cay(\Ga, S)$ on $\Ga$ is integral over $K$ if and only if $S$ is the union of some orbits under the above-mentioned action of $\Gal(\QQQ(\om_n)/K)$ on $\Ga$.
\end{thm}

In \cite[Proposition 3.2]{LiFei12}, Li also gave an upper bound for the number of Cayley digraphs on $\Ga$ which are integral over $K$.

In the case when $K=\QQQ$, the theorem above gives a necessary and sufficient condition for a Cayley digraph on an abelian group to be integral in the usual sense. The reader is invited to compare Theorem \ref{thm:LiFei12} with Theorem \ref{thm:Alperin12}.


\section{Cospectral Cayley graphs}\label{sec:CoCayG}

Which graphs are determined by their spectrum? When two non-isomorphic graphs are cospectral? These broad questions with a long history have been studied extensively in the literature. See \cite{vanDamH03,vanDamH09} for two surveys on these topics for general graphs. In this section we review some known results on these questions for Cayley graphs.


\subsection{Cospectral Cayley graphs}
\label{subsec:CosCay}

Using Theorem \ref{BabaiS2}, Babai proved the following result.

\begin{thm}
\emph{(\cite[Theorem 5.2]{Babai79})}
\label{thm:babai79}
Given an integer $k \ge 2$ and a prime $p > 64k$, there exist $k$ pairwise non-isomorphic cospectral Cayley graphs on the dihedral group $D_{2p}$.
\end{thm}

In \cite[Theorem 1.3]{AbdollahiJJ16}, Abdollahi, Janbaz and Jazaeri proved that for any prime $p \ge 13$ there exist two non-isomorphic cospectral $6$-regular Cayley graphs on $D_{2p}$. Generalizing the construction in this paper, Abdollahi, Janbaz and Ghahramani \cite{AbdollahiJG17} proved the following two results.

\begin{thm}
\emph{(\cite[Theorem 1.1]{AbdollahiJG17})}
Let $p \ge 23$ be a prime. Then for every integer $k$ between $6$ and $2p-7$ there exist at least two non-isomorphic cospectral $k$-regular Cayley graphs on $D_{2p}$.
\end{thm}

\begin{thm}
\emph{(\cite[Theorem 1.2]{AbdollahiJG17})}
\label{thm:AbdollahiJG17}
Let $p \ge 23$ be a prime. Then for every integer $k$ between $6$ and $p+6$ there exist at least $\choose{\frac{p-1}{2}}{\lfloor k/2 \rfloor - 3}$ pairwise non-isomorphic cospectral $k$-regular ($(2p-k-1)$-regular) Cayley graphs on $D_{2p}$.
\end{thm}

Thus the number of pairwise non-isomorphic cospectral Cayley graphs on $D_{2p}$ is exponential in terms of $p$. Moreover, for any $k \ge 6$, there exist arbitrarily large families of $k$-regular Cayley graphs that are cospectral and pairwise non-isomorphic. Recently, Tang, Cheng, Liu and Feng \cite{TangCLF21} proved the following result in the same spirit for dicyclic groups $\Dic_{4p}$.

\begin{thm}
\emph{(\cite[Theorem 3.7]{TangCLF21})}
Let $p \ge 23$ be a prime and $k$ an integer between $12$ and $2p+13$. If $k \in \{12, 13, 2p+12, 2p+13\}$, then there exist a pair of non-isomorphic cospectral $k$-regular ($(4p-k-1)$-regular) Cayley graphs on $\Dic_{4p}$; if
$14 \le k \le 2p+11$, then there exist $\choose{p}{\lfloor k/2 \rfloor - 6}$ pairwise non-isomorphic cospectral $k$-regular ($(4p-k-1)$-regular) Cayley graphs on $\Dic_{4p}$.
\end{thm}

The Gaussian coefficient ${d \brack i}_q$ gives the number of subspaces of dimension $i$ over $\FFF_q$ in the vector space $\FFF_q^d$. In \cite{LubotzkySV06}, Lubotzky, Samuels and Vishne proved the following result. (Note that in \cite{LubotzkySV06} Cayley graphs are not required to be undirected.)

\begin{thm}
\emph{(\cite[Theorem 1]{LubotzkySV06})}
For every integer $d\ge 5$ ($d \not= 6$), every prime power $q$, and every integer $e \ge 1$ such that $q^e > 4d^2 + 1$, there are two systems $A, B$ of generators of the group $G = \PSL_{d}(\FFF_{q^e})$ such that $\Cay(G, A)$ and $\Cay(G, B)$ are cospectral but not isomorphic. Moreover, the number of generators in $A$ and $B$ can be chosen to be either $k=2(q^d-1)/(q-1)$ or $k = \sum_{i=1}^{d-1} {d \brack i}_q$.
\end{thm}

In particular, this implies that for fixed $d$ and $q$ there are infinitely many cospectral pairs of Cayley graphs which are $k$-regular with the same $k$.

\begin{thm}
\label{CosCayC}
\emph{(\cite[Proposition 7]{LubotzkySV06})}
Let $\Gamma$ be a finite group and $\Gamma'$ a proper subgroup of $\Gamma$.
Suppose that $\Gamma'$ has two inverse-closed generating sets $A'$ and $B'$ of size $k$ such that $\Cay(\Gamma',A')$ and $\Cay(\Gamma',B')$ are cospectral but non-isomorphic. Then $\Gamma$ has two generating sets $A$ and $B$ of size $|\Gamma|-k-1$ such that $\Cay(\Gamma, A)$ and $\Cay(\Gamma, B)$ are cospectral but non-isomorphic.
\end{thm}

Combining this with Theorem \ref{thm:babai79}, the following result was obtained in \cite[Corollary 8]{LubotzkySV06}: For large enough $n$, each of the groups $G = S_n$ and $G = \PSL_{n}(\FFF_{q})$ has two subsets $A$ and $B$ of size $|G| - 3$ such that $\Cay(G, A)$ and $\Cay(G, B)$ are cospectral but non-isomorphic.

A few explicit constructions of cospectral but non-isomorphic Cayley graphs were also given in \cite{LubotzkySV06}.


\subsection{Cospectrality and isomorphism}
\label{subsec:cosp-isom}

Consider a Cayley graph $\Cay(\Gamma,S)$ on a group $\Gamma$ with connection set $S$. It can be verified that every $\sigma \in\Aut(\Gamma)$ induces an isomorphism from $\Cay(\Gamma,S)$ to $\Cay(\Gamma, \sigma(S))$, called a \emph{Cayley isomorphism} (CI). A Cayley graph $\Cay(\Gamma, S)$ is called a \emph{CI-graph} of $\Gamma$ if for any Cayley graph $\Cay(\Gamma, T)$ isomorphic to $\Cay(\Gamma, S)$ there exists $\sigma\in\Aut(\Gamma)$ such that $T = \sigma(S)$. \'{A}d\'{a}m \cite{Adam67,Djokovic70} conjectured that all circulant graphs are CI-graphs of the corresponding cyclic groups. Though disproved in \cite{ElspasT70}, this conjecture inspired a lot of interest on CI-graphs in more than four decades, and a complete solution to the problem of determining when two given circulant graphs or digraphs are isomorphic was given by Muzychuk in \cite{Muzychuk04}. See \cite{Li02} for a survey of CI-graphs.

Motivated by \'{A}d\'{a}m's conjecture, Mans, Pappalardi and Shparlinski \cite{MansPS02} proposed to study the following problem: ``Which circulant graphs $\Cay(\mathbb{Z}_n,S)$ satisfy the \emph{spectral \'{A}d\'{a}m property}, that is, for any $T\subseteq \mathbb{Z}_n$ satisfying $0 \not \in T=-T$, $\Spec(\Cay(\mathbb{Z}_n, S)) = \Spec(\Cay(\mathbb{Z}_n, T))$ implies that there exists an integer $\alpha$ coprime with $n$ such that $T =\alpha S=\{\alpha s: s \in S\}$?'' Obviously, if $T =\alpha S$ for an integer $\alpha$ coprime to $n$, then $\Cay(\mathbb{Z}_n, S)\cong \Cay(\mathbb{Z}_n, T)$. So a circulant graph $\Cay(\mathbb{Z}_n,S)$ satisfies the spectral \'{A}d\'{a}m property if and only if, for any circulant graph $\Cay(\mathbb{Z}_n,T)$, $\Spec(\Cay(\mathbb{Z}_n, S)) = \Spec(\Cay(\mathbb{Z}_n, T))$ implies $\Cay(\mathbb{Z}_n, S)  \cong \Cay(\mathbb{Z}_n, T)$.

A graph $G$ is said to be \emph{determined by its spectrum} (DS for short) if every graph cospectral with $G$ is isomorphic to $G$. In particular, a Cayley graph $\Cay(\Gamma, S)$ is called \emph{Cay-DS} if, for any Cayley graph $\Cay(\Gamma, T)$, $\Spec(\Cay(\Gamma, S)) = \Spec(\Cay(\Gamma, T))$ implies $\Cay(\Gamma, S) \cong \Cay(\Gamma, T)$. If all Cayley graphs on $\Gamma$ are Cay-DS, then the group $\Gamma$ is said to be \emph{Cay-DS}. These concepts can be extended to Cayley digraphs in an obvious way.

It is natural to ask whether all circulant (di)graphs are Cay-DS. The answer is negative, as shown in \cite{ElspasT70, MansPS02}. On the positive side, Turner proved in \cite[Theorem 2]{Turner67} that any cyclic group of prime order is Cay-DS. In \cite{MansPS02}, Mans, Pappalardi and Shparlinski proved the following result for circulants of degree $4$.

\begin{thm}
\emph{(\cite[Theorem 9]{MansPS02})}
Let $n \ge 3$ and let $S=\{\pm a, \pm b\}\subseteq \mathbb{Z}_n$ be such that $\gcd(a,b,n)=1$. Then $\Cay(\mathbb{Z}_n, S)$ is Cay-DS.
\end{thm}

This generalizes an earlier result \cite{LitowM98} which asserts that $\Cay(\mathbb{Z}_n, \{\pm 1, \pm d\})$ is Cay-DS provided that $2 \le d < \min\{n/4, \varphi(n)/2\}$, where $\varphi(n)$ is Euler's function. In \cite{MansPS02}, it was also shown that for any fixed $m$ the probability that a random $m$-element subset $S \subseteq \mathbb{Z}_n$ does not have the spectral \'{A}d\'{a}m property is $O(n^{-1})$.

It was proved in \cite[Corollary 3]{ElspasT70} that the circulant digraph $\Cay(\mathbb{Z}_p, S)$ is Cay-DS for any prime $p$. In \cite{HuangC01}, Huang and Chang proved the following result.

\begin{thm}\label{HCRes}
\emph{(\cite[Theorems 1, 2 and 3]{HuangC01})}
The circulant digraph $\Cay(\mathbb{Z}_n, S)$, $n \ge 3$, is Cay-DS if one of the following conditions holds:
\begin{itemize}
  \item[\rm (a)] $n=p^a$ is a prime power, and $S\subset\mathbb{Z}_n$ does not contain any coset of $\langle p^{a-1}\rangle$ in $\mathbb{Z}_n$, where $\langle p^{a-1}\rangle$ is the subgroup of  $\mathbb{Z}_n$ generated by $p^{a-1}$;
  \item[\rm (b)] $n=p^aq^b$, where $p$ and $q$ are distinct odd primes and $a$ and $b$ are positive integers, and $S\subset \mathbb{Z}_n$ satisfies $\max\{s: s\in S\} \le \varphi(n)$;
  \item[\rm (c)] $n=2^aq^b$, where $q$ is an odd prime and $a$ and $b$ are positive integers, and $S$ is a generating subset of $\mathbb{Z}_n$ that satisfies $\max\{s: s\in S\}\le \varphi(n)$.
\end{itemize}
\end{thm}

In \cite[Section 7]{So06}, So proposed the following conjecture.

\begin{conj}\emph{(\cite[Conjecture 7.3]{So06})}
\label{SoConj}
Two integral circulant graphs are isomorphic if and only if they are cospectral; that is, all integral circulant graphs are Cay-DS.
\end{conj}

So \cite[Section 7]{So06} proved that this conjecture is true if $n$ is a prime power or a product of two distinct primes. This conjecture has also been confirmed \cite[Section 5]{Ilic11} when $n$ is square-free and the set $D$ in the integral circulant graph $\ICG(n,D)$ contains at most two prime factors of $n$, but it is still open in its general form.

Given non-empty sets $A_1, \ldots , A_s$ of integers, define
\begin{equation}
\label{eq:prod-set}
\prod_{i=1}^sA_i:=\{a_1 \cdots a_s: a_i\in A_i, \, 1\le i \le s\}.
\end{equation}
Denote by $\mathbb{P}$ the set of primes. Let $n$ be a positive integer and $p\in \mathbb{P}$. Recall that $e_p(n)$ denotes the exponent of $p$ in $n$. For $\emptyset \ne X \subseteq \NNN$ and $p \in \mathbb{P}$, define
\be
\label{eq:Xp}
X_{p} =  \left\{p^{e_{p}(x)}: x \in  X\right\}.
\ee
If $X = \prod_{p\in \mathbb{P}} X_{p}$, then $X$ is called a \emph{multiplicative set}.

Define the \emph{spectral vector} of an integral circulant graph $\ICG(n,D)$ to be
$$
\overrightarrow{\lambda}(n,D) :=(\lambda_1(n,D),\ldots,\lambda_{n-1}(n,D),\lambda_0(n,D)),
$$
where $\lambda_i(n,D)$ for $i=0,1,\ldots,n-1$ are the eigenvalues of $\ICG(n,D)$ as shown in (\ref{gcdEig1}).
In \cite{SanderS15}, Sander and Sander proved the following weaker form of Conjecture \ref{SoConj}.

\begin{thm}\emph{(\cite[Theorem 1.2]{SanderS15})}
Let $n \ge 3$ be an integer.
Let $D$ and $E$ be multiplicative divisor sets of $n$. Then the integral circulant graphs $\ICG(n,D)$ and $\ICG(n, E)$ are isomorphic if and only if $\overrightarrow{\lambda}(n,D) = \overrightarrow{\lambda}(n, E)$.
\end{thm}

In addition, an explicit formula for the eigenvalues of any integral circulant graph with a multiplicative divisor set was given in \cite[Theorem 1.1]{SanderS15}.

The next two theorems were proved by Abdollahi, Janbaz and Jazaeri in \cite{AbdollahiJJ16}.

\begin{thm}
\emph{(\cite[Theorem 1.2]{AbdollahiJJ16})}
The following hold:
\begin{itemize}
\item[\rm (a)] for any prime $p > 5$, every Sylow $p$-subgroup of any finite Cay-DS group is cyclic; and every Sylow $5$-subgroup of any finite DS group is cyclic;
\item[\rm (b)] every Sylow $2$-subgroup of any finite Cay-DS group is of order at most $16$;
\item[\rm (c)] every Sylow $3$-subgroup of any finite Cay-DS group is either cyclic or is isomorphic to $\mathbb{Z}_3 \times \mathbb{Z}_3$.
\end{itemize}
\end{thm}

\begin{thm}
\label{thm:CayDSD2p}
\emph{(\cite[Theorem 1.3]{AbdollahiJJ16})}
Let $p$ be a prime. The dihedral group $D_{2p}$ is Cay-DS if and only if $p\in\{2,3,5,7,11\}$. Moreover, there exist two cospectral non-isomorphic $6$-regular Cayley graphs on $D_{2p}$ for every prime $p \ge 13$.
\end{thm}

In \cite{HuangHL17}, Huang, Huang and Lu classified all cubic Cayley graphs on $D_{2p}$ and enumerated them up to isomorphism by means of the spectral method, where $p$ is an odd prime. They proved that two such graphs are isomorphic if and only if they are cospectral. In other words, they proved:

\begin{thm}
\label{thm:CayDSD2p1}
\emph{(\cite[Corollary 3.8]{HuangHL17})}
Let $p$ be an odd prime. All cubic Cayley graphs on $D_{2p}$ are Cay-DS.
\end{thm}

Huang, Huang and Lu \cite{HuangHL17} also posed the question of classifying and enumerating connected cubic Cayley graphs on general dihedral groups and determining which of them have the Cay-DS property.

In \cite{TangCLF21}, Tang, Cheng, Liu and Feng proved the following result for dicyclic groups.

\begin{thm}
\emph{(\cite[Proposition 3.8]{TangCLF21})}
Let $p \ge 3$ be a prime. The dicyclic group $\Dic_{4p}$ is Cay-DS if and only if $p = 3$.
\end{thm}

We say that a group $\Ga$ is \emph{DS} (determined by spectrum) if for any Cayley graph $G = \Cay(\Ga, S)$ on $\Ga$ and any graph $G'$ cospectral with $G$ we have $G \cong G'$. Obviously, a DS group is Cay-DS, but the converse is not true.

\begin{thm}
\label{thm:CayDSSol}
\emph{(\cite[Theorem 1.4]{AbdollahiJJ16})}
Every finite DS group is solvable.
\end{thm}


\section{Cayley graphs on finite rings}
\label{sec:ComRing}

In this section we focus on Cayley graphs on the additive groups of finite rings. As we will see shortly, every finite commutative ring is associated with two specific Cayley graphs, the unitary and quadratic unitary Cayley graphs, of which the former as a generalization of the unitary Cayley graph of $\ZZZ_n$ is always integral when the ring is commutative. We will also see that every chain ring gives rise to a family of integral Cayley graphs whose eigenvalues can be computed explicitly. A family of Cayley graphs on the additive group of a finite field and their eigenvalues will also be discussed in this section.

Denote by $R^\times$ the set of units of a ring $R$. A \emph{local ring} \cite{Atiyah69} is a commutative ring with a unique maximal ideal. It is readily seen \cite{Atiyah69,Dummit04} that the set of units of a local ring $R$ with maximal ideal $M$ is given by $R^\times=R\setminus M$.
It is well known \cite{Atiyah69,Dummit04} that every finite commutative ring can be expressed as a direct product of finite local rings, and this decomposition is unique up to permutations of such local rings. We make the following assumption throughout this section.

\begin{assump}
\label{as:1}
{\em Whenever we consider a finite commutative ring $R$ with unit element $1 \neq 0$, we assume that it is decomposed as
$$
R=R_1\times R_2\times\cdots\times R_s
$$
such that
$$
\frac{|R_1|}{m_1} \leq \frac{|R_2|}{m_2} \leq \cdots \leq \frac{|R_s|}{m_s},
$$
where each $R_i$, $1 \le i \le s$, is a local ring with maximal ideal $M_i$ of order $m_i$.}
\end{assump}


\subsection{Unitary Cayley graphs of finite rings}
\label{subsec:UCayComRing}

Let $R$ be a finite ring. The \emph{unitary Cayley graph} of $R$ is defined as the Cayley graph $G_R=\Cay(R,R^\times)$ of the additive group of $R$ with respect to $R^\times$. That is, $G_R$ has vertex set $R$ such that $x, y \in R$ are adjacent if and only if $x-y\in R^\times$. (Note that $G_R$ is the edgeless graph on $R$ if $R$ has no units.) This notion is a generalization of the unitary Cayley graph of $\mathbb{Z}_n$ because
$$
G_{\mathbb{Z}_n} \cong \Cay(\mathbb{Z}_n, \ZZZ_n^{\times}) \cong \ICG(n,\{1\}).
$$

In the case when $R$ is a local ring with maximal ideal $M$, $G_R$ is the complete $(|R|/|M|)$-partite graph with each part containing $|M|$ vertices, and so its spectrum is known (see \cite[Proposition 2.1]{Kiani11}). On the other hand, for a general ring $R$ as decomposed in Assumption \ref{as:1}, we have $G_R = G_{R_1} \otimes \cdots \otimes G_{R_s}$. Furthermore, it is well known that if $G$ and $H$ have eigenvalues $\l_i, 1 \le i \le n$ and $\mu_j, 1 \le j \le m$, respectively, where repetitions are allowed, then their tensor product $G \otimes H$ has eigenvalues $\l_i \mu_j: 1 \le i \le n, 1 \le j \le m$. Combining all these facts, the spectrum of $G_R$ can be easily determined as shown in \cite{Kiani11}. Define
$$
\lambda_C = (-1)^{|C|}\dfrac{|R^\times|}{\prod_{j\in C} \frac{|R^\times_j|}{m_j}}
$$
for every subset $C$ of $\{1,2,\ldots,s\}$.

\begin{thm}
\label{ringspectrum}
\emph{(\cite[Lemma 2.3]{Kiani11}; see also \cite{Akhtar09})}
Let $R$ be a finite commutative ring. The eigenvalues of $G_R$ are
\begin{itemize}
\item[\rm (a)] $\lambda_C$, repeated $\prod_{j\in C} \frac{|R^\times_j|}{m_j}$ times, for $C \subseteq \{1,2,\ldots,s\}$; and
\item[\rm (b)] $0$ with multiplicity $|R|-\prod_{i=1}^s\left(1+\dfrac{|R^\times_i|}{m_i}\right)$.
\end{itemize}
In particular, if $R$ is a finite local ring and $m$ is the order of its unique maximal ideal, then
\begin{equation}
\label{eq:GRlocal}
\Spec (G_R) = \left(|R|-m, (-m)^{\frac{|R|-m}{m}}, 0^{\frac{|R|(m-1)}{m}}\right).
\end{equation}
\end{thm}

In \cite{Liu12}, {Liu} and {Zhou} determined the eigenvalues of the complement $\overline{G}_R$ and the line graph $L(G_R)$ of $G_R$.

\begin{thm}
\emph{(\cite[Corollary 6]{Liu12})}
\label{Corringspectrum}
Let $R$ be a finite commutative ring. The eigenvalues of $\overline{G}_R$ are
\begin{itemize}
\item[\rm (a)] $|R|-1-|R^\times|$;
\item[\rm (b)] $-\lambda_C - 1$, repeated $\prod_{j\in C} \frac{|R^\times_j|}{m_j}$ times, for $\emptyset \ne C \subseteq \{1,2,\ldots,s\}$; and
\item[\rm (c)] $-1$ with multiplicity $|R|-\prod_{i=1}^s\left(1+\dfrac{|R^\times_i|}{m_i}\right)$.
\end{itemize}
In particular, if $R$ is a finite local ring and $m$ is the order of its unique maximal ideal, then
$$
\Spec(\overline{G}_R)= \left((m-1)^{\frac{|R|}{m}}, (-1)^{\frac{|R|}{m}(m-1)}\right).
$$
\end{thm}

\begin{thm}
\label{lineeigenvalue}
\emph{(\cite[Corollary 7]{Liu12})}
Let $R$ be a finite commutative ring. The eigenvalues of $L(G_R)$ are
\begin{itemize}
\item[\rm (a)] $\lambda_C+|R^\times|-2$, repeated $\prod_{j\in C} \frac{|R^\times_j|}{m_j}$ times, for $C \subseteq \{1,2,\ldots,s\}$;
\item[\rm (b)] $|R^\times|-2$ with multiplicity $|R|-\prod_{i=1}^s\left(1+\dfrac{|R^\times_i|}{m_i}\right)$; and
\item[\rm (c)] $-2$, repeated $|R|\left(|R^\times|-2\right)/2$ times.
\end{itemize}
In particular, if $R$ is a finite local ring and $m$ is the order of its unique maximal ideal, then
$$
\Spec (L(G_R))=\left(2|R^\times|-2, (|R^\times|-m-2)^{\frac{|R^\times|}{m}}, (|R^\times|-2)^{(\frac{|R^\times|}{m}+1)(m-1)}, (-2)^{\frac{|R|(|R^\times|-2)}{2}}\right).
$$
\end{thm}

Theorems \ref{ringspectrum}--\ref{lineeigenvalue} together imply the following result.

\begin{cor}
Let $R$ be a finite commutative ring. Then $G_R$, $\overline{G}_R$ and $L(G_R)$ are integral. In particular, $G_R$ is an integral Cayley graph.
\end{cor}

Let $R$ be a finite ring, and let $M_{n}(R)$ be the ring of $n \times n$ matrices over $R$, where $n \ge 1$. Then the set of units of $M_{n}(R)$ is the group $\GL_{n}(R)$ of $n \times n$ invertible matrices over $R$. So the unitary Cayley graph of $M_{n}(R)$ is the graph $G_{M_{n}(R)} = \Cay(M_{n}(R), \GL_{n}(R))$. In \cite{RattaM20}, Rattanakangwanwong and Meemark studied various spectral properties of $G_{M_{n}(R)}$, with a focus on the case when $R$ is a finite field or local ring. Here we present the following two results from this paper on the spectrum of $G_{M_{n}(\FFF_q)}$ when $n = 2$ or $3$. Other spectral properties of $G_{M_{n}(R)}$, including when it is Ramanujan or strongly regular, will be discussed in related later sections.

\begin{thm}
\emph{(\cite[Theorem 2.2]{RattaM20})}
\label{thm:mat-rings-1}
For any prime power $q$, we have
\begin{eqnarray*}
\Spec (G_{M_{2}(\FFF_q)}) = \pmat{(q^2 - 1)(q^2 - q) & -(q^2 - q) & q \\
1  & (q - 1)(q + 1)^2  & (q^2 - 1)(q^2 - q)}.
\end{eqnarray*}
\end{thm}

\begin{thm}
\emph{(\cite[Theorem 2.3]{RattaM20})}
\label{thm:mat-rings-2}
For any prime power $q$, the spectrum of $G_{M_{3}(\FFF_q)}$ is
$\Spec (G_{M_{3}(\FFF_q)}) = (((q^3 - 1)(q^3 - q)(q^3 - q^2))^1, (-(q^3 - q)(q^3 - q^2))^{(q^3 - 1)(q^2 + q + 1)},$ $(q(q^3 - q^2))^{(q^3 - 1)(q^3 - q)(q^2 + q + 1)}, (-q^3)^{(q^3 - 1)(q^3 - q)(q^3 - q^2)})$.
\end{thm}


\subsection{Quadratic unitary Cayley graphs of finite rings}
\label{subsec:QUCayComRing}

Let $R$ be a finite ring. Let $Q_R=\{u^2: u\in R^\times\}$ and set $T_R=Q_R\cup(-Q_R)$. The \emph{quadratic unitary Cayley graph} of $R$, denoted by $\mathcal{G}_R$, is defined as the Cayley graph $\Cay(R,T_R)$ on the additive group of $R$ with respect to $T_R$. That is, $\mathcal{G}_R$ has vertex set $R$ such that $x, y \in R$ are adjacent if and only if $x-y\in T_R$. This notion introduced by Liu and Zhou \cite{Liu15} is a generalization of the quadratic unitary Cayley graph $\mathcal{G}_{\mathbb{Z}_n}$ of $\mathbb{Z}_n$ introduced in \cite{Beaudrap10} as well as the well known Paley graphs. (The \emph{Paley graph} $P(q)$, where $q \equiv 1~\mod~4$ is a prime power, is the Cayley graph on the additive group of the finite field $\mathbb{F}_q$ with respect to the set of nonzero squares.) In \cite{Liu15}, {Liu} and {Zhou} determined the spectra of $\mathcal{G}_R$ for finite commutative rings in the case when all finite fields $R_i/M_i$ are of odd order and at most one of them has order congruent to $3$ modulo $4$.

\begin{thm}
\label{SPECLocalQUCG}
\emph{(\cite[Theorem 2.4]{Liu15})}
Let $R$ be a local ring with maximal ideal $M$ of order $m$.
\begin{itemize}
\item[\rm (a)] If $|R|/m\equiv 1~(\mod\,4)$, then the spectrum of $\mathcal{G}_R$ is
\begin{eqnarray*}
\Spec (\mathcal{G}_R) = \left(\begin{array}{cccc}
\frac{|R|-m}{2} & \frac{1}{2}m\left(-1+\sqrt{\frac{|R|}{m}}\right) & \frac{1}{2}m\left(-1-\sqrt{\frac{|R|}{m}}\right) & 0 \\[0.1cm]
1 & \frac{1}{2}\left(\frac{|R|}{m}-1\right) & \frac{1}{2}\left(\frac{|R|}{m}-1\right) & |R|-\frac{|R|}{m}
\end{array}
\right).
\end{eqnarray*}
\item[\rm (b)] If $|R|/m\equiv 3~(\mod\,4)$, then the spectrum of $\mathcal{G}_R$ is
\begin{eqnarray*}
\Spec (\mathcal{G}_R) = \left(\begin{array}{ccc}
|R|-m & -m & 0 \\ [0.1cm]
1 & \frac{|R|}{m}-1 & |R|-\frac{|R|}{m}
\end{array}
\right).
\end{eqnarray*}
In particular, $\mathcal{G}_R$ is an integral Cayley graph.
\end{itemize}
\end{thm}

Define
$$
\lambda_{A,B} = (-1)^{|B|}\dfrac{|R^\times|}{2^s\prod_{i\in A}\left(\sqrt{\frac{|R_i|}{m_i}}+1\right) \prod_{j\in B}\left(\sqrt{\frac{|R_j|}{m_j}}-1\right)}
$$
for disjoint subsets $A, B$ of $\{1,2,\ldots,s\}$.

\begin{thm}
\label{Spec1mod4}
\emph{(\cite[Theorem 2.6]{Liu15})}
Let $R$ be a finite commutative ring such that $|R_i|/m_i\equiv 1~(\mod\,4)$ for $1 \le i \le s$.  Then the eigenvalues of $\mathcal{G}_R$ are
\begin{itemize}
\item[\rm (a)] $\lambda_{A,B}$, repeated $\dfrac{1}{2^{|A|+|B|}}\prod\limits_{k\in A\cup B}\left(\frac{|R_k|}{m_k}-1\right)$ times, for all pairs of disjoint subsets $A, B$ of $\{1,2,\ldots,s\}$; and
\item[\rm (b)] $0$ with multiplicity $|R|-\sum\limits_{\substack{A, B\subseteq\{1,\ldots,s\}\\ A\cap B=\emptyset}}\left(\dfrac{1}{2^{|A|+|B|}}\prod\limits_{k\in A\cup B}\left(\frac{|R_k|}{m_k}-1\right)\right)$.
\end{itemize}
\end{thm}

\begin{thm}
\label{Spec3mod4}
\emph{(\cite[Theorem 2.7]{Liu15})}
Let $R$ be a finite commutative ring such that $|R_i|/m_i\equiv 1~(\mod\,4)$ for $1 \le i \le s$.  Let $R_0$ be a local ring with maximal ideal $M_0$ of order $m_0$ such that $|R_0|/m_0\equiv 3~(\mod\,4)$. Then the eigenvalues of $\mathcal{G}_{R_0\times R}$ are
\begin{itemize}
\item[\rm (a)] $|R_0^\times|\cdot\lambda_{A,B}$, repeated $\dfrac{1}{2^{|A|+|B|}}\prod\limits_{k\in A\cup B}\left(\frac{|R_k|}{m_k}-1\right)$ times, for all pairs of disjoint subsets $A, B$ of $\{1,2,\ldots,s\}$;
\item[\rm (b)] $-\dfrac{|R_0^\times|}{|R_0|/m_0-1}\cdot\lambda_{A,B}$, repeated $\dfrac{1}{2^{|A|+|B|}}\left(\frac{|R_0|}{m_0}-1\right) \prod\limits_{k\in A\cup B}\left(\frac{|R_k|}{m_k}-1\right)$ times, for all pairs of disjoint subsets $A, B$ of $\{1,2,\ldots,s\}$; and
\item[\rm (c)] $0$ with multiplicity $|R|-\sum\limits_{\substack{A,B\subseteq\{1,\ldots,s\}\\ A\cap B=\emptyset}}\left(\dfrac{1}{2^{|A|+|B|}} \frac{|R_0|}{m_0} \prod\limits_{k\in A\cup B}\left(\frac{|R_k|}{m_k}-1\right)\right)$.
\end{itemize}
\end{thm}

Theorems \ref{SPECLocalQUCG}, \ref{Spec1mod4} and \ref{Spec3mod4} can be specified to obtain the eigenvalues of ${\cal G}_{\mathbb{Z}_n}$.

\begin{thm}
\label{G2pspec}
\emph{(\cite[Corollary 2.8]{Liu15})}
Let $p \ge 5$ be an odd prime and $\alpha\ge1$ an integer.
\begin{itemize}
\item[\rm (a)] If $p\equiv1~(\mod\,4)$, then
\begin{eqnarray*}
\Spec (\mathcal {G}_{\mathbb{Z}_{p^\alpha}})=\left(\begin{array}{cccc}
p^{\alpha-1}(p-1)/2 & p^{\alpha-1}\left(-1+\sqrt{p}\right)/2 & 0      &  p^{\alpha-1}\left(-1-\sqrt{p}\right)/2 \\[0.1cm]
1              & (p-1)/2                & p^\alpha-p  &  (p-1)/2
\end{array}
\right).
\end{eqnarray*}
\item[\rm (b)] If $p\equiv3~(\mod\,4)$, then
\begin{eqnarray*}
\Spec (\mathcal {G}_{\mathbb{Z}_{p^\alpha}})=\left(\begin{array}{ccc}
p^{\alpha-1}(p-1)   & -p^{\alpha-1}        &  0       \\[0.1cm]
1              & p-1         & p^\alpha-p
\end{array}
\right).
\end{eqnarray*}
\end{itemize}
\end{thm}

\begin{thm}
\label{EigenQUCG}
\emph{(\cite[Corollary 2.9]{Liu15})}
Let $n=p_1^{\alpha_1}\ldots p_s^{\alpha_s}$ be an integer in canonical factorization such that $p_i \equiv1~(\mod\,4)$ for each $i$. Then the eigenvalues of ${\cal G}_{\mathbb{Z}_n}$ are
\begin{itemize}
\item[\rm (a)] $(-1)^{|B|}\cdot\dfrac{\varphi(n)}{2^s\prod_{i\in A}\left(\sqrt{p_i}+1\right)\prod_{j\in B}\left(\sqrt{p_j}-1\right)}$, repeated $\dfrac{1}{2^{|A|+|B|}}\prod\limits_{k\in A\cup B}(p_k-1)$ times, for all pairs of disjoint subsets $A, B$ of $\{1,2,\ldots,s\}$;
\item[\rm (b)] $0$ with multiplicity $n-\sum\limits_{\substack{A,\,B\subseteq\{1,\ldots,s\}\\A\cap B=\emptyset}}\left(\dfrac{1}{2^{|A|+|B|}}\prod\limits_{k\in A\cup B}(p_k-1)\right)$.
\end{itemize}
\end{thm}

\begin{thm} \label{EigenZ_n121}
\emph{(\cite[Corollary 2.10]{Liu15})}
Let $n=p^\alpha p_1^{\alpha_1}\ldots p_s^{\alpha_s}$ be an integer in canonical factorization such that $p\equiv3~(\mod\,4)$ and $p_i \equiv1\,(\mod\,4)$ for each $i$. Then the eigenvalues of ${\cal G}_{\mathbb{Z}_n}$ are
\begin{itemize}
\item[\rm (a)] $(-1)^{|B|}\cdot\dfrac{\varphi(n)}{2^s\prod_{i\in A}\left(\sqrt{p_i}+1\right)\prod_{j\in B}\left(\sqrt{p_j}-1\right)}$, repeated $\dfrac{1}{2^{|A|+|B|}}\prod\limits_{k\in A\cup B}(p_k-1)$ times, for all pairs of disjoint subsets $A, B$ of $\{1,2,\ldots,s\}$;
\item[\rm (b)] $(-1)^{|B|+1}\cdot\dfrac{\varphi(n)}{2^s(p-1)\prod_{i\in A}\left(\sqrt{p_i}+1\right)\prod_{j\in B}\left(\sqrt{p_j}-1\right)}$, repeated $\dfrac{p-1}{2^{|A|+|B|}}\prod\limits_{k\in A\cup B}(p_k-1)$ times, for all pairs of disjoint subsets $A, B$ of $\{1,2,\ldots,s\}$; and
\item[\rm (c)] $0$ with multiplicity $n-\sum\limits_{\substack{A,\,B\subseteq\{1,\ldots,s\}\\A\cap B=\emptyset}}\left(\dfrac{p}{2^{|A|+|B|}}\prod\limits_{k\in A\cup B}(p_k-1)\right)$.
\end{itemize}
\end{thm}


\subsection{A family of Cayley graphs on finite chain rings}
\label{subsec:ring-gcd}

A finite \emph{chain ring} is a finite local ring $R$ such that for any ideals $I, J$ of $R$ we have either $I \subseteq J$ or $J \subseteq I$. Consider a finite chain ring $R$ with unique maximal ideal $M$ and residue field of $q$ elements. The \emph{nilpotency} $s$ of $R$ is the smallest positive integer such that $M^s = \{0\}$. Then $\{0\} = M^s \subset M^{s-1} \subset \cdots \subset M^2 \subset M \subset M^0 = R$ and $|M^i| = q^{s-i}$ for $0 \le i \le s$. For integers $a_1, a_2, \ldots, a_r$ with $0 \le a_1 < a_2 < \cdots < a_r \le s-1$, set
\begin{equation}
\label{eq:calC}
D = \cup_{i=1}^{r} (M^{a_i} \setminus M^{a_i + 1}).
\end{equation}
The Cayley graph $\Cay(R, D)$ was studied in \cite{SuntornpochM16}. In the special case when $R = \ZZZ_{p^s}$ and $D = \{p^{a_i}: 1 \le i \le r\}$ is a set of proper divisors of $p^s$, this graph is exactly the gcd graph of
$\ZZZ_{p^s}$ with respect to $D$.

\begin{thm}
\label{thm:sun}
\emph{(\cite[Section 2]{SuntornpochM16})}
Let $R$ be a finite chain ring and let $D$ be as above. Then $\Cay(R, D)$ is integral. Moreover, the following hold:
\begin{itemize}
\item[\rm (a)] if $a_r = s-1$, then the eigenvalues of $\Cay(R, D)$ are
\begin{itemize}
\item[\rm (i)] $(q-1)\sum_{i=1}^{r} q^{s - a_i - 1}$, with multiplicity $q^{a_1}$;
\item[\rm (ii)] $-q^{s-a_{k-1}-1}+(q-1)\sum_{i=k}^{r} q^{s - a_i - 1}$, with multiplicity $q^{a_{k-1}} (q-1)$, for $2 \le k \le r$;
\item[\rm (iii)] $(q-1)\sum_{i=k}^{r} q^{s - a_i - 1}$, with multiplicity $q^{a_k - a_{k-1} - 1} - q^{a_{k-1}+1}$, for $2 \le k \le r$; and
\item[\rm (iv)] $-1$, with multiplicity $q^{a_{r}} (q-1)$;
\end{itemize}
\item[\rm (b)] if $a_r \ne s-1$, then the eigenvalues of $\Cay(R, D)$ are
\begin{itemize}
\item[\rm (i)] $(q-1)\sum_{i=1}^{r} q^{s - a_i - 1}$, with multiplicity $q^{a_1}$;
\item[\rm (ii)] $-q^{s-a_{k-1}-1}+(q-1)\sum_{i=k}^{r} q^{s - a_i - 1}$, with multiplicity $q^{a_{k-1}} (q-1)$, for $2 \le k \le r$;
\item[\rm (iii)] $(q-1)\sum_{i=k}^{r} q^{s - a_i - 1}$, with multiplicity $q^{a_k} - q^{a_{k-1}+1}$, for $2 \le k \le r$;
\item[\rm (iv)] $-q^{s-a_r -1}$, with multiplicity $q^{a_{r}} (q-1)$; and
\item[\rm (v)] $0$, with multiplicity $q^{a_{r}+1} (q^{s - a_r -1}-1)$.
\end{itemize}
\end{itemize}
\end{thm}


\subsection{Generalized Paley graphs}
\label{subsec:GPaley}

Let $q = p^r$ be a prime power and $k \ge 1$ an integer. The \emph{generalized Paley graph} $X_{q}^k$ is defined as the Cayley digraph on the additive group of $\FFF_q$ with connection set $\{a^k: a \in \FFF_q^*\}$ (see \cite{JohnsonSS16, PodestaV20a}). Since $\{a^k: a \in \FFF_q^*\} = \langle \om^k \rangle = \langle \om^{k'} \rangle$, where $\om$ is a primitive element of $\FFF_q$ and $k' = \gcd(q-1,k)$, we have $X_{q}^k = X_{q}^{k'}$. Hence we may always assume that $k$ is a divisor of $q-1$ in $X_{q}^k$. Note that $X_{q}^k$ is an undirected graph if and only if either $q$ or $(q-1)/k$ is even, and in this case $X_{q}^k$ is exactly the generalized Paley graph $\mathrm{GPaley}(q, (q-1)/k)$ which was defined in \cite{LimP09} as the Cayley graph on the additive group of $\FFF_q$ with connection set the unique subgroup of the multiplicative group $\FFF_q^*$ of order $(q-1)/k$. Note also that $X_{q}^1$ is the complete graph $K_{q}$ and $X_{q}^2$ is exactly the Paley graph $P(q)$ when $q \equiv 1~(\mod~4)$.

In \cite{PodestaV20a}, Podest\'{a} and Videla proved that the eigenvalues of the (undirected) generalized Paley graphs $X_{q}^k$ and their complements $\overline{X_{q}^k}$ can be expressed in terms of the Gaussian periods
\begin{equation}
\label{eq:gas-period}
\eta_{i}^{(k,q)} = \sum_{x \in \om^i \langle \om^k \rangle} \om_{p}^{\Tr(x)},\;\, 0 \le i \le k-1,
\end{equation}
where
\begin{equation}
\label{eq:trx}
\Tr(x) = \Tr_{\FFF_q/\FFF_p}(x) = x+x^p+ \cdots + x^{p^{r-1}}
\end{equation}
is the trace of $x \in \FFF_q$ and $\om^i \langle \om^k \rangle$ is the coset of the subgroup $\langle \om^k \rangle$ of $\FFF_q^*$ containing $\om^i$.

\begin{thm}
\label{PV20a}
\emph{(\cite[Theorem 2.1]{PodestaV20a})}
Let $q = p^r$ be a prime power and $k \ge 1$ a divisor of $q-1$. Suppose that either $q$ or $n = (q-1)/k$ is even. Let $\{\eta_1, \ldots, \eta_s\}$ be the set of distinct Gaussian periods other than $n$ among $\eta_{0}^{(k,q)}, \eta_{1}^{(k,q)}, \ldots, \eta_{k-1}^{(k,q)}$, and let $m_i$ be the multiplicity of $\eta_i$ for $1 \le i \le s$. Let $m$ be the multiplicity (possibly $0$) that $n$ occurs among these Gaussian periods. Then
\begin{eqnarray*}
\Spec (X_{q}^k)=\left(\begin{array}{cccc}
n  &  \eta_1  &  \ldots  & \eta_s \\[0.1cm]
mn+1  & m_1 n & \ldots & m_s n
\end{array}
\right)
\end{eqnarray*}
and
\begin{eqnarray*}
\Spec (\overline{X_{q}^k})=\left(\begin{array}{cccc}
(k-1)n  &  -\eta_1 - 1  &  \ldots  & -\eta_s - 1 \\[0.1cm]
mn+1  & m_1 n & \ldots & m_s n
\end{array}
\right).
\end{eqnarray*}
Moreover, if $k$ divides $(q-1)/(p-1)$, then both $X_{q}^k$ and $\overline{X_{q}^k}$ are integral.
\end{thm}

In particular, this implies the well-known fact that the Paley graph $P(q) = X_{q}^2$ has eigenvalues $(q-1)/2$ (with multiplicity $1$) and $(-1 \pm \sqrt{q})/2$ (both with multiplicity $(q-1)/2$). For $k=3, 4$, the eigenvalues of $X_{q}^k$ have been worked out explicitly in \cite{PodestaV20a} when $k$ divides $(q-1)/(p-1)$.

In the case when $k=2$ and $q \equiv 1~(\mod~4)$ or $k > 2$ and $k$ divides $p^{t} + 1$ for some proper divisor $t$ of $r/2$ (where $r$ is even), $(k, q)$ is called a \emph{semiprimitive pair} of integers and the generalized Paley graph $X_{q}^k$ is said to be \emph{semiprimitive}. In this case the eigenvalues of $X_{q}^k$ and $\overline{X_{q}^k}$ can be explicitly computed as shown in \cite{PodestaV20a}.

\begin{thm}
\label{PV20b}
\emph{(\cite[Theorem 3.3]{PodestaV20a})}
Let $q = p^r$ be a prime power with $r$ even, and let $k \ge 1$ be a divisor of $q-1$. Suppose that $(k, q)$ is a semiprimitive pair, with $t$ the smallest divisor of $r/2$ such that $k$ divides $p^t + 1$. Set $n = (q-1)/k$. Then
\begin{eqnarray*}
\Spec (X_{q}^k)=\left(\begin{array}{ccc}
n  & -\frac{1}{k}\left((k-1)(-1)^{\frac{r}{2t}}\sqrt{q}+1\right) & \frac{1}{k}\left((-1)^{\frac{r}{2t}}\sqrt{q}-1\right) \\[0.1cm]
1  & n  & (k-1)n
\end{array}
\right)
\end{eqnarray*}
and
\begin{eqnarray*}
\Spec (\overline{X_{q}^k})=\left(\begin{array}{ccc}
(k-1)n & \frac{k-1}{k}\left((-1)^{\frac{r}{2t}}\sqrt{q}-1\right) & -1 - \frac{1}{k}\left((-1)^{\frac{r}{2t}}\sqrt{q}-1\right) \\[0.1cm]
1  & n  & (k-1)n
\end{array}
\right).
\end{eqnarray*}
\end{thm}


\section{Energies of Cayley graphs}
\label{sec:EnCay}

The \emph{energy} of a graph $G$ with $n$ vertices and eigenvalues $\lambda_1,\lambda_2,\ldots, \lambda_n$ is defined as
$$
\E(G) = \sum_{i=1}^n |\lambda_i|.
$$
This concept was introduced by Gutman \cite{Gutman78} in the study of mathematical chemistry. In the past four decades a number of results on energies of various families of graphs have been obtained; see \cite{Li12} for a monograph on this topic. It is well known \cite{KoolenM01} that, for any graph $G$ with $n$ vertices,
\begin{equation}
\label{eq:KoolenM01}
\E(G) \le \frac{n}{2}\left(\sqrt{n}+1\right),
\end{equation}
and the bound can be achieved by infinitely many graphs. It can be easily verified that $\E(K_n) = 2(n-1)$. A graph $G$ with $n$ vertices is called \emph{hyperenergetic} \cite{Gutman99} if $\E(G) > 2(n-1)$.


\subsection{Energies of integral circulant graphs}
\label{Sec-ener}

In this section we assume that
\be
\label{eq:n}
n=p_1^{\alpha_1}p_2^{\alpha_2} \ldots p_s^{\alpha_s} \ge 2
\ee
is an integer in canonical factorization into prime powers, where $p_1<p_2<\ldots<p_s$ are primes and each $\alpha_i\geq1$ is an integer. Recall that the circulant integral graphs are precisely the gcd graphs $\ICG(n,D)$ of cyclic groups (Corollary \ref{cor:intICG}) and that $\Cay(\mathbb{Z}_n, \ZZZ_n^{\times}) = \ICG(n,\{1\})$ is the unitary Cayley graph of $\mathbb{Z}_n$. As before $\varphi(n)$ is Euler's totient function and $D(n)$ denotes the set of positive divisors of $n$.

The first two parts of the following result were proved in \cite[Theorems 2.3 and 2.4]{Ilic09} and \cite[Theorems 3.7 and 3.10]{Ramaswamy09}, and the third part was proved in \cite[Theorem 3.11]{Ramaswamy09}. The first part of this result was also obtained in \cite{Balakrishnan04} in the special case when $n$ is a prime power.

\begin{thm}
\label{thm:energy}
\emph{(\cite{Ilic09,Ramaswamy09})}
Let $n \ge 3$ be an integer.
\begin{itemize}
\item[\rm (a)]
The energy of $\ICG(n,\{1\})$ is given by
$$
\E(\ICG(n,\{1\})) = 2^s\varphi(n).
$$
\item[\rm (b)]
$\ICG(n,\{1\})$ is hyperenergetic if and only if $n$ has at least two distinct odd prime factors.
\item[\rm (c)]
The energy of $\ICG(n,\{1\})$ satisfies
$$
\E(\ICG(n,\{1\})) > \frac{2^{s-1} (n-1)}{s}.
$$
\end{itemize}
\end{thm}

\begin{thm}
	\label{thm:energy1}
	\emph{(\cite[Theorems 3.1 and 3.2]{Ilic09})}
	Let $n \ge 3$ be an integer.
	\begin{itemize}
		\item[\rm (a)]
		The energy of the complement of $\ICG(n,\{1\})$ is equal to
		$$
		2n-2+(2^s-2)\varphi(n)-\prod_{i=1}^sp_i+\prod_{i=1}^s(2-p_i).
		$$
		\item[\rm (b)]
The complement of $\ICG(n,\{1\})$ is hyperenergetic if and only if $s\ge2$ and $n\neq 2p$, where $p$ is a prime.
	\end{itemize}
\end{thm}


The next two results were obtained by Ili\'{c} and Ba\v{s}i\'{c} in \cite{Ilic11}. Denote $p_i^\a \parallel n$ if $p_i^\a ~|~ n$ but $p_i^{\a+1}\nmid  n$.

\begin{thm}
	\label{thm:energy2}
\emph{(\cite[Theorem 4.1]{Ilic11})}
If $n \ge 4$, then for each $i$ and any integer $\g$ with $1\le \gamma\le \alpha_i$,
$$
\E(\ICG(n,\{1,p_i^\gamma\})) =
\begin{cases}
                        2^{s-1}(\varphi(n)+\varphi(n/p_i)), & \text{if $p_i\parallel n$} \\
                        2^{s-1}(2\varphi(n)+(p_i^\gamma-2p+2)\varphi(n/p_i)), & \text{if $p_i^\gamma\parallel n$ with $\gamma\ge2$} \\
                        2^{s}(\varphi(n)+(p_i^\gamma-p_i+1)\varphi(n/p_i)),  & \text{if $p_i^\gamma\nparallel n$.}
\end{cases}
$$
\end{thm}

\begin{thm}
\label{thm:energy3}
\emph{(\cite[Theorem 4.2]{Ilic11})}
If $n \ge 4$, then for $1\le i < j \le s$,
$$
\E(\ICG(n,\{p_i,p_j\}))=
\begin{cases}
                        2^{s}\varphi(n), & \text{if $p_i\parallel n$ and $p_j\parallel n$} \\
                        3\cdot 2^{s-1}\varphi(n), & \text{if $2\parallel n$ and $p_j^2| n$} \\
                        2^{s-1}(2\varphi(n)+\varphi(n/p_j)\varphi(p_j)),  & \text{if $p_i\parallel n$ with $p_i\neq 2$}\\
                        &\text{and $p_j^2| n$}\\
                        2^{s-1}(2\varphi(n)+\varphi(n/p_i)\varphi(p_i)),  & \text{if $p_i^2| n$ and $p_j\parallel n$}\\
                        2^{s-1}(2\varphi(n)+\varphi(n/p_i)\varphi(p_i)+\varphi(n/p_j)\varphi(p_j)),  & \text{if $p_i^2| n$ and $p_j^2| n$}.
\end{cases}
$$
\end{thm}

The first formula in Theorem \ref{thm:energy2} was obtained earlier in \cite[Theorem 4.1]{Ilic09} (see also \cite[Corollary 3.8]{Mollahajiaghaei12}), and in the special case when $n=p_1p_2 \ldots p_s$ the first formula in Theorem \ref{thm:energy3} was also obtained in \cite[Theorem 4.2]{Ilic09}.

Theorem \ref{thm:energy3} implies that there are a large number of non-cospectral regular hyperenergetic graphs of the same order and same energy. See \cite[Section 5]{Ilic11} for more details.

The next three results were obtained by Mollahajiaghaei in \cite{Mollahajiaghaei12}.

\begin{thm}
	\emph{(\cite[Theorem 3.10]{Mollahajiaghaei12})}
	Let $n = p^km$, where $p$ is a prime, $k > 1$ and $\gcd(m,p) = 1$. Then
	$$
	\E(\ICG(n;\{1,p\})) = 2(p^k-p^{k-2})\E(\ICG(m,\{1\})).
	$$
\end{thm}

\begin{thm}
	\emph{(\cite[Theorem 3.11]{Mollahajiaghaei12})}
	Let $n = p^km$, where $p$ is a prime, $k > 1$ and $\gcd(m,p) = 1$. Then
	$\ICG(n;\{1,p\})$ is hyperenergetic if and only if $m$ has at least two distinct prime factors or $m$ is an
	odd integer.
\end{thm}

\begin{thm}
\emph{(\cite[Theorems 4.2 and 4.7]{Mollahajiaghaei12})}
Let $n=p_1p_2 \ldots p_sm \ge 3$ and $D=$ $\{p_1,p_2,\ldots,p_s\}$, where $p_1, p_2,\ldots, p_s$ are distinct primes each of which is greater than $s$, and $m$ is a positive integer with $\gcd(p_1p_2\ldots p_s, m) = 1$. Let $t$ be the number of prime factors of $m$. Then
$$
\E(\ICG(n,D))=2^t\left(s2^{s-1}-(2^{s-1}-2)\sum_{i=1}^s\frac{1}{p_i-1}\right)\varphi(n).
$$
Moreover, $\ICG(n,D)$ is hyperenergetic.
\end{thm}

The next result tells us when the energy of an integral circulant graph is divisible by $4$.

\begin{thm}
\emph{(\cite[Theorems 3.2 and 3.3]{Ilic11})}
Let $n \ge 3$ be an integer and $D \subseteq D(n) \setminus \{n\}$. If $n$ is odd, then $\E(\ICG(n, D))$ is divisible by $4$; if $n$ is even, then $\E(\ICG(n, D))$ is not divisible by $4$ if and only if $n/2 \notin D$ and $\sum_{d \in D} (-1)^d \varphi(n/d) < 0$.
\end{thm}

In fact, when $n$ is odd, $\sum_{d \in D} (-1)^d \varphi(n/d)$ is an eigenvalue of $\ICG(n, D)$.

\delete
{
Given a divisor $d$ of $n$, let
$$
M_n(d)=\{d,2d,\ldots, n-d\}
$$
be the set of positive multiples of $d$ less than $n$. In \cite[Theorems 3.1 and 3.2]{Chelvam12} it was proved that $\Cay(\mathbb{Z}_n, \mathbb{Z}_n^{*} \setminus M_n(d))$ is an integral graph with energy $2n\left(1-\frac{1}{d}\right)$. It was further proved in \cite[Theorem 4.1]{Chelvam12} that this graph is strongly regular with parameters $\left(n, n-\frac{n}{d}, n-\frac{2n}{d}, n-\frac{n}{d}\right)$. \sm{I think we should remove these results in \cite{Chelvam12} (and also delete this paper from the bibliography) since all of them are trivial. In fact, as shown in \cite[Theorem 3.3]{Chelvam12} $\Cay(\mathbb{Z}_n, \mathbb{Z}_n^{*} \setminus M_n(d))$ is a complete $d$-partite graph. So everything about this graph is trivial. Double check.}
}

In \cite{Le12}, {Le} and {Sander} found a connection between $A$-convolutions satisfying a weak form of regularity and the spectra of integral circulant graphs. Using this, they obtained a multiplicative decomposition of the energies of integral circulant graphs with multiplicative divisor sets. A few definitions are in order before presenting their results. Call $A = (A(n))_{n\in\mathbb{N}}$ a \emph{divisor system} if $\emptyset\neq A(n)\subseteq D(n)$ for each $n\in \mathbb{N}$. Denote by $\mathbb{C}^{\mathbb{N}}$ the set of functions from $\mathbb{N}$ to $\mathbb{C}$. The \emph{A-convolution} of two functions $f, g \in \mathbb{C}^{\mathbb{N}}$ is the function $f*_Ag \in \mathbb{C}^{\mathbb{N}}$ defined by
$$
(f*_Ag)(n) = \sum_{d\in A(n)}f(d)g\left(\frac{n}{d}\right),\, n \in \mathbb{N}.
$$
The $A$-convolution is \emph{regular} if it satisfies the following conditions:
\begin{itemize}
\item $(\mathbb{C}^{\mathbb{N}},+,*_A)$ is a commutative ring with unity;
\item if $f, g \in \mathbb{C}^{\mathbb{N}}$ are both multiplicative, so is $f*_Ag$;
\item the constant function $1$ has an inverse $\mu_A$ with respect to the operation $*_A$ such that $\mu_A(p^s)\in\{-1, 0\}$ for all prime powers $p^s$ (generalization of M\"{o}bius' $\mu$ function).
\end{itemize}
In \cite{Narkiewicz63}, {Narkiewicz} proved that any regular $A$-convolution has the multiplicative property in the sense that
$$
A(mn)=A(m)A(n):=\{ab:a\in A(m), b\in A(n)\}
$$
for all coprime integers $m,n\in\mathbb{N}$. This implies that
 \begin{equation}
 \label{RdsMulti}
A\left(\prod_{i=1}^{l} p_i^{s_i}\right) = \prod_{i=1}^{l} A(p_i^{s_i}) = \left\{\prod_{i=1}^{l} a_i: a_i\in A(p_i^{s_i})\, \text{ for } 1\le i\le l\right\}
 \end{equation}
for any distinct primes $p_i$ and nonnegative integers $s_i$.

A divisor system $A = (A(n))_{n\in\mathbb{N}}$ is called \emph{multiplicative} if for all $n\in \mathbb{N}$,
$$
A(n)=\prod_{p\in\mathbb{P},\,p|n}A(p^{e_p(n)}),
$$
where $e_p(n)$ is the exponent of $p$ in $n$, $\PPP$ is the set of primes, and the product on the right-hand side is defined as in (\ref{RdsMulti}). It follows that $A$ is multiplicative if and only if each $A(n)$ is multiplicative. Recall that the Ramanujan sum $c(k,n)$ was defined in  \eqref{eq:RamSum}.

\begin{thm}
\emph{(\cite[Theorem 4.1]{Le12})}
\label{MdsEner}
Let $A = (A(n))_{n\in\mathbb{N}}$ be a multiplicative divisor system. Then for any integer $n \ge 3$,
$$
\E(\ICG(n, A(n))) = \sum_{k=1}^n|\lambda_k(n, A(n))|=\prod_{p\in\mathbb{P},\,p|n}\E(\ICG(p^{e_p(n)}, A(p^{e_p(n)}))),
$$
where $\lambda_k(n, A(n)):=(1 *_{A}c(k,\cdot))(n)$ with $1$ the constant function taking value $1$.
\end{thm}

Recall from \eqref{eq:Xp} that $X_{p} =  \left\{p^{e_{p}(x)}: x \in  X\right\}$ for any prime $p$ and non-empty set $X$ of positive integers. Using Theorem \ref{MdsEner}, {Le} and {Sander} obtained the following result.

\begin{thm}
\emph{(\cite[Corollary 4.1]{Le12})}
\label{MdsEnPow}
Let $A = (A(n))_{n\in\mathbb{N}}$ be a multiplicative divisor system. Then, for any $n = p_1^{\alpha_1}p_2^{\alpha_2} \ldots p_s^{\alpha_s} \ge 3$ as in \eqref{eq:n}, we have $A(p_j^{\alpha_j}) = A(n)_{p_j}$ for $1 \le  j \le s$, and moreover
$$
\E(\ICG(n, A(n))) = \prod_{j=1}^s \E(\ICG(p_j^{\alpha_j}, A(n)_{p_j})).
$$
\end{thm}

This implies the following result.

\begin{thm}
\emph{(\cite[Theorem 4.2]{Le12})}
Let $n \ge 3$ and $D \subseteq D(n)$. Suppose that $D$ can be factorized as $D = \{g\} \cdot X$ for some positive integer $g$ and multiplicative set $X \subseteq D(n/g)$. Then
$$
\E(\ICG(n,D)) = g \prod_{p \in \PPP,\ p \mid n} \E(\ICG(p^{e_p(n/g)}, X_p)).
$$
\end{thm}

Given $D\subseteq D(n)$, set $\overline{D}=D(n)\setminus D$ and
\be
\label{eq:Phi}
\Phi(n,D)=\sum_{d\in D}\varphi\left(\frac{n}{d}\right).
\ee

\begin{thm}
\emph{(\cite[Theorem 6.1]{Le12-1})}
Let $n \ge 3$ and let $D \subset D(n)$ be a multiplicative set containing $n$. Then
$$
\E(\ICG(n,\overline{D})) = \prod_{p\in\mathbb{P},\,p|n}\E(\ICG(p^{e_p(n)}, D_{p}))+n-2\Phi(n,D).
$$
\end{thm}

As mentioned in \cite[Section 5]{Le12}, Theorem \ref{MdsEnPow} reduces the computation of the energies of all integral circulant graphs with respect to multiplicative divisor sets to that of $\ICG(p^a, D)$ for an arbitrary prime power $p^a$ and any divisor set $D\subseteq D(p^a)$. The latter was achieved in \cite{Sander11},  where the following result was proved.

\begin{thm}
\emph{(\cite[Theorem 2.1; Corollary 2.2]{Sander11})}
\label{simplecase}
Let $p$ be a prime and let $a$ be a positive integer. Let $D = \{p^{a_1},p^{a_2},\ldots,p^{a_r}\}$, where $0 \le a_1< a_2<\cdots< a_r < a$. Then
$$
\E(\ICG(p^a, D))=2(p-1)\left(p^{a-1}r-(p-1)\sum_{k=1}^{r-1}\sum_{i=k+1}^rp^{a-a_i+a_k-1}\right).
$$
Moreover, $\ICG(p^a,D)$ is hyperenergetic if and only if
$$
\sum_{k=1}^{r-1}\sum_{i=k+1}^r\frac{1}{p^{a_i-a_k}}<\frac{1}{p-1}\left(r-\frac{p^a-1}{p^{a-1}(p-1)}\right).
$$
\end{thm}

We will see in Theorem \ref{thm:sun1} that this result can be generalized to the Cayley graph $\Cay(R, D)$ on a finite chain ring $R$ as in Section \ref{subsec:ring-gcd}.

The condition $a_r < a$ in Theorem \ref{simplecase} is required for otherwise $\ICG(p^a,D)$ would have a loop. In the case when $a_r = a$, we have the following result.

\begin{thm}
\emph{(\cite[Proposition 5.1]{Le12})}
\label{loopcase}
Let $p$ be a prime and let $a$ be a positive integer. Let $D = \{p^{a_1},p^{a_2},\ldots,p^{a_r}\}$, where $0\le a_1< a_2<\cdots< a_r = a$. Then
$$
\E(\ICG(p^a,D))=2(p-1)\left(p^{a-1}(r-1)-(p-1)\sum_{k=1}^{r-2}\sum_{i=k+1}^{r-1}p^{a-a_i+a_k-1} -\sum_{k=1}^{r-1}p^{a_k}\right)+p^a.
$$
\end{thm}

Given $D=\{d_1,d_2,\ldots,d_r\}\subseteq D(p^a)$ with $d_1<d_2<\cdots<d_r$, define
$$
\eta(p^a, D):=\left\{\begin{array}{ll}
                      2(p-1)p^{a-1}\left(r-(p-1)\sum\limits_{k=1}^{r-1}\sum\limits_{i=k+1}^{r}\dfrac{d_k}{d_i}\right), & \text{if $d_r<p^a$},\\[0.6cm]
                       2(p-1)p^{a-1}\left(r-1-(p-1)\sum\limits_{k=1}^{r-2}\sum\limits_{i=k+1}^{r-1}\dfrac{d_k}{d_i} -\dfrac{1}{p^{a-1}}\sum\limits_{k=1}^{r-1}d_k\right)+p^a, &\text{if $d_r=p^a$}.
                    \end{array}
\right.
$$

Theorems \ref{MdsEnPow}, \ref{simplecase} and \ref{loopcase} together imply the following result.

\begin{thm}
\emph{(\cite[Theorem 5.1]{Le12})}
Let $n \ge 3$ and $D \subseteq D(n)$. Suppose that $D$ can be factorized into $D =
\{g\}\cdot X$ for some positive integer $g$ and multiplicative set $X \subseteq D(n/g)$. Then
$$
\E(\ICG(n, D))=g\prod_{p\in\mathbb{P},\,p|n}\eta(p^{e_p(n/g)}, X_p).
$$
\end{thm}

In view of Theorem \ref{DisPow1}, all distance powers of the unitary Cayley graph $G_{\mathbb{Z}_n}$ ($= \Cay(\mathbb{Z}_n, \ZZZ_n^{\times}) $ $= \ICG(n, \{1\})$) are integral circulant graphs. The energies of some distance powers of $G_{\mathbb{Z}_n}$ were computed by Liu and Li in \cite{LiuL16} and as a consequence when such a distance power is hyperenergetic was determined.

\begin{thm}\emph{(\cite[Theorem 1.1]{LiuL16})}
Let $n=p_1^{\alpha_1}p_2^{\alpha_2}\cdots p_s^{\alpha_s} \ge 2$, where $p_1<p_2<\cdots<p_s$ are distinct primes and $\alpha_i\ge1$.
\begin{itemize}
\item[\rm (a)] If $n$ is a prime, then $\E(G_{\mathbb{Z}_n}^{\{1\}})=2(n-1)$.
\item[\rm (b)] If $n=2^{\alpha}$ with $\alpha>1$ or $n$ is an odd composite number, then
\begin{align*}
    &\E\left(G_{\mathbb{Z}_n}^{\{1\}}\right)=2^s\varphi(n),\quad \E\left(G_{\mathbb{Z}_n}^{\{2\}}\right)=2n-2+(2^s-2)\varphi(n)-\prod_{i=1}^sp_i+\prod_{i=1}^s(2-p_i), \\
    &\E\left(G_{\mathbb{Z}_n}^{\{1,2\}}\right)=2(n-1).
  \end{align*}
\item[\rm (c)] If $n$ is even but has an odd prime factor, then
\begin{align*}
     &\E\left(G_{\mathbb{Z}_n}^{\{1\}}\right)=2^s\varphi(n),\quad \E\left(G_{\mathbb{Z}_n}^{\{2\}}\right)=2(n-2),\quad \E\left(G_{\mathbb{Z}_n}^{\{3\}}\right) = n-4\varphi(n)+ 2^s\varphi(n),\\
     &\E\left(G_{\mathbb{Z}_n}^{\{1,2\}}\right) = 2n-2+(2^s-2)\varphi(n)-\prod_{i=1}^sp_i, \quad \E\left(G_{\mathbb{Z}_n}^{\{1,3\}}\right)=n, \\
     &\E\left(G_{\mathbb{Z}_n}^{\{2,3\}}\right)=2n-2+(2^s-2)\varphi(n)-\prod_{i=1}^sp_i, \quad \E\left(G_{\mathbb{Z}_n}^{\{1,2,3\}}\right)=2(n-1).
  \end{align*}
\end{itemize}
\end{thm}

\begin{thm}\emph{(\cite[Corollary 2.10]{LiuL16})}
Let $n=p_1^{\alpha_1}p_2^{\alpha_2}\cdots p_s^{\alpha_s} \ge 2$, where $p_1<p_2<\cdots<p_s$ are distinct primes and $\alpha_i\ge1$. Then the following hold:
\begin{itemize}
\item[\rm (a)] $G_{\mathbb{Z}_n}^{\{1\}}$ is hyperenergetic if and only if $s\ge3$ or $s=2$ and $n$ is odd;
\item[\rm (b)] $G_{\mathbb{Z}_n}^{\{2\}}$ is hyperenergetic if and only if $n$ is an odd composite number;
\item[\rm (c)] $G_{\mathbb{Z}_n}^{\{3\}}$ is hyperenergetic if and only if $p_1=2$ and $s\ge3$;
\item[\rm (d)] $G_{\mathbb{Z}_n}^{\{1,2\}}$ is hyperenergetic if and only if $p_1=2$ and $s\ge2$ except when $n=2p_2$;
\item[\rm (e)] $G_{\mathbb{Z}_n}^{\{2,3\}}$ is hyperenergetic if and only if $p_1=2$ and $s\ge2$ except when $n=2p_2$;
\item[\rm (f)] $G_{\mathbb{Z}_n}^{\{1,3\}}$ is not hyperenergetic;
\item[\rm (g)] $G_{\mathbb{Z}_n}^{\{1,2,3\}}$ is not hyperenergetic.
\end{itemize}
\end{thm}

In \cite{ZhouXu15}, hyperenergetic integral circulant graphs were discussed and a method for constructing hyperenergetic integral circulant graphs using the Cartesian product of graphs was given.


\subsection{Extremal energies of integral circulant graphs}

In this section we survey some known results on the characterization of those graphs having maximal or minimal energy among all integral circulant graphs with a given order. We adopt the following notation from \cite{Sander11}:
$$
\E_{\min}(n):=\min\{\E(\ICG(n,D)): D\subseteq D(n)\setminus\{n\}\}
$$
$$
\E_{\max}(n):=\max\{\E(\ICG(n,D)): D\subseteq D(n)\setminus\{n\}\}.
$$
A subset $D\subseteq D(n)\setminus\{n\}$ is called
\emph{$n$-minimal} or \emph{$n$-maximal} if $\E(\ICG(n,D))= \E_{\min}(n)$ or $\E(\ICG(n,D))= \E_{\max}(n)$, respectively.

\begin{thm}
\emph{(\cite[Theorem 3.1]{Sander11})}
Let $p$ be a prime and $s$ a positive integer. Then
$$
\E_{\min}(p^s) = 2(p-1)p^{s-1}.
$$
Moreover, the $p^s$-minimal sets are exactly the sets $\{p^t\}$ for $0 \le t \le s-1$.
\end{thm}

\begin{thm}
\emph{(\cite[Theorem 3.2]{Sander11})}
Let $p$ be a prime. Then the following hold:
\begin{itemize}
\item[\rm (a)] $\E_{\max}(p) = 2(p-1)$, the only $p$-maximal set being $D = \{1\}$;
\item[\rm (b)] $\E_{\max}(p^2) = 2(p-1)(p+1)$, the only $p^2$-maximal set being $D = \{1, p\}$;
\item[\rm (c)] $\E_{\max}(p^3) = 2(p-1)(2p^2-p+1)$, the only $p^3$-maximal set being $D = \{1, p^2\}$, except when $p=2$ for which $D = \{1,2,4\}$ is also $2^3$-maximal;
\item[\rm (d)] $\E_{\max}(p^4) = 2(p-1)(2p^3+1)$, the only $p^4$-maximal sets being $D = \{1, p, p^3\}$ and $D = \{1, p^2, p^3\}$.
\end{itemize}
\end{thm}

Using tools from convex optimization, it was proved in \cite[Theorem 4.2]{Sander11-1} that $\E_{\max}(p^s)$ lies between $s(p-1)p^{s-1}$ and $2s(p-1)p^{s-1}$ approximately.

In \cite[Theorem 2.1]{Sander11}, it was proved that, for $0 \le a_1<a_2<\cdots<a_{r-1}< a_r \le s-1$,
$$
\E(\ICG(p^s, \{p^{a_1},\ldots,p^{a_r}\})) = 2(p-1)p^{s-1}\left(r-(p-1)h_{p}(a_1, \ldots, a_r)\right),
$$
where
$$
h_{p}(a_1, \ldots, a_r) = \sum_{k=1}^{r-1} \sum_{i=k+1}^{r} \frac{1}{p^{a_i - a_k}}.
$$
Thus $\E_{\max}(p^s)$ can be obtained by minimizing $h_{p}(a_1, \ldots, a_r)$. As observed in \cite{Sander13}, the minimum of $h_{p}(a_1, \ldots, a_r)$ occurs only when $a_1 = 0$ and $a_r = s-1$. This approach was used in \cite[Theorems 2.1 and 2.2]{Sander12}, where the minimum value of $h_{p}(0, a_2, \ldots, a_{r-1}, s-1)$ was determined when $s \equiv 1\ \mod\ (r-1)$ or $s \equiv 0\ \mod\ (r-1)$. Finally, in \cite[Theorem 1.1]{Sander13}, the following formula for $\E_{\max}(p^s)$ was obtained and all $p^s$-maximal sets were determined.

\begin{thm}
\emph{(\cite[Theorem 1.1]{Sander13})}
Let $p$ be a prime and let $r$ and $s$ be positive integers.
\begin{itemize}
\item[\rm (a)] If $s$ is odd, then
$$
\E_{\max}(p^s)=\frac{1}{(p+1)^2}\left((s+1)(p^2-1)p^s+2(p^{s+1}-1)\right).
$$
Moreover, if $p\ge3$ then $\{1, p^{2}, p^{4}, \ldots, p^{s-3}, p^{s-1}\}$ is the only $p^s$-maximal subset, and if $p = 2$ then $\{1, p^{2}, p^{4}, \ldots, p^{s-3}, p^{s-1}\}$ and $\{1, p^{1}, p^{3}, p^{5}, \ldots, p^{s-4}, p^{s-2}, p^{s-1}\}$ are the only $p^s$-maximal subsets.
\item[\rm (b)] If  $s$ is even, then
$$
\E_{\max}(p^s)=\frac{1}{(p+1)^2}\left(s(p^2-1)p^s+2(2p^{s+1}-p^{s-1}+p^2-p-1)\right)
$$
and $\{1, p^{2}, p^{4}, \ldots, p^{s-2}, p^{s-1}\}$ and $\{1, p^{1}, p^{3}, p^5, \ldots, p^{s-3}, p^{s-1}\}$ are the only $p^s$-maximal subsets.
\end{itemize}
\end{thm}

In general, it seems difficult to obtain the exact values of $\E_{\min}(n)$ and $\E_{\max}(n)$ for arbitrary integers $n \ge 3$ due to the lack of explicit formulas for $\E(\ICG(n,D))$. The following relaxed forms of $\E_{\min}(n)$ and $\E_{\max}(n)$ were introduced in \cite{Le12-1}:
$$
\widetilde{\E}_{\min}(n):=\min\{\E(\ICG(n,D)): D\subseteq D(n)\setminus\{n\} \text{ multiplicative}\}
$$
$$
\widetilde{\E}_{\max}(n):=\max\{\E(\ICG(n,D)): D\subseteq D(n)\setminus\{n\}  \text{ multiplicative}\}.
$$
Of course,
$$
\E_{\min}(n)\le\widetilde{\E}_{\min}(n)\le\widetilde{\E}_{\max}(n)\le \E_{\max}(n).
$$
For any prime power $p^s$, set
$$
\theta(p^s)=\left\{\begin{array}{ll}
                      \dfrac{1}{(p+1)^2}\left((s+1)(p^2-1)p^s+2(p^{s+1}-1)\right), &   \text{if $2\nmid s$},\\[0.3cm]
                     \dfrac{1}{(p+1)^2}\left(s(p^2-1)p^s+2(2p^{s+1}-p^{s-1}+p^2-p-1)\right),& \text{if $2|s$}.
                   \end{array}\right.
$$

\begin{thm}
\emph{(\cite[Theorem 2.1]{Le12-1})}
\label{MdsPow1}
Let $n=p_1^{\alpha_1}p_2^{\alpha_2} \ldots p_s^{\alpha_s} \ge 3$ be an integer in canonical factorization. Then
$$
\widetilde{\E}_{\max}(n)= \prod_{i=1}^s\theta(p_i^{\alpha_i}).
$$
Moreover, a multiplicative set $D\subseteq D(n)\setminus\{n\}$ satisfies $\E(\ICG(n,D)) = \widetilde{\E}_{\max}(n)$ if and only if $D=\prod_{i=1}^sD^{(i)}$, where
$$
D^{(i)}=\left\{\begin{array}{lll}
                      \{1,p_i^2,p_i^4,\ldots,p_i^{\alpha_i-3},p_i^{\alpha_i-1}\}, &   \text{if $2\nmid \alpha_i$, $p_i\ge3$},\\[0.2cm]
                     \{1,2^2,2^4,\ldots,2^{\alpha_i-3}, 2^{\alpha_i-1}\} \text{ or }  \{1,2,2^3,\ldots,2^{\alpha_i-4}, 2^{\alpha_i-2}, 2^{\alpha_i-1}\}, & \text{if $2\nmid \alpha_i$, $p_i=2$},\\[0.2cm]
                      \{1,p_i^2,p_i^4,\ldots, p_i^{\alpha_i-4}, p_i^{\alpha_i-2}, p_i^{\alpha_i-1}\} \text{ or }  \{1,p_i,p_i^3,\ldots, p_i^{\alpha_i-3}, p_i^{\alpha_i-1}\},&\text{if $2|\alpha_i$}.
                   \end{array}\right.
$$
\end{thm}

A set $D\subseteq D(p^s)$ is called \emph{uni-regular} if $D=\{p^i, p^{i+1},\ldots,p^{j-1},p^j\}$ for some integers $i, j$ with $0\le i \le j \le s$.

\begin{thm}
\emph{(\cite[Theorem 2.2]{Le12-1})}
\label{MdsPow2}
Let $n=p_1^{\alpha_1}p_2^{\alpha_2}\ldots p_s^{\alpha_s} \ge 3$ be an integer in canonical factorization, where $p_1<p_2<\ldots<p_s$. Then
$$
\widetilde{\E}_{\min}(n)= 2n\left(1-\frac{1}{p_1}\right).
$$
Moreover, a multiplicative set $D\subseteq D(n)\setminus\{n\}$ satisfies $\E(\ICG(n,D)) = \widetilde{\E}_{\min}(n)$ if and only if $D=\prod_{i=1}^sD^{(i)}$, where $D^{(1)}=\{p_1^u\}$ for some $u\in \{0,1,\ldots,\alpha_1-1\}$,  and $D^{(i)}$ is an arbitrary uni-regular set with $p_i^{\alpha_i}\in D^{(i)}$ for $2\le i \le s$.
\end{thm}

Theorems \ref{MdsPow1} and \ref{MdsPow2} together imply the following result.

\begin{thm}
\emph{(\cite[Corollary 2.1]{Le12-1})}
Let $n=p_1^{\alpha_1}p_2^{\alpha_2} \ldots p_s^{\alpha_s} \ge 3$ be an integer in canonical factorization, where $p_1<p_2< \ldots < p_s$. Then
\begin{itemize}
\item[\rm (a)]
$\E_{\max}(n) \ge \widetilde{\E}_{\max}(n)= \prod_{i=1}^s\theta(p_i^{\alpha_i})$; and
\item[\rm (b)]
$\E_{\min}(n)\le \widetilde{\E}_{\min}(n)= 2n\left(1-\dfrac{1}{p_1}\right)$.
\end{itemize}
\end{thm}

Denote by $\tau(n)$ the number of positive divisors of $n$ and by $\omega(n)$ the number of distinct prime factors of $n$.

\begin{thm}
\emph{(\cite[Theorem 2.3]{Le12-1})}
Let $n \ge 3$ be an integer. Then
\begin{itemize}
\item[\rm (a)]
$
\E_{\max}(n) \le n\sum\limits_{d|n}\dfrac{\varphi(n)\tau(d)}{d}=n\prod\limits_{p\in \mathbb{P},\, p|n}\left(\dfrac{1}{2}\left(1-\dfrac{1}{p}\right)(e_p(n)+1)(e_p(n)+2)+\dfrac{1}{p}\right)
$;
\item[\rm (b)]
$
\E_{\max}(n) < \left(\dfrac{3}{4}\right)^{\omega(n)}n\tau(n)^2$; and
\item[\rm (c)]
$
\E_{\max}(n) \le \widetilde{\E}_{\max}(n)\tau(n)$.
\end{itemize}
\end{thm}

We finish this section by mentioning three conjectures from \cite{Le12-1}. As far as we know, all these conjectures are still open.

\begin{conj}
\label{conj:61}
\emph{(\cite[Conjecture 6.1]{Le12-1})}
For any integer $n \ge 3$,
$$
\E_{\min}(n) = 2n\left(1-\dfrac{1}{p}\right),
 $$
where $p$ is the smallest prime factor of $n$.
\end{conj}

\begin{conj}
\label{conj:62}
\emph{(\cite[Conjecture 6.2]{Le12-1})}
For any integer $n \ge 3$ and $D \subseteq D(n)\setminus\{n\}$, if $\E(\ICG(n,D))=\E_{\min}(n)$, then $D$ is a multiplicative divisor set.
\end{conj}

\begin{conj}
\label{conj:63}
\emph{(\cite[Conjecture 6.3]{Le12-1})}
Let $n \ge 3$ be an integer, and let $D_1, D_2\subseteq D(n)\setminus\{n\}$ be multiplicative sets. If $\ICG(n,D_1)$ and $\ICG(n,D_2)$ are cospectral, then $D_1=D_2$.
\end{conj}


\subsection{Energies of circulant graphs}

In \cite{Shparlinski06}, Shparlinski gave a construction of circulant graphs of very high energy using Gauss sums. Given a prime $p$, let $Q_p$ denote the set of all quadratic residues modulo $p$ in the set $\{1, \ldots, (p-1)/2\}$.

\begin{thm}
\emph{(\cite[Theorem 1]{Shparlinski06})}
For any prime $p \equiv 1$ ($\mod~4$),
$$
\E(\Cay(\ZZZ_p, Q_p \cup (-Q_p))) \ge \frac{(p-1)(\sqrt{p}+1)}{2}.
$$
\end{thm}

Thus $\Cay(\ZZZ_p, Q_p \cup (-Q_p))$ has energy close to the upper bound in \eqref{eq:KoolenM01}.

Let $n$ and $d$ be integers with $1 \le d \le n-1$. The average energy of all circulant graphs of order $n$ and degree $d$ is given by
$$
\E(n,d) = \frac{1}{\choose{\lc n/2 \rc - 1}{\lf d/2 \rf}} \sum \E(\Cay(\ZZZ_n, S)),
$$
where the sum is running over all subsets $S$ of $\ZZZ_n \setminus \{0\}$ with $S = -S$ and $|S| = d$. Using \eqref{eq:KoolenM01}, it can be verified \cite{BlackburnS08} that
$$
\E(n,d) \le d + \sqrt{d(n-1)(n-d)}.
$$
On the other hand, we have the following asymptotic lower bound due to Blackburn and Shparlinski \cite{BlackburnS08}.

\begin{thm}
\emph{(\cite[Theorem 5]{BlackburnS08})}
As $n \rightarrow \infty$, for any integer $d$ such that $4 \le d = o(n^{1/2})$ and $dn$ is even, we have
$$
\E(n,d) \ge
\begin{cases}
\frac{dn}{\sqrt{3d-3+o(1)}}, & \text{ if $d$ is even}\\
\frac{dn}{\sqrt{3d-3+(1/d)+o(1)}}, & \text{ if $d$ is odd.}
\end{cases}
$$
\end{thm}


\subsection{Energies of Cayley graphs on finite rings}
\label{subsec:en-uni-qua}

As seen in Section \ref{Sec-ener}, the energies of the unitary Cayley graph $G_{\mathbb{Z}_n}$ ($= \Cay(\mathbb{Z}_n, \ZZZ_n^{\times}) = \ICG(n, \{1\})$) and its complement have been determined. These results were generalized by Kiani, Aghaei, Meemark and Suntornpoch \cite{Kiani11} to the unitary Cayley graph of any finite commutative ring.

\begin{thm}
\emph{(\cite[Theorems 2.4, 2.5 and 4.1]{Kiani11})}
Let $R$ be a finite commutative ring as in Assumption \ref{as:1}. Then
$$
\E(G_R)=2^s|R^{\times}|
$$
and
$$
\E(\overline{G}_R)=2|R|-2+(2^s-2)|R^{\times}|-\prod_{i=1}^s|R_{i}|/m_{i} + \prod_{i=1}^s(2-|R_{i}|/m_{i}).
$$
Moreover,
\begin{itemize}
\item[\rm (a)] if $s=1$, then $G_R$ is not hyperenergetic;

\item[\rm (b)] if $s = 2$, then $G_R$ is hyperenergetic if and only if $|R_1|/m_1 \ge 3$ and $|R_2|/m_2 \ge 4$;

\item[\rm (c)] if $s \ge 3$, then $G_R$ is hyperenergetic if and only if $|R_{s-2}|/m_{s-2} \ge 3$, or $|R_{s-1}|/m_{s-1} \ge 3$ and $|R_{s}|/m_{s} \ge 4$.
\end{itemize}
\end{thm}

In \cite[Theorem 18]{Liu12},  {Liu} and {Zhou} determined the energy of the line graph $L(G_{R})$ of $G_R$.

\begin{thm}
\emph{(\cite[Theorem 18]{Liu12})}
\label{thmE}
Let $R$ be a finite commutative ring as in Assumption \ref{as:1}. Then
$$
\E(L(G_{R})) =
\left\{\begin{array}{ll}
2^{s+1}(|R^\times|-1)^2,  & \text{if $2=|R_1|/m_1=\cdots=|R_s|/m_s$,}\\
                          & \text{or $R=\underbrace{\mathbb{F}_2\times\cdots\times\mathbb{F}_2}_{s-1}\times\mathbb{F}_3$;}\\ [0.2cm]
2^{t+1} +2|R|\left(|R^\times|-2\right), & \text{if $2=|R_1|/m_1=\cdots=|R_t|/m_t<|R_{t+1}|/m_{t+1}$}\\
                  & \text{with $1\leq t<s$ and $R\ncong\underbrace{\mathbb{F}_2\times\cdots\times\mathbb{F}_2}_{s-1}\times\mathbb{F}_3$};\\ [0.2cm]
2|R|\left(|R^\times|-2\right), & \text{if $3\leq|R_1|/m_1\leq\cdots\leq|R_s|/m_s$ and $R\ncong\mathbb{F}_3$}.\\
                                \end{array}\right.
$$
\end{thm}

In the special case where $R = \mathbb{Z}_n$, Theorem \ref{thmE} yields the following result.

\begin{cor}
\emph{(\cite[Corollary 19]{Liu12})}
\label{corTE}
Let $n=p_1^{\alpha_1}p_2^{\alpha_2}\ldots p_s^{\alpha_s} \ge 3$ be an integer in canonical factorization, where $p_1<p_2<\ldots<p_s$ are primes and each $\a_i \ge 1$ is an integer. Then
$$
\E(L(G_{\mathbb{Z}_n})) =
\left\{\begin{array}{lcl}
           4,               & & \text{if $n=3$;}\\ [0.2cm]
           8,               & & \text{if $n=6$;}\\ [0.2cm]
4\left(2^{\alpha_1-1}-1\right)^2,  & &   \text{if $n=2^{\alpha_1}$;}\\[0.2cm]
4 +2n\left(\left(\prod_{i=1}^sp_i^{\alpha_1-1}(p_i-1)\right)-2\right), & & \text{if $2=p_1$ and $n\ne6$};\\ [0.2cm]
2n\left(\left(\prod_{i=1}^sp_i^{\alpha_1-1}(p_i-1)\right)-2\right), & & \text{if $3\le p_1$ and $n\neq3$}.\\
                                \end{array}\right.
$$
\end{cor}

By Theorem \ref{thmE}, we know exactly when $L(G_{R})$ is hyperenergetic.

\begin{cor}
\emph{(\cite[Corollary 20]{Liu12})}
Let $R$ be as in Assumption \ref{as:1}. Then $L(G_{R})$ is hyperenergetic if and only if one of the following holds:
\begin{itemize}
\item[\rm (a)] $|R^\times| \ge 4$;
\item[\rm (b)] $s = 1$ and $|R| = 2m \ge 8$;
\item[\rm (c)] $s \ge 2$, $2=|R_1|/m_1=\cdots=|R_s|/m_s$, and $|R^\times| \ge 2$.
\end{itemize}
\end{cor}

Recall that $M_{n}(R)$ is the ring of $n \times n$ matrices over a finite ring $R$, and $G_{M_{n}(R)} = \Cay(M_{n}(R), \GL_{n}(R))$ is the unitary Cayley graph of $M_{n}(R)$. In \cite{RattaM20}, Rattanakangwanwong and Meemark obtained the following result.

\begin{thm}
\emph{(\cite[Theorem 4.7]{RattaM20})}
For any finite commutative ring $R$ and any integer $n \ge 2$, $G_{M_{n}(R)}$ is hyperenergetic.
\end{thm}

Recall from Section \ref{subsec:QUCayComRing} that $\mathcal{G}_{R}$ is the quadratic unitary Cayley graph of a finite commutative ring $R$. The next two results were obtained by Liu and Zhou in \cite{Liu15}.

\begin{thm}
\emph{(\cite[Theorem 3.1]{Liu15})}
\label{EnergyLQUCG}
Let $R$ be a local ring with maximal ideal $M$ of order $m$. Then the following hold:
\begin{itemize}
\item[\rm (a)] if $|R|/m\equiv 1\,(\mod\,4)$, then $\E(\mathcal{G}_{R})=\left(\sqrt{|R|/m}+1\right)|R^\times|\big/2$;
\item[\rm (b)] if $|R|/m\equiv 3\,(\mod\,4)$, then $\E(\mathcal{G}_{R})=2|R^\times|$.
\end{itemize}
\end{thm}

\begin{thm}
\emph{(\cite[Theorem 3.2]{Liu15})}
\label{Energy13mod4}
Let $R$ be as in Assumption \ref{as:1} such that $|R_i|/m_i\equiv 1\,(\mod\,4)$ for $1 \le i \le s$, and let $R_0$ be a local ring with maximal ideal $M_0$ of order $m_0$ such that $|R_0|/m_0\equiv 3\,(\mod\,4)$. Then the following hold:
\begin{itemize}
\item[\rm (a)] $\E(\mathcal{G}_{R})=\dfrac{|R^\times|}{2^s}\prod\limits_{i=1}^s\left(\sqrt{|R_i|/m_i}+1\right)$;
\item[\rm (b)] $\E(\mathcal{G}_{R_0\times R})=\dfrac{|R_0^\times||R^\times|}{2^{s-1}}\prod\limits_{i=1}^s\left(\sqrt{|R_i|/m_i}+1\right)$.
\end{itemize}
\end{thm}

Theorem \ref{Energy13mod4} implies the following result.

\begin{cor}
\emph{(\cite[Corollary 3.3]{Liu15})}
\label{HyEnerCor1}
Let $R$ and $R_0$ be as in Theorem \ref{Energy13mod4}. Then the following hold:
\begin{itemize}
\item[\rm (a)]  $\mathcal{G}_{R}$ is hyperenergetic except when $R=R_1$ with $|R_1|/m_1=5$ or $R=R_1\times R_2$ with $|R_1|/m_1=|R_2|/m_2=5$;
\item[\rm (b)]  $\mathcal{G}_{R_0\times R}$ is hyperenergetic except when $|R_0|/m_0=3$ and $R=R_1$ with $|R_1|/m_1=5$.
\end{itemize}
\end{cor}

In the special case where $R = \mathbb{Z}_n$, Theorems \ref{EnergyLQUCG} and \ref{Energy13mod4} and Corollary \ref{HyEnerCor1} together imply the following result.

\begin{cor}
\emph{(\cite[Corollary 3.4]{Liu15})}
Let $n=p_1^{\alpha_1} p_2^{\alpha_2} \ldots p_s^{\alpha_s}$ be an integer in canonical factorization such that $p_i \equiv1\,(\mod\,4)$ for $1 \le i \le s$. Let $p\equiv3\,(\mod\,4)$ be a prime and $\alpha\ge1$ an integer. Then the followng hold:
\begin{itemize}
\item[\rm (a)]  $\E(\mathcal {G}_{\mathbb{Z}_{p^\alpha}})=2\varphi(p^\alpha)$;
\item[\rm (b)]  $\E(\mathcal {G}_{\mathbb{Z}_n}) =\dfrac{\varphi(n)}{2^s}\prod\limits_{i=1}^s\left(\sqrt{p_i}+1\right)$;
\item[\rm (c)]  $\E(\mathcal {G}_{\mathbb{Z}_{np^\alpha}}) =\dfrac{\varphi(n)\varphi(p^\alpha)}{2^{s-1}}\prod\limits_{i=1}^s\left(\sqrt{p_i}+1\right)$;
\item[\rm (d)]  $\mathcal {G}_{\mathbb{Z}_{n}}$ is hyperenergetic except when $n=5^{\alpha_1}$;
\item[\rm (e)]  $\mathcal {G}_{\mathbb{Z}_{np^\alpha}}$ is hyperenergetic except when $np^\alpha=3^\alpha\cdot5^{\alpha_1}$.
\end{itemize}
\end{cor}

Finally, for any semiprimitive pair $(k, q)$ as described in Theorem \ref{PV20b}, the energies of the semiprimitive generalized Paley graph $X_{q}^k$ and its complement $\overline{X_{q}^k}$ can be easily computed (see \cite{PodestaV20, PodestaV20a}): If $r/2t$ is odd, then
$$
\E(X_{q}^k) = \E(\overline{X_{q}^k}) = \frac{2(k-1)(q-1)(\sqrt{q} + 1)}{k^2};
$$
if $r/2t$ is even, then
$$
\E(X_{q}^k) = \frac{2(q-1)\left((k-1)\sqrt{q} + 1\right)}{k^2}
$$
and
$$
\E(\overline{X_{q}^k}) = \frac{2(q-1)\left((k-1)\sqrt{q} + (k - 1)^2\right)}{k^2}.
$$


\subsection{Energies of Cayley graphs on finite chain rings and gcd graphs of unique factorization domains}
\label{subsec:eng-ring-gcd}

The following result due to Suntornpoch and Meemark \cite{SuntornpochM16} gives a formula for the energy of the Cayley graph $\Cay(R, D)$ on a finite chain ring $R$ (see Section \ref{subsec:ring-gcd} for related definitions). In the special case when $R = \ZZZ_{p^s}$, it gives exactly Theorem \ref{simplecase}.

\begin{thm}
\label{thm:sun1}
\emph{(\cite[Section 2]{SuntornpochM16})}
Let $R$ be a finite chain ring and $D$ be as in \eqref{eq:calC}. Then
$$
\EE(\Cay(R, D)) = 2(q-1)\left(q^{s-1}r - (q-1)\sum_{k=1}^{r-1} \sum_{i=k+1}^{r} q^{s - a_i + a_k - 1}\right).
$$
\end{thm}

Let $R$ be a unique factorization domain (UFD). Let $c$ be a nonzero nonunit element of $R$. Then the quotient $R/(c) = \{x + (c): x \in R\}$ is a commutative ring. Assume that this ring is finite. Let $D$ be a set of proper divisors of $c$. In \cite{Kiani11}, the \emph{gcd graph} of $R$ with respect to $D$ was defined as the graph with vertex set $R/(c)$ such that $x + (c)$ and $y + (c)$ are adjacent if and only if $\gcd(x-y, c) \in D$ (the gcd considered here is unique up to associate). We denote this graph by $\ICG(R/(c), D)$. (Note that another notation was used in \cite{Kiani11}.) This notion is a generalization of the gcd graphs in Section \ref{CirUgcd}: If $R = \mathbb{Z}$, $n \ge 3$ is an integer and $D$ is a set of positive proper divisors of $n$, then $\ICG(R/(n), D)$ is exactly the gcd graph $\ICG(n, D)$ of $\ZZZ_n$. Moreover, in the special case when $D = \{1\}$ (where $1$ is the multiplicative identity of $R$), $\ICG(R/(c), \{1\}) = G_{R/(c)}$ is the unitary Cayley graph of the finite commutative ring $R/(c)$. In \cite{Kiani11}, Kiani, Aghaei, Meemark and Suntornpoch obtained the following result.

\begin{thm}
\emph{(\cite[Theorems 3.1 and 3.4]{Kiani11})}
Let $R$ be a unique factorization domain.
\begin{itemize}
\item[\rm (a)] Let $c = p_{1}^{a_1} \ldots p_{n}^{a_n} \in R$ be factorized into a product of irreducible elements. Assume that $R/(c)$ is finite. For $1 \le i \le n$, if $a_i = 1$, then
$$
\E(\ICG(R/(c), \{1, p_i\})) = 2^{n-1} |R/(p_i)| |R/(c/p_i)^{\times}|.
$$
\item[\rm (b)] Let $c = p_{1} \ldots p_{k} p_{k+1}^{a_{k+1}} \ldots p_{n}^{a_n} \in R$ be factorized into a product of irreducible elements, where $a_i > 1$ for $k+1 \le i \le n$. Assume that $R/(c)$ is finite. Then
$$
\E(\ICG(R/(c), \{p_i, p_j\})) = 2^{n-1} |R/(c)^{\times}|
$$
for $1 \le i < j \le k$.
\end{itemize}
\end{thm}

The following result was proved by Suntornpoch and Meemark with the help of the NEPS construction.

\begin{thm}
\label{thm:sun2}
\emph{(\cite[Theorem 3.4]{SuntornpochM16})}
Let $R$ be a unique factorization domain and $c = p_1^{s_1} \ldots p_k^{s_k}$ a nonzero nonunit element of $R$ which is factorized as a product of irreducible elements. Assume that $R/(c)$ is finite. Assume further that for some $l$ with $1 \le l \le k$ and each $i$ with $1 \le i \le l$ there exists a set $D_i =\{p_{i}^{a_{i1}}, p_{i}^{a_{i2}}, \ldots, p_{i}^{a_{ir_i}}\}$ such that $0 \le a_{i1} < a_{i2} < \cdots < a_{ir_i} \le s_i - 1$. Set
$$
D = \{p_1^{a_{1t_1}} \cdots p_l^{a_{l t_l}}p_{l+1}^{s_{l+1}} \cdots p_k^{s_k}: t_i \in \{1, 2, \ldots, r_i\} \text{ for } 1 \le i \le l\}.
$$
Then
$$
\EE(\ICG(R/(c), D)) = \EE(\ICG(R/(p_1^{s_1}), D_1) \cdots \EE(\ICG(R/(p_l^{s_l}), D_l)  \prod_{j=l+1}^{k} |R/(p_j^{s_j})|.
$$
\end{thm}


\subsection{Skew energy of orientations of hypercubes}

Let $G$ be a digraph of order $n$ with no loops or parallel arcs such that at most one of $(u, v)$ and $(v, u)$ occurs as an arc for each pair of distinct vertices $u, v \in V(G)$. The \emph{skew-adjacency matrix} \cite{CaversCFGHKMT12} of $G$ is the $n \times n$ matrix $A_S(G)$ whose $(u, v)$-entry is $1$ if $(u, v)$ is an arc of $G$, $-1$ if $(v, u)$ is an arc of $G$, and $0$ if there is no arc between $u$ and $v$.  Since $A_S(G)$ is skew-symmetric, its eigenvalues are all purely imaginary numbers, and the collection of them with multiplicities is called the \emph{skew spectrum} of $G$. The \emph{skew energy} $\E_{S}(G)$ of $G$ is defined as the sum of the moduli of the eigenvalues of $G$. It is known that $\E_{S}(G) \le n \sqrt{\Delta}$, where $\Delta$ is the maximum degree of the underlying graph of $G$. In particular, the skew energy of every orientation of any $k$-regular graph of order $n$ is at most $n \sqrt{k}$. In \cite{Tian11}, it was shown that for any $d \ge 1$ there exists an orientation of the hypercube $H(d, 2)$ which achieves this bound, that is, with skew energy $2^d \sqrt{d}$. It was also proved \cite{Tian11} that there exists an orientation of $H(d,2)$ whose skew spectrum is obtained from the spectrum of $H(d,2)$ by multiplying each eigenvalue by the imaginary unit $\mathrm{i}$. Both orientations were given algorithmically in \cite{Tian11}.


\subsection{Others}

In \cite{Ghorbani15a}, estimates of the energy and Estrada index of Cayley graphs $\Cay(\Ga, S)$ were obtained, where $\Ga$ is $D_{2n}$ or $U_{6n}$ (see Section \ref{subsec:order6n}) and $S$ is a normal generating set of $\Ga$.

A proper colouring of a graph $G$ is an assignment of colours to its vertices such that adjacent vertices receive distinct colours. Let $G$ be a graph and $c$ a proper colouring of $G$. Define $A_{c}(G)$ to be the matrix with rows and columns indexed by the vertices of $G$ such that the $(u,v)$-entry is equal to $1$ if $u$ and $v$ are adjacent (so that $c(u) \neq c(v)$), $-1$ if $u$ and $v$ are non-adjacent and $c(u) = c(v)$, and $0$ if $u=v$ or $u$ and $v$ are non-adjacent and $c(u) \neq c(v)$. The eigenvalues of $A_c(G)$ are called the \emph{colour eigenvalues} of $(G, c)$. The \emph{colour energy} $E_c(G)$ of $G$ with respect to $c$ is the sum of the absolute values of the colour eigenvalues of $(G, c)$. In \cite{AdigaSS14}, Adiga, Sampathkumar and Sriraj derived formulas for the colour energies of the unitary Cayley graph $\Cay(\ZZZ_n, \ZZZ_n^{\times})$ and its complement. They also obtained explicit formulas for the colour energies of the gcd graphs $\ICG(n, D)$ of $\ZZZ_n$ in the case when $n=p_1^{\alpha_1}p_2^{\alpha_2}\ldots p_{i-1}^{\alpha_{i-1}}p_i^1$ $p_{i+1}^{\alpha_{i+1}} \ldots p_s^{\alpha_s}$ and $D=\{1, p_i\}$ or $n=p_1p_2\ldots p_s$ is a square-free integer and $D=\{p_i, p_j\}$.


\section{Ramanujan Cayley graphs}
\label{sec:Ram}

The \emph{isoperimetric number} (or Cheeger constant) of a graph $G$ is defined as
$$
h(G) = \min \left\{\frac{|\partial(S)|}{|S|}: S \subset V(G),\ 0 < |S| \le \frac{1}{2} |V(G)|\right\},
$$
where $\partial(S)$ is the set of edges of $G$ with one end-vertex in $S$ and the other end-vertex in $V(G) \setminus S$. This important parameter measures how well the graph expands and how strong the graph is connected in some sense. An infinite sequence of finite regular graphs with fixed degree but increasing orders is called a family of expanders if the isoperimetric numbers of such graphs are all bounded from below by a common positive constant. A well-known result of Alon-Milman \cite{Alon86, Alon85} and Dodziuk \cite{Dodziuk84} asserts that, for any connected $k$-regular graph $G$,
$$
\frac{k-\l_{2}}{2} \le h(G) \le \sqrt{2k(k-\l_{2})},
$$
where $\l_2$ is the second largest eigenvalue of $G$. (A better upper bound due to Mohar \cite{Mohar89} says that $h(G) \le \sqrt{k^2 - \l_{2}^2}$.) So we may also define expander families using the \emph{spectral gap} $k-\l_{2}$: A sequence of graphs $G_{1}, G_2, G_3, \ldots$ is called a family of \emph{$\varepsilon$-expander graphs} \cite{Hoory06}, where $\varepsilon> 0$ is a fixed constant, if (i) all these graphs are $k$-regular for a fixed integer $k \ge 3$; (ii) $k-\lambda_2(G_i) \ge \varepsilon$ for each $i$; and (iii) $|V(G_i)| \rightarrow \infty$ as $i \rightarrow \infty$.

The well known Alon-Boppana bound (see, for example, \cite[Theorem 0.8.8]{Davidoff03}) asserts that, for any family of finite connected $k$-regular graphs $\{G_i\}_{i \ge 1}$ such that $|V(G_i)| \rightarrow \infty$ as $i \rightarrow \infty$, we have $\liminf_{i \rightarrow \infty} \lambda(G_i) \ge 2\sqrt{k-1}$, where $\lambda(G_i)$ is the maximum in absolute value of an eigenvalue of $G_i$ other than $\pm k$. This motivated the following definition: A finite $k$-regular graph $G$ is called \emph{Ramanujan} \cite{Hoory06, Murty03} if $\lambda(G) \leq 2\sqrt{k-1}$.

Over the years a significant amount of work has been done in constructing expander families as well as Ramanujan graphs of a fixed degree, and Cayley graphs play an important role in this area of research. The reader is referred to three survey papers \cite{Murty03, Hoory06, Lubotzky12} and two books \cite{Davidoff03, Krebs11} on expanders and Ramanujan graphs.

In this section we give an account of known results on the characterization of Ramanujan Cayley graphs. Note that, since abelian groups can never yield expander families (see \cite{Cioa06, FMT06}), all Ramanujan Cayley graphs on abelian groups surveyed in this section cannot form expander families.


\subsection{Ramanujan integral circulant graphs}
\label{sec:RamCayAbelian}

In \cite{Droll10}, Droll gave the following characterization of Ramanujan unitary Cayley graphs of $\ZZZ_n$, $n \ge 2$.

\begin{thm}
\label{XnCor1}
\emph{(\cite[Theorem 1.2]{Droll10})}
Let $n=p_1^{\alpha_1}p_2^{\alpha_2}\ldots p_s^{\alpha_s} \ge 2$ with $p_1<p_2< \ldots <p_s$ be an integer in canonical factorization. The unitary Cayley graph $G_{\mathbb{Z}_n}$ is Ramanujan if and only if one of the following holds:
\begin{itemize}
\item[\rm (a)] $n=2^{\alpha_1}$  with $\alpha_1\geq1$;
\item[\rm (b)] $n=p_1^{\alpha_1}$ with $p_1$ odd and $\alpha_1=1,2$;
\item[\rm (c)] $n=4p_2p_3$ with $p_2 < p_3\leq 2p_2-3$;
\item[\rm (d)] $n=p_1p_2$ with $3 \le p_1 < p_2\leq 4p_1-5$, or $n=2p_2p_3$ with $3 \le p_2 < p_3\leq 4p_2-5$;
\item[\rm (e)] $n=2p_2^2,~4p_2^2$ with $p_2$ odd, or $n=2^{\alpha_1}p_2$ with $p_2>2^{\alpha_1-3}+1$.
\end{itemize}
\end{thm}

Since $G_{\mathbb{Z}_n} = \Cay(\ZZZ_n, \ZZZ_n^{\times}) = \ICG(n,\{1\})$ is a special unitary Cayley graph as well as a special integral circulant graph, it is natural to ask the following questions: Besides the graphs in Theorem \ref{XnCor1}, which unitary Cayley graphs of finite commutative rings are Ramanujan? And which integral circulant graphs $\ICG(n, D)$ are Ramanujan? The former question was answered by Liu and Zhou in \cite{Liu12}, and their result will be presented in Theorem \ref{XnThm1}. The latter question was studied by {Le} and {Sander} in \cite{Le13} when $n$ is a prime power but $D$ is arbitrary, and by {Sander} in \cite{Sander15} when $D$ is multiplicative but $n$ is arbitrary. Set
$$
D(n;m):=\{d\in D(n): d\ge m\}
$$
for any positive integers $m, n$ with $m \le n$.

\begin{thm}
\emph{(\cite[Theorem 1.1]{Le13})}
\label{PowRam1}
Let $p^s \ge 3$ be a prime power and $D\subseteq D(p^s)\setminus \{p^s\}$. Then $\ICG(p^s,D)$ is Ramanujan if and only if one of the following holds:
\begin{itemize}
\item[\rm (a)] $D=D\left(p^{\lceil\frac{s}{2}\rceil-1}\right)\cup D'$ for some $D'\subseteq D\left(p^{s-1};p^{\lceil\frac{s}{2}\rceil}\right)$;
\item[\rm (b)]  $D=\{1\}$ when $p=2$ and $s\ge3$;
\item[\rm (c)]  $D=D\left(p^{\frac{s-3}{2}}\right)\cup D'$ such that $|D|\ge2$ for some $D'\subseteq D\left(p^{s-1};p^{\frac{s+1}{2}}\right)$ when $p\in\{2,3\}$ and $s\ge3$ is odd;
\item[\rm (d)] $D=D\left(2^{\frac{s-4}{2}}\right)\cup D'$ for some $D'\subseteq D\left(2^{s-1};2^{\frac{s}{2}}\right)$ satisfying $\emptyset\neq D'\neq \left\{2^{s-1}\right\}$ when $p=2$ and $s\ge4$ is even;
\item[\rm (e)]  $D=\left\{1,2^2,2^3,2^4\right\}$ when $p=2$ and $s=5$;
\item[\rm (f)] $D=D\left(5^{\frac{s-3}{2}}\right)\cup\left\{5^{\frac{s+1}{2}}\right\}\cup D'$ for some $D'\subseteq D\left(5^{s-1};5^{\frac{s+3}{2}}\right)$ when $p=5$ and $s\ge5$ is odd;
\item[\rm (g)] $D=D\left(2^{\frac{s-5}{2}}\right)\cup\left\{2^{\frac{s-1}{2}}\right\}\cup D'$ for some $D'\subseteq D\left(2^{s-1};2^{\frac{s+1}{2}}\right)$ satisfying
\[3-2\sqrt{2}+\frac{1}{2^{\frac{s-3}{2}}}\le 2^{\frac{s-1}{2}}\sum_{d'\in D'}\frac{1}{d'}\] when $p=2$ and $s\ge5$ is odd.
\end{itemize}
\end{thm}

Part (a) of Theorem \ref{PowRam1} implies the following result.

\begin{cor}
\emph{(\cite[Corollary 1.1]{Le13})}
Let $p^s \ge 3$ be a prime power and $D = \{1,p,\ldots,p^{r-1}\}$, where $s/2\le r \le s$. Then $\ICG(p^s,D)$ is Ramanujan. In particular, there is a Ramanujan integral circulant graph $\ICG(p^s,D)$ for every prime power $p^s \ge 3$.
\end{cor}

The next four results were proved by Sander in \cite{Sander15}.

\begin{thm}\emph{(\cite[Theorem 2.1]{Sander15})}
$\ICG(n,D(n)\setminus\{n\})$ is Ramanujan for every integer $n \ge 3$.
\end{thm}

Recall from \eqref{eq:Phi} that $\Phi(n,D) = \sum_{d\in D}\varphi\left(\frac{n}{d}\right)$ for $D\subseteq D(n)$.

\begin{thm}
\emph{(\cite[Theorem 2.2]{Sander15})}
Let $n \ge 3$ be an integer and let $D\subseteq D(n)\setminus\{n\}$ be multiplicative. Write $D=\prod_{p\in\mathbb{P}, p|n} D^*(p)$ for suitable $D^*(p)\subseteq D(p^{e_p(n)})$. Then $\ICG(n,D)$ is Ramanujan if and only if $1\in D$ and
$$
\max_{p|n, p\in\mathbb{P}}^*\frac{p^{e_p(n)}}{\Phi(p^{e_p(n)}, D^*(p))}\le 1+ \frac{2\sqrt{\Phi(n,D)-1}}{\Phi(n,D)},
$$
where the operator $\max\limits^*$ is defined by
$$
\max_{p|n, p\in\mathbb{P}}^*:=\left\{\begin{array}{cl}
                                                \max\limits_{p|n, p\in\mathbb{P}}, & \text{if $4|n$ and $D^*(2) \neq \{1\}$,} \\[0.3cm]
                                                \max\limits_{p|n, p\in\mathbb{P},p>2},& \text{otherwise.}
                                              \end{array}
\right.
$$
\end{thm}

\begin{thm}
\emph{(\cite[Theorem 2.3]{Sander15})}
For every even integer $n \ge 4$, there exists a multiplicative divisor set $D\subseteq D(n)\setminus\{n\}$ such that $\ICG(n, D)$ is a Ramanujan graph.
\end{thm}

\begin{thm}
\emph{(\cite[Theorem 2.4]{Sander15})}
Let $n \ge3$ be an odd integer.
\begin{itemize}
\item[\rm (a)]
If $n$ satisfies $\max\limits_{p\in \mathbb{P},\ p|n}p^{e_p(n)} \ge\dfrac{n^{3/2}+n}{2(n-1)}$, then there is a multiplicative divisor set $D\subseteq D(n)\setminus\{n\}$ such that $\ICG(n, D)$ is a Ramanujan graph.
\item[\rm (b)] If $n>8295$ satisfies $\max\limits_{p\in \mathbb{P},\ p|n}p^{e_p(n)} < \dfrac{n^{3/2}+n}{2(n-1)}$, then there is no multiplicative divisor set $D\subseteq D(n)\setminus\{n\}$ such that $\ICG(n, D)$ is a Ramanujan graph.
\end{itemize}
\end{thm}

Recently, Liu and Li \cite{LiuL16} considered some distance powers of the unitary Cayley graph $G_{\mathbb{Z}_n}$ and determined when they are Ramanujan.

\begin{thm}\emph{(\cite[Theorem 1.2]{LiuL16})}
Let $n \ge 2$ be an integer.
\begin{itemize}
\item[\rm (a)] $G_{\mathbb{Z}_n}^{\{1\}}$ is Ramanujan if and only if one of the following holds:
(i) $n$ is a prime; (ii) $n=2^{\alpha}$ for some $\alpha \ge 1$; (iii) $n=p^{2}$ for an odd prime $p$; (iv) $n=pq$ for some primes $p, q$ with $3 \le p < q \le 4p-5$; (v) $n=2pq$ for some primes $p, q$ with $3 \le p < q \le 4p-5$; (vi) $n=4pq$ for some primes $p, q$ with $p < q \le 2p-3$; (vii) $n=2p^2$ or $4p^2$ for an odd prime $p$, or $n=2^{\alpha}p$ for some integer $\alpha$ and prime $p > 2^{\alpha-3}+1$.
\item[\rm (b)] $G_{\mathbb{Z}_n}^{\{2\}}$ is Ramanujan if and only if one of the following holds:
(i) $n=2^{\alpha}$ for some $\alpha \ge 1$; (ii) $n=p^{\alpha}$ for some odd prime $p$ and integer $\alpha \ge 2$; (iii) $n = 15$, $21$ or $35$; (iv) $n$ is even but has an odd prime factor.
\item[\rm (c)] $G_{\mathbb{Z}_n}^{\{3\}}$ is Ramanujan if and only if $n=30$, $70$ or $4p$ for a prime $p$.
\item[\rm (d)] $G_{\mathbb{Z}_n}^{\{1,2\}}$ is Ramanujan if and only if one of the following holds: (i) $n=2^{\alpha}$ for some $\alpha \ge 1$; (ii) $n=p_1^{\alpha_1}p_2^{\alpha_2}\cdots p_s^{\alpha_s}$ for $s\ge2$ odd primes $p_1, p_2, \ldots, p_s$ and positive integers $\alpha_1, \alpha_2, \ldots, \alpha_s$; (iii) $n=2^{\alpha}p$ for some $\alpha$ and odd prime $p \ge 2^{\alpha-4}+1$; (iv) $n=2p^2$, $4p^2$ or $8p^2$ for some odd prime $p$; (v) $n=54$, $250$ or $686$; (vi) $n=2^{\alpha_1}p_2^{\alpha_2}\cdots p_s^{\alpha_s}$ with $s \ge 3$ for odd primes $p_2, \ldots, p_s$ and positive integers $\alpha_1, \alpha_2, \ldots, \alpha_s$ such that $\varphi(n)\ge n/2+1-2\sqrt{n-2}$.
\item[\rm (e)] $G_{\mathbb{Z}_n}^{\{2,3\}}$ is Ramanujan if and only if  $n=6$, $10$, $12$, $18$, $24$ or $30$.
\item[\rm (f)] $G_{\mathbb{Z}_n}^{\{1,3\}}$ is Ramanujan.
\item[\rm (g)] $G_{\mathbb{Z}_n}^{\{1,2,3\}}$ is Ramanujan.
\end{itemize}
\end{thm}


\subsection{Ramanujan Cayley graphs on dihedral groups}
\label{sec:RCGNon}

Recall from \eqref{eq:gnd} that $S_n(1)=\{k: 1\le k\le n-1, \gcd(k,n)=1\}$. Recall also that $D_{2n}= \langle a,b\ |\ a^n=b^2=1, b^{-1}ab=a^{-1} \rangle$ is the dihedral group of order $2n$. The following result due to Liu and Zhou identifies all Ramanujan graphs in two special families of Cayley graphs on dihedral groups.

\begin{thm}\emph{(\cite{LiuZhou})}\label{RamnujanDH1}
Let $n \ge 2$ be an integer.
\begin{itemize}
  \item[\rm (a)] Let $S=\{a^k, a^kb: k\in S_n(1)\}$. Then $\Cay(D_{2n},S)$ is Ramanujan if and only if one of the following holds:
  \begin{itemize}
  \item[\rm (i)] $n=2^\alpha$ for some integer $\alpha\ge1$;
  \item[\rm (ii)] $n=p q$ or $2pq$ for some primes $p, q$ with $p < q \le 2p-3$;
  \item[\rm (iii)] $n=p^2$ or $2p^2$ for some prime $p$;
  \item[\rm (iv)] $n=p$ for some prime $p$;
  \item[\rm (v)] $n=2^\alpha p$ for some integer $\alpha \ge 1$ and prime $p$ with $p > 2^{\alpha-2}+1$.
  \end{itemize}
  \item[\rm (b)] Let $S=\{a^k: k\in S_n(1)\}\cup \{a^jb\}$ for some $j\in \mathbb{Z}_n$. Then $\Cay(D_{2n},S)$ is Ramanujan if and only if $n \in \{1,2,3,4,5,6,8,10,12\}$.
\end{itemize}
\end{thm}


\subsection{Ramanujan Cayley graphs on finite rings}

Recall that $G_R$ is the unitary Cayley graph of a finite ring $R$. In \cite{Liu12}, {Liu} and {Zhou} determined when $G_R$ or its complement $\overline{G}_R$ is Ramanujan using knowledge of the spectra of these graphs (see Section \ref{subsec:UCayComRing}).

\begin{thm}
\emph{(\cite[Theorem 11]{Liu12})}
Let $R$ be a finite local ring with maximal ideal $M$ of order $m$. Then $G_R$ is Ramanujan if and only if either $|R|=2m$ or $|R|\geq\left(\dfrac{m}{2}+1\right)^2$ and $m \neq 2$.
\end{thm}

\begin{thm}
\label{XnThm1}
\emph{(\cite[Theorem 12]{Liu12})}
Let $R$ be as in Assumption \ref{as:1} with $s \ge 2$. Then $G_R$ is Ramanujan if and only if $R$ satisfies one of the following conditions:
\begin{itemize}
\item[\rm (a)] $R_i/M_i\cong  \mathbb{F}_2$ for $i = 1, 2, \ldots, s$;

\item[\rm (b)] $R_i \cong  \mathbb{F}_2$ for $i = 1, 2, \ldots, s-3$, and  $R_i \cong  \mathbb{F}_3$ for $i = s-2, s-1, s$;

\item[\rm (c)] $R_i \cong  \mathbb{F}_2$ for $i = 1, 2, \ldots, s-3$, $R_i \cong  \mathbb{F}_3$ for $i = s-2, s-1$, and $R_s\cong\mathbb{F}_4$;

\item[\rm (d)] $R_i \cong  \mathbb{F}_2$ for $i = 1, 2, \ldots, s-3$, and  $R_i \cong  \mathbb{F}_4$ for $i = s-2, s-1, s$;

\item[\rm (e)] $R_i \cong  \mathbb{F}_2$ for $i = 1, 2, \ldots, s-2$, $R_{s-1}\cong  \mathbb{F}_3$, and $R_{s}\cong  \mathbb{Z}_9 $ or $\mathbb{Z}_3[X]/(X^2)$;

\item[\rm (f)] $R_1\cong\mathbb{Z}_4$ or $\mathbb{Z}_2[X]/(X^2)$, $R_i \cong  \mathbb{F}_2$ for $i = 2, 3, \ldots, s-2$, and $R_{s-1}\cong  \mathbb{F}_{q_1}$ and $R_{s}\cong  \mathbb{F}_{q_2}$ for some prime powers $q_1, q_2 \ge 3$ such that
$q_1 \le q_2\leq q_1+\sqrt{(q_1-2)q_1}$;

\item[\rm (g)] $R_i \cong  \mathbb{F}_2$ for $i = 1, 2, \ldots, s-2$, and $R_{s-1}\cong  \mathbb{F}_{q_1}$ and  $R_{s}\cong  \mathbb{F}_{q_2}$ for some prime powers $q_1, q_2 \geq3$ such that $q_1 \le q_2\leq2\left(q_1+\sqrt{(q_1-2)q_1}\right)-1$;

\item[\rm (h)] $R_i/M_i\cong  \mathbb{F}_2$ for $i = 1, 2, \ldots, s-1$, $R_{s}/M_{s}\cong  \mathbb{F}_q$ for some prime power $q\geq3$, and $\prod_{i=1}^s m_i \leq2\left(q-1+\sqrt{(q-2)q}\right)$.
\end{itemize}
\end{thm}

In the special case when $R = \mathbb{Z}_n$, Theorem \ref{XnThm1} gives rise to Theorem \ref{XnCor1}.

\begin{thm}
\label{ComSXnThm1}
\emph{(\cite[Theorem 15]{Liu12})}
Let $R$ be a finite local ring. Then $\overline{G}_R$ is Ramanujan.
\end{thm}

\begin{thm}
\label{ComXnThm222}
\emph{(\cite[Theorem 16]{Liu12})}
Let $R$ be as in Assumption \ref{as:1} with $s \ge 2$. Then $\overline{G}_R$ is Ramanujan if and only if $R$ satisfies one of the following conditions:
\begin{itemize}
\item[\rm (a)] $|R_i|/m_i=2$ for $i = 1, 2, \ldots, s$, and $\prod_{i=1}^s m_i \leq 2^{s+1}-3+2\sqrt{2^s(2^s-3)}$;
\item[\rm (b)] $2=|R_1|/m_1=\cdots=|R_t|/m_t < |R_{t+1}|/m_{t+1}$ for some $t$ with $2 \le t < s$, and $|R^\times|\leq2\sqrt{|R|}-3$;
\item[\rm (c)] $2=|R_1|/m_1 < |R_2|/m_2$ and $|R^\times|\leq2\sqrt{|R|-2}-1$;
\item[\rm (d)] $3\leq|R_1|/m_1$ and
$\frac{|R^\times|}{(|R_1|/m_1)-1}\leq-\left(2(|R_1|/m_1)-3\right)+\sqrt{\left(2(|R_1|/m_1)-3\right)^2+(4|R|-9)}$.
\end{itemize}
\end{thm}

Applying Theorems \ref{ComSXnThm1} and \ref{ComXnThm222} to $\mathbb{Z}_n$, we obtain the following result on the complement of the unitary Cayley graph $G_{\mathbb{Z}_n}$ of $\ZZZ_n$.

\begin{cor}
\emph{(\cite[Corollary 17]{Liu12})}
Let $n \ge 2$ be an integer. Then $\overline{G}_{\mathbb{Z}_n}$ is Ramanujan if and only if either $n$ is a prime power or $n \in \{6, 10, 12, 15, 18, 21, 24, 30, 35\}$.
\end{cor}

The following results were proved by Rattanakangwanwong and Meemark in \cite{RattaM20}, where $G_{M_{n}(R)}$ is the unitary Cayley graph of the ring $M_{n}(R)$ of $n \times n$ matrices over a finite ring $R$.

\begin{thm}
\emph{(\cite[Theorem 3.2]{RattaM20})}
Let $n \ge 2$ be an integer and $q$ a prime power. Then $G_{M_{n}(\FFF_q)}$ is Ramanujan if and only if $n = 2$ or $(n, q) = (3, 2)$.
\end{thm}

\begin{thm}
\emph{(\cite[Theorem 4.5]{RattaM20})}
Let $R$ be a local ring other than a field. Then $G_{M_{n}(R)}$ is not Ramanujan for any $n \ge 2$.
\end{thm}

Ramanujan quadratic unitary Cayley graphs were characterized by {Liu} and {Zhou} in \cite{Liu15} using knowledge of the spectra of quadratic unitary Cayley graphs (see Section \ref{subsec:QUCayComRing}).

\begin{thm}
\label{RamQUCG13mod4}
\emph{(\cite[Theorem 5.1]{Liu15})}
Let $R$ be as in Assumption \ref{as:1} such that $|R_i|/m_i\equiv 1\,(\mod\,4)$ for $1 \le i \le s$, and let $R_0$ be a local ring with maximal ideal $M_0$ of order $m_0$ such that $|R_0|/m_0\equiv 3~(\mod\,4)$. Then the following hold:
\begin{itemize}
\item[\rm (a)] $\mathcal{G}_{R_0}$ is Ramanujan if and only if $|R_0|\ge(m_0+2)^2/4$;
\item[\rm (b)] $\mathcal{G}_{R}$ is Ramanujan if and only if $R$ is isomorphic to $\mathbb{F}_5\times\mathbb{F}_5$ or $\mathbb{F}_q$ for a prime power $q\equiv 1\,(\mod\,4)$;
\item[\rm (c)] $\mathcal{G}_{R_0\times R}$ is Ramanujan if and only if $R_{0}\times R$ is isomorphic to $\mathbb{F}_3\times\mathbb{F}_5$, $\mathbb{F}_3\times\mathbb{F}_9$ or $\mathbb{F}_3\times\mathbb{F}_{13}$.
\end{itemize}
\end{thm}

In the special case where $R = \mathbb{Z}_n$, Theorem \ref{RamQUCG13mod4} yields the following result.

\begin{cor}
\emph{(\cite[Corollary 5.2]{Liu15})}
Let $n \ge 5$ be an integer such that each of its prime factors is congruent to $1$ modulo $4$. Let $p\equiv3~(\mod\,4)$ be a prime and $\alpha \ge 1$ an integer. Then the following hold:
\begin{itemize}
\item[\rm (a)] $\mathcal {G}_{\mathbb{Z}_{n}}$ is Ramanujan if and only if $n$ is a prime;
\item[\rm (b)] $\mathcal {G}_{\mathbb{Z}_{p^\alpha}}$ is Ramanujan if and only if $p^\alpha=p$ or $p^2$;
\item[\rm (c)] $\mathcal {G}_{\mathbb{Z}_{np^\alpha}}$ is Ramanujan if and only if $np^\alpha = 15$ or $39$.
\end{itemize}
\end{cor}

Recently, a family of Cayley graphs from Galois rings have been studied by Satake in \cite{Satake19a}. Let $p$ be a prime, and let $e$ and $r$ be positive integers. Let $GR(p^e, r)$ be the Galois ring of characteristic $p^e$ and order $p^{er}$. It is well known that $GR(p^e, r)$ is a finite, commutative and local ring, and its multiplicative group of units is isomorphic to the direct product $\Ga_1 \times \Ga_2$, where $\Ga_1$ is a cyclic group of order $p^r - 1$ and $G_2$ is a group of order $p^{(e-1)r}$. Define $G_{p^{e}, r}$ to be the Cayley graph on the additive group of $GR(p^e, r)$ with respect to the connection set $\Ga_1 \cup (-\Ga_1)$. This is an undirected graph with degree $2^{r+1}-2$ when $p=2$ and $p^{r}-1$ when $p$ is odd. The eigenvalues of this graph were studied in \cite{Satake19a} with the help of some known estimates of character sums for Galois rings. In particular, the following result was proved in \cite{Satake19a}.

\begin{thm}
\emph{(\cite[Theorem 4.1]{Satake19a})}
Let $r \ge 4$ be an integer. Then $G_{2^{2}, r}$ is Ramanujan. Moreover, all eigenvalues of $G_{2^{2}, r}$  except $2^{r+1}-2$ are in the interval $[-2^{\frac{r}{2}+1}-2, 2^{\frac{r}{2}+1}+2]$.
\end{thm}

We finish this section with a result about when a semiprimitive generalized Paley graph or its complement is Ramanujan (see Section \ref{subsec:GPaley} for related definitions).

\begin{thm}
\label{PV20c}
\emph{(\cite[Theorem 4.1]{PodestaV20a})}
Let $q = p^r$ be a prime power with $r$ even, and let $k \ge 1$ be a divisor of $q-1$. Suppose that $(k, q)$ is a semiprimitive pair. Then $\overline{X_{q}^k}$ is Ramanujan. Moreover, $X_{q}^k$ is Ramanujan if and only if one of the following holds:
\begin{itemize}
\item[\rm (a)] $k = 2$ and $q \equiv 1~(\mod~4)$ (in which case $X_{q}^2$ is the Paley graph $P(q)$);
\item[\rm (b)] $k = 3$, $p = 2$ and $r \ge 4$;
\item[\rm (c)] $k = 3$, $2 < p \equiv 2~(\mod~3)$ and $r \ge 2$;
\item[\rm (d)] $k = 4$, $p = 3$ and $r \ge 4$;
\item[\rm (e)] $k = 4$, $3 < p \equiv 3~(\mod~4)$ and $r \ge 2$;
\item[\rm (f)] $k = 5$, $p = 2$ and $r \ge 8$ is a multiple of $4$;
\item[\rm (g)] $k = 5$, $2 < p \equiv 2, 3~(\mod~5)$ and $r \ge 4$ is a multiple of $4$;
\item[\rm (h)] $k = 5$, $p \equiv 4~(\mod~5)$ and $r \ge 2$.
\end{itemize}
\end{thm}


\subsection{Ramanujan Euclidean graphs}
\label{FiEuGr}

Let $q = p^r$ be an odd prime power and $n \ge 1$ an integer. Define
$$
d(\mathbf{x}, \mathbf{y}) = \sum_{i=1}^n(x_i-y_i)^2
$$
for (column) vectors $\mathbf{x}, \mathbf{y}$ in the linear space $\mathbb{F}_q^n$. Note that this is not a metric in the sense of analysis as $d(\mathbf{x}, \mathbf{y})$ is not real-valued but takes values in $\mathbb{F}_q$. Given $a\in \mathbb{F}_q$, the \emph{Euclidean graph} $E_q(n,a)$, introduced by Medrano, Myers, Stark and Terras in \cite{MedranoMST96}, is defined to be the graph with vertex set $\mathbb{F}_q^n$ such that two vertices $\mathbf{x}, \mathbf{y} \in \mathbb{F}_q^n$ are adjacent if and only if $d(\mathbf{x}, \mathbf{y}) = a$. (Loops are allowed when $a=0$.) Assume that $(q, n, a) \not= (q, 2, 0)$ and $-1$ is not a square in $\mathbb{F}_q$. Then
$$
E_q(n,a) = \Cay(\mathbb{F}_q^n, S_q(n,a))
$$
is a connected Cayley graph on the additive group of $\mathbb{F}_q^n$, where
$$
S_q(n,a)=\{\mathbf{x} \in \mathbb{F}_q^n:  d(\mathbf{x}, \mathbf{0}) = a\}.
$$
As shown in \cite[Theorem 1]{MedranoMST96}, the degree $|S_q(n,a)|$ of this graph is approximately $q^{n-1}$, with an error depending on the value of the quadratic character at $(-1)^{(n-1)/2}a$ or $(-1)^{n/2}$ and the parity of $n$. Define
$$
e_{\mathbf{b}}(\mathbf{x}) = \exp\{2 \pi i \Tr(\mathbf{b}^t \mathbf{x})/p\},\, \mathbf{b}, \mathbf{x} \in \mathbb{F}_q^n,
$$
where $\mathbf{b}^t$ is the transpose of $\mathbf{b}$ and as in \eqref{eq:trx}, $\Tr(u) = \Tr_{\FFF_q/\FFF_p}(u) = u+u^p+ \cdots + u^{p^{r-1}}$ is the trace of $u \in \FFF_q$. As shown in \cite[Proposition 2]{MedranoMST96}, the eigenvalues of $E_q(n,a)$ are given by
$$
\l_{\mathbf{b}} = \sum_{\mathbf{y} \in \mathbb{F}_q^n,\ d(\mathbf{y}, \mathbf{0})=1} e_{\mathbf{b}}(\mathbf{y}),\, \mathbf{b} \in \mathbb{F}_q^n,
$$
and moreover $e_{\mathbf{b}}$ is an eigenfunction corresponding to $\l_{\mathbf{b}}$. Based on this the following was proved by Medrano, Myers, Stark and Terras in \cite{MedranoMST96}.

\begin{thm}
\label{thm:Eqna}
\emph{(\cite[Theorem 3]{MedranoMST96})}
The eigenvalues $\l_{\mathbf{b}}$ of $E_q(n,a)$ for $\mathbf{b} \in \mathbb{F}_q^n \setminus \{0\}$ satisfy
$$
|\l_{\mathbf{b}}| \le 2 q^{(n-1)/2}
$$
and can be expressed as some generalized Kloosterman sums.
\end{thm}

Since the degree of $E_q(n,a)$ is approximately $q^{n-1}$, the bound in Theorem \ref{thm:Eqna} is asymptotic to the Ramanujan bound as $q$ approaches infinity. So Theorem \ref{thm:Eqna} implies that $E_q(n,a)$ is asymptotically Ramanujan for large $q$ (see also \cite[Theorem 1, Chapter 5]{Terras99}).

It was proved in \cite[Proposition 4]{MedranoMST96} that, for fixed $q$ and $n$, all graphs $E_q(n,a)$ for squares $a \ne 0$ are isomorphic, and all graphs $E_q(n,a)$ for non-squares $a \ne 0$ are isomorphic. It was further proved in \cite[Theorem 5]{MedranoMST96} that when $n$ is even these graphs for $a \ne 0$ are isomorphic to each other, and so there are exactly two Euclidean graphs up to isomorphism, namely $E_q(n,0)$ and $E_q(n,1)$ for even $n$.

Similar to $E_q(n,a)$, for a prime power $q$ and a positive integer $n$, the \emph{Euclidean graph} associated with $\mathbb{Z}_q^n$ and $a \in \mathbb{Z}_q$ is defined \cite{MedranoMST98} to be the Cayley graph
$$
X_q(n,a) = \Cay(\mathbb{Z}_q^n, S'_q(n,a)),
$$
where
$$
S'_q(n,a)=\left\{\mathbf{x} \in \mathbb{Z}_q^n: \sum_{i=1}^nx_i^2=a\right\}.
$$
In \cite{MedranoMST98}, Medrano, Myers, Stark and Terras determined the spectrum of $X_q(n,a)$ and proved that for $n \ge 2$ and $q = p^r$ with $p \ge 3$ and $r \ge 2$, $X_q(n,a)$ is not Ramanujan unless $(p, r, n) = (3, 2, 2)$. In the case when $q = 2^r$, the Euclidean graphs $X_q(n,a)$ were studied in \cite{DeDeo03} and \cite{Vinh09}. The following result is a combination of \cite[Theorem 2.5]{MedranoMST98}, \cite[Theorem 11]{DeDeo03} and \cite[Theorems 8 and 10]{Vinh09}.

\begin{thm}\emph{(\cite{DeDeo03,MedranoMST98,Vinh09})}
Let $r, d\ge 2$ and $a\in\mathbb{Z}$ be integers, and let $p$ be a prime not dividing $a$. Then $X_{p^r}(n, a)$ is not Ramanujan except when $(p, r, n) = (2,2,2), (2,2,3)$ or $(3, 2, 2)$.
\end{thm}

In 2009, Bannai, Shimabukuro and Tanaka \cite{BannaiST09} introduced finite Euclidean graphs in a more general setting. Let $q$ be a prime power, $n \ge 1$ an integer, $a$ an elment of $\mathbb{F}_q$, and $Q$ a non-degenerate quadratic form on $\mathbb{F}_q^n$. Define $E_q(n,Q,a)$ to be the graph with vertex set $\mathbb{F}_q^n$ in which distinct $\mathbf{x}, \mathbf{y} \in \mathbb{F}_q^n$ are adjacent if and only if $Q(\mathbf{x}-\mathbf{y})=a$. In the special case where $Q(\mathbf{x})=\sum_{i=1}^n x_i^2$, $E_q(n,Q,a)$ is exactly $E_q(n,a)$.

Consider a non-degenerate quadratic form $Q$ on $V = V_{n}(q) = \mathbb{F}_q^n$. The \emph{orthogonal group} associated with $Q$ is the group of all linear transformations that fix $Q$; that is,
$$
O(V, Q) = \{\s \in \GL(V): Q(\s(\mathbf{x})) = Q(\mathbf{x}) \text{ for all } \mathbf{x} \in V\}.
$$
Natually, the semi-direct product of $V$ (translations) by $O(V, Q)$, $V \rtimes O(V, Q)$, acts transitively on $V$. Let $\OO_0 = \{(\mathbf{x}, \mathbf{x}): \mathbf{x} \in V\}, \OO_1, \ldots, \OO_d$ be the orbital of $V \rtimes O(V, Q)$ on $V$. Treating $\OO_i$ as a relation on $V$, we obtain an association scheme $\XX(O(V, Q), V) = (V, \{\OO_i\}_{0 \le i \le d})$.

Let $\rho$ be a primitive element of $\FFF_q$. If $n=2m$ is even, then there are two inequivalent non-degenerate quadratic forms, denoted by $Q^{+}$ and $Q^{-}$, respectively. The corresponding orthogonal groups are $\GO_{2m}^{+}(q) = O(V, Q^{+})$ and $\GO_{2m}^{-}(q) = O(V, Q^{-})$. So we have the Euclidean graphs $E_q(2m,Q^{+},\rho^i)$ and $E_q(2m,Q^{-},\rho^i)$ for $1 \le i \le q$. The arc sets of $E_q(2m,Q^{+},0)$ and $E_q(2m,Q^{+},\rho^i)$ ($1 \le i \le q-1$), together with the identity relation $R_0 = \{(\mathbf{x}, \mathbf{x}): \mathbf{x} \in V_{2m}(q)\}$, are the relations of $\XX(\GO_{2m}^{+}(q), V_{2m}(q))$. Similarly, the arc sets of $E_q(2m,Q^{-},0)$ and $E_q(2m,Q^{-},\rho^i)$ ($1 \le i \le q-1$) together with $R_0$ are the relations of $\XX(\GO_{2m}^{-}(q), V_{2m}(q))$.

If $n = 2m+1$ and $q$ is odd, then there are two inequivalent non-degenerate quadratic forms $Q$ and $Q'$, but the groups $O(V, Q)$ and $O(V, Q')$ are isomorphic. This group is denoted by $\GO_{2m+1}(q)$. In this case we have the Euclidean graphs $E_q(2m+1,Q,\rho^{2i})$ and $E_q(2m+1,Q,\rho^{2i-1})$ for $1 \le i \le (q-1)/2$, as well as $E_q(2m+1,Q,0)$, whose arc sets together with the identity relation are the relations of $\XX(\GO_{2m+1}(q), V_{2m+1}(q))$.

If $n = 2m+1$ and $q$ is even, then there is exactly one inequivalent non-degenerate quadratic form $Q$. In this case we have the Euclidean graphs $E_q(2m+1,Q,0)$ and $E_q(2m+1,Q,\rho^{i})$ for $1 \le i \le q-1$, and the arc set of $E_q(2m+1,Q,\rho^{i})$ is the union of two relations $R_i$ and $R_{q+i}$ of $\XX(\GO_{2m+1}(q), V_{2m+1}(q))$, where $\GO_{2m+1}(q) = O(V, Q)$.

The next four theorems were proved by Bannai, Shimabukuro and Tanaka in \cite{BannaiST09}.

\begin{thm}
\emph{(\cite[Theorem 3.1]{BannaiST09})}
The graphs $E_q(2m,Q^{-},\rho^i)$ ($1 \le i \le q-1$) are Ramanujan.
\end{thm}

\begin{thm}
\emph{(\cite[Theorem 3.2]{BannaiST09})}
The graphs $E_q(2m,Q^{+},\rho^i)$ ($1 \le i \le q-1$) are Ramanujan if $m$ is sufficiently large (that is, larger than a certain integer determined by $q$).
\end{thm}

\begin{thm}
\emph{(\cite[Theorem 3.3]{BannaiST09})}
Let $q$ be an odd prime power.
\begin{itemize}
\item[\rm (a)] The graphs $E_q(2m+1,Q,\rho^{2i})$ ($1 \le i \le (q-1)/2$) and $E_q(2m+1,Q,0)$ are Ramanujan.
\item[\rm (b)] If $q$ is a prime, then the graphs $E_q(2m+1,Q,\rho^{2i-1})$ ($1 \le i \le (q-1)/2$) are Ramanujan if $m$ is sufficiently large (that is, larger than a certain integer determined by $q$).
\item[\rm (c)] If $q$ is not a prime, then the graphs $E_q(2m+1,Q,\rho^{2i-1})$ ($1 \le i \le (q-1)/2$) are not Ramanujan.
\end{itemize}
\end{thm}

\begin{thm}
\emph{(\cite[Theorem 3.4]{BannaiST09})}
Let $q$ be an even prime power.
\begin{itemize}
\item[\rm (a)] The graphs $(V, R_{q+i})$ ($1 \le i \le q-1$) and $E_q(2m+1,Q,0)$  are Ramanujan unless $q=2$ and $m=1$.
\item[\rm (b)] For $a=\rho^i \ne 0$, the graph $E_q(2m+1,Q,a) = (V, R_i \cup R_{q+i})$ is Ramanujan if and only if $m \ge 1$.
\end{itemize}
\end{thm}

In a companion of \cite{BannaiST09}, Bannai, Shimabukuro and Tanaka \cite{BannaiST04} constructed many other Ramanujan graphs from association schemes obtained from the following actions of orthogonal groups over finite fields: (i)  $O_{2m+1}(q)$ ($q$ odd) acting on the set of non-square-type non-isotropic $1$-dimensional subspaces of $V_{2m+1}(q)$ with respect to the non-degenerate quadratic form $Q(\mathbf{x}) = 2(x_1 x_{m+1} + \cdots + x_m x_{2m}) + x_{2m+1}^2 $; (ii)  $O_{2m+1}(q)$ ($q$ odd) acting on the set of square-type non-isotropic $1$-dimensional subspaces of $V_{2m+1}(q)$ with respect to the same quadratic form; (iii)  $O_{2m+1}(q)$ ($q$ even) acting on the set of negative-type hyperplanes of $V_{2m+1}(q)$ with respect to the non-degenerate quadratic form $Q(\mathbf{x}) = x_1 x_{m+1} + \cdots + x_m x_{2m} + x_{2m+1}^2 $; (iv)  $O_{2m+1}(q)$ ($q$ even) acting on the set of positive-type hyperplanes of $V_{2m+1}(q)$ with respect to the same quadratic form; (v) $O^{\pm}_{2m}(q)$ acting on the set of non-isotropic points of $V_{2m}(q)$; (vi) $U_{n}(q)$ acting on the set of $1$-dimensional non-isotropic subspaces of $V_{n}(q^2)$ with respect to the canonical non-singular Hermitian form. The reader is referred to \cite{BannaiST04} for details.


\subsection{Ramanujan finite upper half plane graphs}
\label{FiUpHaP}

Let $q = p^r$ be an odd prime power and $\delta$ a non-square of the finite field $\mathbb{F}_q$. The \emph{finite upper half plane} over $\mathbb{F}_q$ is defined as
$$
\mathbf{H}_q=\{z=x+y\sqrt{\delta}: x,y\in \mathbb{F}_q, ~y\not=0\}.
$$
For $z=x+y\sqrt{\delta}\in \mathbf{H}_q$ and $w=u+v\sqrt{\delta}\in \mathbf{H}_q$, define
$$
d(z,w)=\frac{(x-u)^2-\delta(y-v)^2}{yv}.
$$
Let $a\in \mathbf{H}_q$. Define $X_q(\delta,a)$ to be the graph with vertex set $\mathbf{H}_q$ such that $z, w \in \mathbf{H}_q$ are adjacent if and only if $d(z,w)=a$. This graph is called \cite{Terras99} a \emph{finite upper half plane graph} on $\mathbb{F}_q$. As shown in \cite[Theorem 2(4), Chapter 19]{Terras99}, this is a Cayley graph on the affine group
$$
\mathrm{Aff}(q)
= \left\{\pmat{y & x \\ 0 & 1}: x, y\in\mathbb{F}_q,~y\not=0 \right\},
$$
namely,
$$
X_q(\delta,a) \cong \Cay(\mathrm{Aff}(q), S_q(\delta,a)),
$$
where
$$
S_q(\delta,a)=\left\{\pmat{y & x \\0 & 1}: x, y\in\mathbb{F}_q,~y\not=0,~x^2=ay+\delta(y-1)^2\right\}.
$$

\begin{thm}\emph{(\cite[Theorem 4, Chapter 20]{Terras99})}
The finite upper half plane graphs $X_q(\delta, a)$, for $a\not=0$, $4\delta$, are Ramanujan.
\end{thm}

The reader is referred to \cite{ChaiLi04} and \cite[Theorem 2, Chapter 19]{Terras99} for more properties of $X_q(\delta,a)$ (mostly due to Angel \cite{Angel93, Angel96}, Angel, Celniker, Poulos, Terras, Trimble and Velasquez \cite{AngelCPTTV92}, Angel, Poulos, Terras, Trimble and Velasquez \cite{AngelCPTTV92, AngelPTTV94}, Celniker \cite{Celniker94}, Celniker, Poulos, Terras, Trimble and Velasquez \cite{CelnikerPTTV93}, Poulos \cite{Poulos91}, and Terras \cite{Terras91, Terras96}).

Let $q = p^r$ be an odd prime power and $\delta$ a positive integral generator of the group $\mathbb{Z}_q^{\times}$. The \emph{finite upper half plane} over $\mathbb{Z}_q$ is defined as
$$
\mathbf{H}'_q=\{z=x+y\sqrt{\delta}: x\in \mathbb{Z}_q,\, y\in \mathbb{Z}_q^{\times}\}.
$$
For $z=x+y\sqrt{\delta}\in \mathbf{H}'_q$ and $w=u+v\sqrt{\delta}\in \mathbf{H}'_q$, define
$$
d'(z,w)=\frac{(x-u)^2-\delta(y-v)^2}{yv}.
$$
Fix $a\in \mathbb{Z}_q$. In \cite{AngelTST95}, Angel, Trimble, Shook and Terras defined the \emph{finite upper half plane graph} $X'_{q}(\delta,a)$ on $\mathbb{Z}_q$ to be the graph with vertex set $\mathbf{H}'_q$ such that $z, w\in \mathbf{H}'_q$ are adjacent if and only if $d'(z,w)=a$.
This graph was proved (see \cite[Theorem 2]{AngelTST95} or \cite[Proposition 3.2]{Angel93}) to be a Cayley graph on the affine group
$$
\mathrm{Aff}(\mathbb{Z}_q)
=\left\{\pmat{y & x \\ 0 & 1}: x\in \mathbb{Z}_q,\, y\in \mathbb{Z}_q^{\times} \right\},
$$
namely,
$$
X'_{q}(\delta,a) \cong \Cay(\mathrm{Aff}(\mathbb{Z}_q), S_{\mathbb{Z}_q}(\delta,a)),
$$
where
$$
S_{\mathbb{Z}_q}(\delta,a)=\left\{
\pmat{y & x \\ 0 & 1}: x\in \mathbb{Z}_q,\, y\in \mathbb{Z}_q^{\times},\, x+y\sqrt{\delta}\in \mathbb{Z}_q,\, d'\left(x+y\sqrt{\delta},\sqrt{\delta}\right)=a\right\}.
$$
If $a \not \equiv 0, 4\delta~\mod~p$, then $X'_{q}(\delta,a)$ is a connected $(p^r + p^{r-1})$-regular graph with $p^r (p^r - p^{r-1})$ vertices (see \cite[Theorems 1 and 2]{AngelTST95}). The following result was proved by Angel, Trimble, Shook and Terras in \cite{AngelTST95}.

\begin{thm}\emph{(\cite[Theorem 6]{AngelTST95})}
Let $p$ be an odd prime, $\d$ a generator of $\mathbb{Z}_{p^2}^{\times}$, and $a \in \mathbb{Z}_{p^2}$. If $p = 3$, then $X'_{p^2}(\delta,a)$ is Ramanujan; if $p \ge 5$, then $X'_{p^2}(\delta,a)$ is not Ramunujan.
\end{thm}

The following result due to Bell and Minei \cite{BellM06} gives a class of Ramanujan Cayley graphs on $\mathrm{Aff}(\mathbb{Z}_p)$.

\begin{thm}
\emph{(\cite[Theorem 13]{BellM06})}
Let $p$ be a prime and let $a$ and $d$ be integers with $\gcd(p-1, d) = 1$
and $\gcd(p-1, d-1) = 2$. Let $g$ be a primitive root of $p$. Set
$$
S=\left\{
\pmat{i^d & g^a i\\ 0 & 1},
\pmat{i^d & -g^a i^{d-1} \\ 0 & 1}: i=1,2,\ldots,p-1\right\}.
$$
Then $\Cay(\mathrm{Aff}(\mathbb{Z}_p), S)$ is a Ramanujan graph.
\end{thm}

Let $p \ge 5$ be an odd prime. Consider the following subgroup of $\GL(3,p)$:
$$
\Gamma'=\left\{\pmat{y& x&z \\
                               0 & 1 &0\\
                               0&0&1}:
                               x,z\in\mathbb{F}_p,~y\in \mathbb{F}_p^{\times}\right\}.
$$
In \cite{Allen98}, Allen studied the Cayley graph $X_p(\delta,a,c)$ on $\Ga'$ with connection set
$$
S_p(\delta,a,c)=\left\{
\pmat{y& x&z \\ 0 & 1 &0\\ 0&0&1}
\in \Gamma': x,z\in\mathbb{F}_p,~y,\delta,a,c\in\mathbb{F}_p^{\times}, ~x^2+cz^2=ay+\delta(y-1)^2\right\}.
$$
It was proved \cite{Allen98} that this is a connected non-bipartite graph with girth $3$ or $4$ and degree $p^2 - p + p \left(\frac{-c}{p}\right) \left(\frac{a(a-4\delta}{p}\right) + \left(\frac{-c}{p}\right)$, where $\left(\frac{x}{y}\right)$ denotes the Legendre symbol. Allen \cite{Allen98} made the following conjectures:

\begin{conj}
\label{AllenConj1}
\emph{(\cite[Conjecture 1]{Allen98})}\label{conj:xpdac}
If $\left(\frac{-c}{p}\right)=\left(\frac{a(a-4\delta}{p}\right)=1$, then $X_p(\delta,a,c)$ is Ramanujan.
\end{conj}

\begin{conj}
\label{AllenConj2}
\emph{(\cite[Conjecture 2]{Allen98})}
If $a=4\delta$, then $X_p(\delta,a,c)$ is Ramanujan.
\end{conj}

Conjecture \ref{conj:xpdac} has been confirmed by Allen for all odd primes less than 225 (for a total of over 2000 graphs).

In a similar fashion, Allen also studied a family of Cayley graphs on $\GL(n, p)$ for $n \ge 4$. In general, these graphs are not Ramanujan as shown in \cite[Section 3]{Allen98}.


\subsection{Heisenberg graphs}

The \emph{Heisenberg group} $H(R)$ over a ring $R$ is the multiplicative group of $3\times 3$ upper triangular matrices with entries in $R$ and ones on the diagonal. Denote such a matrix $\pmat{1 & x & z\\0 & 1 & y\\0 & 0 & 1}$ by $(x, y, z)$. Call
$$
G(p^n) = \Cay(H(\ZZZ_{p^n}), \{(\pm 1, 0, 0), (0, \pm 1, 0)\})
$$
and
$$
G(2^n)' = \Cay(H(\ZZZ_{2^n}), \{(\pm 1, 0, 0), (1, 1, 0), (1, 1, 0)^{-1}\})
$$
the \emph{Heisenberg graphs} \cite{DeDeoMMMST02} of $\ZZZ_{p^n}$. Obviously, these are $4$-regular graphs. In \cite{DeDeoMMMST02}, DeDeo, Mart\'{i}nez, Medrano, Minei, Stark and Terras proved the following result.

\begin{thm}
\emph{(\cite[Theorems 2-3]{DeDeoMMMST02})}
Let $p$ be a prime and $n \ge 1$ an integer. If $p$ is odd, then $G(p^n)$ is the only connected Cayley graph on $H(\ZZZ_{p^n})$ up to isomorphism; if $p=2$, then $G(2^n)$ and $G(2^n)'$ are the only connected Cayley graphs on $H(\ZZZ_{2^n})$ up to isomorphism. Moreover, the spectra of these Heisenberg graphs approach a continuous interval $[-4, 4]$ as $p^n \rightarrow \infty$.
\end{thm}


\subsection{Platonic graphs and beyond}

Let $q = p^r$ be an odd prime power. Define $G^*(q)$ to be the graph with vertex set $V(G^*(q)) = (\FFF_q^2 \setminus \{(0,0)\})/\la \pm 1 \ra$ in which $(a, b), (c, d) \in V(G^*(q))$ are adjacent if and only if the determinant of $\pmat{a & b \\ c & d}$ is equal to $\pm 1$. In particular, for an odd prime $p$, $G^*(p)$ is called the $p$th \emph{Platonic graph}. This name was coined as $G^*(3)$ and $G^*(5)$ are isomorphic to the $1$-skeletons of the tetrahedron and icosahedron, respectively.

\begin{thm}
\label{thm:Pl}
\emph{(\cite[Theorem 1]{DLM07}; also \cite{LM05} and \cite{Nica17})}
Let $q = p^r$ be an odd prime power. Then
\begin{eqnarray*}
\Spec (G^*(q)) = \pmat{q & -1 & \sqrt{q} & -\sqrt{q} \\
1  & q  & (q + 1)(q - 3)/4 & (q + 1)(q - 3)/4}.
\end{eqnarray*}
In particular, $G^*(q)$ is a Ramanujan graph.
\end{thm}

The proof of this result in \cite{LM05} by Li and Meemark used representation-theoretic methods, while the proof of it in \cite{DLM07} by DeDeo, Lanphier and Minei was done by relating $G^*(q)$ to a certain quotient graph of a Cayley graph $G(q)$ on $\PSL(2,q)$ with degree $4$ (when $p \equiv 1\ \mod~4$) or $5$ (when $p \equiv 3\ \mod~4$). Lower and upper bounds for the isoperimetric number of $G^*(q)$ were also given in \cite{DLM07}.

The notion of Platonic graphs was generalized to finite rings by Nica \cite{Nica17}. Given a finite ring $R$ and a positive integer $n$, let $R^{n, u}$ be the set of \emph{unimodular $n$-tuples}, namely those $n$-tuples of elements of $R$ whose entries generate $R$ as an ideal. The \emph{Platonic graph} over $R$, denoted by $\mathrm{Pl}(R)$, is defined to have vertex set $R^{2, u}/\{\pm 1\}$, the set of unimodular pairs taken up to sign, such that two such pairs $[a, b], [c, d]$ are adjacent if and only if $ad-bc = \pm 1$. Using certain Eisenstein character sums, for any finite chain ring $R$ the eigenvalues of $\mathrm{Pl}(R)$ were computed in \cite[Theorem 1.12]{Nica17} in terms of the parameters $q, s$ defined in Section \ref{subsec:ring-gcd}. In particular, in the case when $R$ is a finite field, this result gives rise to Theorem \ref{thm:Pl}. In \cite{Nica17}, two other families of unimodular graphs over finite chain rings $R$, both of which are bipartite Cayley graphs on two copies of $R^{n, u}$, were also introduced and their eigenvalues were determined using Eisenstein sums associated with the unramified extension of $R$.


\subsection{A family of Ramanujan Cayley graphs on $\ZZZ_p[i]$}

Solving a conjecture in \cite[Conjecture 31]{CM16}, Bibak, Kapron and Srinivasan \cite{BKS16} proved the following result using results from number theory.

\begin{thm}
\emph{(\cite[Theorem 9]{BKS16})}
Let $p$ be a prime with $p \equiv 3$ ($\mod~4$) and $S$ the set of units of ring $\ZZZ_p[i]$. Then the Cayley graph $\Cay(\ZZZ_p[i], S)$ on the additive group of $\ZZZ_p[i]$ is a $(p+1)$-regular Ramanujan graph.
\end{thm}


\subsection{Bounds for the degrees of Ramanujan Cayley graphs}\label{sec:BDRCG}

Given an odd integer $m \ge 1$, define $\hat{l}_m$ to be the maximum odd positive integer $l$ such that for every even integer $k$ between $m-l$ and $m-1$, all $k$-regular circulant graphs of order $m$ are Ramanujan. The following result is due to Hirano, Katata and Yamasaki \cite{HiranoKY13}.

\begin{thm}\emph{(\cite[Theorem 1.1]{HiranoKY13})}
Let $m \ge 15$ be an odd integer. Then
$$
\hat{l}_m = 2\left\lfloor\sqrt{m}-\frac{3}{2}\right\rfloor + \varepsilon_m + 1,
$$
where $\varepsilon_m \in \{0, 2\}$. Moreover, the case $\varepsilon_m = 2$ occurs only if $m$ is represented by one of the quadratic polynomials $k^2 +5k +c$ for some $c \in \{\pm1,\pm3,\pm5\}$ and is either a prime or the product of two primes $p, q$ with $p < q < 4p$.
\end{thm}

A parameter $\hat{l}$ similar to $\hat{l}_m$ can be defined \cite{HiranoKY13} for any finite abelian group of odd order $m$. It was noted in \cite{HiranoKY13} that this parameter satisfies $\hat{l} \ge 2\left\lfloor\sqrt{m}-\frac{3}{2}\right\rfloor + 1$. Moreover, for most finite abelian groups of odd order the equality holds, as proved in \cite{HiranoKY13}.

\begin{thm}
\emph{(\cite[Theorem 7.1]{HiranoKY13})}
Let $\Ga$ be a finite abelian group of odd order other than a cyclic group. Then
$$
\hat{l} = 2\left\lfloor\sqrt{m}-\frac{3}{2}\right\rfloor + 1
$$
except when $\Ga = \ZZZ_p \oplus \ZZZ_p$ for an odd prime $p$ between $3$ and $17$, and in these exceptional cases we have
$$
\hat{l} =
\begin{cases}
2\left\lfloor\sqrt{m}-\frac{3}{2}\right\rfloor + 3,\, p = 7,11,13,17\\
2\left\lfloor\sqrt{m}-\frac{3}{2}\right\rfloor + 5,\, p = 5.\\
\end{cases}
$$
\end{thm}

In general, let $\Ga$ be a finite group and $\SS$ a family of connection sets of $\Ga$ (that is, inverse-closed subsets of $\Ga \setminus \{1\}$). Define $\hat{l}(\Ga, \SS)$ to be the maximum positive integer $l$ such that for every $S \in \SS$ such that $1 \le |\Ga| - |S| \le l$ the Cayley graph $\Cay(\Ga, S)$ is Ramanujan (that is, every Cayley graph $\Cay(\Ga, S)$ with $S \in \SS$ and degree between $|\Ga|-l$ and $|\Ga| - 1$ is Ramanujan). In \cite{HiranoKY16}, this parameter was introduced and studied for two families $\SS$ of connection sets of $D_{2p}$, where $p$ is an odd prime. As mentioned above, earlier this parameter was studied in \cite{HiranoKY13} for cyclic groups of odd order and finite abelian groups of odd order, with $\SS$ the set of all connection sets.

The following families were considered by Hirano, Katata and Yamasaki in \cite{HiranoKY16}: $\SS_A$, the family of all connection sets of $\Ga$; $\SS_N$, the family of all normal connection sets of $\Ga$, where as before a connection set is called normal if it is the union of some conjugacy classes. Define
$$
l_0(\Ga) = \max\{|\Ga|-|S|: S \in \SS_N, |\Ga|-|S| \le 2(\sqrt{|\Ga|}-1)\}.
$$
Then $l_0(\Ga) \le \hat{l}(\Ga, \SS_N)$ as shown in \cite[Lemma 2.2]{HiranoKY16}.

A {\em Frobenius group} is a transitive permutation group which is not regular but only the identity element can fix two points. It is well known that any finite Frobenius group can be expressed as $N \rtimes H$ with $N$ a nilpotent normal subgroup, where $N \rtimes H$ acts on $N$ in such a way that $N$ acts on itself by right multiplication and $H$ acts on $N$ by conjugation.

The next three results are due to Hirano, Katata and Yamasaki.

\begin{thm}
\label{thm:HKY16}
\emph{(\cite[Theorem 3.3]{HiranoKY16})}
Let $\Ga = N \rtimes H$ be a finite Frobenius group with Frobenius kernel $N$ and complement $H$. Suppose that $(|N|-1)/|H| \ge 4$. Then
$$
\hat{l}(\Ga, \SS_N) = l_0(\Ga) < |N|.
$$
\end{thm}

It was noted that the condition $(|N|-1)/|H| \ge 4$ cannot be removed for otherwise the result may not be true. Applying Theorem \ref{thm:HKY16} to $D_{2p} = \ZZZ_p \rtimes \ZZZ_2$, where $p$ is a prime, the following result for normal dihedrants (that is, Cayley graphs on dihedral groups with respect to normal connection sets) was obtained.

\begin{thm}
\emph{(\cite[Corollary 3.4]{HiranoKY16})}
Let $p \ge 11$ be a prime. Then
$$
\hat{l}(D_{2p}, \SS_N) = l_0(D_{2p}) = 2 \left\lf \sqrt{2p} - \frac{1}{2} \right\rf - 1.
$$
\end{thm}

\begin{thm}
\emph{(\cite[Theorem 4.3]{HiranoKY16})}
Let $p \ge 29$ be a prime. Then the following hold:
\begin{itemize}
\item[\rm (a)] if $\lf 2 \sqrt{2p}\rf$ is even, then $\hat{l}(D_{2p}, \SS_A) = 2 \left\lf \sqrt{2p} - \frac{1}{2} \right\rf$;
\item[\rm (b)] if $\lf 2 \sqrt{2p}\rf$ is odd, then $\hat{l}(D_{2p}, \SS_A) = 2 \left\lf \sqrt{2p} - \frac{1}{2} \right\rf - 1$ or $2 \left\lf \sqrt{2p} - \frac{1}{2} \right\rf$.
\end{itemize}
\end{thm}

In the case when $\lf 2 \sqrt{2p}\rf$ is odd and $\hat{l}(D_{2p}, \SS_A) = 2 \left\lf \sqrt{2p} - \frac{1}{2} \right\rf$, the prime $p$ is called \emph{exceptional} \cite{HiranoKY16}. A characterization of exceptional primes was given in \cite[Theorem 4.5]{HiranoKY16} and connections with the well-known Hardy-Littlewood conjecture was discussed in \cite[Corollary 4.6]{HiranoKY16}.

In \cite{Yamasaki18}, Yamasaki applied the approach above to the generalized quaternion groups
\begin{equation}
\label{eq:Q4m}
Q_{4m} = \la x, y\ |\ x^{2m} = 1, x^m = y^2, y^{-1}xy = x^{-1}\ra
\end{equation}
and proved the following result.

\begin{thm}
\emph{(\cite[Theorem 4.3]{Yamasaki18})}
Let $m \ge 1$ be an integer. Then
$$
\hat{l}(Q_{4m}, \SS_A) = l_0(Q_{4m}) = \lf 4 \sqrt{m} \rf - 2.
$$
\end{thm}

Denote by $\SS'$ the family of connection sets $S$ of $Q_{4m}$ such that $\la x \ra y \not \subseteq S$. In \cite{Yamasaki18}, the parameter $\hat{l}(Q_{4m}, \SS')$ was also studied and its relation to the Hardy-Littlewood conjecture was discussed.


\section{Second largest eigenvalue of Cayley graphs}
\label{sec:2ndCay}

In this section we use $\l_2$ to denote the second largest eigenvalue of the graph under consideration. This invariant has been a focus of research in spectral graph theory over many years. In particular, the second largest eigenvalue of Cayley graphs has been studied extensively in the context of expander graphs, owning to the basic result (see Section \ref{sec:Ram}) that a $k$-regular graph is a good expander if and only if its spectral gap $k-\l_2$ (that is, the algebraic connectivity) is large. In this section we review some of the known results on the second largest eigenvalue of several families of Cayley graphs.


\subsection{Cayley graphs on Coxeter groups}
\label{subsec:cox}

A Coxeter group \cite{Davis2008} generated by a set $\{s_1, s_2, \ldots, s_{l}\}$ is defined by integers $m_{ij}$ with $m_{ii} = 1$ and $m_{ij} = m_{ji} \ge 2$ for $i \ne j$ and has a presentation $\Ga \cong \langle s_1, s_2, \ldots, s_{l} \mid (s_i s_j)^{m_{ij}} = 1\text{ for all } i, j\rangle$. In this section we focus on the second largest eigenvalue of Cayley graphs on $S_n$ and other Coxeter groups. As before, we use $S_n$ to denote the symmetric group on $[n] = \{1, 2, \ldots, n\}$, where $n \ge 2$. We begin with the following result due to Flatto, Odlyzko and Wales \cite{FOW85}.

\begin{thm}
\label{thm:FOW85}
\emph{(\cite{FOW85})}
Let $T = \{(1,n), (2,n), \ldots, (n-1, n)\}$. Then $\l_2(\Cay(S_n, T)) = n-2$.
\end{thm}

Note that this graph $\Cay(S_n, T)$ is exactly the star graph of degree $n-1$ discussed in Section \ref{subsec:sn}. Hence Theorem \ref{thm:FOW85} can be derived from Theorem \ref{thm:star}.

The following result due to Bacher gives the second largest eigenvalue of another well-known Cayley graph on $S_n$, namely the \emph{bubble sort graph}.

\begin{thm}
\label{thm:Bacher94}
\emph{(\cite{Bacher94})}
Let $T_1 = \{(1,2), (2,3), \ldots, (n-1, n)\}$. Then $\l_2(\Cay(S_n, T_1)) = n-3+2\cos(\pi/n)$ and its multiplicity is equal to $n-1$.
\end{thm}

The proof of Theorem \ref{thm:Bacher94} as given in \cite{Bacher94} made use of the irreducible representations of $S_n$. A similar approach was used by Friedman \cite{Friedman00} to prove the following result in which the extremal graphs achieving the lower bound $n-2$ are isomorphic to the star graph of degree $n-1$.

\begin{thm}
\label{thm:Fri}
\emph{(\cite[Theorem 1.1]{Friedman00})}
Let $T$ be a set of $n-1$ transpositions in $S_n$. Then $\l_2(\Cay(S_n, T)) \ge n-2$, with equality if and only if $T= \{(i, j): j \ne i\}$ for some fixed $i$.
\end{thm}

Given a set $T$ of transpositions in $S_n$, define $G_T$ to be the graph with vertex set $[n]$ such that $i, j \in [n]$ are adjacent if and only if $(i, j) \in T$. Theorem \ref{thm:Fri} was proved through establishing the following result.

\begin{thm}
\emph{(\cite[Theorem 1.2]{Friedman00})}
Let $T$ be a set of $n-1$ transpositions in $S_n$. If $G_T$ is bipartite, then $\l_2(G_T)$ occurs as an eigenvalue of $\Cay(S_n, T)$ with multiplicity at least $n-1$ times its multiplicity in $G_T$.
\end{thm}

Friedman conjectured that $\Cay(S_n, T)$ and $G_T$ have the same algebraic connectivity as long as $G_T$ is bipartite (see \cite[Conjecture 1.1]{Friedman00}). It turns out that this is a special case of the following more general conjecture of Aldous, known as Aldous' spectral gap conjecture.

\begin{conj}
\label{conj:Aldous}
For any set $T$ of transpositions in $S_n$, $\Cay(S_n, T)$ and $G_T$ have the same algebraic connectivity.
\end{conj}

In the language of probability theory, this conjecture asserts that the random walk and the interchange process on the graph have the same spectral gap. See \cite{Cesi10} for more information about Aldous's conjecture.

\begin{thm}
\label{Aldous-conj-spec}
Aldous' spectral gap conjecture is true if $G_T$ is
\begin{itemize}
\item[\rm (a)] a star (\cite{FOW85}),
\item[\rm (b)] a complete graph (\cite{Diaconis81}), or
\item[\rm (c)] a complete multipartite graph (\cite[Theorem 3.1]{Cesi10}).
\end{itemize}
\end{thm}

Major progresses on Aldous' spectral gap conjecture were made in \cite{HJ96} and \cite{KN97}, where irreducible representations of $S_n$ were heavily used. The conjecture in its general form was finally proved by Caputo, Liggett and Richthammer in \cite{CLR09}.

\begin{thm}
\emph{(\cite[Theorem 1.1]{CLR09}, unweighted version)}
\label{thm:AldThm}
Aldous' spectral gap conjecture is true for any set $T$ of transpositions in $S_n$.
\end{thm}

In fact, what they proved is a stronger result which will be presented in Theorem \ref{thm:AldThm1} below. A key ingredient in their proof was an inequality called the octopus inequality (see \cite[Theorem 2.3]{CLR09} for details). In \cite[Theorem 4.2]{Cesi16}, Cesi gave a simpler and more transparent proof of this inequality by looking at Aldous' spectral gap conjecture from an algebraic perspective. In the same paper Cesi also gave a self-contained algebraic proof of Aldous' spectral gap conjecture, showing how the conjecture follows from the octopus inequality.

Let $\tau$ be an element of $S_n$. Write $\tau$ as a string $\tau = [\tau_1, \tau_2, \ldots, \tau_n]$, where $\tau_i = \tau (i)$ for $1 \le i \le n$. A \emph{substring} in $\tau$ is a sequence of the form $[\tau_i, \ldots, \tau_j]$ for some $1 \le i < j \le n$, and \emph{reversing} this substring yields $[\tau_j, \ldots, \tau_i]$. A \emph{substring reversal} of $\tau$ is any permutation obtained from $\tau$ by reversing a substring in $\tau$. The \emph{reversal graph} $R_n$ is the graph with vertex set $S_n$ in which two vertices are adjacent if and only if they are substring reversals of each other. It is evident that $R_n$ is a Cayley graph on $S_n$ with degree ${\choose{n}{2}}$. (More precisely, $R_n$ is a right Cayley graph as $\sigma, \tau \in S_n$ are adjacent in $R_n$ if and only if $\tau = \sigma \alpha$ for some permutation $\alpha$ which can be obtained by reversing a substring in the identity element $[1, 2, \ldots, n]$.) A \emph{prefix reversal} of $\tau$ is any permutation obtained from $\tau$ by reversing a prefix substring $[\tau_1, \ldots, \tau_j]$ (where $2 \le j \le n$) in $\tau$. The \emph{pancake graph} $P_n$ is the Cayley graph on $S_n$ in which two vertices are adjacent if and only if they are prefix reversals of each other. Obviously, $P_n$ is a spanning subgraph of $R_n$ with degree $n-1$. Answering a question posed in \cite{GunnellsSW07}, the following result was proved by Cesi in \cite{Cesi09}.

\begin{thm}
\emph{(\cite{Cesi09})}
Let $P_n$ be the pancake graph. Then $\l_{2}(P_n) = n-2$.
\end{thm}

In \cite{ChungT17}, Chung and Tobin generalized this result as follows.

\begin{thm}
\emph{(\cite[Theorem 2]{ChungT17})}
Let ${\cal G}_n$ be the family of graphs with vertex set $S_n$ such that for each $2 \le i \le n$, every vertex $[\tau_1, \tau_2, \ldots, \tau_n]$ is adjacent to exactly one vertex of the form $[\tau_i, \a_2, \ldots, \a_{i-1}, \tau_1, \tau_{i+1}, \ldots, \tau_n]$, where $[\a_2, \ldots, \a_{i-1}]$ is a permutation of $\{\tau_2, \ldots, \tau_{i-1}\}$. Then $\l_{2}(G) = n-2$ for any $G \in {\cal G}_n$.
\end{thm}

Note that many graphs in ${\cal G}_n$ are Cayley graphs, including the pancake graph $P_n$ and the star graph $\Cay(S_n, T)$ where $G_T = K_{1, n}$. However, not every graph in this family is a Cayley graph. In \cite{ChungT17}, Chung and Tobin also proved the following result.

\begin{thm}
\emph{(\cite[Theorem 1]{ChungT17})}
Let $R_n$ be the reversal graph. Then $\l_{2}(R_n) = {\choose{n}{2}} - n$.
\end{thm}

There is a signed version of $R_n$, namely, the \emph{signed reversal graph} $SR_n$, whose vertices are the \emph{signed permutations} which are defined as $n$-tuples $[\tau_1, \tau_2, \ldots, \tau_n]$ with all $\tau_i \in \{\pm 1, \pm 2, \ldots, \pm n\}$ such that $[|\tau_1|, |\tau_2|, \ldots, |\tau_n|]$ is a permutation in $S_n$. Two vertices $[\s_1, \s_2, \ldots, \s_n]$ and $[\tau_1, \tau_2, \ldots, \tau_n]$ are adjacent in $SR_n$ if and only if $[\tau_1, \tau_2, \ldots, \tau_n]$ can be obtained from $[\s_1, \s_2, \ldots, \s_n]$ by reversing a substring $[\s_{i}, \ldots, \s_j]$, where $1 \le i \le j \le n$, and negating each reversed entry. In \cite[Section 5]{CioabaRT20}, Cioab\v{a}, Royle and Tan determined part of the spectrum of $SR_n$ and showed that for any eigenvalue $\mu$ of $R_n$, $\mu + n$ is an eigenvalue of $SR_n$. In \cite[Conjecture 6.1]{CioabaRT20}, they conjectured that for any $n \ge 2$ the second largest eigenvalue of $SR_n$ is ${\choose{n}{2}}$.

In recent years, Aldous' spectral gap conjecture has been generalized in different directions. In \cite{ParzanchevskiP20}, Parzanchevski and Puder explored its counterpart for normal Cayley graphs on $S_n$. To explain their results, let us first mention an equivalent version of Aldous' spectral gap conjecture. It is well known that the left regular representation $\l_{\text{reg}}$ of a group $\Ga$ can be decomposed into a direct sum of all irreducible representations of $\Ga$. This implies that the adjacency matrix of any Cayley graph $\Cay(\Ga, S)$ can be decomposed into the direct sum $d_1 \rho_1(S) \oplus d_2 \rho_2(S) \oplus \cdots \oplus d_k \rho_k(S)$, where $\{\rho_1, \rho_2, \ldots, \rho_k\}$ is a complete set of inequivalent irreducible matrix representations of $\Ga$, and for each $i$, $d_i$ is the dimension of $\rho_i$ and $\rho_i(S) = \sum_{s \in S} \rho_{i}(s)$. (See, for example, \cite[Proposition 7.1]{Krebs11}.) Thus the multiset of eigenvalues of $\Cay(\Ga, S)$ is the union of the multisets of eigenvalues of $\rho_i(S)$, for $1 \le i \le k$. In particular, the second largest eigenvalue $\l_{2}(\Cay(\Ga, S))$ is obtained by some $\rho_i$ other than the trivial representation of $\Ga$, meaning that $\l_{2}(\Cay(\Ga, S))$ is equal to the largest eigenvalue $\l_{1}(\rho_i(S))$ of $\rho_i(S)$. Recall that the irreducible representations of $S_n$ are in one-to-one correspondence with the partitions of integer $n$, and so we can index the former by the latter. In particular, the standard representation of $S_n$ corresponds to partition $(n-1,1)$ and hence is given by $\rho_{(n-1,1)}$. Since the defining representation of $S_n$ can be decomposed into the direct sum of one standard representation and one trivial representation, for any set $T$ of transpositions in $S_n$, the algebraic connectivity of $G_T$ is given by $\l_{1}(\rho_{(n-1,1)}(T))$. Therefore, Aldous' spectral gap theorem can be restated as follows.

\begin{thm}
\emph{(\cite[Theorem 1.1]{CLR09}, unweighted version)}
\label{thm:AldThm2}
For any set $T$ of transpositions in $S_n$,
$$
\l_{2}(\Cay(S_n, T)) = \l_{1}(\rho_{(n-1,1)}(T)).
$$
\end{thm}

Denote by ${\mathsf{eight}}_n$ the set of the following eight partitions of $n$: $(n-1, 1), (n-2, 2), (n-3, 3), (n-4, 4), (n-3, 2, 1), (n-1, 1)^t, (n-2, 2)^t, (n-2, 1, 1)^t$, where $t$ indicates the conjugate of a partition obtained from transposing the corresponding Young diagram. In \cite{ParzanchevskiP20}, Parzanchevski and Puder proved that for any conjugacy class $T$ of $S_n$, the largest eigenvalue of $\Cay(S_n, T)$ other than its degree $|T|$, written $\l(\Cay(S_n, T))$, is obtained by one of the representations indexed by these partitions. Note that $\l(\Cay(S_n, T)) = \l_{2}(\Cay(S_n, T))$ if and only if $\Cay(S_n, T)$ is connected.

\begin{thm}
\emph{(\cite[Theorem 1.7]{ParzanchevskiP20})}
\label{thm:SingleClass}
There exists a positive integer $N_0$ such that for every $n \ge N_0$, if $T$ is a conjugacy class of $S_n$, then
$$
\l(\Cay(S_n, T)) = \max_{\gamma \in {\mathsf{eight}}_n} \l_{1}(\rho_{\gamma}(T)).
$$
\end{thm}

Parzanchevski and Puder also obtained an asymptotic upper bound on the second largest eigenvalue of a general weighted normal Cayley graph on $S_n$. See \cite[Theorem 1.10]{ParzanchevskiP20} for details.

Recently, Li, Xia and Zhou \cite{LiXZ21} proved the following three results, in which the Cayley graphs involved are all normal. A Cayley graph $\Cay(S_n,S)$ on $S_n$ is said to have the \emph{Aldous property} \cite{LiXZ21} if its largest eigenvalue $\l(\Cay(S_n,S))$ strictly smaller than $|S|$ is attained by the standard representation of $S_n$ (that is, $\l(\Cay(S_n,S)) = \lambda_1(\rho_{(n-1,1)}(S))$). Recall that the support of a permutation $\s \in S_n$, $\mathrm{supp}(\s)$, is the set of elements of $[n]$ not fixed by $\s$. For $\emptyset \ne I \subseteq \{2,3,\ldots,n-1, n\}$ and $2\leq k\leq n$, set
$$
T(n,I) = \{\sigma\in S_n: |\mathrm{supp}(\sigma)|\in I\}
$$
and
$$
T(n,k)=\{\sigma \in S_n: 2\leq |\mathrm{supp}(\sigma)|\leq k\}.
$$

\begin{thm}
\emph{(\cite{LiXZ21})}
\label{thm:LiXZ21b}
There exists a positive integer $N$ such that for every $n\ge N$ and any conjugacy class $S$ of $S_n$, the normal Cayley graph $\Cay(S_n, S)$ has the Aldous property if and only if $2 \le |\mathrm{supp}(\sigma)| \le n-2$ for some (and hence all) $\sigma \in S$.
\end{thm}

\begin{thm}
\emph{(\cite{LiXZ21})}
\label{thm:LiXZ21a}
There exists a positive integer $N$ such that for every $n\ge N$ and any $\emptyset \ne I \subset \{2,3,\ldots,n-1, n\}$ with $|I \cap \{n-1, n\}| \ne 1$, the normal Cayley graph $\Cay(S_n,T(n,I))$ has the Aldous property if and only if $I \cap \{n-1, n\} = \emptyset$.
\end{thm}

\begin{thm}
\emph{(\cite{LiXZ21})}
\label{thm:LiXZ21}
There exists a positive integer $N$ such that for every $n\ge N$ and any $2\leq k\leq n$, the connected normal Cayley graph $\Cay(S_n,T(n,k))$ has the Aldous property.
\end{thm}

Note that the case $2\leq k\leq n-2$ in Theorem \ref{thm:LiXZ21} is covered by the sufficiency part of Theorem \ref{thm:LiXZ21a}. The fact that $\mathrm{Cay}(S_n,T(n,n-1))$ has the Aldous property can be derived from \cite[Theorem 7.1]{Renteln07}. So Theorem \ref{thm:LiXZ21} can be viewed as a corollary of Theorem \ref{thm:LiXZ21a} and \cite[Theorem 7.1]{Renteln07}.

As noticed in \cite{LiXZ21}, Theorem \ref{thm:LiXZ21a} together with some known results from \cite{Deng11, KuLW16a, Renteln07} implies that, for sufficiently large $n$, among the $k$-point-fixing graphs $\Cay(S_n, S(n, k))$ ($0\leq k\leq n-2$), $\Cay(S_n, S(n, 0))$ (that is, the derangement graph $\Cay(S_n, \mathscr{D}_n)$) and $\Cay(S_n, S(n, 1))$ are the only graphs without the Aldous property, where, as in Section \ref{subsec:sn}, $S(n, k) = T(n,\{n-k\})$ is the set of permutations of $S_n$ fixing exactly $k$ points. Combining this with Theorems \ref{thm:LiXZ21a} and \ref{thm:LiXZ21}, we see that, as far as the Aldous property of $\mathrm{Cay}(S_n,T(n,I))$ for sufficiently large $n$ is concerned, the only unsettled case is the one in which $|I| \ge 2$, $|I \cap \{n-1, n\}| = 1$ and $I \neq \{2,3,\ldots,n-1\}$. Two problems arising from this case were posed in \cite{LiXZ21}; see Problems \ref{prob:LiXZa} and \ref{prob:LiXZb} in Section \ref{sec:open}.

It is natural to consider whether Aldous' spectral gap theorem holds for some groups other than the symmetric groups. This line of research was undertaken by Cesi \cite{Cesi20} for the Weyl group $W(B_n)$ and by Kassabov \cite{Kassabov11} in a study of the Kazhdan property T of groups. A restatement of \cite[Theorem 1.1]{CLR09} is needed before we can state their results. Let $\Ga$ be a finite group, and let $w = \sum_{g \in \Ga} w_{g} g \in \mathbb{C}\Ga$ be such that all coefficients $w_g$ are real, nonnegative, and symmetric in the sense that $w_{g^{-1}} = w_g$ for $g \in G$. The weighted Cayley graph $\Cay(\Ga, w)$ is defined as the undirected graph with vertex set $\Ga$ and edges $\{g, hg\}$, $g, h \in \Ga$, such that each edge $\{g, hg\}$ carries weight $w_h$. If we identify $w$ with its support $S(w) = \{h \in \Ga: w_h \ne 0\}$, then $\Cay(\Ga, w)$ can be identified with the weighted Cayley graph $\Cay(\Ga, S(w))$ where each edge $\{g, hg\}$ ($g \in \Ga$ and $h \in S(w)$) carries weight $w_h$. The \emph{Laplacian} of $\Cay(\Ga, w)$ is the operator $\De_{\Cay(\Ga, w)}$ acting on the set $\mathbb{C}^{\Ga}$ of functions $f: \Ga \rightarrow \mathbb{C}$ in such a way that $(\De_{\Cay(\Ga, w)}f)(g) = \sum_{h \in \Ga} w_{h} (f(g)-f(hg))$, $g \in \Ga$. It is readily seen that $\De_{\Cay(\Ga, w)} = D - A$, where $D = \left(\sum_{g \in \Ga} w_{g}\right)I$ with $I$ the $|\Ga| \times |\Ga|$ identity matrix, and $A$ is the weighted adjacency matrix of $\Cay(\Ga, w)$. It is also easy to see that $\De_{\Cay(\Ga, w)} = \sum_{h \in \Ga} w_{h} (I - \l_{\text{reg}}(h))$, where $\l_{\text{reg}}$ is the left regular representation of $\Ga$ acting on $\mathbb{C}^{\Ga}$ defined by $(\l_{\text{reg}}(h)f)(g) = f(h^{-1}g)$, $g, h \in \Ga$. Denote by $\mu_{\Ga}(w)$ the second smallest eigenvalue of $\De_{\Cay(\Ga, w)}$, with multiple eigenvalues counted separately. Then $\mu_{\Ga}(w) > 0$ if and only if $\Cay(\Ga, w)$ is connected (that is, $S(w)$ generates $\Ga$). In general, for any representation $\rho$ of $\Ga$ on a complex vector space $V$ and any real, nonnegative and symmetric $w \in \mathbb{C}\Ga$, the representation Laplacian $\De_{\Ga}(w, \rho)$ is defined \cite{Cesi20} as the linear operator $\De_{\Ga}(w, \rho) = \sum_{h \in \Ga} w_{h} (I_V - \rho(h))$ on $V$, where $I_V$ is the identity on $V$. Denote by $\mu_{\Ga}(w, \rho)$ the second smallest eigenvalue of $\De_{\Ga}(w, \rho)$. In particular, we have $\De_{\Cay(\Ga, w)} = \De_{\Ga}(w, \l_{\text{reg}})$ and $\mu_{\Ga}(w) = \mu_{\Ga}(w, \l_{\text{reg}})$. The weighted version of Aldous' spectral gap theorem can be restated as follows.

\begin{thm}
\emph{(\cite{CLR09}; see \cite[Theorem 1.1]{Cesi20})}
\label{thm:AldThm1}
Let $w \in \mathbb{C}S_n$ be such that $S(w) \subseteq \{(i,j): 1 \le i < j \le n\}$ and all coefficients are nonnegative and symmetric. (That is, $w$ is of the form $w = \sum_{(i, j)} b_{ij} (i, j) \in \mathbb{C}S_n$, where $b_{ij} = b_{ji} \ge 0$.) Then
$$
\mu_{S_n}(w) = \mu_{S_n}(w, \rho_{n}^{0}),
$$
where $\rho_{n}^{0}$ is the defining representation of $S_n$ associated with the natural
action of $S_n$ on $\{1, 2, \ldots, n\}$.
\end{thm}

This version of \cite[Theorem 1.1]{CLR09} inspires the question of whether a similar result holds for some other group-representation pairs $(\Ga, \rho)$. In \cite{Cesi20}, Cesi proved that this is indeed the case for the Weyl group $W_n = W(B_n)$ and a certain $2n$-dimensional permutation representation $\mathbf{P}_n$ of $W_n$ associated with a specific action of $W_n$ with degree $2n$. One may define $W_n$ as the subgroup of $\GL(n,\mathbb{C})$ consisting of all $n \times n$ matrices which have exactly one non-zero entry in each row and each column, and this non-zero entry is either $1$ or $-1$. Then $S_n$ is embedded into $W_n$ as the subgroup of all matrices with nonnegative entries, and $W_n$ is embedded into $S_{2n}$ via $\mathbf{P}_n$ as a subgroup of signed permutations. For $1 \le i \le n$, let $s_{\{i\}} = \diag(1, \ldots, 1, -1, 1, \ldots, 1)$ where the unique $-1$ occurs as the $i$th diagonal entry.

\begin{thm}
\emph{(\cite[Theorem 1.2]{Cesi20})}
\label{thm:Weyl}
Let $w \in \mathbb{C}W_n$ be given by
$$
w = \sum_{i=1}^n a_{i} s_{\{i\}} + \sum_{(i, j)} b_{ij} (i, j),\;\, a_{i} \ge 0, b_{ij} \ge 0,
$$
where the second sum is running over all transpositions $(i, j)$ in $S_n$. Then
$$
\mu_{W_n}(w) = \mu_{W_n}(w, \mathbf{P}_n).
$$
\end{thm}

A conjecture generalizing this result was posed by Cesi in \cite[Conjecture 5.2]{Cesi20}. See Problem \ref{pb:Weyl} in Section \ref{sec:open}.

It is well known that any Coxeter group generated by $S = \{s_1, s_2, \ldots, s_{l}\}$ has a defining representation on an $l$-dimensional vector space such that each $s_i$ acts as a reflection with respect to a hyperplane. The following result was obtained by Kassabov \cite{Kassabov11} as an application of a result on the Kazhdan property T to spectral gaps of random walks on finite Coxeter groups. It shows that Aldous' spectral gap conjecture holds for every finite Coxeter group $\Ga$ generated by $S$ and its defining representation in the special case when $w = \sum_{s \in S} s \in \mathbb{C}\Ga$.

\begin{thm}
\emph{(\cite[Theorem 1.3]{Kassabov11}; see also \cite[Theorem 5.1]{Cesi16})}
\label{thm:Weyl}
Let $\Ga$ be a finite Coxeter group with Coxeter generating set $S$, and let $\rho$ be the defining representation of $\Ga$. Let $w = \sum_{s \in S} s \in \mathbb{C}\Ga$. Then
$$
\mu_{\Ga}(w) = \mu_{\Ga}(w, \rho).
$$
\end{thm}

The following result for general Coxeter groups was obtained by Akhiezer in \cite{Akhiezer07}.

\begin{thm}
\emph{(\cite[Theorem]{Akhiezer07})}
\label{thm:Akhiezer07}
Let $(W, S)$ be a finite Coxeter system, $S = \{s_1, s_2, \ldots, s_l\}$ the set of Coxeter generators, $h$ the Coxeter number of $(W, S)$, and $m_1, m_2, \ldots, m_l$ the exponents, $0 < m_1 \le m_2 \le \cdots \le m_l < h$. Then the numbers
$$
l - 2 + 2\cos \frac{\pi m_i}{h},\, i = 1, 2, \ldots, l
$$
are eigenvalues of $\Cay(W, S)$. Moreover, if $W$ is irreducible, then each of these eigenvalues has multiplicity at least $l$.
\end{thm}

The following result follows from Theorem \ref{thm:Akhiezer07} immediately.

\begin{cor}
\emph{(\cite[Corollary]{Akhiezer07})}
Let $(W, S)$ be a finite irreducible Coxeter system with $l = |S|$ generators and Coxeter number $h$. Then
$$
\l_{2}(\Cay(W, S)) \ge l - 2 + 2\cos \frac{\pi}{h}.
$$
\end{cor}

Eigenvalues of Cayley graphs on Coxeter groups have also been studied in \cite{IvrissimtzisP13} in the context of spectral representations associated with random walks on vertex-transitive graphs.

Recall from Section \ref{subsec:sn} that the complete transposition graph $T_n$ is the Cayley graph on the symmetric group $S_n$ with connection set consisting of all transpositions in $S_n$. The eigenvalues of $T_n$ were obtained in \cite{KalpakisYesha97}; see Theorem \ref{thm:KalpakisYesha97}. The second largest eigenvalue of $T_n$ was used to obtain the exact value of the bisection width of $T_n$, answering an open question posed by F. T. Leighton.


\subsection{Normal Cayley graphs on highly transitive groups}

A partition $\PP = \{V_1, \ldots, V_m\}$ of the vertex set of a graph $G$ is called an \emph{equitable partition} of $G$ if for each pair $i, j$ with $1 \le i, j \le m$, every vertex of $V_i$ has the same number (say, $b_{ij}$) of neighbours in $V_j$; the matrix $B_{\PP} = (b_{ij})$ is called the \emph{quotient matrix} of $G$ with respect to $\PP$. It is known that the eigenvalues of $B_{\PP}$ are eigenvalues of $G$  (see, for example, \cite{GodsilG2001}).

Let $\Ga$ be a group acting transitively on $[n] = \{1, 2, \ldots, n\}$. Denote by $\Ga_i$ the stabilizer of $i \in [n]$ in $\Ga$. Then $|\Ga:\Ga_i| = n$ as $\Ga$ is transitive. Let $\Ga = \Ga_{1,i} \cup \cdots \cup \Ga_{n, i}$ be the left coset decomposition of $\Ga$, where each $\Ga_{j,i}$ is a left coset of $\Ga_i$ in $\Ga$. Then for any Cayley graph $G = \Cay(\Ga,  S)$ on $\Ga$, $\PP_i = \{\Ga_{1,i}, \ldots, \Ga_{n, i}\}$ is an equitable partition of $G$ with corresponding quotient matrix $B_{\PP_i} = (b_{ij})$ given by $b_{ij} = |S \cap \Ga_{j,i}|$. Fix $k \in [n]$. It is not difficult to see that for all $i \in [n]$ the subgraph $G[\Ga_{k,i}]$ of $G$ induced by $\Ga_{k,i}$ is isomorphic to the Cayley graph $\Cay(\Ga_k, S \cap \Ga_k)$, and the spanning subgraph of $G$ induced by those edges with end-vertices in distinct  $\Ga_{k,i}, \Ga_{k,i'}$, $i \ne i'$, is isomorphic to the Cayley graph $\Cay(\Ga, S \setminus (S \cap \Ga_k))$ (see \cite[Lemma 2.5]{HuangHC18}). Based on these observations the following result was proved by Huang, Huang and Cioab\u{a} in \cite{HuangHC18}.

\begin{thm}
\emph{(\cite[Theorem 2.6]{HuangHC18})}
\label{thm:HandHC}
Let $\Ga$ be a group acting transitively on $[n] = \{1, 2, \ldots, n\}$, and let $G = \Cay(\Ga,  S)$ be a Cayley graph on $\Ga$. Then $\PP_i$ is an equitable  partition of $G$ and the corresponding quotient matrix $B_{\PP_i}$ is symmetric and independent of the choice of $i \in [n]$. Moreover, if $\l$ is an eigenvalue of $G$ which is not an eigenvalue of $B_{\PP_i}$, then for each $k \in [n]$ we have
$$
\l \le \l_2(\Cay(\Ga_k, S \cap \Ga_k)) + \l_2(\Cay(\Ga, S \setminus (S \cap \Ga_k))).
$$
\end{thm}

Using this result, Huang, Huang and Cioab\u{a} computed explicitly the second largest eigenvalue of a normal Cayley graph on a highly transitive group under a certain condition; see \cite[Theorem 3.5]{HuangHC18} for details. This result was then used to obtain the second largest eigenvalues of a majority of normal Cayley graphs on $S_n$ such that the connection set contains no element with support of size greater than $5$; see \cite[Theorem 4.1]{HuangHC18}. As corolloaries, the following known results for $n \ge 3$ were recovered (see parts (a) and (b) of Theorem \ref{Aldous-conj-spec}): If $G_T$ is a complete graph, then the spectral gap of $\Cay(S_n, T)$ is equal to $n$ \cite{Diaconis81}; if $G_T$ is a star, then the spectral gap of $\Cay(S_n, T)$ is equal to $1$ \cite{FOW85}.

In \cite{Dai18}, Dai introduced a variant of the family of Johnson graphs, called the full-flag Johnson graphs, and discussed their combinatorial properties. He proved among other things that for $n > d \ge 1$ the full-flag Johnson graph $FJ(n,d)$ is isomorphic to the Cayley graph on $S_n$ with connection set the set of $(n-d)$-reducible permutations of $S_n$, where for $1 \le m \le n-1$ a permutation $\sigma \in S_n$ is called \emph{$m$-reducible} if $m$ is the largest integer such that $[n]$ can be partitioned into $m$ parts $I_1, \ldots, I_m$ such that each part $I_i$ consists of consecutive integers in $[n]$ and is fixed setwise by $\sigma$. The eigenvalues of $FJ(n,1)$ were studied in \cite{Dai18}, where it was conjectured that the second largest eigenvalue of $FJ(n,1)$ is equal to that of a specific tridiagonal matrix of order $n$. In \cite{HuangHC18}, Huang, Huang and Cioab\u{a} noted that this conjecture is true and follows easily from Theorem \ref{thm:AldThm} and an argument using Theorem \ref{thm:HandHC}.


\subsection{Three Cayley graphs on alternating groups}

In \cite{HuangH18}, Huang and Huang determined the second largest eigenvalues of three Cayley graphs on the alternating group $A_n$, where $n \ge 3$. Set $T_1 = \{(1, 2, i), (1, i, 2): 3 \le i \le n\}$, $T_2 = \{(1, i, j), (1, j, i): 2 \le i < j \le n\}$ and $T_3 = \{(i, j, k), (i, k, j): 1 \le i < j < k \le n\}$. Call $AG_n = \Cay(A_n, T_1)$, $EAG_n = \Cay(A_n, T_2)$ and $CAG_n = \Cay(A_n, T_3)$ the \emph{alternating group graph}, \emph{extended alternating group graph} and \emph{complete alternating group graph}, respectively. Since $T_1 \subseteq T_2 \subseteq T_3$, $AG_n$ is a spanning subgraph of $EAG_n$ and in turn $EAG_n$ is a spanning subgraph of $CAG_n$, and all three graphs are connected as $T_1$ is a generating set of $A_n$.

\begin{thm}
\emph{(\cite[Theorems 2.3, 3.2 and 3.4]{HuangH18})}
\label{thm:HandH}
Let $n \ge 4$. Then $\l_{2}(AG_n) = 2n-6$, $\l_{2}(EAG_n) = n^2 - 5n + 5$ and $\l_{2}(CAG_n) = n(n-2)(n-4)/3$.
\end{thm}


\subsection{Cayley graphs on abelian groups}

In \cite{FMT06}, it was proved that for any abelian group $\Ga$ with order $n$ and any $k$-regular Cayley graph $\Cay(\Ga, S)$ on $\Ga$ we have $\lambda_2(\Cay(\Ga,S)) \geq k-O(kn^{-\frac{4}{k}})$, where the constant in the big $O$ term does not depend on $k$ and $n$. More explicitly, the following result was proved by Friedman, Murty and Tillich in \cite{FMT06}.

\begin{thm}
\emph{(\cite[Theorem 6]{FMT06})}
Let $k \ge 1$ be a fixed integer and let $\mu \in [0, 1]$. Then there exists a constant $C_k$ satisfying
$$
C_k \le \frac{1}{2} (1-\mu)^{-1} \mu^{-\frac{4}{k}} \pi^2
$$
such that for any abelian group $\Ga$ of order $n$ and any $k$-regular Cayley graph $\Cay(\Ga, S)$ on $\Ga$ we have
$$
\lambda_2(\Cay(\Ga,S)) \geq k - C_k k n^{-\frac{4}{k}} + o(n^{-\frac{4}{k}}).
$$
In particular, as $k \rightarrow \infty$ we can take $C_k \preccurlyeq \pi^2/2$.
\end{thm}

It was also shown in \cite{FMT06} that for any fixed $k$, the lower bound above cannot be improved for large odd primes $n$; that is, if $n$ is an odd prime, then most $k$-regular graphs on $n$ vertices have the second largest eigenvalue at most $k - \Omega(k n^{-\frac{4}{k}})$.

A well-known result of J.-P. Serre \cite{Davidoff03} asserts that for any $\ve > 0$ and integer $k \ge 1$ there exists a constant $c = c(\ve, k) > 0$ such that any $k$-regular graph $G$ of order $n$ has at least $cn$ eigenvalues no less than $(2-\ve) \sqrt{k-1}$. Elementary proofs of this result were given in \cite{Nilli04} and \cite{Cioaba06}, and a result of this type for Cayley graphs on abelian groups was established by Cioab\v{a} in \cite{Cioa06}.

\begin{thm}
\emph{(\cite[Theorem 1.2]{Cioa06})}
For any $\ve > 0$ and integer $k \ge 3$, there exists a constant $C = C(\ve, k) > 0$ such that any Cayley graph $\Cay(\Ga, S)$ on any abelian group $\Ga$ with order $n$ has at least $Cn$ eigenvalues no less than $k - \ve$.
\end{thm}


\subsection{A family of Cayley graphs with small second largest eigenvalues}
\label{subsec:ssle}

In \cite{Fri95}, Friedman studied the second largest eigenvalue of two families of directed Cayley graphs and one family of Cayley graphs. We only mention the third family here and leave the other two to Section \ref{subsec:dissle}. Let $p \equiv 3$ ($\mod~4$) be a prime. Let $\AGL(1,p)$ be the group of affine linear transformations $t_{a, b}$ of $\ZZZ_p$, where $t_{a, b}: \ZZZ_p \rightarrow \ZZZ_p, x \mapsto ax+b$ for $a \in \ZZZ_p^*$ and $b \in \ZZZ_p$. In \cite{Fri95}, $\mathrm{SQRT}(p)$ was defined to be the Cayley graph on $\AGL(1,p)$ with respect to the connection set $\{t_{r^2, r}, t_{-r^2, r}: r \in \ZZZ_p^*\}$. Clearly, this is an undirected graph of order $p(p-1)$ and degree $2(p-1)$. In \cite{Fri95}, Friedman proved the following result.

\begin{thm}
\emph{(\cite[Theorem 1.3]{Fri95})}
Let $p \equiv 3$ ($\mod~4$) be a prime. The graph $\mathrm{SQRT}(p)$ above has second largest eigenvalue in absolute value at most $2\sqrt{p}$.
\end{thm}


\subsection{First type Frobenius graphs}
\label{subsec:inter-net}

Let $N \rtimes H$ be a finite Frobenius group (see Section \ref{sec:BDRCG}). A Cayley graph $\Cay(N, S)$ on $N$ is called \cite{FLP} a {\em first type Frobenius graph} if $S = a^H$ for some $a \in N$ such that $\la a^H \ra = N$ and either $|H|$ is even or $a$ is an involution, where $x^H = \{h^{-1}xh: h \in H\}$ for $x \in N$. In \cite{Zhou2009}, it was shown that first type Frobenius graphs exhibit ``perfect" routing properties in some sense. To be precise, let us define a \emph{shortest path routing} of a graph $G$ to be a set of oriented shortest paths in $G$ which contains exactly one oriented shortest path from $u$ to $v$ for each ordered pair of distinct vertices $(u, v)$ of $G$. The \emph{load} of an arc under such a routing is the number of paths in the routing that traverse the arc in its direction, and the maximum load among all arcs is called the maximum-arc-load of the routing. The \emph{minimal arc-forwarding index} \cite{Hey97} of a graph $G$, $\overrightarrow{\pi}_{m}(G))$, is the minimum of the maximum-arc-loads over all shortest path routings of $G$. Obviously,
$$
\overrightarrow{\pi}_{m}(G)) \ge \frac{\sum_{u, v \in V(G)} d(u,v)}{2|E(G)|},
$$
where $d(u, v)$ is the distance in $G$ between $u$ and $v$. In \cite[Theorem 6.1]{Zhou2009}, it was proved that any first type Frobenius graph attains this lower bound and thus has the smallest possible minimal arc-forwarding index.

On the other hand, it is known \cite[Corollary 1]{DS} that the second largest eigenvalue of random walks on a connected graph $G$ is bounded from above by
$$
1 - \frac{2|E(G)|}{\Delta(G)^2 \diam(G)  \overrightarrow{\pi}_{m}(G)},
$$
where $\diam(G)$ and $\Delta(G)$ are the diameter and maximum degree of $G$, respectively. In the case when $G$ is regular, this can be translated into an upper bound on the second largest eigenvalue $\l_2(G)$ of $G$. This upper bound, together with \cite[Theorem 6.1]{Zhou2009} and \cite[Theorem 1.6]{FLP}, implies the following result.

\begin{thm}
\emph{(\cite[Corollary 6.4]{Zhou2009})}
\label{cor:eigen}
Let $N \rtimes H$ be a finite Frobenius group. Let $\Cay(N, S)$ be a first kind Frobenius graph of $N \rtimes H$. Let $d$ be the diameter of $\Cay(N, S)$ and $n_i$ the number of $H$-orbits on $N$ at distance $i$ in $\Cay(N, S)$ to the identity element of $N$, $1 \le i \le d$. Then
$$
\l_{2}(\Cay(N, S)) \le |H| - \frac{|N|}{d\sum_{i=1}^{d}in_i}.
$$
\end{thm}


\subsection{Distance-$j$ Hamming graphs and distance-$j$ Johnson graphs}
\label{subsec:Hamming-Johnson}

In this section we restrict our attention to two families of graphs studied by Brouwer, Cioab\u{a}, Ihringer and McGinnis in \cite{BrouwerCIM18}. One of the motivations of their study was that computing the smallest eigenvalues of these graphs implies that the approximation ratio of the max-cut semidefinite programming algorithm of Goemans and Williamson \cite{GoemansW95} is best possible.

Let $q \ge 2$ and $d \ge 1$ be integers. The \emph{distance-$j$ Hamming graph} $H(d,q,j)$, where $0 \le j \le d$, is the graph with vertex set $\ZZZ_q^d$ in which two vectors are adjacent if and only if they differ in exactly $j$ coordinates. Obviously, $H(d,q,j)$ is a Cayley graph for each $j$ and $H(d,q,1)$ is the Hamming graph $H(d, q)$ in the usual sense. Using Corollary \ref{BabaiS1a} and Krawtchouk or Kravchuk polynomials, it can be shown (see \cite{BrouwerCIM18} and \cite{Martin09}) that the eigenvalues of $H(d,q,j)$ are given by
$$
K_{j}(i) = \sum_{h=0}^j (-1)^h (q-1)^{j-h} {\choose{i}{h}} {\choose{d-i}{j-h}},\;\, 0 \le i \le d.
$$
Solving a conjecture of van Dam and Sotirov \cite{vanDamS16}, Brouwer, Cioab\u{a}, Ihringer and McGinnis \cite{BrouwerCIM18} proved the following result.

\begin{thm}
\emph{(\cite{BrouwerCIM18})}
\label{dist-j-hamming}
Let $q \ge 2$, $d \ge 1$ and $j \ge d - \frac{d-1}{q}$, where $j$ is even when $q = 2$. Then the smallest eigenvalue of $H(d,q,j)$ is $K_j(1)$.
\end{thm}

In \cite{BrouwerCIM18}, Brouwer, Cioab\u{a}, Ihringer and McGinnis also proved that in most cases $K_j(1)$ is not only the smallest eigenvalue, but also the second largest eigenvalue in absolute value, the only exception being the case when $(d, q) = (4, 3)$. Earlier van Dam and Sotirov's conjecture was confirmed for $q=2$ by other researchers (see \cite{BrouwerCIM18} for details).

Let $n \ge 2$ and $d \ge 1$ be integers. The \emph{distance-$j$ Johnson graph} $J(n,d,j)$, where $0 \le j \le d$, is the graph with vertices the $d$-subsets of a fixed $n$-set such that two $d$-subsets are adjacent if and only if they meet at a $(d-j)$-subset. In particular, $J(n,d,1)$ is the Johnson graph $J(n,d)$ in the usual sense and $J(n,d,d)$ is commonly known as the Kneser graph. The eigenvalues of $J(n,d,j)$ are given by
$$
E_{j}(i) = \sum_{h=0}^i (-1)^{i-h} {\choose{i}{h}} {\choose{d-h}{j}} {\choose{n-d-i+h}{n-d-j}},\;\, 0 \le i \le d.
$$
In \cite{BrouwerCIM18}, Brouwer, Cioab\u{a}, Ihringer and McGinnis proved the following result.

\begin{thm}
\emph{(\cite[Theorem 3.10]{BrouwerCIM18})}
\label{dist-j-johnson}
Let $j > 0$. Then $E_j(1)$ is the smallest eigenvalue of $J(n,d,j)$ if and only if $j(n - 1) \ge d(n-d)$. In this case $E_j(1)$ is also the second largest eigenvalue in absolute value among the eigenvalues of $J(n,d,j)$.
\end{thm}

In particular, if $n = 2d$ and $j > \frac{d}{2}$, then the smallest eigenvalue of $J(n,d,j)$, and the second largest eigenvalue of $J(n,d,j)$ in absolute value, is $E_j(1)$, solving a conjecture of Karloff \cite{Karloff99}.

Theorem \ref{dist-j-johnson} is relevant to this survey as some (but not all) distance-$j$ Johnson graphs are Cayley graphs. See \cite{JonesR16} for detailed information about when a distance-$j$ Johnson graph is a Cayley graph.


\subsection{Algebraically defined graphs}
\label{subsec:alg-def}

Let $q = p^e$ be a prime power, where $p$ is a prime and $e \ge 1$ an integer. Let $k \ge 2$ be an integer. For $2 \le i \le k$, take an arbitrary polynomial $h_i$ in $2(i-1)$ indeterminants over $\FFF_q$. Define $B\Ga_k = B\Ga(q; h_2, \ldots, h_k)$ to be the bipartite graph with biparition $\{P_k, L_k\}$, where each of $P_k$ and $L_k$ is a copy of $\FFF_q^k$, such that $(p_1, p_2, \ldots, p_k) \in P_k$ and $(l_1, l_2, \ldots, l_k) \in L_k$ are adjacent if and only if their coordinates satisfy
$$
p_i + l_i = h_{i}(p_1, l_1, p_2, l_2, \ldots, p_{i-1}, l_{i-1}),\;\, 2 \le i \le k.
$$
One can see that $B\Ga_k$ is a $q$-regular graph with $2q^k$ vertices. It was introduced by Lazebnik and Woldar as a generalization of some graphs introduced earlier by Lazebnik and Ustimenko (see the survey papers \cite{LazebnikSW17, LazebnikW01} for details).

Three specifications of graphs $B\Ga_k$ have received special attention in recent years. In the case when $h_i = p_1 l_1^{i-1}$, $2 \le i \le k+1$, $B\Ga(q; h_2, \ldots, h_{k+1})$ is called the \emph{Wenger graph}, denoted by $W_k(q)$; the spectra of such graphs have been determined completely by Cioab\u{a}, Lazebnik and Li in \cite{CioabaLL14}. In the case when $h_i = p_1^{p^{i-2}} l_1$, $2 \le i \le k+1$, $B\Ga(q; h_2, \ldots, h_{k+1})$ is called the \emph{linearized Wenger graph}, denoted by $L_k(q)$; the eigenvalues of such graphs have been determined by Cao, Lu, Wan, Wang and Wang \cite{CaoLWWW15} and their multiplicities have been determined by Yan and Liu \cite{YanL17}. In the case when $h_2 = p_1 l_1, h_3 = p_1 l_2$ and for $4 \le i \le k$, $h_i = -p_{i-2}l_1$ when $i \equiv 0, 1~\mod~4$  and $h_i = p_{1}l_{i-2}$ when $i \equiv 2, 3~\mod~4$, the graph $B\Ga(q; h_2, \ldots, h_k)$ is a bipartite edge-transitive graph, denoted by $D(k,q)$. It was conjectured in \cite{Ustimenko03} (see also \cite{LazebnikSW17}) that the second largest eigenvalue of $D(k,q)$ is at most $2\sqrt{q}$. In the special case when $k = 4$, this conjecture was proved by Moorhouse, Sun and Williford in \cite{MoorhouseSW17}. It turns out that the distance-two graph $G^{(2)}$ of each of $G = W_{k}(q), L_{k}(q), D(4, q)$ is a Cayley graph, and this fact was used to bound the second largest eigenvalue of $G$. (In general, given a bipartite graph $G$ with bipartition $\{V_1, V_2\}$, the \emph{distance-two graph} $G^{(2)}$ of $G$ on $V_1$ is the graph with vertex set $V_1$ in which two vertices are adjacent if and only if they are at distance two in $G$.)

In a recent paper \cite{CioabaLS17}, Cioab\u{a}, Lazebnik and Sun introduced a family of Cayley graphs on certain abelian groups by directly generalizing the defining system of equations for $W_{k}(q)^{(2)}$ and $L_{k}(q)^{(2)}$. Let $k \ge 3$ and let $f_i, g_i \in \FFF_{q}[x]$, $3 \le i \le k$, be polynomials of degrees at most $q-1$ such that $g_{i}(-x) = -g_{i}(x)$ for each $i$. Define $S(k, q) = S(k, q; f_3, g_3, \ldots, f_k, g_k)$ to be the graph with vertex set $\FFF_q^k$ such that $(a_1, a_2, \ldots, a_k) \in \FFF_q^k$ and $(b_1, b_2, \ldots, b_k) \in \FFF_q^k$ are adjacent if and only if their coordinates satisfy
$$
b_i - a_i = g_{i}(b_1 - a_1) f_{i}\left(\frac{b_2 - a_2}{b_1 - a_1}\right),\;\, 3 \le i \le k.
$$
The requirement $g_{i}(-x) = -g_{i}(x)$ ensures that this is an undirected graph. It was noticed \cite{CioabaLS17} that $S(k, q)$ is a Cayley graph on the additive group of the vector space $\FFF_q^k$ with connection set $\{(a, au, g_{3}(a)f_{3}(u), \ldots,  g_{k}(a)f_{k}(u)): a \in \FFF_{q}^*, u \in \FFF_q\}$. Thus $S(k, q)$ is a $q(q-1)$-regular graph with $q^k$ vertices. With the help of Theorem \ref{BabaiS2}, Cioab\u{a}, Lazebnik and Sun \cite[Theorem 2.1]{CioabaLS17} determined the spectrum of $S(k, q)$. Using this they obtained the following sufficient condition for the graphs $S(k, q)$ to form a family of expander graphs.

\begin{thm}
\emph{(\cite[Theorem 2.2]{CioabaLS17})}
Let $k \ge 3$ and let $(q_{m})_{m \ge 1}$ be an increasing sequence of prime powers, where $q_m = p_m^{e_m}$ with $p_m$ a prime and $e_m \ge 1$ an integer. Let
$$
S(k, q_m) = S(k, q_m; f_{3,m}, g_{3,m}, \ldots, f_{k,m}, g_{k,m}).
$$
Let $d_{f}^{(m)}$ be the maximum degree of polynomials $f_{3,m}, \ldots, f_{k,m}$ and $d_{g}^{(m)}$ the maximum degree of polynomials $g_{3,m}, \ldots, g_{k,m}$. Suppose that $1 \le d_{f}^{(m)} = o_{m}(q_m)$, $1 \le d_{g}^{(m)} = o_{m}(\sqrt{q_m})$, $d_{g}^{(m)} < p_m$, and for all $m$ at least one of the following conditions is satisfied:
\begin{itemize}
\item[\rm (a)] the polynomials $1, x, f_{3,m}, \ldots, f_{k,m}$ are $\FFF_{q}$-linearly independent, and $g_{i,m}$ has linear term for all $i$, $3 \le i \le k$;
\item[\rm (b)] the polynomials $f_{3,m}, \ldots, f_{k,m}$ are $\FFF_{q}$-linearly independent, and there exists some $j$, $2 \le j \le d_{g}^{(m)}$, such that each $g_{i,m}$, $3 \le i \le k$, contains a term $c_{i,j}^{(m)}x^j$ with $c_{i,j}^{(m)} \ne 0$.
\end{itemize}
Then $S(k, q_m)$ is connected and
$$
\l_2(S(k, q_m)) = o_{m}(q_m^2).
$$
\end{thm}

In \cite{CioabaLS17}, Cioab\u{a}, Lazebnik and Sun also obtained better upper bounds on the second largest eigenvalue of $S(k, q)$ for two specifications of polynomials $f_i, g_i$. Define the \emph{exponential sum} of $f \in \FFF_{q}[x]$ by $\ve_f = \sum_{x \in \FFF_q} \om_{p}^{\Tr(f(x))}$, where as before $\Tr(x) = x + x^p + \cdots + x^{p^{e-1}}$ is the trace of $x \in \FFF_q$ over $\FFF_p$.

\begin{thm}
\emph{(\cite[Theorem 2.3]{CioabaLS17})}
Let $q$ be an odd prime power with $q \equiv 2~\mod~3$ and $4 \le k \le q+1$. Let $g_{i}(x) = x^3$ and $f_{i}(x) = x^{i-1}$ for each $i$, $3 \le i \le k$. Then $S(k,q)$ is connected and
$$
\l_{2}(S(k,q)) = \max\{q(k-3), (q-1)M_q\},
$$
where $M_q = \max_{a, b \in \FFF_q^*} \ve_{ax^3 + bx} \le 2 \sqrt{q}$.
\end{thm}

In particular, if $(q-1)M_q \le q(k-3)$, then $\l_{2}(S(k,q)) = q(k-3) < q(k-2) = \l_{2}(W_{k-1}(q)^{(2)})$.

\begin{thm}
\emph{(\cite[Theorem 2.4]{CioabaLS17})}
Let $q$ be an odd prime power with $q \equiv 2~\mod~3$ and $3 \le k \le e+2$. Let $g_{i}(x) = x^3$ and $f_{i}(x) = x^{p^{i-2}}$ for each $i$ with $3 \le i \le k$. Then $S(k,q)$ is connected and
$$
\l_{2}(S(k,q)) \le \max\{q(p^{k-3} - 1), (q-1)M_q\},
$$
where $M_q = \max_{a, b \in \FFF_q^*} \ve_{ax^3 + bx} \le 2 \sqrt{q}$.
\end{thm}

In particular, if $(q-1)M_q \le q(p^{k-3} - 1)$, then $\l_{2}(S(k,q)) = q(p^{k-3} - 1) < q(p^{k-2} - 1) = \l_{2}(L_{k-1}(q)^{(2)})$.


\section{Perfect state transfer in Cayley graphs}
\label{sec:pst}

Let $G$ be a graph. Define
\be
\label{eq:HGt}
H_G(t)=\exp(itA(G))=\sum_{k=0}^{\infty}\frac{\mathrm{i}^kt^kA(G)^k}{k!}
\ee
and call it the \emph{transition matrix} of $G$, where $A(G)$ is the adjacency matrix of $G$ and $\mathrm{i}=\sqrt{-1}$. We say that \emph{perfect state transfer} occurs from a vertex $u$ to another vertex $v$ in $G$ if there exists a time $\tau$ such that $|H_G(\tau)_{u,v}|=1$. If there exists a time $\tau$  such that $|H_G(\tau)_{u,u}|=1$, then $G$ is \emph{periodic} at vertex $u$ with \emph{period} $\tau$.  A graph is called \emph{periodic} if it is periodic at every vertex with the same period. As seen in \cite{Christandl04, Christandl05, Kay10}, these two concepts arise from quantum computing, in which a continuous-time quantum walk on a graph $G$ with transition matrix $H_G(t)$ plays a significant role. The reader is referred to \cite{Coutinho15, Godsil11, Godsil12} and the surveys \cite{Godsil12-1, Kay10, Stevanovic11} for background information and results on perfect state transfer and periodicity in graphs.

It is known that if $G$ has perfect state transfer from $u$ to $v$ at time $\tau$ then $G$ is periodic at both $u$ and $v$ with period $2\tau$ (see \cite{Christandl05} or \cite[Lemma 2.1]{Godsil11}). Moreover, there is an interesting relation between the eigenvalues of a graph $G$ and its periodicity (see \cite[Corollary 3.3]{Godsil11} or \cite[Theorem 2.3]{Zhou11}), which asserts that $G$ is periodic if and only if the eigenvalues of $G$ are all  integers or all rational multiples of $\sqrt{\Delta}$ for some square-free integer $\Delta$. Furthermore, if the second alternative holds, then $G$ is bipartite. This criterion shows that perfect state transfer can occur in integral graphs but not always (see \cite{Christandl04}). In addition, if $G$ is regular, then its largest eigenvalue is an integer and so the second alternative above cannot occur. Thus a regular graph is periodic if and only if it is integral (see \cite[Section 5]{Godsil12-1}). Hence all integral Cayley graphs are periodic. In particular, all integral Cayley graphs mentioned in Section \ref{sec:ICG} are periodic. Moreover, from the discussion above one can see that integral Cayley graphs are good candidates if one wants to construct regular graphs having perfect state transfer. In fact, as noted in \cite[Lemma 2.2]{TanFC17} and \cite[Theorem 3.5]{PalB16}, if a Cayley graph on an abelian group has perfect state transfer, then it must be integral.

In \cite[Theorem 4.3]{Coutinho15}, necessary and sufficient conditions for a graph belonging to an association scheme to admit perfect state transfer were obtained. In particular, this result applies to distance-regular graphs (see \cite[Corollaries 4.5 and 5.1]{Coutinho15}) and was used in \cite{Coutinho15} to determine which distance-regular graphs listed in \cite[Chapter 14]{BrouwerCN89} admit perfect state transfer. Note that many distance-regular graphs are Cayley graphs, and \cite[Corollaries 4.5 and 5.1]{Coutinho15} can be used to test when a distance-regular Cayley graph admits perfect state transfer.


\subsection{Perfect state transfer in Cayley graphs on abelian groups}
\label{PSTinACGs}

Let $\Gamma=\mathbb{Z}_{n_1}\oplus\cdots\oplus\mathbb{Z}_{n_d}$ be an abelian group, where each $n_i \ge 2$. It is well known that for any $\mathbf{x} = (x_1, \ldots,x_d)\in\Gamma$ the mapping $\chi_{\mathbf{x}}: \Gamma \rightarrow \mathbb{C}$ defined by
\begin{equation}
\label{eq:chixg}
\chi_{\mathbf{x}}(\mathbf{g}) = \exp\left(2\pi \mathrm{i} \sum_{j=1}^{d} \frac{x_j g_j}{n_j}\right),\, \mathbf{g} = (g_1, g_2, \ldots, g_d) \in \Ga
\end{equation}
is a character of $\Gamma$. Let $S$ be a subset of $\Gamma$ not containing the identity element of $\Gamma$ such that $S^{-1} = S$. By Corollary \ref{BabaiS1a},
the eigenvalues of $\Cay(\Gamma, S)$ are
$$
\alpha_{\mathbf{x}} = \sum_{\mathbf{g} \in S} \chi_{\mathbf{x}}(\mathbf{g}),\; \mathbf{x} \in \Gamma.
$$
Thus $\Cay(\Gamma, S)$ is integral if and only if $\alpha_{\mathbf{x}} \in \mathbb{Z}$ for all $\mathbf{x}\in\Gamma$.

The following result due to Tan, Feng and Cao tells us the period of an integral Cayley
graph on an abelian group at every vertex.

\begin{thm}
\emph{(\cite[Theorem 2.3]{TanFC17})}
Let $\Gamma = \mathbb{Z}_{n_1}\oplus\cdots\oplus\mathbb{Z}_{n_d}$ be an abelian group and $\Cay(\Gamma, S)$ an integral Cayley graph on $\Ga$, where each $n_i \ge 2$. Set $c=\gcd(|S|-\alpha_{\mathbf{x}}: \mathbf{x} \in \Gamma)$. Then $\Cay(\Gamma, S)$ is periodic at every vertex with period $2\pi l/c$ for each $l=1,2,\ldots$.
\end{thm}

The \emph{$2$-adic valuation} of rational numbers is the mapping $\nu_{2}:\mathbb{Q} \rightarrow \mathbb{Z}\cup \{\infty\}$ defined by
\begin{equation}
\label{eq:2adic}
\nu_{2}(0)=\infty,\ \nu_{2}\left(\frac{2^l a}{b}\right)=l,\ \text{where $a,b,l\in \mathbb{Z}$ with $ab$ odd}.
\end{equation}
Note that if $\Gamma$ contains involutions then $|\Gamma|$ is even. Moreover, $\chi_{\mathbf{a}}(\mathbf{g})=\pm1$ for any involution $\mathbf{a}$ and any element $\mathbf{g}$ of $\Gamma$. The following result gives a necessary and sufficient condition for a Cayley graph on an abelian group to admit perfect state transfer.

\begin{thm}\label{PSTAbelian11}
\emph{(\cite[Theorem 2.4]{TanFC17})}
Let $\Gamma = \mathbb{Z}_{n_1}\oplus\cdots\oplus\mathbb{Z}_{n_d}$ be an abelian group of order at least $3$ and $\Cay(\Gamma, S)$ a connected Cayley graph on $\Gamma$, where each $n_i \ge 2$. Then for any distinct $\mathbf{g}, \mathbf{h} \in \Gamma$, $\Cay(\Gamma, S)$ has perfect state transfer from $\mathbf{g}$ to $\mathbf{h}$ if and only if  the following conditions are satisfied:
\begin{itemize}
  \item[\rm (a)] $\Cay(\Gamma, S)$ is an integral graph (that is, $\alpha_{\mathbf{x}} \in \mathbb{Z}$ for all $\mathbf{x} \in \Gamma$);
  \item[\rm (b)] the element $\mathbf{a} := \mathbf{g}-\mathbf{h}$ is an involution of $\Gamma$;
  \item[\rm (c)] the values of $\nu_{2}(|S|-\alpha_{\mathbf{x}})$ for all $\mathbf{x} \in \chi_{\mathbf{a}}^{-1}(-1)$ are the same, say, $\rho$, and moreover $\nu_{2}(c_{\mathbf{a}}) \ge\rho+1$, where $c_{\mathbf{a}} = \gcd(|S|-\alpha_{\mathbf{x}}: \mathbf{x} \in \chi_{\mathbf{a}}^{-1}(1))$.
\end{itemize}
Moreover, if conditions {\rm (a)-(c)} are satisfied, then $\Cay(\Gamma, S)$ has perfect state transfer from $\mathbf{g}$ to $\mathbf{h}$ at time $(\pi/c)+(2\pi l/c)$, for $l=1,2,\ldots$, where $c=\gcd(|S|-\alpha_{\mathbf{x}}: \mathbf{x} \in \Gamma)$.
\end{thm}

By condition (b) in Theorem \ref{PSTAbelian11}, if  $\Cay(\Gamma, S)$ admits perfect state transfer, then $|\Gamma|$ must be even. It was proved in \cite{Petkovic11} that for circulant graphs the stronger necessary condition that $n$ is a multiple of $4$ should be satisfied. Recently, Tan, Feng and Cao proved that the same condition should be respected for all integral Cayley graphs on abelian groups.

\begin{thm}
\emph{(\cite[Theorem 3.5]{TanFC17})}
Let $\Cay(\Gamma, S)$ be an integral Cayley graph on an abelian group $\Gamma$ with $|\Gamma| \ge 6$ and $|\Gamma| \equiv 2~(\mod~4)$. Then $\Cay(\Gamma, S)$ has no perfect state transfer between any two distinct vertices.
\end{thm}


\subsection{A few families of Cayley graphs on abelian groups admitting perfect state transfer}


\subsubsection{Perfect state transfer in integral circulant graphs}

Perfect state transfer in integral circulant graphs has been studied extensively. The main problem, as raised by Angeles-Canul, Norton, Opperman, Paribello, Russell and Tamon in \cite{Angeles10}, is to characterize those integral circulant graphs which admit perfect state transfer. Partial results on this problem were obtained in several papers, including \cite{Angeles10, Basic09-2, Basic10, Basic09-1, Basic15, Petkovic11, Saxena07}, \cite[Section 7]{Godsil12-1} and \cite[Section 9]{Stevanovic11}. A complete characterization, stated below, was finally obtained by Ba\v{s}i\'{c} in \cite{Basic13}. Recall from Corollary \ref{cor:intICG} that a circulant graph is integral if and only if it is isomorphic to a gcd graph $\ICG(n, D)$ as defined in \eqref{eq:ICG}, where $D \subseteq D(n) \setminus \{n\}$ with $D(n)$ as defined in \eqref{eq:Dn}. Denote $kD = \{kd: d \in D\}$ for any positive integer $k$.


\begin{thm}\label{CIMainPST}\emph{(\cite[Theorem 22]{Basic13})}
The integral circulant graph $\ICG(n,D)$ has perfect state transfer if and only if $n\in 4\mathbb{N}$ and $D=D_2\cup 2D_2\cup 4D_2 \cup D_3 \cup\{n/2^a\}$, where $D_2=\{d\in D: n/d\in8\mathbb{N}+4\}\setminus\{n/4\}$, $D_3=\{d\in D: n/d\in8\mathbb{N}\}$, and $a\in\{1,2\}$.
\end{thm}

A weighted graph whose weighted adjacency matrix is circulant is called a \emph{weighted circulant graph}. In \cite{Basic14IS}, Ba\v{s}i\'c studied perfect state transfer and periodicity in weighted circulant graphs. He proved that a weighted circulant graph is periodic if and only if it is integral. He also gave a criterion for the existence of perfect state transfer in a weighted circulant graph in terms of its eigenvalues (see \cite[Theorem 8]{Basic14IS}). As applications he found some classes of weighted circulant graphs having perfect state transfer (see \cite[Theorems 10]{Basic14IS}) as well as ones admitting no perfect state transfer (see \cite[Theorems 13-15]{Basic14IS}).


\subsubsection{Perfect state transfer in NEPS of complete graphs}

Recall from Section \ref{HammingE} that cubelike graphs are precisely Cayley graphs on elementary abelian $2$-groups and are exactly NEPS graphs of copies of $K_2$. Specifically, every subset $S$ of $\mathbb{Z}_2^d \setminus \{\mathbf{0}\}$ gives rise to a cubelike graph $\Cay(\mathbb{Z}_2^d, S)$, where $\mathbf{0}$ is the zero element of $\mathbb{Z}_2^d$. The existence of perfect state transfer in cubelike graphs was studied in \cite{Bernasconi08, Cheung11, Christandl04}. See also \cite[Section 7]{Godsil12-1} and \cite[Section 8]{Stevanovic11}.

\begin{thm}
\emph{(\cite[Theorem 1]{Bernasconi08}, or \cite[Theorem 2.3]{Cheung11})}
Let $d \ge 1$ be an integer.
Let $S$ be a subset of $\mathbb{Z}_2^d \setminus \{\mathbf{0}\}$ and let $\mathbf{a}$ be the sum of the elements of $S$. If $\mathbf{a}\not= 0$, then perfect
state transfer occurs in $\Cay(\mathbb{Z}_2^d, S)$ from $\mathbf{u}$ to $\mathbf{u}+\mathbf{a}$ at time $\pi/2$, for any $\mathbf{u} \in \mathbb{Z}_2^d$. If $\mathbf{a}=0$, then $\Cay(\mathbb{Z}_2^d, S)$ is periodic with period $\pi/2$.
\end{thm}

The elements of $\mathbb{Z}_2^d$ can be treated as column vectors of the linear space $\mathbb{F}_2^d$. Given $S \subseteq \mathbb{Z}_2^d \setminus \{\mathbf{0}\}$, let $M(S)$ be the matrix with vectors in $S$ as its columns. The \emph{code} of $S$, $C(S)$, is defined as the row space of $M(S)$ in the linear space $\mathbb{F}_2^d$. The \emph{support} of $\mathbf{x} \in C(S)$ is defined as $\mathrm{supp}(\mathbf{x})=\{i: \text{ the $i$th entry of $\mathbf{x}$ is nonzero}\}$.

\begin{thm}\label{PSTCub121}
\emph{(\cite[Theorem 4.1]{Cheung11})}
Let $d \ge 1$ be an integer.
Let $S \subseteq \mathbb{Z}_2^d \setminus \{\mathbf{0}\}$ and $\mathbf{u} \in \mathbb{Z}_2^d \setminus \{\mathbf{0}\}$. Then the following conditions are equivalent:
\begin{itemize}
  \item[\rm (a)] there is perfect state transfer from $\mathbf{0}$ to $\mathbf{u}$ at time $\pi/2\Delta$ in the cubelike graph $\Cay(\mathbb{Z}_2^d, S)$;
  \item[\rm (b)] all codewords in the code of $S$ have weight divisible by $\Delta$, and $\Delta^{-1}w(\mathbf{a}^T M(S))$ and $\mathbf{a}^T\mathbf{u}$ have the same parity for all vectors $\mathbf{a} \in  \mathbb{Z}_2^d$, where $w(\cdot)$ denotes the Hamming weight;
  \item[\rm (c)] $\Delta$ divides $|\mathrm{supp}(\mathbf{x})\cap\mathrm{supp}(\mathbf{y})|$ for any two codewords $\mathbf{x}$ and $\mathbf{y}$ in the code of $S$.
\end{itemize}
\end{thm}

In particular, if $\Cay(\mathbb{Z}_2^d, S)$ admits perfect state transfer from $\mathbf{0}$ to $\mathbf{u}$ at time $\pi/2\Delta$, then $\Delta$ must be the greatest common divisor of the weights of the codewords in the code of $S$ (see \cite[Corollary 4.2]{Cheung11}). As an application of Theorem \ref{PSTCub121}, Cheung and Godsil \cite{Cheung11} presented a cubelike graph over $\mathbb{Z}_2^5$ having perfect state transfer at time $\pi/4$. They asked \cite{Cheung11} further whether the minimum time can be less than $\pi/4$.

In \cite{Godsil12-1}, Godsil asked the following: ``Are there cubelike graphs having perfect state transfer at time $\tau$, where $\tau$ is arbitrarily small?'' Recently, this question has been answered positively by Tan, Feng and Cao in \cite{TanFC17}.

\begin{thm}\label{PSTTimeB1}
\emph{(\cite[Theorem 4.3]{TanFC17})}
Let $d \ge 2$ be an integer. Let $S$ be a subset of $\mathbb{Z}_2^d \setminus \{\mathbf{0}\}$ that generates $\mathbb{Z}_2^d$, and let $0\not= \mathbf{a} \in \mathbb{Z}_2^d$. If the cubelike graph $\Cay(\mathbb{Z}_2^d, S)$ has perfect state transfer from $\mathbf{u}$ to $\mathbf{u}+\mathbf{a}$ at time $t$, for some $\mathbf{u} \in \mathbb{Z}_2^d$, then the minimum time $t$ is $\pi/2^l$ for some $l$ with $1 \le l \le \lfloor d/2 \rfloor$.
\end{thm}

Given a subset $S$ of $\mathbb{Z}_2^d$, define the Boolean function $f:\mathbb{Z}_2^d\rightarrow \mathbb{Z}_2$ with respect to $S$ by
$$
f(\mathbf{x})=\left\{\begin{array}{rl}
              1, & \text{if $\mathbf{x}\in S$},\\[0.2cm]
              0, &  \text{otherwise}.
            \end{array}
\right.
$$
We call $S$ the \emph{support} of $f$ and write $S = \mathrm{supp}(f)$.

Let $B_d$ be the ring of Boolean functions with $d$ variables. For $f \in B_d$, we have the function
$(-1)^f: \mathbb{Z}_2^d \rightarrow \{\pm1\}$ which maps $\mathbf{x}$ to $(-1)^{f(\mathbf{x})}$. The \emph{Fourier transformation of $(-1)^f$} over $(\mathbb{Z}_2^d,+)$, also called the \emph{Walsh transformation of $f$}, is the function $W_f :\mathbb{Z}_2^d \rightarrow \mathbb{Z}$ defined by
$$
W_f(\mathbf{y})=\sum_{\mathbf{x}\in\mathbb{Z}_2^d}(-1)^{f(\mathbf{x})+\mathbf{x}\cdot \mathbf{y}},\; \mathbf{y}\in \mathbb{Z}_2^d,
$$
where $\mathbf{x}\cdot \mathbf{y}=\sum_{i=1}^d x_iy_i\in \mathbb{Z}_2$ for $\mathbf{x}=(x_1,\ldots,x_d)  \in \mathbb{Z}_2^d$ and $\mathbf{y}=(y_1,\ldots,y_d) \in \mathbb{Z}_2^d$. For $d = 2m~(m \ge 2)$, $f \in B_d$ is called a \emph{bent function} in $B_d$ if $|W_f(\mathbf{y})| = 2^m$ for all $\mathbf{y}\in \mathbb{Z}_2^d$.

The next two results due to Tan, Feng and Cao tell us when the lower bound given in Theorem \ref{PSTTimeB1} can be achieved.

\begin{thm}
\emph{(\cite[Theorem 4.5]{TanFC17})}
Let $d = 2m+1$, where $m \ge 2$. Let $f$ be a bent function in $B_{d-1}$ with $0\notin \mathrm{supp}(f) =\{\mathbf{z} \in \mathbb{Z}_2^{d-1}: f(\mathbf{z}) = 1\}$ and $S = \{(0, \mathbf{z}), (1, \mathbf{z}): \mathbf{z} \in \mathrm{supp}(f)\} \subset \mathbb{Z}_2^{d}$. Then the following hold:
\begin{itemize}
  \item[\rm (a)] the cubelike graph $\Cay(\mathbb{Z}_2^d, S)$ is connected;
  \item[\rm (b)] for $\mathbf{a} = (1, 0,0,\ldots,0) \in \mathbb{Z}_2^{d}$, $\Cay(\mathbb{Z}_2^d, S)$ has perfect state transfer from $\mathbf{u}$ to $\mathbf{u} + \mathbf{a}$ at time $\pi/2^m$, for any $\mathbf{u} \in \mathbb{Z}_2^d$;
  \item[\rm (c)] the minimum period with which $\Cay(\mathbb{Z}_2^d,S)$ is periodic at any vertex is $\pi/2^m$.
\end{itemize}
\end{thm}

\begin{thm}
\emph{(\cite[Theorem 4.5]{TanFC17})}
Let $d = 2m$, where $m \ge 2$. Let $f$ be a bent function in $B_{d}$ and $S= \mathrm{supp}(f) =\{\mathbf{z}\in \mathbb{Z}_2^{d}: f(\mathbf{z}) = 1\}$. Then the cubelike graph $\Cay(\mathbb{Z}_2^d, S)$ is connected and the minimum period with which $\Cay(\mathbb{Z}_2^d,S)$ is periodic at any vertex is equal to $\pi/2^m$.
\end{thm}

In \cite{ZhengLZ18}, Zheng, Liu and Zhang proved that $G_{d, k} := \NEPS(H(k,2), \ldots, H(k,2); \mathcal{B})$ is  isomorphic to a cubelike graph, where $d \ge 1$, $k \ge 2$, $\emptyset \ne \mathcal{B} \subseteq \mathbb{Z}_2^{d}\setminus \{\mathbf{0}\}$ and $H(k,2)$ is the hypercube of dimension $k$. In the same paper they also proved the following statements: If $k \ge 3$ is odd and the sum of the elements of $\mathcal{B}$ is not equal to $\mathbf{0}$, then $G_{d, k}$ admits perfect state transfer at time $\pi/2$; if $k$ is even, then $G_{d, k}$ admits perfect state transfer if and only if $\mathcal{B}$ contains at least one element $\beta$ with Hamming weight $w(\beta)=1$, and in this case the perfect state transfer occurs at time $\pi/2$.

Perfect state transfer in general NEPS of complete graphs was considered by Li, Liu, Zhang and Zhou \cite{LiLZZ18} who proved the next three results, where
$$
\NEPS_{n_1, \ldots, n_d, r \cdot 2}(\mathcal{B}) := \NEPS(K_{n_1},\ldots,K_{n_{d}}, \underbrace{K_2,\ldots, K_2}_{r}; \mathcal{B})
$$
for $d \ge 1, r \ge 1, n_1, \ldots, n_d \ge 3$ and $\emptyset \ne \mathcal{B} \subseteq \mathbb{Z}_2^{d+r}\setminus \{\mathbf{0}\}$.

\begin{thm}
\emph{(\cite[Theorem 1.1]{LiLZZ18})}
Suppose that the last $r$ coordinates of each element of $\mathcal{B}$ are all equal to $0$. Then $\NEPS_{n_1, \ldots, n_d, r \cdot 2}(\mathcal{B})$ is periodic with period $2\pi/\gcd(n_1,\ldots,n_d)$.
\end{thm}

\begin{thm}
\emph{(\cite[Theorem 1.2]{LiLZZ18})}
Suppose that the first $d$ coordinates of each element of $\mathcal{B}$ are all equal to $0$.
\begin{itemize}
\item[\rm (a)]  If $\sum_{\beta \in \mathcal{B}} \beta \neq \mathbf{0}$, then for each vertex $\mathbf{u} \in \mathbb{Z}_2^{d+r}$, $\NEPS_{n_1, \ldots, n_d, r \cdot 2}(\mathcal{B})$ admits perfect state transfer from $\mathbf{u}$ to $\mathbf{u} + \sum_{\beta \in \mathcal{B}} \beta$ at time $\pi/2$;
\item[\rm (b)]
If $\sum_{\beta \in \mathcal{B}} \beta = \mathbf{0}$, then $\NEPS_{n_1, \ldots, n_d, r \cdot 2}(\mathcal{B})$ is periodic with period $\pi/2$.
\end{itemize}
\end{thm}

In each case the transition matrix at $\pi/2$ was also given in \cite[Theorem 1.2]{LiLZZ18}. In general, the sufficient condition in (a) above is not necessary.

\begin{thm}
\emph{(\cite[Theorem 1.3]{LiLZZ18})}
\label{thm:LiLZZ18}
Suppose that $\mathcal{B}$ contains at least one element whose last $r$ coordinates are not all $0$ and at least one element whose first $d$ coordinates are not all $0$. Then $\NEPS_{n_1, \ldots, n_d, r \cdot 2}(\mathcal{B})$ is periodic  with period $2\pi$ or $\pi$ depending on whether $\gcd(n_1,\ldots,n_d)$ is odd or even.
\end{thm}

Moreover, sufficient conditions for $\NEPS_{n_1, \ldots, n_d, r \cdot 2}(\mathcal{B})$ to admit perfect state transfer or is periodic with period $\pi/2$ were also given in \cite[Theorem 1.3]{LiLZZ18}. In the special case when $\mathcal{B}$ is the standard basis of $\mathbb{Z}_2^{d+r}$, Theorem \ref{thm:LiLZZ18} says that the Hamming graph $H(n_1, \ldots, n_d, r \cdot 2) := H(n_1, \ldots, n_d, \underbrace{2,\ldots, 2}_{r})$ is periodic with period $2\pi$ or $\pi$ depending on whether $\gcd(n_1,\ldots,n_d)$ is odd or even. Furthermore, if $\gcd(n_1,\ldots,n_d)$ is divisible by $4$, then for each $\mathbf{u} \in \mathbb{Z}_2^{d+r}$, $H(n_1, \ldots, n_d, r \cdot 2)$ admits perfect state transfer from $\mathbf{u}$ to $\mathbf{u}+(0,\ldots,0,1,\ldots,1)$ at time $\pi/2$.


\subsubsection{Perfect state transfer in gcd graphs of abelian groups}

Let $\Gamma=\mathbb{Z}_{n_1}\oplus\cdots\oplus\mathbb{Z}_{n_d}$ be an abelian group, where each $n_i \ge 2$. Set $\mathbf{n} = (n_1, \ldots, n_d)$ and let $\mathbf{D} \subseteq \mathbf{D}(\mathbf{n}) \setminus \{\mathbf{n}\}$, where $\mathbf{D}(\mathbf{n})$ is as defined in \eqref{eq:DDn}. Recall from Section \ref{HammingE} that the gcd graph $\ICG(\mathbf{n}, \mathbf{D})$ of $\Gamma$ with respect to $\mathbf{D}$ is the Cayley graph $\Cay(\Gamma, S_{\Ga}(\mathbf{D}))$, where $S_{\Ga}(\mathbf{D})$ is the set of elements $\mathbf{x} = (x_1, \ldots, x_d)$ of $\Gamma$ such that $\gcd(\mathbf{x}, \mathbf{n}) \in \mathbf{D}$. As shown in (\ref{eq:NEPSICG}), the gcd graphs of abelian groups are precisely the NEPS of complete graphs. So the results to be presented below are essentially about perfect state transfer in NEPS of complete graphs, and this subsection can be considered as a continuation of the previous one.

It is known \cite[Corollary 2.7]{KlotzS13} that every gcd graph of any abelian group with order $n$ is isomorphic to a gcd graph of $\Gamma=\mathbb{Z}_{p_1}\oplus\cdots\oplus\mathbb{Z}_{p_k}$, where $n=p_1p_2\ldots p_k$ is the factorization of $n$ into the product of primes but $p_1, p_2, \ldots, p_k$ are not necessarily distinct. In particular, if $n$ is a power of $2$, then the gcd graph is actually a cubelike graph \cite[Theorem 4.5]{PalB17}. Recall that $G \otimes H$ denotes the tensor product of graphs $G$ and $H$.

\begin{thm}\emph{(\cite[Lemma 4.6]{PalB17})}\label{tensorCubl}
Let $\Gamma = \Gamma_1 \oplus \Gamma_2$ be an abelian group, where $\Gamma_1 = \mathbb{Z}_{2^{a_1}} \oplus\cdots \oplus \mathbb{Z}_{2^{a_r}}$ and $\Gamma_2 = \mathbb{Z}_{p_1^{b_1}} \oplus\cdots \oplus \mathbb{Z}_{p_s^{b_s}}$, with $p_1, \ldots, p_s$ odd primes. Set $\mathbf{n} = (2^{a_1}, \ldots, 2^{a_r}, p_1^{b_1},\ldots, p_s^{b_s})$ and $\mathbf{n}_{2} = (p_1^{b_1},\ldots, p_s^{b_s})$. Assume that $d_{r+1}, \ldots , d_{r+s}$ are fixed divisors of $p_1^{b_1}, \ldots, p_s^{b_s}$, respectively. Let $\mathbf{D} \subseteq \mathbf{D}(\mathbf{n}) \setminus \{\mathbf{n}\}$ be such that the last $s$ components of each member of $\mathbf{D}$ are $d_{r+1}, \ldots , d_{r+s}$, respectively. Set $\mathbf{D}_2 = \{(d_{r+1},\ldots,d_{r+s})\}$. Then there exists a cubelike graph $\Cay(\Gamma_1, S_1)$ such that
$$
\ICG(\mathbf{n}, \mathbf{D}) \cong \Cay(\Gamma_1, S_1) \otimes \ICG(\mathbf{n}_{2}, \mathbf{D}_2).
$$
\end{thm}

As shown in \cite{PalB17}, the connection set $S_1$ (which relies on $\mathbf{D}$) in Theorem \ref{tensorCubl} can be constructed explicitly from $\mathbf{D}$. Let $\mathscr{D}$ be the collection of all subsets $\mathbf{D} \subseteq \mathbf{D}(\mathbf{n})$ satisfying the conditions of Theorem \ref{tensorCubl} as well as the
following two conditions: (i) the sum of the elements in $S_1$ is equal to $\mathbf{0}$; (ii) $|S_1| \equiv 0 ~(\mod~4)$ if there exists $i$ with $1 \leq i \leq s$ such that $d_{r+i} < p_i^{b_i}$. Define $\mathscr{D}_1$ to be the collection of all possible unions of pairwise disjoint members of $\mathscr{D}$. That is, each member of $\mathscr{D}_1$ is of the form $\cup_{i=1}^{t} \mathbf{D}_{i}$ for some pairwise disjoint members $\mathbf{D}_{1}, \ldots, \mathbf{D}_{t}$ of $\mathscr{D}$. Define $\mathscr{D}_2$ to be the collection of all subsets $\mathbf{D} \subseteq \mathbf{D}(\mathbf{n})$ such that the last $s$ components of each member of $\mathbf{D}$ are $p_{1}^{b_{1}}, \ldots, p_{s}^{b_{s}}$, respectively, and the sum of the elements in the corresponding $S_1$ is not equal to $\mathbf{0}$.

The following result due to Pal and Bhattacharjya gives a sufficient condition for a gcd graph of an abelian group to admit perfect state transfer.

\begin{thm}\emph{(\cite[Theorem 4.8]{PalB17})}
\label{thm:PalB17}
Let $\Gamma = \Gamma_1 \oplus \Gamma_2$ be an abelian group, where $\Gamma_1 = \mathbb{Z}_{2^{a_1}} \oplus\cdots \oplus \mathbb{Z}_{2^{a_r}}$ and $\Gamma_2 = \mathbb{Z}_{p_1^{b_1}} \oplus\cdots \oplus \mathbb{Z}_{p_s^{b_s}}$, with $p_1, \ldots, p_s$ odd primes. Let $\mathbf{D}_1\in \mathscr{D}_1$ and  $\mathbf{D}_2\in \mathscr{D}_2$. If $\mathbf{D}_1\cap \mathbf{D}_2=\emptyset$ and $\mathbf{D}=\mathbf{D}_1\cup \mathbf{D}_2$ generates $\Gamma$, then $\ICG(\mathbf{n}, \mathbf{D})$ is connected and has perfect state transfer at time $\pi/2$, where $\mathbf{n} = (2^{a_1}, \ldots, 2^{a_r}, p_1^{b_1},\ldots, p_s^{b_s})$.
\end{thm}

In \cite{PalB17}, Pal and Bhattacharjya also gave the following sufficient condition for the existence of gcd graphs of abelian groups having perfect state transfer at time $\pi/2$ (see also \cite{PalB16}).

\begin{thm}\emph{(\cite[Theorem 4.12]{PalB17})}
\label{thm:PalB17a}
Let $\Gamma$ be a finite abelian group with order $|\Gamma|\equiv 0~ (\mod ~4)$. Then there exists a connected gcd graph of $\Gamma$ that admits perfect state transfer at time $\pi/2$.
\end{thm}


\subsection{Perfect state transfer in Cayley graphs on finite commutative (chain) rings and gcd graphs of unique factorization domains}
\label{subsec:pst-rings}

In \cite{ThongsomnukM17}, {Thongsomnuk} and {Meemark} determined when the unitary Cayley graph $G_R$ of a finite commutative ring $R$ admits perfect state transfer.

\begin{thm}
\label{PSTFcRs}
\emph{(\cite[Theorem 2.5]{ThongsomnukM17})}
Let $R\cong R_0\times R_1\times\cdots\times R_s$ be a finite commutative ring with $1\not=0$, where $s \ge 0$ and each $R_i$, $0 \le i \le s$, is a local ring with maximal ideal $M_i$ of size $m_i$. Set $m = m_0 m_1 \ldots m_s$. The  unitary Cayley graph $G_R$ has perfect state transfer if and only if $R$ satisfies one of the following conditions:
\begin{itemize}
\item[\rm (a)] $m = 1$ and $R\cong
\mathbb{F}_{2}\times\mathbb{F}_{2^{a_1}}\times\mathbb{F}_{2^{a_2}}\times \cdots\times\mathbb{F}_{2^{a_s}}$ for some integers $a_1, a_2, \ldots, a_s \ge 1$;
\item[\rm (b)]  $m = 2$ and $R\cong
R_{0}\times\mathbb{F}_{2^{a_1}}\times\mathbb{F}_{2^{a_2}}\times \cdots\times\mathbb{F}_{2^{a_s}}$ for some integers $a_1, a_2, \ldots, a_s \ge 1$, where $R_0$ is $\mathbb{Z}_4$ or $\mathbb{Z}_2[x]/(x^2)$.
\end{itemize}
\end{thm}

This result yields \cite[Corollary 2.6]{MeemarkS14} in the special case when $R$ is a finite local ring.

The following result characterizes when perfect state transfer occurs in unitary Cayley graphs of finite chain rings.

\begin{thm}
\label{PSTFchain}
\emph{(\cite[Theorem 4.1]{ThongsomnukM17})}
Let $R$ be a finite chain ring with nilpotency $s$. Let $\Cay(R, D)$ be the graph defined in Section \ref{subsec:ring-gcd}, where $D$ depending on a sequence of integers $0 \le a_1 < a_2 < \cdots < a_r \le s-1$ is as defined in \eqref{eq:calC}. Then $\Cay(R, D)$ has perfect state transfer if and only
if $q = 2$ and one of the following holds:
\begin{itemize}
\item[\rm (a)]  $s = 1$ or $2$, and $a_1 =1$;
\item[\rm (b)]  $s \ge 2$, and $a_r = s-2$.
\end{itemize}
\end{thm}

Let $R$ be a unique factorization domain (UFD) and $c$ a nonzero nonunit element of $R$. Assume that the commutative ring $R/(c)$ is finite. Let $D$ be a set of proper divisors of $c$. Recall from Section \ref{subsec:eng-ring-gcd} that the gcd graph $\ICG(R/(c), D)$ is the graph with vertex set the quotient ring $R/(c)$ such that $x+(c), y + (c) \in R/(c)$ are adjacent if and only if $\gcd(x-y, c)$ belongs to $D$ up to associate. As applications of Theorem \ref{PSTFcRs}, {Thongsomnuk} and {Meemark} \cite{ThongsomnukM17} determined the existence of perfect state transfer in some gcd graphs of a unique factorization domain.

\begin{thm}
\label{gcdnew1-1}
\emph{(\cite[Theorem 3.1]{ThongsomnukM17})}
Let $R$ be a unique factorization domain. Let $$c = p_1p_2^{a_2} \ldots p_n^{a_n}$$ be an element of $R$ which is factorized into a product of non-associate irreducible elements, where $a_i > 1$ for $2 \le i \le n$. Then $\ICG(R/(c), \{1, p_1\})$ has perfect state transfer if
and only if there exists $i \in \{2, \ldots , n\}$ such that $a_j = 1$ and $R/(p_j)$ is a finite field of characteristic $2$ for each $j \in\{1, 2, \ldots , n\}\setminus\{i\}$, and either \emph{(i)} $a_i = 1$ and $R/(p_i) \cong \mathbb{F}_2$, or \emph{(ii)} $a_i = 2$ and $R/(p_i^2)$ is isomorphic to $\mathbb{Z}_4$ or $\mathbb{Z}_2[x]/(x^2)$.
\end{thm}

\begin{thm}
\label{gcdnew1-2}
\emph{(\cite[Theorem 3.2]{ThongsomnukM17})}
Let $R$ be a unique factorization domain. Let $$c = p_1p_2\ldots p_kq_1^{a_1} \ldots q_l^{a_l}$$ be an element of $R$ which is factorized as a product of non-associate irreducible elements, where $a_i > 1$ for $1 \le i \le l$. Then $\ICG(R/(c), \{p_i, p_j\})$ with $1\le i<j\le k$ has perfect state transfer if
and only if one of the following occurs:
\begin{itemize}
\item[\rm (a)]  $c = p_1p_2\ldots p_k$, $R/(p_s)$ is a finite field of characteristic $2$ for $1 \le s \le k$, and $R/(p_i)\cong R/(p_j)\cong \mathbb{F}_2$ or $R/(p_t)\cong \mathbb{F}_2$ for some $t\not=i,j$;
\item[\rm (b)]  $c = p_1p_2\ldots p_kq_1^2$, $R/(p_s)$ is a finite field of characteristic $2$ for $1 \le s \le k$, and $R/(q_1^2)$ is isomorphic to $\mathbb{Z}_4$ or $\mathbb{Z}_2[x]/(x^2)$.
\end{itemize}
\end{thm}

Let $R$ be a UFD and $c = p_1^{a_1} p_2^{a_2} \ldots p_k^{a_k}$ a nonzero nonunit element of $R$ factorized as a product of irreducible elements. Assume that $R/(c)$ is finite. Assume further that for some $k\ge2$ and each $i$ with $1 \le i \le k$ there exists a set $$
D_i =\{p_{i}^{a_{i1}}, p_{i}^{a_{i2}}, \ldots, p_{i}^{a_{ir_i}}\}
$$
such that $0 \le a_{i1} < a_{i2} < \cdots < a_{ir_i} \le a_i - 1$. Set
$$
D = \{p_1^{a_{1t_1}} \ldots p_k^{a_{k t_k}}: t_i \in \{1, 2, \ldots, r_i\} \text{ for } 1 \le i \le k\}.
$$
Then
$$
\ICG(R/(c), D) \cong \Cay(R/(p_1^{a_1}), D_1)\otimes \cdots\otimes \Cay(R/(p_k^{a_k}), D_k),
$$
where each factor on the right hand side is the Cayley graph $\Cay(R/(p_i^{a_i}), D_i)$ over the finite chain ring $R/(p^{s_i}_i )$ with respect to the divisor set $D_i$.
Using this isomorphism, {Thongsomnuk} and {Meemark} established the following result.

\begin{thm}
\label{gcdnew12}
\emph{(\cite[Theorem 4.2]{ThongsomnukM17})}
Under the assumption above, the following hold:
\begin{itemize}
\item[\rm (a)]  if $\ICG(R/(c), D)$ has perfect state transfer, then $\Cay(R/(p_i^{a_i}), D_i)$ has perfect state transfer for some $i \in\{1, 2, \ldots, k\}$;
\item[\rm (b)] if $\Cay(R/(p_1^{a_1}), D_1)$ has perfect state transfer, and for all $i=2, \ldots, k$, $R/(p_i^{a_i})$ is of even characteristic and $a_{ir_i} = a_i-1$, then $\ICG(R/(c), D)$ admits perfect state transfer.
\end{itemize}
\end{thm}

\subsection{Perfect state transfer in Cayley graphs on non-abelian groups}
\label{subsec:PSTdih}

In the literature there are relatively few results on perfect state transfer in Cayley graphs on non-abelian groups. We review these results in this section.

In \cite{CaoCL20, CaoF19}, necessary and sufficient conditions for a Cayley graph on a dihedral group (that is, a dihedrant) to admit perfect state transfer were obtained and explicit constructions were given. As one may expect, these results rely on the irreducible representations of dihedral groups. As before, in this section $D_{2n} = \langle a,b \mid a^n=b^2=1, b^{-1}ab=a^{-1} \rangle$ is the dihedral group of order $2n \ge 4$ and $S$ is a subset of $D_{2n} \setminus \{1\}$ with $S^{-1}=S$. The following result handles the case when $n$ is odd.

\begin{thm}
\label{thm:CaoCL}
\emph{(\cite[Theorem 8]{CaoCL20})}
Let $n \ge 3$ be an odd integer. Then any connected Cayley graph $\Cay(D_{2n}, S)$ on $D_{2n}$ has no perfect state transfer between any two distinct vertices. Moreover, $\Cay(D_{2n}, S)$ is periodic if and only if it is integral and $S \cap b \la a \ra = \emptyset$ or $b \la a \ra$. Furthermore, in this case the minimum period at all vertices is $2\pi/m$, where $m$ is the greatest common divisor of $\l - |S|$ for $\l$ running over all eigenvalues of $\Cay(D_{2n}, S)$ other than $|S|$.
\end{thm}

In the case when $n = 2m$ is even, following \cite{CaoCL20} let $\psi_{1}$ denote the trivial representation of $D_{2n}$ and $\psi_{4}$ the one-dimensional irreducible representation of $D_{2n}$ defined by
$$
\psi_{4}(a^i) = (-1)^i,\ \psi_{4}(ba^i) = (-1)^{i+1},\ 0 \le i \le n-1.
$$
If $S$ is closed under conjugation, then $\Cay(D_{2n}, S)$ has four (not necessarily distinct) eigenvalues
$$
\l_1 = |S|, \l_2, \l_3, \l_4
$$
which correspond to $\psi_{1}$ and $\psi_{4}$, respectively, and some eigenvalues $\mu_j$ corresponding to the two-dimensional representations $\rho_j$, $1 \le j \le m-1$, where $\rho_j$ is defined by
$$
\rho_j(a^i) = \pmat{\om_{n}^{ij} & 0\\0 & \om_{n}^{-ij}},\ \rho_j(ba^i) = \pmat{0 & \om_{n}^{-ij}\\\om_{n}^{ij} & 0},\ 0 \le i \le n-1.
$$
Identify the elements of $D_{2n}$ with integers $0, 1, \ldots, 2n-1$ in the follow way: for $0 \le u \le n-1$, $a^u$ corresponds to $u$, and for $n \le u \le 2n-1$, $ba^u$ corresponds to $u$. Recall from \eqref{eq:2adic} the $2$-adic valuation $\nu_{2}: \QQQ \rightarrow \ZZZ \cup \{\infty\}$ of rational numbers. With these notations we now present the result from \cite{CaoCL20} for even $n$ which covers \cite[Theorem 3.2]{CaoF19} in the special case when $S$ is closed under conjugation.

\begin{thm}
\label{thm:CaoCLa}
\emph{(\cite[Theorem 11]{CaoCL20})}
Let $n = 2m$ and let $G = \Cay(D_{2n}, S)$ be a connected Cayley graph on the dihedral group $D_{2n}$. Then $G$ has no perfect state transfer between any two distinct vertices unless $S$ is closed under conjugation. Moreover, $G$ is periodic if and only if it is integral, and in this case the minimum period at all vertices is $2\pi/m$, where $m$ is the greatest common divisor of $\l - \l_1$ for $\l$ running over all eigenvalues of $G$ other than $\l_1$. Furthermore, using the notation above, the following hold:
\begin{itemize}
\item[\rm (a)] if $m$ is even, then $G$ has perfect state transfer from $u$ to $v$ if and only if $G$ is integral, $v=u+m$, and there exists a constant $c$ such that for any eigenvalue $\l$ of $G$ we have $\nu_{2}(\l - \l_1) = c$ if $\l = \mu_{2j-1}$ for some $1 \le j \le m/2$ and $\nu_{2}(\l - \l_1) > c$ otherwise;
\item[\rm (b)] if $m$ is odd, then $G$ has perfect state transfer from $u$ to $v$ if and only if $G$ is integral, $v=u+m$, and there exists a constant $d$ such that $\nu_{2}(\l - \l_1) = d$ if $\l = \l_3, \l_4$ or $\mu_{2j-1}$ for some $1 \le j \le (m-1)/2$, and $\nu_{2}(\l - \l_1) > d$ if $\l = \l_2$ or $\mu_{2j}$ for some $1 \le j \le (m-1)/2$.
\end{itemize}
In addition, in each case when the conditions above hold, the minimum time at which $G$ has perfect state transfer from $u$ to $v$ is $\pi/m$.
\end{thm}

Several families of dihedrants which admit perfect state transfer or are periodic were  constructed in \cite[Section 5]{CaoF19} explicitly.

Recall from \eqref{eq:SD} that $SD_{8n}$ denotes the semi-dihedral group of order $8n$. Recently, Luo, Cao, Wang and Wu gave a necessary and sufficient condition for a connected normal Cayley graph on $SD_{8n}$ to admit perfect state transfer. See \cite[Theorems 3.2 and 3.3]{LuoCWW21} for details of these technical results.

In \cite{SinJ21}, Sin and Sorci obtained a necessary and sufficient condition for a normal Cayley graph on an extraspecial $2$-group to admit perfect state transfer. Their choice of extraspecial $2$-groups was based on a result in \cite{ChanCTVZ20} (see also \cite[Theorem 4.1]{Godsil11}) which says that perfect state transfer can only occur in Cayley graphs on groups containing a central involution, and the fact that each extraspecial $2$-group contains a unique central involution. A finite group $\Ga$ is called an \emph{extraspecial $p$-group} if it is a $p$-group whose center $Z$ has order $p$ such that $\Ga/Z$ is isomorphic to $\ZZZ_p^m$ for some $m \ge 1$. It is well known that an extraspecial $p$-group must have order $p^{2n+1}$ for some $n \ge 1$. Let $\Ga$ be an extraspecial $2$-group with order $2^{2n+1}$ and center $Z$, and let $z$ be the unique central involution of $\Ga$. Then each conjugacy class of $\Ga$ with size greater than $1$ must have size $2$. Let $S \subseteq \Ga \setminus \{1\}$ be a union of some conjugacy classes of $\Ga$ that generates $\Ga$. Then $S \setminus \{z\} = \cup_{i=1}^{\ell} K_i$ for some distinct conjugacy classes $K_1, K_2, \ldots, K_{\ell}$ of $\Ga$ with size $2$, where $\ell \ge 1$ and $z$ may or may not be in $S$. Take a representative $x_i \in K_i$ for each $i$, and for each $y \in \Ga/Z \cong \ZZZ_2^{2n}$ let $e_y$ be the number of elements $x_i$ such that $y \cdot \bar{x}_i = 0$ module $2$, where for each $g \in \Ga$, $\bar{g}$ is the image of $g$ in $\Ga/Z$. Recall from \eqref{eq:2adic} that $\nu_{2}(x)$ denotes the $2$-adic valuation of $x \in \mathbb{Q}$.

\begin{thm}
\emph{(\cite[Theorem 4.7]{SinJ21})}
\label{thm:SinJ21}
Let $\Ga$ be an extraspecial $2$-group with center $Z$, and let $z$, $S$ and $\ell$ be as above. Then the normal Cayley graph $\Cay(\Ga, S)$ admits perfect state transfer at $\tau$ if and only if one of the following holds:
\begin{itemize}
\item[\rm (a)] $z \in S$ and $\nu_{2}(\ell - e_y) \ge \nu_{2}(\ell + 1)$ for all $y \in \Ga/Z$;
\item[\rm (b)] $z \notin S$ and $\nu_{2}(\ell - e_y) \ge \nu_{2}(\ell)$ for all $y \in \Ga/Z$.
\end{itemize}
Moreover, in case (a), $\tau$ is an odd multiple of $\pi/2^{\nu_{2}(2\ell + 2)}$, and in case (b), $\tau$ is an odd multiple of $\pi/2^{\nu_{2}(2\ell)}$.
\end{thm}

More information about $\tau$ can be found in \cite[Theorems 4.4 and 4.6]{SinJ21}. Constructions of normal Cayley graphs on some extraspecial $2$-groups that admit perfect state transfer were also given in this paper. In \cite[Theorem 3.1]{SinJ21}, the Cayley graph in Theorem \ref{thm:SinJ21} was shown to be integral and its spectrum was given explicitly.

\subsection{Pretty good state transfer in Cayley graphs}
\label{subsec:pgst}

Let $G$ be a graph, and let $u, v \in V(G)$. We say that a \emph{pretty good state transfer} (PGST) occurs from $u$ to $v$ in $G$ if for any $\ve > 0$ there exists a time $\tau$ such that $|\bfe_{v}^T H_{G}(\tau) \bfe_{u}| \ge 1 - \ve$, where $H_{G}(t)$ is the transition matrix of $G$ as defined in \eqref{eq:HGt} and $\bfe_{v}$ is the characteristic vector of the subset $\{v\}$ of $V(G)$. This condition is equivalent to that there exists a sequence of real numbers $\{t_k\}$ and a complex number $\gamma$ with $|\gamma| = 1$ such that $\lim_{k \rightarrow \infty} H_{G}(t_k) \bfe_{u} = \gamma \bfe_{v}$. The notion of PGST was introduced by Godsil \cite{Godsil12-1} and Vinet and Zhedanov \cite{VinetZ12} independently in 2012 as a relaxation of perfect state transfer, where in the latter paper the term ``almost state transfer" was used for PGST. Since then this concept has attracted considerable attention from both mathematicians and physicists as one can see in, for example, \cite{KemptonLY17} and the references therein. The reader is referred to \cite[Section 1]{CaoWF20} for an account of known results on the existence of PGST in several families of graphs. In the rest of this section we restriction our attention to PGST in Cayley graphs.

Recall that a regular graph is periodic if and only if it is integral. As observed in \cite{PalB2017}, a periodic graph exhibits PGST if and only if it admits perfect state transfer. Since all integral circulant graphs admitting perfect state transfer have been characterized, it follows that integral circulant graphs exhibiting PGST are all known and are given in Theorem \ref{CIMainPST}. So the existence of PGST in circulant graphs is reduced to that in non-integral circulant graphs. Unfortunately, as far as we are aware there is no known necessary and sufficient condition for a non-integral circulant graph to admit PGST. In \cite[Lemma 5]{PalB2017}, Pal and Bhattacharjya observed that in a circulant graph $\Cay(\ZZZ_n, S)$, PGST occurs only from $u$ to $u + (n/2)$ for $u \in \ZZZ_n$, and of course this requires $n/2 \in S$ and $n$ is even. They also proved that the cycle $C_n$ of length $n$ admits PGST if and only if its complement admits PGST, and this occurs if and only if $n = 2^k$ for some $k \ge 2$ (see \cite[Theorem 13]{PalB2017}). In \cite[Theorem 7]{PalB2017}, the same authors proved further that, if $n = 2^k$ for some $k \ge 3$, then $\Cay(\ZZZ_n, S_{n}(D) \cup \{-1, 1\})$ with $1 \not \in D$ and $S_{n}(D)$ as defined in \eqref{eq:SnD} (that is, the edge-disjoint union of $C_n$ and a gcd-graph $\ICG(n,D)$) as well as its complement both exhibit PGST. This result was generalized in \cite{Pal18} where the following result was proved. Recall from \eqref{eq:gnd} the definition of $S_{n}(d)$.

\begin{thm}
\emph{(\cite[Theorem 2.4]{Pal18})}
Let $\Cay(\ZZZ_n, S)$ be a non-integral circulant graph with order $n = 2^k$, where $k \ge 2$. Let $d$ be the smallest divisor $n$ such that $\emptyset \ne S \cap S_{n}(d) \ne S_{n}(d)$. If $|S \cap S_{n}(d)| \equiv 2~(\mod ~4)$, then $\Cay(\ZZZ_n, S)$ admits PGST.
\end{thm}

It is unknown whether the sufficient condition in this theorem is also necessary. On the other hand, the following necessary condition was obtained by Pal in \cite{Pal18} in the case when $n$ is not a power of $2$.

\begin{thm}
\emph{(\cite[Theorem 2.6]{Pal18})}
Let $n \ge 3$ be an integer which is not a power of $2$. Then $\Cay(\ZZZ_n, S)$ does not admit PGST unless each odd prime factor of $n$ divides at least one element of $S$.
\end{thm}

Regarding normal Cayley graphs on dihedral groups, the next two results were proved by Cao, Wang and Feng in \cite{CaoWF20}.

\begin{thm}
\emph{(\cite[Theorem 3]{CaoWF20})}
Let $n \ge 3$ be an odd integer. Then any connected normal Cayley graph on the dihedral group $D_{2n}$ does not admit PGST between any pair of distinct vertices.
\end{thm}

\begin{thm}
\emph{(\cite[Theorem 4]{CaoWF20})}
Let $n = 2^k$, where $k \ge 3$. Then a normal Cayley graph $\Cay(D_{2n}, S)$ on the dihedral group $D_{2n} = \langle a,b\left|a^n=b^2=1, b^{-1}ab=a^{-1}\right \rangle$ admits PGST with respect to a time sequence in $2 \pi \ZZZ$ provided that $b\langle a \rangle \subseteq S$ and either $S \cap \langle a \rangle = \{a^{h}, a^{-h}\}$ for some odd integer $h$ between $1$ and $2^{k-1} - 1$, or $S \cap \langle a \rangle = \{a^{h}, a^{-h}\} \cup \left(\cup_{j=1}^r \{a^{2^{m_j}l}: 1 \le l \le 2^{k-m_j}-1, l \text{ is odd }\}\right)$ for some odd integer $h$ between $1$ and $2^{k-1} - 1$ and pairwise distinct integers $m_1, \ldots, m_r$ between $1$ and $k-1$.
\end{thm}

In \cite{CaoWF20}, it is also proved that, except when $S \cap \langle a \rangle$ satisfies a certain condition (the so-called ``power-two case"), a normal Cayley graph $\Cay(D_{2n}, S)$ on $D_{2n}$ admits PGST if and only if $n$ is a power of $2$. As noted in the same paper, in the ``power-two case" the normal Cayley graph $\Cay(D_{2n}, S)$ may admit PGST.

Pretty good state transfer in normal Cayley graphs on semi-dihedral groups $SD_{8n}$ was studied by Wang and Cao in \cite{WangC21}. Due to the structure of $SD_{8n}$, we may assume without loss of generality that any normal Cayley graph $\Cay(SD_{8n}, S)$ satisfies $b \la a \ra \subseteq S$ and $S \cap \la a \ra = \{a^{h_j}, a^{-h_j}, a^{(2n-1)h_j}, a^{-(2n-1)h_j}: 1 \le j \le r\}$ for some pairwise distinct positive integers $h_1, \ldots, h_r$. If $2^{\nu_{2}(n) - \nu_{2}(h_j)} h_j = n$ for every $1 \le j \le r$, then $S \cap \la a \ra$ is said to be in the ``power-two" case.

\begin{thm}
\emph{(\cite[Theorems 3.5 and 3.6]{WangC21})}
Let $SD_{8n}$ be the semi-dihedral group as shown in \eqref{eq:SD}, and let $\Cay(SD_{8n}, S)$ be a non-integral normal Cayley graph with $S$ as above. Then the following hold:
\begin{itemize}
\item[\rm (a)] if $n$ is odd, then $\Cay(SD_{8n}, S)$ does not admit PGST;
\item[\rm (b)] if $n$ is even but not a power of $2$, then $\Cay(SD_{8n}, S)$ does not admit PGST unless $S \cap \la a \ra$ is in the power-two case.
\end{itemize}
\end{thm}

In \cite{WangC21}, it was observed that, indeed, the study of the ``power-two" case can be reduced to the case where $n$ is a power of $2$, and moreover the following result was obtained in this case.

\begin{thm}
\emph{(\cite[Theorem 3.3]{WangC21})}
Let $SD_{8n}$ be the semi-dihedral group as shown in \eqref{eq:SD}. Assume that $n = 2^k$, $k \ge 2$. Then a normal Cayley graph $\Cay(SD_{8n}, S)$ admits PGST with respect to a sequence in $2\pi \ZZZ$ provided that one of the following holds:
\begin{itemize}
\item[\rm (a)] $b \la a \ra \subseteq S$ and $S \cap \la a \ra = \{a^{h}, a^{-h}, a^{(2n-1)h}, a^{-(2n-1)h}\}$ for some odd integer $h$ between $1$ and $2n-1$;
\item[\rm (b)] $b \la a \ra \subseteq S$ and $
S \cap \la a \ra = \{a^{h}, a^{-h}, a^{(2n-1)h}, a^{-(2n-1)h}\} \cup S_0$, where
$$
S_0 = \cup_{j=1}^r \{a^{2^{m_j}l}, a^{2^{m_j}(2n-1)l}: 1 \le l \le 2^{k+2-m_j}-1, l \text{ is odd}\},
$$
for some odd integer $h$ between $1$ and $2n - 1$ and pairwise distinct integers $m_1, \ldots, m_r$ between $1$ and $k+1$.
\end{itemize}
\end{thm}

The following result for general normal Cayley graphs is an analog of \cite[Theorem 4.1]{Godsil11}.

\begin{thm}
\emph{(\cite[Theorem 1]{CaoWF20})}
Let $G = \Cay(\Ga, S)$ be a normal Cayley graph on a finite group $\Ga$. Let $u, v$ be a pair of distinct elements  of $\Ga$. If $G$ admits PGST (respectively, perfect state transfer) from $u$ to $v$, then $vu^{-1}$ is an involution lying in the center of $\Ga$, and moreover there exist a sequence of real numbers $\{t_k\}$, a complex number $\gamma$ with $|\gamma| = 1$ and a permutation matrix $P$ of order $2$ with no fixed points such that $\lim_{k \rightarrow \infty} H_{G}(t_k) = \gamma P$ (respectively, $H_{G}(t_k) = \gamma P$ and $t_k = t$ for all $k$).
\end{thm}


\section{Distance-regular Cayley graphs}
\label{sec:drCay}

A connected graph $G$ with diameter $d$ is called \emph{distance-regular} (see \cite[Section 3.7]{Cvetkovic10}) if there exist nonnegative integers $b_0, b_1, \ldots , b_{d-1}$ and
$c_1,c_2,\ldots,c_{d}$ such that for any pair of vertices $u,v$ at distance $i$, we have
$$
b_i = |N_{i+1}(u)\cap N_1(v)|,\, 0 \le i \le d-1;\quad
c_i = |N_{i-1}(u)\cap N_1(v)|,\, 1 \le i \le d,
$$
where $N_i(u)$ is the set of vertices of $G$ at distance $i$ from $u$. It is well known (see  \cite[Theorem 3.7.3]{Cvetkovic10}) that the eigenvalues of a distance-regular graph are determined by its  \emph{intersection array} $\{b_0,b_1, \ldots, b_{d-1}; c_1,c_2, \ldots,c_d\}$. A distance-regular graph of diameter two is a connected strongly regular graph. More explicitly, a \emph{strongly regular graph} with parameters $(v, k, \l, \mu)$ is a $k$-regular graph with $v$ vertices in which any two adjacent vertices have exactly $\l$ common neighbours and any two non-adjacent vertices have exactly $\mu$ common neighbours (see \cite[Section 1.2]{Cvetkovic10}). Obviously, if a graph $G$ is strongly regular, then so is its complement $\overline{G}$. A strongly regular graph $G$ is \emph{nontrivial} if both $G$ and $\overline{G}$ are connected and \emph{trivial} otherwise. (Equivalently, $G$ is trivial if either $\mu = 0$ or $\mu = k$.) It is well known (see \cite[Theorem 3.6.5]{Cvetkovic10}) that the eigenvalues of a strongly regular graph with parameters $(v, k, \l, \mu)$ are
$$
k, \quad \frac{1}{2}\left((\l-\mu)+\sqrt{\Delta}\right), \quad \frac{1}{2}\left((\l-\mu)-\sqrt{\Delta}\right),
$$
with multiplicities
$$
1, \quad \frac{1}{2}\left((v-1)-\frac{2k+(v-1)(\l-\mu)}{\sqrt{\Delta}}\right), \quad \frac{1}{2}\left((v-1)+\frac{2k+(v-1)(\l-\mu)}{\sqrt{\Delta}}\right),
$$
respectively, where $\Delta = (\l-\mu)^2+4(k-\mu)$.
The reader is referred to the monograph \cite{BrouwerCN89} and the survey paper \cite{vanDamKT16} for comprehensive treatments of distance-regular graphs.

It is well known that strongly regular Cayley graphs are essentially partial difference sets.
A $k$-subset $S$ of a group $\Gamma$ with order $v$ is called an $(v,k,\l,\mu)$-\emph{partial difference set} (PDS) in $\Gamma$ if the expressions $gh^{-1}$, for $g, h \in S$ with $g\not=h$, represent each nonidentity element in $S$ exactly $\l$ times and each nonidentity element not in $S$ exactly $\mu$ times. It follows immediately that a Cayley graph $\Cay(\Ga, S)$ is strongly regular if and only if $S$ is a PDS in $\Gamma$ such that $1 \notin S$ and $S^{-1}=S$.

In this section we restrict ourselves to distance-regular Cayley graphs. We barely touch strongly regular Cayley graphs since we believe that a survey on partial difference sets should be a separate treatise written by experts in the area (see in \cite{Ma94} for such a survey published in 1994).


\subsection{Distance-regular Cayley graphs}
\label{subsec:dis-reg-cay}

A well known strongly regular graph with parameters $\left(q, \frac{q-1}{2}, \frac{q-5}{4}, \frac{q-1}{4}\right)$ is the Paley graph $P(q)$, which is the Cayley graph on the additive group of $\mathbb{F}_q$ with connection set the set of nonzero squares of $\mathbb{F}_q$, where $q$ is a prime power with $q \equiv 1 \pmod{4}$. The following result, obtained independently by Bridges and Mena \cite{BridgesM79}, Hughes, van Lint and Wilson \cite{HughesLW79}, Ma \cite{Ma84}, and partially by Maru\v{s}i\v{c} \cite{Marusic89}, gives a complete classification of strongly regular circulant graphs.

\begin{thm}
\emph{(\cite{BridgesM79, HughesLW79, Ma84, Marusic89})}
If $G$ is a nontrivial strongly regular circulant graph, then $G$ is
isomorphic to a Paley graph $P(p)$ for some prime $p \equiv 1~(\mod ~4)$.
\end{thm}

In \cite{MiklavicP03}, Miklavi\v{c} and Poto\v{c}nik generalized this result to distance-regular circulant graphs. Denote by $K_{n,n} - nK_2$ the complete bipartite graph $K_{n,n}$ with a perfect matching removed.

\begin{thm}\emph{(\cite[Theorem 1.2]{MiklavicP03})}
Let $G$ be a circulant graph with $n \ge 3$ vertices. Then $G$ is distance-regular
if and only if it is isomorphic to one of the following graphs:
\begin{itemize}
  \item[\rm (a)] cycle $C_n$;
  \item[\rm (b)] complete graph $K_n$;
  \item[\rm (c)] complete $t$-partite graph $K_{m,\ldots, m}$, where $tm = n$;
  \item[\rm (d)] $K_{m,m}-mK_2$, where $2m = n$ and $m$ is odd;
  \item[\rm (e)]  Paley graph $P(n)$, where $n \equiv 1~(\mod~ 4)$ is a prime.
\end{itemize}
\end{thm}

In \cite{MiklavicS14}, Miklavi\v{c} and {\v{S}}parl classified all distance-regular Cayley graphs on abelian groups with respect to ``maximal" inverse-closed generating sets. A few definitions are in order before presenting their result. The \emph{Shrikhande graph} is the Cayley graph on $\mathbb{Z}_{4}\times \mathbb {Z}_{4}$ with connection set $\{\pm (1,0), \pm (0,1), \pm (1,1)\}$; this is a strongly regular graph with parameters $(16, 6, 2, 2)$. The \emph{Doob graph} $D(m,n)$ (where $n, m\ge1$) is the Cartesian product of the Hamming graph $H(n, 4)$ with $m$ copies of the Shrikhande graph. A distance-regular graph $G$ of diameter $d$ is called \emph{antipodal} if the relation $R$ on $V(G)$ defined by $xRy \Leftrightarrow d(x, y) \in \{0, d\}$ is an equivalence relation, and \emph{non-antipodal} otherwise. In the former case the \emph{antipodal quotient} of $G$ is the graph with vertices the equivalence classes of $R$ such that two equivalence classes are adjacent if and only if there is at least one edge between them in $G$. For example, the hypercube $H(d,2)$ is antipodal since each vertex $x$ has a unique antipodal vertex $\bar{x}$ whose distance to $x$ is equal to $d$. The antipodal quotient of $H(d,2)$ is the quotient of $H(d,2)$ with respect to the partition $\{\{x, \bar{x}\}: x \in V(H(d,2))\}$ of $V(H(d,2))$.

\begin{thm}\emph{(\cite[Theorem 1.1]{MiklavicS14})}
Let $\Gamma$ be an abelian group with identity $1$ and let $S$ be an inverse-closed subset of $\Gamma\setminus\{1\}$ which generates $\Gamma$ such that $S\setminus\{s,s^{-1}\}$ does not generate $\Gamma$ for at least one element $s\in S$. Then $\Cay(\Gamma, S)$ is distance-regular if and only if it is isomorphic to one of the following graphs:
\begin{itemize}
  \item[\rm (a)] complete bipartite graph $K_{3,3}$;
  \item[\rm (b)] complete tripartite graph $K_{2,2,2}$;
  \item[\rm (c)] $K_{6,6}-6K_2$;
  \item[\rm (d)] cycle $C_n$ for $n\ge3$;
  \item[\rm (e)] Hamming graph $H(d,n)$, where $d\ge1$ and $n\in \{2,3,4\}$;
  \item[\rm (f)]  Doob graph $D(m,n)$, where $n,m\ge1$;
  \item[\rm (g)] antipodal quotient of the hypercube $H(d,2)$, where $d \ge 2$.
\end{itemize}
\end{thm}

Let $v$, $k$ and $\l$ be nonnegative integers and $\Gamma$ a group of order $v$. A $k$-subset $D$ of $\Gamma$ is called a \emph{$(v, k,\l)$-difference set} if every nonidentity element of $\Ga$ has exactly $\l$ representations as a product $g h^{-1}$ with $g, h \in D$; $D$ is \emph{trivial} if $k \in \{0, 1, v-1, v\}$ and \emph{nontrivial} otherwise. Recall that Cayley graphs on dihedral groups $D_{2n}= \langle a,b\ |\ a^n=b^2=1, b^{-1}ab=a^{-1}\rangle$ are called dihedrants. Given $S,T\subseteq \mathbb{Z}_n$, let $a^S = \{a^i: i \in S\}$
and $a^Tb = \{a^ib: i \in T \}$. Given $A \subseteq \mathbb{Z}_n$ and $i \in \mathbb{Z}_n$, let $i +A = \{i +a: a \in A\}$ and $iA = \{ia: a \in A\}$. Distance-regular dihedrants have been classified by Miklavi\v{c} and Poto\v{c}nik in \cite{MiklavicP07}.

\begin{thm}
\emph{(\cite[Theorem 4.1]{MiklavicP07})}
Let $S,T\subseteq \mathbb{Z}_n$, where $n \ge 2$. Let $G=\Cay(D_{2n},a^S\cup a^Tb)$ be a connected dihedrant other than $C_{2n}$, $K_{2n}$, the complete $t$-partite graph $K_{m,\ldots,m}$ (where $tm=2n$), or $K_{n,n} - nK_2$. Then $G$ is distance-regular if and only if one of the following
holds:
\begin{itemize}
  \item[\rm (a)] $S =\emptyset$ and $T$ is a nontrivial difference set in $\mathbb{Z}_n$;
  \item[\rm (b)] $n$ is even, $S$ is a non-empty subset of $1+ 2\mathbb{Z}_n$, and either
      \begin{itemize}
  \item[\rm (i)] $T \subseteq 1 + 2\mathbb{Z}_n$ and $a^{-1+S} \cup a^{-1+T}b$ is a nontrivial difference set in the dihedral group $\langle a^2, b\rangle$, or
  \item[\rm (ii)] $T \subseteq 2\mathbb{Z}_n$ and $a^{-1+S} \cup a^{T}b$ is a nontrivial difference set in the dihedral group $\langle a^2, b\rangle$.
\end{itemize}
\end{itemize}
Moreover, if either (a) or (b) occurs, then $G$ is bipartite and non-antipodal with diameter $3$.
\end{thm}

In \cite{Momihara18}, Momihara gave a construction of strongly regular Cayley graphs on the additive groups of finite fields based on three-valued Gauss periods. As corollaries, he obtained two infinite families and one sporadic example of new strongly regular Cayley graphs. The construction in \cite{Momihara18} can be viewed as a generalization of the construction of strongly regular Cayley graphs given by Bamberg, Lee, Momihara and Xiang \cite{BambergLMX17}.

\begin{thm}
\emph{(\cite[Theorem 1.1]{Momihara18})}
Let $q$ be a prime power. There exists a strongly regular Cayley graph on the additive group of $\FFF_{q^6}$ with negative Latin square type parameters $(q^6, r(q^3 + 1), q^3 + r^2 - 3r, r^2 - r)$, where $r = M(q^2 - 1)/2$, in the following cases:
\begin{itemize}
  \item[\rm (a)] $M = 3$ and $q \equiv 7~(\mod~24)$;
  \item[\rm (b)] $M = 7$ and $q \equiv 11, 51~(\mod~56)$.
  \end{itemize}
\end{thm}

\begin{thm}
\emph{(\cite[Theorem 15]{AalipourA14})}
Let $R$ be a finite commutative ring. Let $Z^{*}(R)$ be the set of nonzero zero divisors of $R$ and $\Cay(R, Z^{*}(R))$ the Cayley graph on the additive group of $R$ with connection set $Z^{*}(R)$. Then the following statements are equivalent:
\begin{itemize}
  \item[\rm (a)] $\Cay(R, Z^{*}(R))$ is edge-transitive;
  \item[\rm (b)] $\Cay(R, Z^{*}(R))$ is strongly regular;
  \item[\rm (c)] $R$ is a local ring, or $R = \ZZZ_2^d$ for some $d \ge 2$, or $R = \FFF_q \times \FFF_q$ for some $q \ge 3$.
  \end{itemize}
Moreover, if $R$ is not a local ring, then each of these statements is equivalent to that $\Cay(R, Z^{*}(R))$ is distance-regular.
\end{thm}

It was also noted in \cite[Corollary 16]{AalipourA14} that $G_R$ is strongly regular if and only if $R$ is a local ring, or $R = \ZZZ_2^d$ for some $d \ge 2$, or $R = \FFF_q \times \FFF_q$ for some $q \ge 3$.

Recall that, for a finite ring $R$ and an integer $n \ge 1$, the unitary Cayley graph of the ring $M_{n}(R)$ of $n \times n$ matrices over $R$ is $G_{M_{n}(R)} = \Cay(M_{n}(R), \GL_{n}(R))$. In \cite[Theorem 2.3]{KianiM15}, Kiani and Mollahajiaghaei proved that $G_{M_{2}(\FFF_{q})}$ is strongly regular with parameters $(q^4, q^4 - q^3 - q^2 + q, q^4 - 2q^3 - q^2 + 3q, q^4 - 2q^3 + q)$. This result also follows from Theorem \ref{thm:mat-rings-1}. In fact, based on Theorems \ref{thm:mat-rings-1} and \ref{thm:mat-rings-2}, Rattanakangwanwong and Meemark proved the following stronger result.

\begin{thm}
\emph{(\cite[Theorem 2.4]{RattaM20})}
\label{thm:mat-rings-3}
For any prime power $q$, $G_{M_{n}(\FFF_q)}$ is strongly regular if and only if $n = 2$.
\end{thm}

The following result was also proved in \cite{RattaM20}.

\begin{thm}
\emph{(\cite[Theorem 4.3]{RattaM20})}
Let $R$ be a local ring other than a field. Then $G_{M_{n}(R)}$ is not strongly regular for any $n \ge 2$.
\end{thm}

All distance-regular graphs with degree $3$ are known \cite[Theorem 7.5.1]{BrouwerCN89}; all intersection arrays for distance-regular graphs with degree $4$ are known \cite{BrouwerK99}; and all intersection arrays for distance-regular graphs with girth $3$ and degree $6$ or $7$ are known. With the help of these results, van Dam and Jazaeri \cite{vanDamJ18} determined all Cayley graphs among such distance-regular graphs with small degree, as well as all Cayley graphs among distance-regular graphs of degree $5$ with one of the known feasible intersection arrays. See \cite{vanDamJ18} for details.

\delete
{
A PDS with parameters $(n, (n-1)/2, (n-5)/4, (n-1)/4)$ for some integer $n > 5$ is called a \emph{Paley} PDF \cite{LeungM95}.

In \cite[Corollary 3.1]{Davis94}, Davis constructed the first infinite family of Paley PDSs such that the underlying group has a prime power order but is not elementary abelian; they are defined on $\ZZZ_{p^2} \times \ZZZ_{p^2}$ with parameters $(p^4, (p^4-1)/2, (p^4-5)/4, (p^4-1)/4)$, where $p$ is a prime. Moreover,  he showed in \cite[Corollary 3.2]{Davis94} the existence of a PDS with Paley parameters on $\ZZZ_{p^2}^{2a} \times \ZZZ_{p}^{4b}$ whenever $a+b$ is a power of $2$. In \cite{LeungM95}, Leung and Ma used local rings to construct examples of Paley PDSs in abelian groups of arbitrary rank. They also obtained some new non-existence results. As an application, they proved that a Paley PDS in an abelian group $G$ of rank $2$ exists if and only if $G = \ZZZ_{p^s} \times \ZZZ_{p^s}$ for some odd prime $p$ and positive integer $s$.

In \cite[Theorem 1.6]{LeifmanM05}, Leifman and Muzychuk classified strongly regular Cayley graphs on $\mathbb{Z}_{p^n}\oplus\mathbb{Z}_{p^n}$ using Schur ring method, where $p$ is a prime and $n \ge 1$ is an integer. The reader is referred to \cite{LeifmanM05} for technical details of this result and its proof. As a consequence they obtained a complete classification of strongly regular Cayley graphs with Paley parameters (that is, $(n, (n-1)/2, (n-5)/4, (n-1)/4)$) over abelian groups of rank $2$.
}


\subsection{Distance-regular Cayley graphs with least eigenvalues $-2$}

In \cite{Seidel68}, Seidel classified all strongly regular graphs with least eigenvalue $-2$. In general, all regular graphs with least eigenvalue $-2$ have been classified (see \cite[Theorem 3.12.4]{BrouwerCN89}). In particular, a distance-regular graph with least eigenvalue $-2$ is strongly regular or the line graph of a regular graph with girth at least five. Using these results, Abdollahi, van Dam and Jazaeri classified in \cite{AbdollahiDJ17} all distance-regular Cayley graphs with least eigenvalue $-2$ and diameter at most three. We present their results in Theorems \ref{SRCG-21} and \ref{SRCG-22}, but before that let us first mention a few well-known strongly regular graphs.

The $5$-regular \emph{Clebsch graph} is the Cayley graph on $\ZZZ_2^4$ with connection set
$$
\{(1,0,0,0), (0,1,0,0), (0,0,1,0), (0,0,0,1), (1,1,1,1)\}.
$$
In other words, the $5$-regular Clebsch graph is the graph obtained from the $4$-dimensional hypercube $H(4, 2)$ by adding an edge between each pair of antipodal vertices. The complement of the $5$-regular Clebsch graph is called the $10$-regular \emph{Clebsch graph}. Both Clebsch graphs are strongly regular, with parameters $(16,5,0,2)$ and $(16, 10, 6, 6)$, respectively. The \emph{Schl\"{a}fli graph} is the graph obtained from $G=L(K_8)$ (the line graph of $K_8$) by applying the following operations (see \cite[Example 1.2.5]{Cvetkovic10}): Select a vertex $v$ of $G$; switch edges and non-edges of $G$ between $N_{G}(v)$ and $V(G) \setminus N_{G}(v)$, where $N_{G}(v)$ is the neighbourhood of $v$ in $G$; and delete $v$ from the resultant graph. This is the unique strongly regular graph with parameters $(27, 16, 10, 8)$ and is a Cayley graph on $\ZZZ_9 \rtimes \ZZZ_3$ or $(\ZZZ_3 \times \ZZZ_3) \rtimes \ZZZ_3$ (see \cite[Proposition 4.2]{AbdollahiDJ17}). The \emph{cocktail party graph} $CP(n)$ is the complete $n$-partite graph $K_{2, \ldots, 2}$ with each part of size two. Thus $CP(n)$ is the Cayley graph on any group $\Ga$ of order $2n$ containing an involution $a$ with connection set $\Ga \setminus \la a \ra$. The \emph{triangular graph} $T(n)$ is the line graph of $K_n$, and the \emph{lattice graph} $L_2(n)$ is the line graph of $K_{n,n}$. Obviously, $T(4)$ is isomorphic to $CP(3)$ and hence is a Cayley graph. As mentioned in \cite[Corollary 4.6]{AbdollahiDJ17}, for $n > 4$, $T(n)$ is a Cayley graph if and only if $n \equiv 3 ~(\mod~4)$ and $n$ is a prime power. The lattice graph $L_2(n)$ is isomorphic to the Hamming graph $H(2,n)$ and hence is a Cayley graph.

The following result relies on Seidel's classification \cite{Seidel68} of strongly regular graphs with least eigenvalue $-2$.

 \begin{thm}\emph{(\cite[Theorem 4.8]{AbdollahiDJ17})}\label{SRCG-21}
 A graph $G$ is a strongly regular Cayley graph with least eigenvalue at least
$-2$ if and only if $G$ is isomorphic to one of the following graphs:
\begin{itemize}
  \item[\rm (a)] cycle $C_5$, Clebsch graph, Shrikhande graph, or Schl\"{a}fli graph;
  \item[\rm (b)] cocktail party graph $CP(n)$ with $n \ge 2$;
  \item[\rm (c)] triangular graph $T(n)$, with $n = 4$, or $n \equiv 3 ~(\mod~4)$ and $n > 4$ a prime power;
  \item[\rm (d)] lattice graph $L_2(n)$, with $n \ge 2$.
\end{itemize}
\end{thm}

The \emph{dual} of an incidence structure $D$ is the incidence structure obtained from $D$ by interchanging the roles of points and blocks but retaining the incidence relation. The \emph{incidence graph} of $D$ is the bipartite graph with points in one part and blocks in the other part such that two vertices are adjacent if and only if they are incident in $D$. An incident point-block pair of $D$ is usually called a \emph{flag} of $D$. An \emph{isomorphism} from an incidence structure $D$ to an incidence structure $D'$ is a bijection that maps points to points, blocks to blocks, and preserves the incidence relation. An isomorphism from $D$ to itself is called an \emph{automorphism} of $D$.

A \emph{projective plane} is a point-line incidence structure such that any two distinct points are joined by exactly one line, any two distinct lines intersect in a unique point, and there exists a set of four points no three of which are on a common line. It is well known that for any finite projective plane $\pi$ there exists an integer $n \ge 1$, called the order of $\pi$, such that each line contains $n + 1$ points, each point is on $n + 1$ lines, and the number of points and the number of lines are both equal to $n^2 + n + 1$. In other words, $\pi$ is a symmetric $2$-$(n^2 + n + 1, n+1, 1)$ design. An automorphism of $\pi$ is usually called a \emph{collineation} of $\pi$. An isomorphism between $\pi$ and its dual is called a \emph{correlation} of $\pi$. A projective plane is Desarguesian or non-Desarguesian depending on whether it satisfies Desargues' Theorem.


The following result relies on the classification of distance-regular graphs with least eigenvalue $-2$ (see \cite[Theorems 3.12.4 and 4.2.16]{BrouwerCN89}).

\begin{thm}\emph{(\cite[Theorem 5.8]{AbdollahiDJ17})}\label{SRCG-22}
Let $G$ be a distance-regular Cayley graph with diameter three and least eigenvalue
at least $-2$. Then $G$ is isomorphic to one of the following graphs:
\begin{itemize}
  \item[\rm (a)] cycle $C_6$ or $C_7$;
  \item[\rm (b)] line graph of the incidence graph of the Desarguesian projective plane of order $2$ or
$8$;
  \item[\rm (c)] line graph of the incidence graph of a non-Desarguesian projective plane of order
$q$, where $q^2 + q + 1$ is prime and $q$ is even and at least $2 \times 10^{11}$;
  \item[\rm (d)] line graph of the incidence graph of a projective plane of odd order with a group of
collineations and correlations acting regularly on its flags.
\end{itemize}
\end{thm}


\section{Generalizations of Cayley graphs}
\label{sec:genCay}

The concept of a Cayley graph can be generalized or varied in several different ways.
In this section we discuss three avenues of generalization and variation, with an emphasis on eigenvalues of such generalized Cayley graphs.


\subsection{Eigenvalues of $n$-Cayley graphs}
\label{subsec:eignCay}

Let $\Ga$ be a finite group and $R, S, T$ be (not necessarily non-empty) subsets of $\Ga$ such that $R^{-1} = R, S^{-1} = S$ and $1 \notin R \cup S$. The \emph{bi-Cayley graph} $\Cay(\Ga; R, S, T)$ on $\Ga$ is the graph with vertex set $\Ga \times \{0, 1\}$ such that $(g, i), (h, j)$ are adjacent if and only if one of the following holds: (i) $i=j=0$ and $g^{-1}h \in R$; (ii) $i=j=1$ and $g^{-1}h \in S$; (iii) $i=0, j=1$ and $g^{-1}h \in T$. Bi-Cayley graphs are also called \emph{semi-Cayley graphs} in the literature (see, for example, \cite{deResminiJ92}), and a bi-Cayley graph on a cyclic group is called a \emph{bicirculant} or more specifically an \emph{$m$-bicirculant} if the cyclic group has order $m$. A bi-Cayley graph $\Cay(\Ga; R, S, T)$ such that $T = \{1\}$ is called a \emph{one-matching bi-Cayley graph} on $\Ga$. The famous Petersen and Hoffman-Singleton graphs are examples of bi-Cayley graphs.

In \cite{Gao10}, Gao and Luo proved the following result, which shows that the computation of eigenvalues of bi-Cayley graphs on abelian groups can be reduced to that of three Cayley graphs on the same group.

\begin{thm}
\label{thm:bi-Cay}
\emph{(\cite[Theorem 3.2]{Gao10})}
Let $\Ga = \ZZZ_{n_1} \oplus \cdots \oplus \ZZZ_{n_d}$ be an abelian group, where each $n_i \ge 2$, and let $G = \Cay(\Ga; R, S, T)$ be a bi-Cayley graph on $\Ga$. Then the eigenvalues of $G$ are given by
$$
\frac{1}{2} \left(\l_{r_1, \ldots, r_d}^R + \l_{r_1, \ldots, r_d}^S \pm \sqrt{(\l_{r_1, \ldots, r_d}^R - \l_{r_1, \ldots, r_d}^S)^2 + 4|\l_{r_1, \ldots, r_d}^T|^2}\right),
$$
for $r_j = 0,1,\ldots,n_j-1$, $1 \le j \le d$, where $\l_{r_1, \ldots, r_d}^R$, $\l_{r_1, \ldots, r_d}^S$ and $\l_{r_1, \ldots, r_d}^T$ are the eigenvalues of $\Cay(\Ga, R)$, $\Cay(\Ga, S)$ and $\Cay(\Ga, T)$, respectively.
\end{thm}

In particular, this implies that for abelian groups $\Ga$, $\Cay(\Ga; R, R, T)$ is integral provided that $\Cay(\Ga, R)$ is integral (see \cite[Corollary 3.5]{Gao10}).

In \cite{GaoLH16}, Gao, L\"u and Hao obtained formulas for the Laplacian and signless Laplacian spectra of bi-Cayley graphs on abelian groups. As applications of their main result, special formulas for the Laplacian and signless Laplacian spectra are also given for two classes of bi-Cayley graphs, namely one-matching bi-Cayley graphs and the join of two Cayley graphs over isomorphic abelian groups. In particular, a method for constructing Laplacian and signless Laplacian integral bi-Cayley graphs was given.

\begin{thm}
\label{thm:GaoLH16}
\emph{(\cite[Theorem 1]{GaoLH16})}
Let $\Ga = \ZZZ_{n_1} \oplus \cdots \oplus \ZZZ_{n_d}$ be an abelian group, where each $n_i \ge 2$, and let $G = \Cay(\Ga; R, S, T)$ be a bi-Cayley graph on $\Ga$. Then the Laplacian eigenvalues (respectively, signless Laplacian eigenvalues) of $G$ are given by
$$
\frac{1}{2} \left(\mu_{r_1, \ldots, r_d}^R + \mu_{r_1, \ldots, r_d}^S + 2|T| \pm \sqrt{(\mu_{r_1, \ldots, r_d}^R - \mu_{r_1, \ldots, r_d}^S)^2 + 4|\l_{r_1, \ldots, r_d}^T|^2}\right),
$$
for $r_j = 0,1,\ldots,n_j-1$, $1 \le j \le d$, where $\l_{r_1, \ldots, r_d}^T$ are the eigenvalues of $\Cay(\Ga, T)$, and $\mu_{r_1, \ldots, r_d}^R$ and $\mu_{r_1, \ldots, r_d}^S$ are the Laplacian eigenvalues (respectively, signless Laplacian eigenvalues) of $\Cay(\Ga, R)$ and $\Cay(\Ga, S)$, respectively.
\end{thm}

In particular, this implies that for abelian groups $\Ga$, $\Cay(\Ga; R, R, T)$ is a Laplacian and signless Laplacian integral graph provided that $\Cay(\Ga, R)$ and $\Cay(\Ga, T)$ are both integral (\cite[Corollary 4.6]{GaoLH16}).

A few special cases and applications of Theorem \ref{thm:GaoLH16} were also given in \cite{GaoLH16}. In particular, the Laplacian and signless Laplacian eigenvalues of the \emph{I-graphs} $I(n, j, k) := \Cay(\ZZZ_n, \{\pm j\}, \{\pm k\}, \{0\})$ were obtained. These bi-Cayley graphs on $\ZZZ_n$ are interesting generalizations of the generalized Petersen graphs, and their Laplacian and signless Laplacian eigenvalues can be obtained using Theorem \ref{thm:GaoLH16}.

\begin{cor}
\emph{(\cite[Corollary 4.2]{GaoLH16})}
The Laplacian eigenvalues (respectively, signless Laplacian eigenvalues) of the I-graph $I(n, j, k)$ are given by
$$
\frac{1}{2} \left(\mu_{r}^j + \mu_{r}^k + 2 \pm \sqrt{(\mu_{r}^j - \mu_{r}^k)^2 + 4}\right),
$$
for $r = 0,1,\ldots,n-1$, where $\mu_{r}^j$ and $\mu_{r}^k$ are the Laplacian eigenvalues (respectively, signless Laplacian eigenvalues) of $\Cay(\Ga, \{\pm j\})$ and $\Cay(\Ga, \{\pm k\})$, respectively.
\end{cor}

In \cite{ZouMeng07}, Zou and Meng determined the eigenvalues of the bi-circulants $\Cay(\ZZZ_n; \emptyset, \emptyset, T)$. This was achieved through investigating connections between the eigenvalues of the Cayley digraph $\Cay(\Ga, T)$ and those of the bi-Cayley graph $\Cay(\Ga; \emptyset, \emptyset, T)$ for any finite abelian group $\Ga$. They also obtained asymptotic results on the number of spanning trees of bi-circulants.

A finite group $\Ga$ is called \emph{bi-Cayley integral} if $\Cay(\Ga; \emptyset, \emptyset, T)$ is an integral graph for any $T \subseteq \Ga$. This concept was introduced by Arezoomand and Taeri in \cite{ArezoomandT15}, where it was proved that a finite group is a bi-Cayley integral group if and only if it is isomorphic to $\ZZZ_3$, $S_3$ or $\ZZZ_2^d$ for some positive integer $d$.

Let $R$ be a finite commutative ring with unit element $1 \ne 0$, and let $R^{\times}$ be the set of units of $R$. The one-matching bi-Cayley graph $\Cay(R; R^{\times}, R^{\times}, \{0\})$ on the additive group of $R$ is called \cite{Liu19} the \emph{unitary one-matching bi-Cayley graph} of $R$. In \cite{Liu19}, Liu computed the energies of this graph and its line graph and determined when these graphs are hyperenergetic. He also proved \cite{Liu19} that neither of these graphs can be \emph{hypoenergetic} (that is, with energy strictly less than the number of vertices). In \cite{LiuYan20}, the bi-Cayley graph $\Cay(R; R^{\times}, R^{\times}, R^{\times})$ on the additive group of $R$ is called the \emph{unitary homogeneous bi-Cayley graph} of $R$. In this paper, Liu and Yan computed the energies of this graph, its complement graph, and its line graph, and determined when these graphs are hyperenergetic. In the same paper they also determined when these graphs are Ramanujan.

In general, a graph (or digraph) $G$ is called an \emph{$n$-Cayley graph} (or $n$-Cayley digraph) if $\Aut(G)$ contains a semiregular subgroup $\Ga$ that has exactly $n$ orbits on $V(G)$; we also say that $G$ is an $n$-Cayley graph (or digraph) on the group $\Ga$. A $1$-Cayley graph is precisely a Cayley graph and a $2$-Cayley graph is exactly a bi-Cayley graph. A $3$-Cayley graph is often called a \emph{tri-Cayley graph}, and in particular a $3$-Cayley graph on a cyclic group is called a \emph{tricirculant} or more specifically \emph{an $m$-tricirculant} if the cyclic group has order $m$. In \cite{ArezoomandT13}, Arezoomand and Taeri gave for any $n \ge 2$ a factorization of the characteristic polynomials of $n$-Cayley digraphs in terms of the irreducible representations of the underlying group. As applications they obtained a formula for the eigenvalues of any $n$-Cayley graph on an abelian group (thus generalizing Theorem \ref{thm:bi-Cay}), determined the spectrum of any Cayley graph on a finite group containing a cyclic subgroup of index $2$, and obtained a few results on the main eigenvalues of an $n$-Cayley digraph. A few results on the eigenvalues of ``quasi-abelian" $n$-Cayley graphs can be found in \cite{Arezoomand20}.

The concept of an $n$-Cayley graph was also studied by Sjogren in \cite{Sjogren94} under a different name. In \cite{Sjogren94}, a \emph{$\Ga$-graph} was defined as a graph $G$ whose automorphism group contains a semiregular subgroup which is isomorphic to $\Ga$. Note that this is exactly an $n$-Cayley graph on $\Ga$, where $n$ is the number of orbits of $\Ga$ on $V(G)$. (The term ``regular" in \cite{Sjogren94} is meant ``semiregular" by our definition in Section \ref{sec:int}.) Among other things it was proved in \cite{Sjogren94} that the Laplacian spectrum of any $n$-Cayley graph $G$ can be expressed in terms of the Laplacian spectrum of the quotient graph of $G$ with respect to the partition of $V(G)$ into such orbits. As a corollary it was proved that the Laplacian spectrum of any Cayley graph of odd order consists of $0$ and a set of eigenvalues each with an even multiplicity.

The reader is referred to \cite{Arezoomand} for a recent survey on eigenvalues, Laplacian eigenvalues and signless Laplacian eigenvalues of $n$-Cayley graphs.


\subsection{Strongly regular $n$-Cayley graphs}
\label{subsec:srnCay}

Throughout this section we set
$$
\De = (\l-\mu)^2+4(k-\mu)
$$
when discussing a strongly regular graph with parameters $(v, k, \l, \mu)$.

In \cite{Marusic88}, Maru\v{s}i\v{c} obtained necessary conditions for a bicirculant or tricirculant to be strongly regular, starting up the study of strongly regular $n$-Cayley graphs. For example, he proved that if a strongly regular $p$-bicirculant exists, where $p$ is a prime, then $p$ must be of the form $p = 2s^2 + 2s + 1$ for some positive integer $s$, and the parameters of this strongly regular bicirculant satisfy $\mu = \l + 1 = s^2$ or $(s+1)^2$. In \cite[Theorem 2.2]{deResminiJ92}, de Resmini and Jungnickel proved that a bi-Cayley graph
$\Cay(\Ga; R, S, T)$ on a group $\Ga$ of order $n$ is strongly regular with parameters $(2n, k, \l, \mu)$ if and only if $k, \l, \mu$ and $R, S, T$ satisfy three equations in the group ring $\ZZZ \Ga$. Any triple $(T, R, S)$ satisfying these equations for suitable parameters $k, \l, \mu$ is called \cite{deResminiJ92} a \emph{partial difference triple}. Thus strongly regular bi-Cayley graphs are essentially partial difference triples, akin to the relation between strongly regular Cayley graphs and partial difference sets. Using this relation, the following result was obtained in \cite{deResminiJ92}.

\begin{thm}
\emph{(\cite[Proposition 2.3]{deResminiJ92})}
Let $\Cay(\Ga; R, S, T)$ be a strongly regular bi-Cayley graph with parameters $(2n, k, \l, \mu)$, where $n$ is the order of $\Ga$. Then
$$
|T| = \frac{2k - (\l-\mu) \pm \sqrt{\De}}{4},\; |R| = |S| = k - |T|,
$$
and moreover $\sqrt{\De}$ is an integer.
\end{thm}

In the case when $\Ga$ is abelian, the following necessary conditions were obtained in \cite{deResminiJ92}.

\begin{thm}
\emph{(\cite[Proposition 3.1 and Corollary 3.3]{deResminiJ92})}
Let $\Cay(\Ga; R, S, T)$ be a strongly regular bi-Cayley graph with parameters $(2n, k, \l, \mu)$, where $\Ga$ is an abelian group with order $n$. Then, for every non-principal character $\chi$ of $\Ga$, either
$$
\chi(T) = 0\ \text{ and }\ \chi(R) = \chi(S) = \frac{(\l-\mu) \pm \sqrt{\De}}{2}
$$
or
$$
\chi(R) + \chi(S) = \l-\mu.
$$
In particular, $\chi(R + S) \in \{(\l-\mu)-\sqrt{\De}, \l-\mu, (\l-\mu)+\sqrt{\De}\}$ for every non-principal character $\chi$ of $\Ga$. Moreover, in the special case when $R = S$, all of $\l$, $\mu$ and $\sqrt{\De}$ must be even.
\end{thm}

In \cite[Theorems 4.2 and 4.4]{deResminiJ92}, de Resmini and Jungnickel proved the following result (which generalizes the above-mentioned result of Maru\v{s}i\v{c} \cite{Marusic88}): Let $G = \Cay(\ZZZ_n; R, S, T)$ be a nontrivial strongly regular $n$-bicirculant with parameters $(2n, k, \l, \mu)$. If $n$ is odd or is not divisible by $\sqrt{\De}$, then up to complementation the parameters of $G$ are given by
$$
(n, |T|, |R|, \l, \mu) = (2s^2 + 2s + 1, s^2, s^2 + s, s^2 - 1, s^2)
$$
for some positive integer $s$, and moreover $G$ is trivial if $R = S$. (Note that $k = |T| + |R|$.) Solving a problem posed in \cite{deResminiJ92}, Leung and Ma proved in \cite{LeungM93} that there are other forms of parameters for nontrivial strongly regular bicirculants. In fact, they determined all possible parameters as stated in the next theorem.

\begin{thm}
\emph{(\cite[Theorem 3.1]{LeungM93})}
Let $G = \Cay(\ZZZ_n; R, S, T)$ be a nontrivial strongly regular bicirculant with parameters $(2n, k, \l, \mu)$. Then up to complementation the parameters of $G$ are in one of the following forms:
\begin{itemize}
\item[\rm (a)] $(n, |T|, |R|, \l, \mu) = (2s^2 + 2s + 1, s^2, s^2 + s, s^2 - 1, s^2)$, where $s \ge 1$;
\item[\rm (b)] $(n, |T|, |R|, \l, \mu) = (2s^2, s^2, s^2 - s, s^2 - s, s^2 - s)$, where $s \ge 2$;
\item[\rm (c)] $(n, |T|, |R|, \l, \mu) = (2s^2, s^2, s^2 + s, s^2 + s, s^2 + s)$, where $s \ge 3$;
\item[\rm (d)] $(n, |T|, |R|, \l, \mu) = (2s^2, s^2 \pm s, s^2, s^2 \pm s, s^2 \pm s)$, where $s \ge 2$.
\end{itemize}
\end{thm}

Another interesting result proved in \cite{LeungM93} is as follows.

\begin{thm}
\emph{(\cite[Theorems 2.1 and 2.2]{LeungM93})}
Let $G = \Cay(\Ga; R, S, T)$ be a nontrivial strongly regular bi-Cayley graph with parameters $(2n, k, \l, \mu)$, where $n$ is the order of $\Ga$. If $R \cup S$ is contained in a proper normal subgroup $\Sigma$ of $\Ga$, then one of the following holds:
\begin{itemize}
\item[\rm (a)] $(n, |T|, |R|, \l, \mu) = (8, 4, 1, 0, 2)$ and $|\Sigma| = 2$ or $4$;
\item[\rm (b)] $(n, |T|, |R|, \l, \mu) = (25, 5, 2, 0, 1)$ and $|\Sigma| = 5$;
\item[\rm (c)] $(n, |T|, |R|, \l, \mu) = (2s^2, 2s, 2s - 2, 2s - 2, 2)$ and $|\Sigma| = s^2$ for some integer $s \ge 2$.
\end{itemize}
Moreover, the third case occurs if and only if there exist subgroups $K, L$ of $\Ga$ and elements $g \in \Sigma$ and $h \in \Ga \setminus \Sigma$ with $|K| = |L| = s$, $K \cap L = \{1\}$ and $g^{-1}Kg \cap h^{-1}Lh = \{1\}$ such that $T = Kg \cup Lh$, $R = (K \cup L) \setminus \{1\}$ and $S = (g^{-1}Kg \cup h^{-1}Lh) \setminus \{1\}$.
\end{thm}

In \cite{deResminiJ92}, some restrictions on strongly regular Cayley graphs with an even number of vertices were obtained using the results in the paper on strongly regular bi-Cayley graphs. In particular, it was proved in \cite[Theorems 5.3 and 5.6]{deResminiJ92} that if $\Ga$ is the dihedral group $D_{2n}$, the generalized quaternion group $Q_{4m}$ (see \eqref{eq:Q4m}), or a $2$-group of order $|\Ga| \ne 16, 64$ with a cyclic subgroup of index $2$, then there is no nontrivial strongly regular Cayley graph on $\Ga$. In \cite[Theorem 4.2]{LeungM93}, it was further proved that in the third case the same statement is also true when $|\Ga| = 16, 64$ unless $\Ga$ is one of two specific groups of order $16$ (which produce isomorphic strongly regular Cayley graph with parameters $(16, 6, 2, 2)$).

In \cite{MalnicMS07}, Malni\v{c}, Maru{\v{s}}i{\v{c}} and {\v{S}}parl gave a necessary condition for the existence of a strongly regular vertex-transitive bicirculant of order twice a prime and constructed three new strongly regular bicirculants having 50, 82 and 122 vertices, respectively. These graphs together with their complements form the first known pairs of complementary strongly regular bicirculants which are vertex-transitive but not edge-transitive.

In line with the work in \cite{deResminiJ92, LeungM93}, a theory of strongly regular tri-Cayley graphs was developed by Kutnar, Maru\v{s}i\v{c}, Miklavi\v{c} and \v{S}parl in \cite{KutnarMMS09}. Let $\Ga$ be a group with order $n$, and let $A, B, C, R, S$ and $T$ be subsets of $\Ga$ such that $A^{-1} = A, B^{-1} = B, C^{-1} = C$ and $1 \not \in A \cup B \cup C$. Define $G = \Cay(\Ga; A, B, C; R, S, T)$ to be the graph with vertex set $\Ga \times \ZZZ_3$ and adjacency relation $\sim$ as follows: for $g \in \Ga$, $(g, 0) \sim (ag, 0)$ for $a \in A$, $(g, 1) \sim (bg, 1)$ for $b \in B$, $(g, 2) \sim (cg, 2)$ for $c \in C$, $(g, 0) \sim (rg, 1)$ for $r \in R$, $(g, 1) \sim (sg, 2)$ for $s \in S$, and $(g, 2) \sim (tg, 0)$ for $t \in T$. It can be verified that $G$ is a tri-Cayley graph and conversely any tri-Cayley graph is of this form. Similar to \cite[Theorem 2.2]{deResminiJ92}, it was shown in \cite[Proposition 3.1]{KutnarMMS09} that $G$ is strongly regular with parameters $(3n, k, \l, \mu)$ if and only if $k, \l, \mu$ and the subsets $A, B, C, R, S, T$ of $\Ga$ satisfy six equations in the group ring $\ZZZ \Ga$. Based on these equations, Kutnar, Maru\v{s}i\v{c}, Miklavi\v{c} and \v{S}parl \cite[Propositions 3.2 and 3.3]{KutnarMMS09} obtained a necessary condition for a tri-Cayley graph to be strongly regular. In particular, they proved that $ \sqrt{\De}$ must be an integer if $G$ is a nontrivial strongly regular graph. In the case when $\Ga$ is abelian, they proved the following result. Recall that the triangular graph $T(n)$ is the line graph of $K_n$.

\begin{thm}
\emph{(\cite[Proposition 4.1 and Corollary 4.2]{KutnarMMS09})}
Let $\Ga$ be an abelian group of order $n$, and let $G = \Cay(\Ga; A, B, C; R, S, T)$ be a tri-Cayley graph on $\Ga$. If $G$ is strongly regular with parameters $(3n, k, \l, \mu)$, then
$$
\chi(A) + \chi(B) + \chi(C) \in \left\{\frac{3((\l-\mu) \pm \sqrt{\De})}{2}, \frac{3(\l-\mu) \pm \sqrt{\De}}{2}\right\}
$$
for every non-principal character $\chi$ of $\Ga$. In particular, if in addition $B = \emptyset$, $|A|=|C|$,  and $G$ is nontrivial, then either $A = \emptyset$ and $G$ is isomorphic to the Paley graph $P(9)$, or $|A|=2$ and $G$ is isomorphic to the complement of the triangular graph $T(6)$.
\end{thm}

In \cite{KutnarMMS09}, Kutnar, Maru\v{s}i\v{c}, Miklavi\v{c} and \v{S}parl also obtained possible parameters for strongly regular tricirculants in two special cases.

\begin{thm}
\emph{(\cite[Theorem 5.3]{KutnarMMS09})}
Let $G = \Cay(\ZZZ_n; A, B, C; R, S, T)$ be a tricirculant with $|A| + |B| + |C| \leq |\Ga \setminus A| + |\Ga \setminus B| + |\Ga \setminus C|$ and $A \cup B \cup C \ne \emptyset$. Suppose that $G$ is a nontrivial strongly regular graph with parameters $(3n, k, \l, \mu)$. Suppose further that either $n$ is prime or $n$ is coprime to $6 \sqrt{\De}$. Then there exists an integer $s$ such that the following hold:
\begin{itemize}
\item[\rm (a)] if $|A|$, $|B|$ and $|C|$ are not all equal, then
$$
(3n, k, \l, \mu) = (3(12s^2 +9s+2), (3s+1)(4s+1), s(4s+3), s(4s+1));
$$
moreover, if in addition $|A| = |C| \ne |B|$, then
$$
|A| = 2s(2s+1),\ |B| = (s + 1)(4s + 1),\ |R| = s(4s + 1),\ |T|=(2s+1)^2;
$$
\item[\rm (b)] if $|A| = |B| = |C|$, then
$$
(3n, k, \l, \mu) = (3(3s^2 - 3s + 1), s(3s-1), s^2 + s - 1, s^2)
$$
and
$$
|A| = s(s-1),\ |R| = |S| = |T| = s^2.
$$
\end{itemize}
\end{thm}

In \cite{KutnarMMM13}, Kutnar, Malni\v{c}, Mart\'{i}nez and Maru\v{s}i\v{c} introduced a new class of graphs, called \emph{quasi $n$-Cayley graphs}. Such graphs admit a group of automorphisms that fixes one vertex of the graph and acts semiregularly with $n$-orbits on the set of the rest vertices. They determined when these graphs are strongly regular, leading to the definition of a quasi-partial difference family (or QPDF for short), and gave several infinite families and sporadic examples of QPDFs. They also studies properties of QPDFs and determined, under several conditions, the form of the parameters of QPDFs when the group involved is cyclic.

A brief survey of strongly regular graphs and digraphs admitting a semiregular cyclic group of automorphisms can be found in the first part of \cite{Martinez15} by Mart\'{i}nez. In the second part of the same paper, Mart\'{i}nez studied some new types of such digraphs. Using partial sum families, he determined the form of the parameters and obtained some directed strongly regular graphs derived from these partial sum families with previously unknown parameters.

Cayley graphs on semigroups are generalizations of Cayley graphs on groups. Given a semigroup $\Ga$ and a subset $S$ of $\Ga$, the Cayley (di)graph $\Cay(\Ga, S)$ is defined to be the digraph with vertex set $\Ga$ such that there is an arc $(x, y)$ from $x$ to $y$, where $x \ne y$, if and only if $y = sx$ for some $s \in S$.
Cayley (di)graphs on semigroups have been studied extensively due to their wide applications in various areas including automata theory, but here we only mention one result which is related to strongly regular graphs. A completely $0$-simple inverse semigroup is called a \emph{Brandt semigroup}. In \cite[Theorem 2.12]{HaoGL11}, Hao, Gao and Luo obtained a necessary and sufficient condition for the connected components of a Cayley graph on a Brandt semigroup to be strongly regular, and as corollaries they recovered two known results on strongly regular Cayley graphs \cite{Ma94} and strongly regular bi-Cayley graphs \cite{deResminiJ92} on groups.


\subsection{Cayley sum graphs}
\label{subsec:CaySum}

Let $\Ga$ be a finite abelian group and $S$ a subset of $\Ga$. The \emph{Cayley sum graph} $\Cay^{+}(\Ga, S)$ on $\Ga$ with respect to $S$ is the graph with vertex set $\Ga$ in which $x, y \in \Ga$ are adjacent if and only if $x+y \in S$. Since $\Ga$ is abelian, this is an undirected graph of degree $|S|$, with a loop at $x$ if $x + x \in S$. If $S$ is required to be \emph{square-free} (that is, $x + x \not \in S$ for any $x \in \Ga$), then $\Cay^{+}(\Ga, S)$ is an undirected graph without loops. Cayley sum graphs are also known as \emph{addition Cayley graphs}, \emph{addition graphs}, and \emph{sum graphs} \cite{Chung89} in the literature.

Let $\Ga = \ZZZ_{n_1} \oplus \cdots \oplus \ZZZ_{n_d}$ be a finite abelian group, where each $n_i \ge 2$. Recall that the characters of $\Ga$ are $\chi_{\mathbf{x}}$ as given in \eqref{eq:chixg}, for $\mathbf{x} = (x_1, \ldots, x_d) \in \Ga$. The set of real-valued characters is $R = \{\chi_{\mathbf{x}}: \mathbf{x} \in \Ga, \mathbf{x} + \mathbf{x} = \mathbf{0}\}$. Let $C$ be a set containing exactly one character from each conjugate pair $\{\chi_{\mathbf{x}}, \chi_{-\mathbf{x}}\}$ (where $\mathbf{x} \in \Ga$ and $\mathbf{x} + \mathbf{x} \ne \mathbf{0}$). Then the set of characters of $\Ga$ is $R \cup \{\chi, \overline{\chi}: \chi \in C\}$. The following result was observed by a few authors in \cite{Alon07,Chung89,DeVosGMS09}.

\begin{thm}
\emph{(\cite[Theorem 2.1]{DeVosGMS09})}
\label{thm:DeVosGMS09}
Let $\Ga = \ZZZ_{n_1} \oplus \cdots \oplus \ZZZ_{n_d}$ be a finite abelian group, where each $n_i \ge 2$, and let $R$ and $C$ be as above. Let $S$ be a subset of $\Ga$. Then the multiset of eigenvalues of the Cayley sum graph $\Cay^{+}(\Ga, S)$ is
$$
\{\chi(S): \chi \in R\} \cup \{\pm |\chi(S)|: \chi \in C\}.
$$
Moreover, the corresponding eigenvectors are $\chi$ (for $\chi \in R$) and the real and imaginary parts of $\a \chi$ (for $\chi \in C$ with a suitable complex scalar $\a$ depending on $\chi(S)$ only).
\end{thm}

A finite abelian group $\Ga$ is called a \emph{Cayley sum integral group} \cite{AmooshahiT16a} if for every subset $S$ of $\Ga$ the Cayley sum graph $\Cay^{+}(\Ga, S)$ is integral. In \cite{AmooshahiT16a}, Amooshahi and Taeri proved that all Cayley sum integral groups are represented by $\ZZZ_3$ and $\ZZZ_2^d$, where $d \ge 1$. In the same paper they also classified all simple connected cubic integral Cayley sum graphs.

A finite abelian group $\Ga$ is called \emph{Cayley sum integral} if all Cayley sum graphs on any subgroup of $\Ga$ with respect to square-free subsets are integral. For an integer $k \ge 2$, let $\AA_k$ be the class of finite abelian groups $\Ga$ such that for any subgroup $\Sigma$ of $\Ga$ all Cayley sum graphs $\Cay^{+}(\Sigma, S)$ on $\Sigma$ with $|S| = k$  are integral.

\begin{thm}
\emph{(\cite[Theorems 1-3]{MaW16})}
The following hold:
\begin{itemize}
\item[\rm (a)] $\AA_2$ consists of the following groups: $\ZZZ_2^n$, $\ZZZ_4$, $\ZZZ_2 \times \ZZZ_3^m$, $n \ge 2, m \ge 1$;
\item[\rm (b)] $\AA_3$ consists of the following groups: $\ZZZ_2^n$, $\ZZZ_6$, $\ZZZ_8$, $n \ge 2$;
\item[\rm (c)] all finite groups which are Cayley sum integral are represented by: $\ZZZ_2^n$, $\ZZZ_4$, $\ZZZ_6$, $n \ge 1$.
\end{itemize}
\end{thm}

The notion of a Cayley sum graph was generalized to all finite groups in \cite{AmooshahiT16}. Let $\Ga$ be a finite group (with multiplicative operation) and $S$ a subset of $\Ga$. The \emph{Cayley sum graph} $\Cay^{+}(\Ga, S)$ on $\Ga$ with respect to $S$ is the digraph with vertex set $\Ga$ such that there is an arc from $x$ to $y$ if and only if $xy \in S$. In the case when $x^2 \in S$ for some $x \in \Ga$, this digraph contains a semi-edge from $x$; such a semi-edge contributes only one to the degree of its end-vertex and the corresponding diagonal entry in the adjacency matrix. We can extend the definition of a Cayley sum graph to the case when $S$ is a multiset of $\Ga$, in which case $\Cay^{+}(\Ga, S)$ has parallel arcs. If $S$ is normal (that is, the union of some conjugacy classes of $\Ga$), then $\Cay^{+}(\Ga, S)$ can be considered as an undirected graph.

\begin{thm}
\emph{(\cite[Theorem 7]{AmooshahiT16})}
Let $\Ga$ be a finite non-abelian group. There exists a square-free normal subset $S$ of $\Ga$ with $|S| = 3$ such that $\Cay^{+}(\Ga, S)$ is connected and integral if and only if $\Ga$ is isomorphic to the dihedral group $D_6$ of order $6$.
\end{thm}

Given a finite group $\Ga$ and a real-valued function $\a: \Ga \rightarrow \RRR$, the \emph{Cayley sum colour graph} $\Cay^{+}(\Ga, \a)$ is the complete directed graph with vertex set $\Ga$ in which each arc $(x, y) \in \Ga \times \Ga$ is associated with \emph{color} $\a(xy)$ (see \cite{AmooshahiT15} and \cite{AmooshahiT16}). In the case when $\a$ is the characteristic function of a subset $S$ of $\Ga$, $\Cay^{+}(\Ga, \a)$ becomes the Cayley sum graph $\Cay^{+}(\Ga, S)$ when arcs of colour $0$ are omitted. In \cite{AmooshahiT15} and \cite[Section 3]{AmooshahiT16}, Cayley sum colour graphs and their eigenvalues via irreducible representations of the underlying groups have been studied. The following result gives Theorem \ref{thm:DeVosGMS09} in the special case when $\a$ is the characteristic function of $S \subseteq \Ga$.

\begin{thm}
\emph{(\cite[Theorem 3.4]{AmooshahiT15})}
\label{thm:AmooshahiT15}
Let $\Cay^{+}(\Ga, \a)$ be a Cayley sum colour graph of a finite abelian group $\Ga$ and $\{\chi_1, \chi_2, \ldots, \chi_n\}$ a complete set of irreducible inequivalent characters of $\Ga$. If $\chi_k$ is real-valued, then $\sum_{g \in \Ga} \a(g) \chi_{k}(g)$ is an eigenvalue of $\Cay^{+}(\Ga, \a)$; if $\chi_k$ is not real-valued, then $\pm |\sum_{g \in \Ga} \a(g) \chi_{k}(g)|$ are two eigenvalues of $\Cay^{+}(\Ga, \a)$.
\end{thm}

The \emph{anti-shift operator} maps every vector $(c_0, c_1, \ldots, c_{n-1})$ to $(c_1, c_2, \ldots, c_{n-1}, c_0)$. An $n \times n$ \emph{anti-circulant matrix} is a matrix whose $i$th row is obtained from the first row by applying anti-shift operator $i-1$ times. A graph is called \emph{anti-circulant} if it has an anti-circulant adjacency matrix. It was observed in \cite[Lemma 1.1]{AmooshahiT15} that a graph is anti-circulant if and only if it is a Cayley sum graph of a cyclic group. Using this and Theorem \ref{thm:AmooshahiT15}, the eigenvalues of any anti-circulant graph were computed by Amooshahi and Taeri in \cite[Theorem 3.6]{AmooshahiT15}.

Sidon sets in groups are important objects of study which have numerous applications. Let $\Ga$ be a finite abelian group and $S$ a subset of $\Ga$. If for every nonzero element $x$ of $\Ga$, there is at most one pair $(a, b)$ of elements of $\Ga$ such that $x = a-b$, then $S$ is called a \emph{Sidon set} in $\Ga$. It can be verified that $|S| < \sqrt{|\Ga|} + (1/2)$ for any Sidon set $S$ of $\Ga$. Sidon sets of size close to this bound, namely those with size $\sqrt{|\Ga|} - \d$ for a small number $\d$, are the most interesting. The following result shows that, if $S$ is a Sidon set of size $\sqrt{|\Ga|} - \d$, then the Cayley sum graph $\Cay^{+}(\Ga, S)$ is pseudo-random. As before, we use $\l(G)$ to denote the maximum absolute value of the eigenvalues of a graph $G$ other than the largest one.

\begin{thm}
\emph{(\cite[Theorem 2.3]{Vinh13})}
Let $\Ga$ be a finite abelian group and $S$ a Sidon set of $\Ga$ with $|S| = \sqrt{|\Ga|} - \d$. Then the Cayley sum graph $\Cay^{+}(\Ga, S)$ satisfies
$$
\l(\Cay^{+}(\Ga, S)) \le \sqrt{2(1+\d) |X|^{1/2}}.
$$
\end{thm}

We now move on to Cayley sum graphs of finite commutative rings. Let $R$ be a finite commutative ring. The \emph{unitary Cayley sum graph} of $R$ is defined as the Cayley sum graph $G^{+}_R=\Cay^{+}(R,R^\times)$ on the additive group of $R$, where as before $R^\times$ is the set of units of $R$. This graph was studied by Podest\'{a} and Videla in \cite{PodestaV20} in the context of equienergetic and non-cospectral graphs. In the case when $R$ is a finite local ring with unique maximal ideal $M$, we have $R^\times = R \setminus M$ and the spectrum of $G^{+}_R$ was shown \cite{PodestaV20} to be
\begin{equation}
\label{eq:GRplus}
\Spec (G^{+}_R) =
\left(|R|-m, m^{\frac{|R|-m}{2m}}, 0^{\frac{|R|(m-1)}{m}}, (-m)^{\frac{|R|-m}{2m}}\right),
\end{equation}
where $m = |M|$. (Note that $0$ does not occur in this spectrum when $R$ is a finite field (that is, $m = 1$).) As noted in \cite{PodestaV20}, this together with \eqref{eq:GRlocal} implies that $G_R$ and $G^{+}_R$ are integral equienergetic non-cospectral graphs for any finite local ring $R$, where $G_R$ is the unitary Cayley graph of $R$ as defined in Section \ref{subsec:UCayComRing} and two graphs are said to be \emph{equienergetic} if they have the same energy. For a general finite commutative ring $R$ as in Assumption \ref{as:1}, we have $G^{+}_R = G^{+}_{R_1} \otimes \cdots \otimes G^{+}_{R_s}$ as observed in \cite{PodestaV20}. Similar to Theorem \ref{ringspectrum}, from this and \eqref{eq:GRplus} one can easily determine the spectrum of $G^{+}_R$. Using this, in \cite{PodestaV20} the following result was obtained and a necessary and sufficient condition for $G^{+}_R$ to be Ramanujan was given.

\begin{thm}
\emph{(\cite[Theorem 3.4]{PodestaV20})}
Let $R$ be a finite commutative ring with odd order. Then $G_R$ and $G^{+}_R$ are integral equienergetic non-cospectral graphs.
\end{thm}

In \cite{PalanivelC19}, Palanivel and Chithra computed explicitly the eigenvalues of the unitary Cayley sum graph $G_{\ZZZ_{n}}^{+} = \Cay^{+}(\ZZZ_{n}, \ZZZ_{n}^{\times})$ and its complement, and based on these they obtained lower bounds on the energies of these two graphs and gave necessary and sufficient conditions for them to be hyperenergetic.

As the counterpart of the generalized Paley graph $X_{q}^k$ (see Section \ref{subsec:GPaley}), the \emph{generalized Paley sum graph} $X^{k,+}_{q}$ is defined \cite{PodestaV20} as the Cayley sum graph on the additive group of $\FFF_q$ with respect to $\{a^k: a \in \FFF_q^*\}$, where $q = p^r$ is a prime power and $k \ge 1$ is an integer. Similar to what we saw in Section \ref{subsec:GPaley}, in studying $X^{k,+}_{q}$ we can always assume that $k$ is a divisor of $q-1$. As in Theorem \ref{PV20a}, let $\{\eta_1, \ldots, \eta_s\}$ be the set of distinct Gaussian periods other than $n$ among $\eta_{0}^{(k,q)}, \eta_{1}^{(k,q)}, \ldots, \eta_{k-1}^{(k,q)}$ (see \eqref{eq:gas-period}), and let $m_i$ be the multiplicity of $\eta_i$ for $1 \le i \le s$. Let $m \ge 0$ be the number of times that $n$ occurs in these Gaussian periods. Rearranging if necessary, we may assume that $\eta_1 = -\eta_2, \ldots, \eta_{2t-1} = -\eta_{2t}$ and $\eta_{i} \ne -\eta_{j}$ for distinct $i, j$ between $2t$ and $s$, where $t$ is some integer between $0$ and $\lfloor s/2 \rfloor$. Set $m'_{2i-1} = m'_{2i} = (m_{2i-1} + m_{2i})/2$ for $1 \le i \le t$. The following result obtained in \cite{PodestaV20} reduces the computation of the eigenvalues of $X^{k,+}_{q}$ to that of Gaussian periods.

\begin{thm}
\label{PV20d}
\emph{(\cite[Theorem 4.1]{PodestaV20})}
Let $q = p^r$ be an prime power with $r \ge 2$. Let $k \ge 1$ be a divisor of $q-1$ and set $n = (q-1)/k$. With the above notation, we have
\begin{eqnarray*}
\Spec (X^{k,+}_{q})=
\left(\begin{array}{cccccccc}
n  &  -n  &  \eta_1 & \ldots & \eta_{2t} & \pm \eta_{2t+1} & \ldots & \pm \eta_s \\ [0.1cm]
1+\frac{1}{2}mn & \frac{1}{2}mn & m'_{1}n & \ldots & m'_{2t}n & \frac{1}{2} m_{2t+1}n & \ldots & \frac{1}{2} m_{s}n
\end{array}
\right).
\end{eqnarray*}
Moreover, if $k$ divides $(q-1)/(p-1)$, then $X^{k,+}_{q}$ is integral.
\end{thm}

In the case when $k=3, 4$ and $k$ divides $(q-1)/(p-1)$, the eigenvalues of $X^{k,+}_{q}$ have been computed explicitly in \cite{PodestaV20}. As a corollary of Theorems \ref{PV20a} and \ref{PV20d}, it was noted that $X_{q}^k$ and $X^{k,+}_{q}$ are equienergetic but non-cospectral for $q = p^r$ with $r \ge 2$.

Similar to Theorem \ref{PV20b}, the eigenvalues of $X^{k,+}_{q}$ can be computed explicitly in the case when $(k, q)$ is a semiprimitive pair.

\begin{thm}
\label{PV20e}
\emph{(\cite[Theorem 4.3]{PodestaV20})}
Let $q = p^r$ be an odd prime power with $r \ge 2$ even. Let $k \ge 1$ be a divisor of $q-1$ and set $n = (q-1)/k$. Suppose that $(k, q)$ is a semiprimitive pair, with $t$ the smallest divisor of $r/2$ such that $k$ divides $p^t + 1$. Then
\begin{eqnarray*}
\Spec (X^{k,+}_{q})=\left(\begin{array}{ccc}
n  & \pm \frac{1}{k}\left((k-1)(-1)^{\frac{r}{2t}}\sqrt{q}+1\right) & \pm \frac{1}{k}\left((-1)^{\frac{r}{2t}}\sqrt{q}-1\right) \\ [0.1cm]
1  & \frac{1}{2}n  & \frac{1}{2}(k-1)n
\end{array}
\right).
\end{eqnarray*}
\end{thm}

All pairs $(k, q)$ such that $X^{k}_{q}$ and $X^{k,+}_{q}$ are equienergetic, non-cospectral and Ramanujan were identified in \cite[Theorem 6.7]{PodestaV20}; this happens if and only if $(k, q)$ is as in Theorem \ref{PV20c} with $p$ odd.

We finish this section by mentioning that the Cayley sum graph $\Cay^{+}(R, Z^{*}(R))$ on the additive group of a finite commutative ring $R$ was studied in \cite{AalipourA14}, where $Z^{*}(R)$ is the set of nonzero zero divisors of $R$.


\subsection{Group-subgroup pair graphs}
\label{subsec:gp-subgp}

The following generalization of Cayley graphs was introduced in \cite{Reyes-Bustos16}. Let $\Ga$ be a group, $\Si$ a subgroup of $\Ga$, and $S$ a subset of $\Ga$ such that $S \cap \Si$ is closed under taking inverse elements. The \emph{group-subgroup pair graph} (or \emph{pair-graph}) $\Cay(\Ga, \Si, S)$ is the graph with vertex set $\Ga$ and edges $\{x, xs\}$ for $x \in \Si$ and $s \in S$. Obviously, the subgraph of $\Cay(\Ga, \Si, S)$ induced by $\Si$ is the Cayley graph $\Cay(\Si, S \cap \Si)$, and in particular $\Cay(\Ga, \Ga, S)$ is the Cayley graph $\Cay(\Ga, S)$. In \cite{Reyes-Bustos16}, Reyes-Bustos investigated several combinatorial properties of pair-graphs and determined some eigenvalues of $\Cay(\Ga, \Si, S)$, including the largest eigenvalue. Among others the following result was proved in \cite{Reyes-Bustos16}.

\begin{thm}
\label{thm:Rey}
\emph{(\cite[Theorem 5.1]{Reyes-Bustos16})}
Let $\Ga$ be a group, and let $\Si$ be a subgroup of $\Ga$ with index $|\Ga:\Si| = k+1 \ge 2$, with $\{x_0 = 1, x_1, \ldots, x_k\}$ a set of right coset representatives. Let $S$ be a subset of $\Ga$ such that $S \cap \Si$ is inverse-closed and $S \setminus \Si \neq \emptyset$, and set $S_i = S \cap \Si x_{i}$ for $0 \le i \le k$. Then
$$
\mu^{\pm} = \frac{1}{2} \left( |S_{0}| \pm \left(|S_{0}|^2 + 4 \sum_{i=1}^k |S_{i}|^2\right)^{\frac{1}{2}} \right)
$$
are eigenvalues of $\Cay(\Ga, \Si, S)$. Moreover, the corresponding eigenfunctions $f^{\pm}$ are given by $f^{\pm}(g) = \mu^{\pm}$ for $g \in \Si$ and $f^{\pm}(g) = |S_{i}|$ for $g \in {\Si} x_{i}$, $1 \le i \le k$.
\end{thm}

A few other results on the spectra of pair-graphs and their applications in constructing Ramanujan graphs can be found in \cite[Sections 5-6]{Reyes-Bustos16}.

In \cite{Kimoto18}, Kimoto proved the following result which determines the eigenvalues of $\Cay(\Ga, \Si, S)$ in terms of the characters of $\Si$ in the case when $\Si$ is abelian.

\begin{thm}
\label{thm:Kim}
\emph{(\cite[Theorem 3.5]{Kimoto18})}
Let $\Ga$ be a group, and let $\Si$ be abelian subgroup of $\Ga$ with order $n = |\Si|$ and index $|\Ga:\Si| = k+1 \ge 2$, with $\{x_0 = 1, x_1, \ldots, x_k\}$ a set of right coset representatives. Let $S$ be a subset of $\Ga$ such that $S \cap \Si$ is inverse-closed, and let $\Si_{i}$ be the subset of $\Si$ defined by the condition $S \cap \Si x_{i} = \Si_{i} x_i$, for $0 \le i \le k$. Then the eigenvalues of $\Cay(\Ga, \Si, S)$ are
$$
\frac{1}{2} \left( \sum_{g \in \Si_{0}}\chi(g) \pm \left(\big(\sum_{g \in \Si_{0}}\chi(g)\big)^2 + 4 \sum_{i=1}^k \big|\sum_{g \in \Si_{i}}\chi(g)\big|^2\right)^{\frac{1}{2}} \right),
$$
with $\chi$ running over all characters of $\Si$, and $0$ with multiplicity at least $(k-1)n$.
\end{thm}

In \cite[Theorem 3.8]{Kimoto18}, Kimoto also obtained a lower bound for the second largest eigenvalue of $\Cay(\Ga, \Si, S)$ in the case when $\Ga$ is abelian and $k=1$.


\section{Directed Cayley graphs}
\label{sec:diCay}

Eigenvalues of directed Cayley graphs appeared at several places before this point. In this section we review some of the known results about them which were not mentioned in previous sections.


\subsection{Eigenvalues of directed Cayley graphs}
\label{subsec:eigdiCay}

Let $\Ga$ be a finite group and $S$ a subset of $\Ga \setminus \{1\}$. Consider the Cayley digraph $\Cay(\Ga, S)$ on $\Ga$ (which is viewed as an undirected graph when $S$ is inverse-closed). Denote by $\mathrm{Irred}(\Ga)$ the set of irreducible complex characters of $\Ga$ and by $A_S$ the adjacency matrix of $\Cay(\Ga, S)$. The following result was first noted by MacWilliams and Mann in \cite{MacWilliamsM68}.

\begin{thm}
\emph{(\cite{MacWilliamsM68}; see also \cite[Theorem 3.1]{DuceyJ14})}
Let $\Ga$ be a finite abelian group and $S$ a subset of $\Ga \setminus \{1\}$. Let $M$ denote the character table of $\Ga$, with rows indexed by $\mathrm{Irred}(\Ga)$ and columns indexed by $\Ga$ in the same order as the one underlying $A_S$. Then
$$
\frac{1}{|\Ga|} M A_S^t \overline{M}^t = \diag\left(\sum_{s \in S} \chi(s) \right)_{\chi \in \mathrm{Irred}(\Ga)}.
$$
Thus the eigenvalues of $\Cay(\Ga, S)$ take the form $\sum_{s \in S} \chi(s)$, for $\chi \in \mathrm{Irred}(\Ga)$.
\end{thm}

A matrix with real entries is called \emph{normal} if it commutes with its transpose. This is equivalent to saying that the sum of the squares of the absolute values of its eigenvalues equals the sum of the squares of the absolute values of its entries. It can be shown (see \cite[Lemma 1]{LyubshinS09}) that for a finite group $\Ga$ and a subset $S \subseteq \Ga \setminus \{1\}$, the adjacency matrix of $\Cay(\Ga, S)$ is normal if and only if $S^{-1} S = S S^{-1}$. The next two results were obtained by Lyubshin and Savchenko.

\begin{thm}
\emph{(\cite[Corollaries 3 and 4]{LyubshinS09})}
Let $\Ga$ be a finite group and $S$ a minimal generating set of $\Ga$. Then the adjacency matrix of $\Cay(\Ga, S)$ is normal if and only if any two elements of $S$ have at least one common out-neighbour in $\Cay(\Ga, S)$. Moreover, if $\Ga$ is of odd order, then the adjacency matrix of $\Cay(\Ga, S)$ is normal if and only if $\Ga$ is abelian.
\end{thm}

\begin{thm}
\emph{(\cite[Theorem 1]{LyubshinS09})}
Let $\Ga$ be a finite group. If the adjacency matrix of every Cayley digraph of degree two on $\Ga$ is normal, then $\Ga$ is abelian or $\Ga \cong Q_8 \times \ZZZ_2^d$ for some nonnegative integer $d$, where $Q_8$ is the quaternion group.
\end{thm}

A digraph is called \emph{strongly connected} if it contains a directed path from any vertex to any other vertex.

\begin{thm}
\emph{(\cite[Theorem 2]{LyubshinS09})}
Let $\Ga$ be a finite group. If there exists a strongly connected Cayley digraph of degree two on $\Ga$ whose adjacency matrix is normal, then one of the following holds:
\begin{itemize}
\item[\rm (a)] $\Ga$ is an abelian group of rank at most two;
\item[\rm (b)] $\Ga \cong \la a, c\mid a^{2k} = c^n = 1, a^{-1}ca = c^{-1}\ra$, for some $k \ge 1$ and $n \ge 3$;
\item[\rm (c)] $\Ga \cong \la a, c \mid a^{4p}=1, c^{2q} = a^{2p}, a^{-1}ca = c^{-1}\ra$, for some $p \ge 1$ and $q \ge 3$.
\end{itemize}
\end{thm}

The next four theorems were proved by Godsil in \cite{Godsil82}.

\begin{thm}
\emph{(\cite[Theorem 3.1]{Godsil82})}
Let $G$ be a digraph with maximum degree greater than one. Then there is a Cayley digraph $H$ such that the minimal polynomial of $G$ divides that of $H$.
\end{thm}

\begin{thm}
\label{thm:Godsil82}
\emph{(\cite[Corollary 3.2]{Godsil82})}
The following hold:
\begin{itemize}
\item[\rm (a)] every algebraic integer is an eigenvalue of some Cayley digraph;
\item[\rm (b)] there exist Cayley digraphs whose adjacency matrices are not diagonalizable.
\end{itemize}
\end{thm}

\begin{thm}
\emph{(\cite[Theorem 3.3]{Godsil82})}
Let $G$ be a graph with minimum degree greater than one. Then there is a Cayley graph $H$ such that the minimal polynomial of $G$ divides that of $H$.
\end{thm}

\begin{thm}
\emph{(\cite[Corollary 3.4]{Godsil82})}
If $\theta$ is an eigenvalue of a symmetric integral matrix, then it is an eigenvalue of some Cayley graph.
\end{thm}

The last result above was sharpened by Babai in the following way. A digraph is called \emph{arc-transitive} (respectively, \emph{$2$-arc-transitive}) if its automorphism group is transitive on its set of arcs (respectively, directed paths of length $2$).

\begin{thm}
\label{thm:Babai85}
\emph{(\cite[Theorem 1.2]{Babai85})}
Let $A$ be an integral matrix. Then there exists an arc-transitive digraph $G$ such that the minimal polynomial of $A$ divides that of the adjacency matrix of  $G$. Moreover, if $A$ is symmetric, then $G$ may be required to be a $2$-arc-transitive graph.
\end{thm}

This theorem implies that there exist arc-transitive digraphs with non-diagonalizable adjacency matrices, answering a question of Cameron who asked whether there exists an arc-transitive digraph whose adjacency matrix is not diagonalizable.

A major tool used in the proof of Theorem \ref{thm:Babai85} is the following result (see \cite[Theorem 1.3]{Babai85}) which is of interest for its own sake: Every finite regular multidigraph has a finite arc-transitive covering digraph, and every finite regular multigraph has a finite $2$-arc-transitive covering graph.

Cospectral Cayley digraphs were studied in \cite{Meng98}.



\subsection{Two families of directed Cayley graphs}
\label{subsec:dissle}

In this section we focus on two families of Cayley digraphs introduced by Friedman in \cite{Fri95} (see Section \ref{subsec:ssle} for a related family of Cayley graphs). Let $p$ be a prime. Define \cite{Fri95} $\mathrm{SUMPROD}(p)$ to be the digraph with vertex set $\ZZZ_{p} \times \ZZZ_{p}^{\times}$, with each vertex $(x, y)$ having an arc to $(x+a, ya)$ for each $a = 1, 2, \ldots, p-1$. In other words, $\mathrm{SUMPROD}(p)$ is the Cayley digraph on $\ZZZ_{p} \times \ZZZ_{p}^{\times}$ with connection set $\{(a, a): a \in \ZZZ_{p}^{\times}\}$. Obviously, this digraph has order $p(p-1)$, with each vertex having in-degree and out-degree $p-1$.

Let  $d \ge 2$ an integer. Define \cite{Fri95} $\mathrm{POWER}(p, d)$ to be the digraph with vertex set $\ZZZ_{p}^d$, with each vertex $(x_1, x_2, \ldots, x_d)$ having an arc to each of $(x_1 + a, x_2 + a^2, \ldots, x_d + a^d)$ for $a = 0, 1, \ldots, p-1$. In other words, $\mathrm{POWER}(p, d)$ is the Cayley digraph on $\ZZZ_{p}^d$ with connection set $\{(a, a^2, \ldots, a^d): a \in \ZZZ_{p}\}$. This graph has order $p^d$, with each vertex having in-degree and out-degree $p$.

The following results show that both digraphs above have relatively small second largest eigenvalue in absolute value.

\begin{thm}
\emph{(\cite[Theorem 1.1]{Fri95})}
The graph $\mathrm{SUMPROD}(p)$ has second largest eigenvalue in absolute value at most $\sqrt{p}$.
\end{thm}

\begin{thm}
\emph{(\cite[Theorem 1.2]{Fri95})}
The graph $\mathrm{POWER}(p, d)$ has second largest eigenvalue in absolute value at most $(d-1)\sqrt{p}$.
\end{thm}


\subsection{Two more families of directed Cayley graphs}
\label{subsec:twofam}

Consider a finite field $\FFF_{q}$ of characteristic $p$. Take an irreducible polynomial $f(x)$ of degree $n \ge 2$ over $\FFF_{q}$. Let $\Ga_f = (\FFF_{q}[x]/((f(x)))^* = (\FFF_{q}[\a])^* = \FFF_{q^n}^*$, where $\a = \bar{x}$. Of course $\Ga_f $ is a cyclic group of order $q^n - 1$. For $1 \le d < n$, let $P_d$ be the set of monic primary polynomials of degree $d$ in $\FFF_{q}[x]$, where a polynomial in $\FFF_{q}[x]$ is called \emph{primary} if it is a power of an irreducible polynomial. Set $E_d = \{g(\a): g \in P_{d}\}$. Then $E_d$ is a proper subset of $\Ga_f $. In \cite{LuWWZ14}, Lu, Wan, Wang and Zhang studied the Cayley digraph $G_{d}(n,q,\a) := \Cay(\Ga_f, E_d)$. In the special case when $d=1$, we have $E_1 = \a + \FFF_{q}$ and $G_{1}(n,q,\a)$ is Chung's difference graph \cite{Chung89}. In general, $G_{d}(n,q,\a)$ is a  regular digraph of order $q^n-1$ with degree $|P_d|$ approximately $q^d/d$. Among other things Lu, Wan, Wang and Zhang \cite{LuWWZ14} proved the following result, in which an expander graph is understood as a regular digraph such that the modulus of every nontrivial eigenvalue is ``much" less than the degree in some sense.

\begin{thm}
\emph{(\cite[Theorem 6]{LuWWZ14})}
Let $\d$ be a constant with $0 < \d < 1$ such that $n+d-1 \le q^{d/2}(1-\d)$. Then every nontrivial eigenvalue $\l$ of the adjacency operator for $G_{d}(n,q,\a)$ satisfies
$$
|\l| \le \frac{q^d}{d} (1-\d) \le |P_d| (1-\d).
$$
In particular, $G_{d}(n,q,\a)$ is an expander graph.
\end{thm}

A weighted digraph $G_{d}^*(n,q,\a)$ obtained from $G_{d}(n,q,\a)$ by assigning a specific weight to each of its arcs was also considered in \cite{LuWWZ14}. It was proved in \cite{LuWWZ14} that $G_{d}^*(n,q,\a)$ is also an expander graph.

Let $R$ be a finite local commutative ring with maximal ideal $M$ and residue field $K = R/M$ equipped with the reduction map $\bar{}: R \rightarrow K$. A polynomial in $R[x]$ is called \emph{regular} if its reduction is not zero in the residue field. Let $f(x)$ be a primary regular non-unit polynomial in $R[x]$. That is, its reduction $\bar{f}(x)$ is a power of an irreducible polynomial in $K[x]$. By Hensel's Lemma,
$$
f(x) = \d(x) \pi(x)^s + \b(x),
$$
where $s \ge 1$ is an integer, $\d(x)$ is a unit, $\b(x) \in M[x]$, and $\pi(x)$ is a monic irreducible polynomial in $R[x]$ of degree $n \ge 1$ such that $\bar{\pi}(x)$ is irreducible in $K[x]$. Let $I = \la f(x) \ra$
be the principal ideal generated by $f(x)$. Then $R[x]/I$ is a local ring of order $|R|^{ns}$. Denote by $\Ga_{f}(R)$ the unit group $(R[x]/I)^{\times}$. For each integer $d$ with $1 \le d < n$, let $P_{d}(R)$ be the set of monic primary polynomials of degree $d$ in $R[x]$. Define $G_{d}(R, f)$ to be the Cayley digraph on $\Ga_{f}(R)$ with respect to the set $P_{d}(R)+I$ of cosets of $I$ with representatives in $P_{d}(R)$. That is, $G_{d}(R, f)$ has vertex set $\Ga_{f}(R)$ and there is an arc from $g_1(x) + I$ to $g_2(x) + I$ if and only if $(g_2(x) + I)(g_1(x) + I)^{-1} \in P_{d}(R)+I$. This graph was introduced by Rasri and Meemark in \cite{RasriM17} as a generalization of the digraph $G_{d}(n,q,\a)$ above, the latter being the digraph obtained when $R$ is a finite field. Among other things the following result was proved in \cite{RasriM17}.

\begin{thm}
\label{thm:RasriM17}
\emph{(\cite[Theorem 3.1]{RasriM17})}
Let $G_{d}(R, f)$, $n$ and $s$ be as above.
Let $\d$ be a constant with $0 < \d < 1$ such that
$$
ns - 1+ d|M|^d \le |M|^{d/2} |R|^{d/2}(1-\d).
$$
Then every nontrivial eigenvalue $\l$ of the adjacency matrix of $G_{d}(R, f )$ satisfies
$$
|\l| \le \frac{|R|^d}{d} (1-\d) \le |P_d(R)| (1-\d).
$$
In particular, $G_{d}(R, f)$ is an expander graph.
\end{thm}

It was also proved in \cite[Theorem 3.2]{RasriM17} that the weighted digraph obtained from $G_{d}(R, f )$ by assigning a specific weight to each arc is also an expander graph. In \cite{RasriM17}, the construction of $G_{d}(R, f)$ was generalized to the case when $f$ is any regular polynomial in $R[x]$ and a counterpart of Theorem \ref{thm:RasriM17} under this general setting was given in \cite[Theorem 5.3]{RasriM17}.


\section{Miscellaneous}
\label{sec:miscell}


\subsection{Random Cayley graphs}

The well-known Alon-Roichman theorem \cite{AlonR94} states that for $\ve > 0$, with high probability, $O(\log|\Ga|/\ve^2)$ elements chosen independently and uniformly from a finite group $\Ga$ give rise to a Cayley graph with second eigenvalue no more than $\ve$. In particular, such a graph is an expander with high probability. Landau and Russell \cite{LandauRussell04}, and independently Loh and Schulman \cite{LohSchulman04}, improved the bounds in the theorem. Following Landau and Russell, in \cite{ChristofidesM08} Christofides and Markstr\"{o}m gave a new proof of the result, improving the bounds even further, and proved a generalization of the Alon-Roichman theorem to random coset graphs. See also \cite{Pak99} for an improvement of the Alon-Roichman theorem.

A $k$-regular graph $G$ of order $n$ is called a \emph{$c$-expander} if the quantity
$$
h^*(G) = \min\left\{\frac{|N(S)|}{|S|(n-|S|)}: \emptyset \ne S \subset V(G)\right\}
$$
is at least $c/n$, where $N(S)$ is the set of neighbours of the vertices of $S$ in $G$. (Since the minimum of $|\partial(S)|/|S|(n-|S|)$ for $\emptyset \ne S \subset V(G)$ is between $h(G)/n$ and $2h(G)/n$ and $k |N(S)| \ge |\partial(S)| \ge |N(S) \setminus S|$, where $h(G)$ is the isoperimetric number of $G$ as defined in Section \ref{sec:Ram}, we see that $h^*(G)$ can be bounded from above as well as from below by $h(G)$.) In \cite{AlonR94}, Alon and Roichman proved that for any $\d$, $0 < \d < 1$, there exists a real number $c(\d) > 0$ such that for every group $\Ga$ of order $n$ and for every set $S$ of $c(\d) \log n$ random elements of $\Ga$, the expected value of the second largest eigenvalue in absolute value of the normalized adjacency matrix of $\Cay(\Ga, S \cup S^{-1})$ is at most $1 - \d$. Moreover, the probability that such graph is a $\d$-expander tends to $1$ as $n \rightarrow \infty$.

In \cite{LovettMR15}, Lovett, Moore and Russell proved that there exists a family of groups $\Ga_n$ and nontrivial irreducible representations $\rho_n$ such that, for any constant $t$, the average of $\rho_n$ over $t$ uniformly random elements $g_1, \ldots, g_t \in \Ga_n$ has operator norm 1 with probability approaching 1 as $n \rightarrow \infty$. More explicitly, settling a conjecture of Wigderson, they proved that there exist families of finite groups $\Ga$ for which $\Omega(\log \log |\Ga|)$ random elements are required to bound the norm of a typical representation below 1.

The following concepts were introduced in the study of quasirandomness and expansion of graphs. A $k$-regular graph $G$ with order $n$ is called \emph{$\ve$-uniform} if, for all $S, T \subseteq V(G)$,
$$
\left| e(S, T) - \frac{k}{n}|S||T|\right| \le \ve k n,
$$
where $e(S, T)$ is the number of edges of $G$ between $S$ and $T$. The graph $G$ is called an \emph{$(n, k, \l)$-graph} if all eigenvalues of $G$ except $k$ are bounded from above in absolute value by $\l$. The well-known expander mixing lemma \cite{AlonC88} (see also \cite[Lemma 2.5]{Hoory06}) says that, if $G$ is an $(n, k, \l)$-graph, then
$$
\left| e(S, T) - \frac{k}{n}|S||T|\right| \le \l \sqrt{|S| |T|}
$$
for all $S, T \subseteq V(G)$, where the left-hand side is called the discrepancy. Thus, if the second largest eigenvalue $\l_{2}(G)$ of a $k$-regular graph $G$ satisfies $|\l_{2}(G)| \le \ve k$, then $G$ is $\ve$-uniform. In \cite{KohayakawaRS16}, Kohayakawa, R\"odl and Schacht proved the following result.

\begin{thm}
\label{thm:KRS}
\emph{(\cite[Theorem 1.6]{KohayakawaRS16})}
Let $\Ga$ be an abelian group. Then every $\ve$-uniform Cayley graph $\Cay(\Ga, S)$ is an $(n, k, \l)$-graph with $n = |\Ga|$, $k = |S|$ and $\l \le C \ve k$ for some absolute constant $C$.
\end{thm}

In other words, having small discrepancy and having large eigenvalue gap are equivalent properties for Cayley graphs on abelian groups, even if they are sparse. (Roughly, the discrepancy of a graph measures how uniformly the edges are distributed among the vertices.) This answers a question of Chung and Graham \cite{ChungG02} for the particular case of Cayley graphs on abelian groups, while in general the answer is negative.

In \cite{ConlonZ17}, Conlon and Zhao generalized Theorem \ref{thm:KRS} to any finite group.

\begin{thm}
\emph{(\cite[Theorem 1.4]{ConlonZ17})}
Let $\Ga$ be a finite group. Then every $\ve$-uniform Cayley graph $\Cay(\Ga, S)$ is an $(n, k, \l)$-graph with $n = |\Ga|$, $k = |S|$ and $\l \le 8 \ve k$.
\end{thm}

As a corollary, Conlon and Zhao \cite{ConlonZ17} also proved that the same result holds for all vertex-transitive graphs. That is, every $k$-regular $\ve$-uniform vertex-transitive graph of order $n$ is an $(n, k, \l)$-graph with $\l \le 8 \ve k$ (see \cite[Corollary 1.5]{ConlonZ17}).

In a study of quasirandom groups and product-free subsets of groups, Gowers \cite[Theorem 4.5]{Gowers08} proved among other appealing results that a finite group $\Ga$ has the property that the Cayley digraph with respect to every subset of $\Ga$ is quasirandom is equivalent to a few other properties of $\Ga$, including the property of having no nontrivial low-dimensional representations. This result is relevant to the present survey as the second largest eigenvalue, in absolute value, of a regular graph is responsible for quasirandomness of the graph. See \cite{Gowers08} for details. The reader is also referred to \cite{KrivelevichS06} for a survey on quasirandom graphs and pseudo-random graphs, where in particular several families of Cayley graphs are discussed.


\subsection{Distance eigenvalues of Cayley graphs}

The \emph{distance matrix} of a graph of order $n$ is the $n \times n$ matrix whose $(u, v)$-entry is equal to the distance between vertices $u$ and $v$ in the graph. The eigenvalues of this matrix are called the \emph{distance eigenvalues} of the graph. A graph whose distance eigenvalues are all integers is called \emph{distance integral}. The \emph{distance energy} of a graph is defined as the sum of the absolute values of its distance eigenvalues.

In \cite{Ilic10}, Ili\'{c} gave a characterization of the distance eigenvalues of integral circulant graphs $\ICG(n, D)$ and proved that these graphs have integral distance eigenvalues. In the same paper he also computed the distance eigenvalues and distance energy of the unitary Cayley graphs $\Cay(\ZZZ_n, \ZZZ_n^{\times})$.

The distance eigenvalues of the Hamming graph $H(d, q)$ were determined by Indulal in \cite{Indulal09} and Cioab\u{a}, Elzinga, Markiewitz, Vander Meulen and Vanderwoerd in \cite{CioabaEMVV17}, where in the latter paper these eigenvalues were used to solve a graph addressing problem \cite{GrahamP71} for Hamming graphs.

Let $(W, S)$ be a finite Coxeter system, and let $T = \{wsw^{-1}: w \in W, s \in S\}$ be the set of all reflections of $W$. The Cayley graphs $\Cay(W, S)$ and $\Cay(W, T)$ are called the \emph{weak order graph} and
\emph{absolute order graph} of $(W, S)$, respectively. It is known that every $w \in W$ can be written as a word in the simple reflections $S$ (respectively, reflections $T$), and the minimum number of such reflections that must be used is the length $\ell_S(w)$ (respectively, absolute length $\ell_T(w)$) of $w$. It is also known that the graph distance between two vertices $u, v \in W$ in these two graphs is equal to $\ell_S(uv^{-1})$ and $\ell_T(uv^{-1})$, respectively. In other words, the distance matrices of $\Cay(W, S)$ and $\Cay(W, T)$ are $(\ell_{S}(uv^{-1}))_{u, v \in \Ga}$ and $(\ell_{T}(uv^{-1}))_{u, v \in \Ga}$, respectively. The following result was obtained by Renteln \cite{Renteln11}.

\begin{thm}
\emph{(\cite[Theorems 1 and 6]{Renteln11})}
\label{thm:Rn11}
Let $(W, S)$ be a finite Coxeter system and $T$ the set of all reflections of $W$. Then the distance eigenvalues of the absolute order graph $\Cay(W, T)$ are given by
$$
\frac{1}{\chi(1)} \sum_{K} |K| \ell_{T}(w_K) \chi(w_K), \text{ with multiplicity } \chi(1)^2,
$$
where $K$ runs over all conjugacy classes of $W$, $w_K$ is any element of $K$, and $\chi$ ranges over all irreducible characters of $W$. Moreover, $\Cay(W, T)$ is distance integral.
\end{thm}

Using this result, Renteln computed explicitly the distance spectra of some absolute order graphs (see \cite[Section 3.3]{Renteln11}) and in particular proved the following result for the dihedral group $D_{2n}$ of order $2n$.

\begin{thm}
\emph{(\cite[Theorem 9]{Renteln11})}
\label{thm:Rn11a}
 The characteristic polynomial of the distance matrix of the absolute order graph $\Cay(D_{2n}, T)$ of $D_{2n}$ is given by
$$
(x-3n+2)(x-n+2)(x+2)^{2n-2}.
$$
\end{thm}

The weak order graph of $D_{2n}$ is simply the cycle of length $2n$ and so its distance eigenvalues can be easily computed (see \cite[Theorem 16]{Renteln11}). However, it was noted in \cite{Renteln11} that in general the computation of the distance spectrum of the weak order graphs is more complicated and that there appear to be far fewer nontrivial eigenvalues in the weak order case than the absolute order case. Computational results in \cite{Renteln11} indicate that in most cases the majority of the distance eigenvalues of weak order graphs are zero, and in some cases the distance spectra are integral.  In \cite{Renteln11}, the distance eigenvalues of the weak order graphs of Coxeter groups of a few types were given, a variant of the distance matrix of $\Cay(W, S)$ was introduced, the eigenvalues of the former were related to that of the latter, and a few open questions were posed.

A real (respectively, complex) \emph{reflection group} is a finite group generated by reflections of a Euclidean (respectively, unitary) vector space. A reflection group is \emph{irreducible} if that reflection representation is irreducible.

Let $W$ be a reflection group on $V$. For $w \in V$, define $\codim(w)$ to be the codimension of the fixed point space $V^w = \{v \in W: wv = v\}$. If $W$ is a real reflection group, then $\ell_T(w) = \codim(w)$ for all $w \in W$, but this is not true in general for complex reflection groups. The \emph{codimension matrix} of $W$ is defined to be the matrix $(\codim(uv^{-1}))_{u, v \in W}$. The spetrum of this matrix is called the \emph{codimension spectrum} of $W$.

Let $\d_T$ denote the characteristic function of $T$. Since $\d_T, \ell_T$ and $\codim$ are all class functions for $W$, Theorem \ref{thm:FG16} can be used to prove the following result.

\begin{thm}
\emph{(\cite[Section 4.3]{Foster-GreenwoodK16})}
Let $W$ be a finite complex reflection group and $T$ the set of all reflections of $W$. Then the eigenvalues of $\Cay(W, T)$, the distance eigenvalues of $\Cay(W, T)$ and the eigenvalues of the codimension matrix of $W$ are given by
$$
\theta_{\chi}(f) = \frac{1}{\chi(1)} \sum_{w \in W} f(w) \chi(w), \text{ with multiplicity } \chi(1)^2,
$$
with $\chi$ ranging over the irredducible characters of $W$, for $f = \d_T, \ell_T, \codim$, respectively.

Furthermore, the largest eigenvalue of $\Cay(W, T)$, the largest distance eigenvalue of $\Cay(W, T)$ and the largest eigenvalue of the codimension matrix of $W$ are $\theta_{1}(f)$ corresponding to the trivial character, for $f = \d_T, \ell_T, \codim$, respectively.
\end{thm}

Note that for $f = \d_T$ the formula above is the same as the one in Theorem \ref{thm:Rn11}.

\begin{thm}
\label{thm:CayWT}
\emph{(\cite[Corollary 4.11 and Theorem 5.4]{Foster-GreenwoodK16})}
Let $W$ be a finite irreducible complex reflection group and $T$ the set of all reflections of $W$. Then the Cayley graph $\Cay(W, T)$ is integral and distance integral, and the codimension spectrum of $W$ is integral as well.
\end{thm}

Note that the distance integrality of $\Cay(W, T)$ in Theorem \ref{thm:CayWT} extends Theorem \ref{thm:Rn11a} from irreducible real reflection groups to irreducible complex reflection groups.

Finally, in \cite[Theorem 4.3]{HuangL21}, Huang and Li obtained a necessary and sufficient condition for a Cayley graph $\Cay(\Ga, S)$ on a generalized dihedral group $\Ga = \Dih(\Si)$ to be distance integral, and in \cite[Theorem 5.1]{HuangL21} they proved that in the special case when $|S \cap \Si| = 1$, $\Cay(\Ga, S)$ is integral if and only if it is distance integral. A few special cases where $\Cay(\Ga, S)$ is distance integral were also identified in this paper.


\subsection{Others}

Fullerenes are of great importance for chemistry as evidenced by the well-known Buckministerfullerene $C_{60}$. A \emph{fullerene} can be represented by a $3$-regular graph on a closed surface with pentagonal and hexagonal faces such that its vertices are carbon atoms of the molecule and two vertices are adjacent if there is a bond between the corresponding atoms. Fullerenes exist in the sphere, torus, projective plane, and the Klein bottle. A fullerene is called \emph{toroidal} if it lies on the torus. In \cite{Kang11}, Kang classified all fullerenes which are Cayley graphs and determined their eigenvalues.

The \emph{power graph} of a group $\Ga$ is the graph with vertex set $\Ga$ in which two distinct elements $x, y$ are adjacent if and only if $x^m = y$ or $y^m = x$ for some positive integer $m$. In \cite{ChattopadhyayP15}, Chattopadhyay and Panigrahi studied relations between the power graph and the unitary Cayley graph $\Cay(\ZZZ_n, \ZZZ_n^{\times})$ of the cyclic group $\ZZZ_n$.

Given a graph $G$ of order $n$, say, with vertex set $\{1, 2, \ldots, n\}$, define $S(G) = \{e_i + e_j: ij \in E(G)\} \subseteq \ZZZ_2^n$, where $\{e_1, e_2, \ldots, e_n\}$ is the standard basis of $\ZZZ_2^n$. In \cite{QinXM09}, Qin, Xiao and Miklavi\v{c} observed that $G$ can be isometrically embedded as a subgraph of the cubelike graph $\Cay(\ZZZ_2^n, S(G))$. Among other things they found some relations between the spectrum of $G$ and that of $\Cay(\ZZZ_2^n, S(G))$.

A Cayley graph $G = \Cay(\Ga, S)$ is called \emph{normal edge-transitive} if the normalizer of $\Ga$ in $\Aut(G)$ is transitive on the set of edges of $G$. In \cite{Ghorbani15}, the eigenvalues of normal edge-transitive Cayley graphs on $D_{2n}$ and $T_{4n}$ were given by Ghorbani, where $D_{2n}$ is the dihedral group of order $2n$ and $T_{4n} = \la a, b\ | \ a^{2n} = 1, b^2 = a^n, bab^{-1} = a^{-1} \ra$.

Denote by $\Ga_n$ the multiplicative group of the upper unitriangular $2 \times 2$ matrices over $\ZZZ~\mod~n$. In \cite{GuptaSL16}, Gupta, Shwetha Shetty and Lokesha studied $\Cay(\Ga_n, S)$ for a certain set $S$ of two generators of $\Ga_n$. Among their findings are the spectra and energies of the adjacency, Laplacian, normalized Laplacian and signless Laplacian matrices of $\Cay(\Ga_n, S)$ when $n$ is odd.

It is known that for any vertex-transitive graph $G$ and any eigenvalue $\l$ of $G$ the multiplicity of $\l$ is decreased by one as an eigenvalue of $G-v$ for any $v \in V(G)$. In \cite{Savchenko05}, Savchenko proved that the same statement holds for vertex-transitive digraphs.

A graph is called \emph{singular} if its adjacency matrix is singular or, equivalently, it has $0$ as an eigenvalue. Obviously, by Theorem \ref{Conjeign}, a connected normal Cayley graph $\Cay(\Ga, S)$ is singular if and only if there is an irreducible character $\chi$ of $\Ga$ such that $\sum_{g \in \Ga} \chi(g) = 0$. Starting from this observation, Siemons and Zalesski obtained in \cite{SiemonsZ19} several results on the singularity of some connected normal Cayley graphs on non-abelian simple groups including alternating groups.

The \emph{$G$-graph} with respect to a non-empty subset $S$ of a group $G$ is the multigraph with vertex set the set of $\la s \ra$-orbits on $G$ under the left regular multiplication of $G$, for $s$ running over $S$, such that there are $|\la s_1 \ra x_1 \cap \la s_2 \ra x_2|$ edges between vertices $\la s_1 \ra x_1$ and $\la s_2 \ra  x_2$. In \cite{BadaouiBM19}, a connection between Cayley graphs and $G$-graphs was established and using this connection the eigenvalues of two specific families of Cayley graphs on the dicyclic and dihedral groups, respectively, were computed explicitly.

Recall that the generalized Paley graph $X_{q}^k$ was defined in Section \ref{subsec:GPaley}. In \cite{JohnsonSS16}, it was shown that if $p$ is an odd prime and $k \ge 2$ is an integer such that $d = \gcd(p-1, k)$ divides $(p-1)/2$ then the isoperimetric number $h(X_{p}^k)$ of $X_{p}^k$ satisfies $(p+(1-d)\sqrt{p})/2d \le h(X_{p}^k) \le (4 \sum_{i=1}^l \g_i)/(p-1)$, where $\{a^k: a \in \FFF_p^*\} = \{\pm \g_1, \ldots, \pm \g_l\}$ with $l = (p-1)/2d$ and $0 \le \g_i \le (p-1)/2$ for each $i$.

A graph $G$ is called a \emph{$t$-existentially closed} (\emph{$t$-e.c.}) \emph{graph}, where $t \ge 1$ is an integer, if for any disjoint (but not necessarily non-empty) subsets $A, B$ of $V(G)$ satisfying $|A \cup B| = t$, there exists a vertex $z \in V(G) \setminus (A \cup B)$ adjacent to all vertices of $A$ but none of $B$. It is known that random graphs are $t$-e.c. for any $t \ge 1$. Another property shared by random graphs is the one of being best pseudo-random. In the case of regular graphs, the expander mixing lemma implies that this property can be defined using eigenvalues: A regular graph of order $n$ and degree $d(n)$ is \emph{best pseudo-random} if $\l(G) = O(\sqrt{d(n)})$ as $n \rightarrow \infty$, where $\l(G)$ is the largest eigenvalue of $G$ other than $d(n)$ in absolute value. Many known $t$-e.c. graphs are best pseudo-random, and it was not clear whether there are $t$-e.c. graphs which are not best pseudo-random. In \cite{Satake19}, Satake answered this question affirmatively using a special family of quadratic unitary Cayley graphs (see Section \ref{subsec:QUCayComRing}). He proved that for any $t \ge 1$, any prime power $q \equiv 1~(\mod~4)$ and any odd integer $e \ge 1$, if $q^e - (t 2^{t-1} - 2^t + 1)q^{e-\frac{1}{2}} - t 2^t q^{e-1} + t 2^{t-1} > 0$, then the quadratic unitary Cayley graph $\GG_{\ZZZ_{q^e}}$ is a $t$-e.c. graph. This together with Theorem \ref{G2pspec} gives an infinite family of $t$-e.c. graphs which are not best pseudo-random for any $t \ge 1$.

A \emph{mixed graph} is a graph in which both edges and arcs may present. The \emph{mixed adjacency matrix} $M(G)$ of a mixed graph $G$ of order $n$ is the matrix whose $(u, v)$-entry is $1$ if there is an edge between $u$ and $v$ or an arc from $u$ to $v$, $-1$ if there is an arc from $v$ to $u$, and $0$ otherwise. In \cite{AdigaRS16}, Adiga, Rakshith and So studied various spectral properties of $M(G)$ and in particular obtained bounds on the mixed energy of $G$. In the same paper they introduced the \emph{mixed unitary Cayley graph} $M_n$: The vertex set of $M_n$ is $\{0, 1, \ldots, n-1\}$; for $0 \le i < j \le n-1$ with $\gcd(j-i, n) = 1$, there is an edge between $i$ and $j$ if $\legendre{j-i}{n} = 1$, there is an arc from $i$ to $j$ if $\legendre{j-i}{n} = -1$ and $j-i < \lf n/2 \rf$, and there is an arc from $j$ to $i$ if $\legendre{j-i}{n} = -1$ and $j-i > \lf n/2 \rf$, where $\legendre{a}{n}$ denotes the Jacobi symbol. With the help of Ramanujan sums, they determined the eigenvalues of $M_n$ in \cite[Theorem 4.15]{AdigaRS16} (see also \cite{AdigaR16}) in the case when $n$ has an even number of prime factors congruent to $3$ modulo $4$. The energies of some mixed unitary Cayley graphs were computed by Adiga and Rakshith in \cite{AdigaR16}.

In \cite{DalfoF20}, Dalf\'{o} and Fiol gave a method for constructing equitable partitions of Cayley digraphs (which can be undirected or mixed graphs) of permutation groups on $n$ letters. They proved that every partition of the integer $n$ gives rise to an equitable partition of the Cayley digraph, and using representation theory they also gave a method for finding all eigenvalues and eigenspaces of the corresponding quotient digraphs.

Let $G = (V, E)$ be a digraph with vertex set $V$ and arc set $E$, where loops and multiple arcs are allowed. Let $\Ga$ be a finite group with a generating set $S$. A \emph{voltage assignment} of $G$ is a mapping $\a: E \rightarrow S$, and the associated \emph{lifted digraph} (or \emph{lift} for short) $G^{\a}$ is the digraph with vertex set $V(G^{\a}) = V \times G$, where there is an arc from vertex $(u, x)$ to vertex $(v, y)$ if and only if $(u, v) \in E$ and $y = x\a (u,v)$. Let $\mathbf{B}$ be the square matrix over $\CCC\Ga$ indexed by the vertices of $G$ with $(u,v)$-entry $\a (u,v)$ if $(u,v) \in E$ and $0$ otherwise. Denote by $r$ the number of vertices of $G$ and $n$ the order of $\Ga$. Let $\{\rho_1, \rho_2, \ldots, \rho_k\}$ be a complete set of irreducible matrix representations of $\Ga$, and let $d_i$ be the dimension of $\rho_i$ for each $i$. Let $\rho_{i}(\mathbf{B})$ be the complex matrix obtained by replacing each element $g \in \Ga$ in the entries of $\mathbf{B}$ by the $d_{i} \times d_{i}$ matrix $\rho_{i}(g)$, and let $\mu_{u, j}$, $u \in V$, $j \in \{1, 2, \ldots, d_i\}$ be the eigenvalues of $\rho_{i}(\mathbf{B})$. In \cite[Theorem 2.1]{DalfoFS19}, Dalf\'{o}, Fiol and \v{S}ir\'{a}\v{n} proved that the $rn$ eigenvalues of the lift $G^{\a}$ are the $rd_{i}$ eigenvalues of $\rho_{i}(\mathbf{B})$, for $i \in \{1, 2, \ldots, k\}$, each repeated $d_i$ times. The same result was also obtained independently by Li \cite{Li21} using a slightly different language. In the special case when $G$ has only one vertex with $|S|$ loops, each with voltage an element of $S$, this result yields Theorem \ref{BabaiS2} because in this case the lift $G^{\a}$ is exactly the Cayley digraph $\Cay(\Ga, S)$. Related results on eigenvalues and eigenspaces of arbitrary lifts of digraphs can be found in \cite{DalfoFPS21}.

A \emph{signed graph} \cite{BelardoCKW18} is a graph in which each edge is labelled as positive or negative. A signed graph is called \emph{balanced} if the number of negative edges in each cycle is even. The \emph{adjacency matrix} of a signed graph $G$ is obtained from that of the underlying graph of $G$ by changing the entries corresponding to negative edges from 1 to $-1$, and the \emph{energy} $\EE(G)$ of $G$ is the sum of the absolute values of the eigenvalues of the adjacency matrix of $G$. A signed graph $G$ of order $n$ is said to be \emph{hyperenergetic} if $\EE(G) > 2n-2$. If $G$ is a signed graph, define its \emph{line signed graph} $L(G)$ to be the signed graph whose underlying graph is that of the line graph of the underlying graph of $G$, such that the sign of an edge $\{e, f\}$ of $L(G)$ is negative if and only if both $e$ and $f$ are negative edges of $G$. Let $R$ be a finite commutative ring. In \cite{MeemarkS14a}, Meemark and Suntornpoch defined the \emph{unitary Cayley signed graph} $G(R)$ of $R$ to be the signed graph whose underlying graph is the unitary graph $\Cay(R, R^{\times})$ such that an edge $xy$ (where $x, y \in R, x-y \in R^{\times}$) is negative if and only if $\{x, y\} \cap R^{\times} = \emptyset$. In the same paper they also determined when $G(R)$ is balanced, when $L(G(R))$ is balanced, when $G(R)$ is hyperenergetic and balanced, and when $L(G(R))$ is hyperenergetic and balanced. The reader is referred to \cite{BelardoCKW18} for a survey of known results and open problems on the adjacency spectra of signed graphs.

The \emph{resistance distance} $r_{uv}$ between two vertices $u, v$ in a connected graph $G$ is the effective electrical resistance between them when unit resistors are placed on every edge of $G$. The \emph{Kirchhoff index} of $G$ is the quantity $Kf(G) = \frac{1}{2} \sum_{u \in V(G)} \sum_{v \in V(G)} r_{uv}$. It is known that $Kf(G)$ is equal to the order of $G$ times the sum of the reciprocals of the nonzero Laplacian eigenvalues of $G$. In \cite{GaoLL11}, Gao, Luo and Liu obtained closed-form formulas for the Kirchhoff index and resistance distances of Cayley graphs on finite abelian groups in terms of the Laplacian eigenvalues and eigenvectors, respectively. In particular, they gave formulas for the Kirchhoff index of the hexagonal torus network, the multidimensional torus and the hypercubes, as well as formulas for the Kirchhoff index and resistance distances of complete multipartite graphs.

In \cite{GrigorchukN12}, Grigorchuk and Nowak explored relations between the diameter, the first positive eigenvalue of the discrete $p$-Laplacian, and the $\ell_p$-distortion of a finite graph. They proved an inequality connecting these three quantities, and apply it to families of Cayley and Schreier graphs. They also showed that the $\ell_p$-distortion of Pascal graphs is bounded, which allows one to obtain estimates for the convergence to zero of the spectral gap.

In \cite{MineiSkogman09}, Minei and Skogman presented a block diagonalization method for the adjacency matrices of two types of covering graphs, using the irreducible representations of the Galois group of the covering graph over the base graph. The first type of covering graph is the Cayley graph over the finite ring $\ZZZ_{p^n}$, and the second one resembles lattices with vertices $\ZZZ_n \times \ZZZ_n$ for large $n$. Using the block diagonalization method, they obtained explicit formulas for the eigenvalues of one lattice and nontrivial bounds on the eigenvalues of another lattice.

In \cite{LaffertyRockmore99}, Lafferty and Rockmore described numerical computation of eigenvalue spacings of 4-regular Cayley graphs on cyclic and symmetric groups, and of two-dimensional special linear groups over prime fields. Denote by $P(s)$ the number of indices where the gap between two consecutive distinct eigenvalues of a Cayley graph is at most $s$. It was observed \cite{LaffertyRockmore99} that for the above-mentioned groups, after a linear transformation the function $P(s)$ approximates $1-e^{-s}$. In contrast, in random 4-regular graphs a linear transformation of $P(s)$ seems to approximate $1-e^{-\pi s^2/4}$.

Consider a continuous-time quantum walk on a graph $G$ with transition matrix $H_{G}(t)$ as defined in \eqref{eq:HGt}. The probability that at time $\tau$ the quantum walk with initial state $u$ is in state $v$ is given by $|H_{G}(\tau)_{u,v}|^2$. We say that $G$ admits \emph{uniform mixing} at time $\tau$ if this probability is the same for all vertices $u$ and $v$. As noticed in \cite{GodsilZ2017}, uniform mixing on graphs is rare, and most known examples are integral Cayley graphs on abelian groups. In \cite{GodsilZ2017}, Godsil and Zhan constructed among other things infinite families of Cayley graphs on $\ZZZ_n^d$ admitting uniform mixing. Several other results about uniform mixing on Cayley graphs can also be found in \cite{GodsilZ2017} and the references therein.

In \cite{ChenG20}, Chen and Godsil introduced the notion of perfect edge state transfer in a graph. A necessary and sufficient condition for the existence of perfect edge state transfer in connected normal Cayley graphs on dihedral groups was obtained in \cite{LuoCXC21}, and several constructions of such Cayley graphs were also given in this paper.

Extending the results in \cite{TanFC17} to the weighted case, Cao, Feng and Tan \cite{CaoFT21} gave a necessary and sufficient condition for a weighted Cayley graph on an abelian group to admit perfect state transfer, and proved that any integral weighted Cayley graph on an abelian group is periodic at every vertex. See \cite[Theorems 2.1 and 2.2]{CaoFT21} for details.

Given integers $k\geq 2$, $N\geq0$ and $M\geq1$, the \emph{spider-web graph} $\SS_{k,N,M}$ is the graph with vertex set $\{0, 1, \ldots, k-1\}^N \times \ZZZ_M$ such that each vertex $((x_1,\ldots,x_N),i)$ is adjacent to $((x_2,\ldots,x_N,y),i+1)$, for $y \in \{0, 1, \ldots, k-1\}$. In \cite{GrigorchukLN16}, Grigorchuk, Leemann and Nagnibeda proved that, as $N, M \rightarrow \infty$, $\SS_{k,N,M}$ converges in some sense to the Cayley graph on the lamplighter group $\ZZZ_k \wr \ZZZ$. In the proof of this result, they realised $\SS_{k,N,M}$ as the tensor product of the de Bruijn graphs $B_{k,N}$ with cycle $C_M$. Since the spectra of the de Bruijn graphs are known, this enabled them to compute the spectra of spider-web graphs by relating to the Laplacian of the Cayley graph on the lamplighter group.

Let $n$ and $t$ be integers with $n \ge 2$ and $1 \le t \le n$. As before, let $S_n$ be the symmetric group over $[n] = \{1,2,\ldots,n\}$. A subset $A \subset S_n$ is called \emph{$t$-set-intersecting} if for any $\sigma, \pi \in A$ there exists some $t$-set $T \subset [n]$ such that $\sigma(T)=\pi(T)$. Generalizing a result of Frankl and Deza (for $t=1$) and settling a conjecture of K\"{o}rner (for $t=2$), Ellis \cite{Ellis12} proved that, if $n$ is sufficiently large depending on $t$, and $A \subset S_n$ is $t$-set-intersecting, then $|A| \leq t!(n-t)!$ and equality holds only if $A$ is a coset of the stabilizer of a $t$-set. In the proof of this result the Cayley graph $G_{(t)}$ on $S_n$ with connection set generated by the set $D_{(t)}$ of $t$-derangements in $S_n$ played a key role, where an element of $S_n$ is called a \emph{$t$-derangement} if it fixes no $t$-set of $[n]$ setwise. This graph is called the \emph{$t$-derangement graph} on $[n]$ and can be defined equivalently as the graph with vertex set $S_n$ in which $\sigma,\pi \in S_n$ are adjacent if and only if $\sigma(T) \ne \pi(T)$ for any $t$-set $T \subset [n]$. In particular, $G_{(1)}$ is the derangement graph on $[n]$ as seen in Section \ref{subsec:sn}. A $t$-set-intersecting family of permutations in $S_n$ is precisely an independent set in $G_{(t)}$. Thus the task is to bound the independence number of $G_{(t)}$ and analyze maximum independent sets in $G_{(t)}$. A key step towards this goal is to construct a `pseudo-adjacency matrix' $A$ for $G_{(t)}$, which is a suitable real linear combination of the adjacency matrices of certain subgraphs of $G_{(t)}$; the subgraphs used are normal Cayley graphs. This enabled Ellis  \cite{Ellis12} to apply a weighted version of Hoffman's theorem to obtain the desired independence number using the spectra of matrix $A$.

Finally, let us mention three families of Cayley graphs from the domain of interconnection networks. Let $d \ge 1$ be an integer. The \emph{cube-connected-cycle graph} $CCC(d)$ is obtained from the hypercube $H(d,2)$ by replacing each vertex of $H(d,2)$ by a cycle of length $d$ and can be defined as a certain Cayley graph on the group $\ZZZ_2^n \rtimes \ZZZ_n$ (see \cite{Hey97, MZ17} for details). The \emph{shuffle-exchange graph} $SE(d)$ is defined to have vertex set $\ZZZ_2^d$ such that $(a_1, a_2, \ldots, a_d)$ and $(b_1, b_2, \ldots, b_d)$ are adjacent if and only if either they differ at only the last coordinate, or $(b_1, b_2, \ldots, b_d) = (a_2, a_3, \ldots, a_d, a_1)$, or $(b_1, b_2, \ldots, b_d) = (a_d, a_1, \ldots, a_{d-1})$. Both $CCC(d)$ and $SE(d)$ are popular networks in parallel computing \cite{Hey97}. The \emph{spectral set} of a graph is the set of its eigenvalues with multiplicities ignored. In \cite{RSW12}, Riess, Strehl and Wanka determined the spectral sets of $CCC(d)$ and $SE(d)$ for $d \ge 3$. It was noted that for odd integers $d \ge 3$ the spectral sets of these two graphs are identical, and for even integers $d \ge 4$ the spectral set of $SE(d)$ is a proper subset of that of $CCC(d)$.

Let $n \ge 2$ be an even integer and $\De$ an integer with $1 \le \De \le \lfloor \log_2 n \rfloor$. The \emph{Kn\"{o}del graph} $W_{\De, n}$ is the graph with vertex set $\{(i, j): i=1,2,\ 0 \le j \le \frac{n}{2} - 1\}$ and edges joining $(1, j)$ and $(2, j+2^k - 1~\mod~\frac{n}{2})$ for $0 \le j \le \frac{n}{2} - 1$ and $0 \le k \le \De - 1$. It is not difficult to see that $W_{\De, n}$ is the Cayley graph on the semidirect product $\ZZZ_{n/2} \rtimes \ZZZ_2$ (with the underlying action given by $y^x = (-1)^x y$, $x \in \ZZZ_2$, $y \in \ZZZ_{n/2}$) with respect to the connection set $\{(1, 2^k - 1): 0 \le k \le \De - 1\}$. Kn\"{o}del graphs have been studied extensively as a topological structure for interconnection networks. Computational issues pertaining to the eigenvalues of Kn\"{o}del graphs were considered by Harutyunyan and Morosan in \cite{HM06}. Recently, Balakrishnan, Paulraja, So and Vinay \cite{BalakrishnanPSV18} noted that the eigenvalues of $W_{\De, n}$ are $\pm \big|\sum_{k=0}^{\De - 1} \om_{n}^{2^{k+1} t}\big|$, $0 \le t \le \frac{n}{2} - 1$, and using this fact they determined the eigenvalues of $W_{\De, 2^{\De}}$  together with their multiplicities, for all $\De \ge 4$.


\section{Open problems and research topics}
\label{sec:open}

In this section we present a list of problems that we are aware of or come to our mind at the time of writing. Needless to say, this list is not meant to be exhaustive. Some problems on the list are actually broad research topics, while other problems are more specific and focused.

\begin{prob}
\label{pb-SachsM}
Prove or disprove Conjecture 1 in \cite{SachsM82} due to Sachs and Stiebitz: For any integer $n = 2^e m$, where $m$ is odd, there exists a vertex-transitive graph (in which loops and parallel edges are allowed) with order $n$ and $2^e$ simple eigenvalues if and only if $m \ge 2^{e-1}$.
\end{prob}

In the special case when $m=m_1 \cdots m_e$ with each $m_i \ge 3$ an odd integer or $e=2$ and $m \ge 3$, the sufficiency in Problem \ref{pb-SachsM} has been proved by Sachs and Stiebitz themselves in \cite[Theorems 8--9]{SachsM82}.

\begin{prob}
Prove or disprove Conjecture 7.3 in \cite{So06} posed by So: There are exactly $2^{\tau(n)-1}$ integral circulant graphs on $n$ vertices, where $\tau(n)$ is the number of divisors of $n$.
\end{prob}

Partial results on So's conjecture can be found in Section \ref{subsec:cosp-isom}.

Theorem \ref{thm:Alperin12} gives a characterization of integral Cayley graphs on abelian groups. It would be natural to ask whether a characterization of integral Cayley graphs on nilpotent groups could be achieved.

\begin{prob}
Characterize integral Cayley graphs on finite nilpotent groups.
\end{prob}

This problem is likely to be challenging or even inaccessible. Nevertheless, one may attempt it for some specific families of non-abelian nilpotent groups.

A special Cayley graph on $U_{6n} = \langle a,b~|~ a^{2n}=b^3=1, a^{-1}ba=b^{-1}\rangle$ was shown to be integral in Theorem \ref{thm:U6n}, but the integrality of other Cayley graphs on this group is unknown. So we have the following problem for the sake of completeness.

\begin{prob}
Characterize integral Cayley graphs on $U_{6n}$ for any $n \ge 1$.
\end{prob}

Theorem \ref{nomalSn} says that all normal Cayley graphs on the symmetric group $S_n$ are integral. So the searching for integral Cayley graphs on symmetric groups is reduced to the non-normal case. Under what conditions is a non-normal Cayley graph on $S_{n}$ integral? It may be challenging to give a uniform answer to this question for all non-normal Cayley graphs on $S_{n}$. As such one may attempt some special non-normal Cayley graphs on $S_{n}$. Recall from Section \ref{subsec:cox} that, for a set $T$ of transpositions in $S_n$, $G_T$ is the graph with vertex set $[n]$ such that $i, j \in [n]$ are adjacent if and only if $(i, j) \in T$. It is well known that $\Cay(S_n, T)$ is connected if and only if $G_T$ is connected. In the case when $G_T$ is the star $K_{1, n-1}$, $\Cay(S_n, T)$ is the star graph which is integral by Theorems \ref{thm:star} and \ref{thm:star-int}. In the case when $G_T$ is the complete graph $K_{n}$, $\Cay(S_n, T)$ is the complete transposition graph which is integral by Theorem \ref{thm:KalpakisYesha97}. The next two problems are motivated by these interesting special cases.

\begin{prob}
Let $T$ be a set of transpositions in $S_n$ such that $G_T$ is a tree, where $n \ge 2$. Give a necessary and sufficient condition for $\Cay(S_n, T)$ to be integral.
\end{prob}

\begin{prob}
Let $T$ be a set of transpositions in $S_n$ such that $G_T$ is connected, where $n \ge 2$. Give a necessary and sufficient condition for $\Cay(S_n, T)$ to be integral.
\end{prob}

The next a few problems are about cospectral Cayley graphs, Cay-DS graphs and Cay-DS groups. Recall that a Cayley graph on a group $\Ga$ is said to be Cay-DS if any Cayley graph on $\Ga$ having the same spectrum as it must be isomorphic to it, and a group is called Cay-DS if every Cayley graph on it is Cay-DS. As seen in Theorem \ref{thm:AbdollahiJG17}, for any $k \ge 6$, there exist arbitrarily large families of $k$-regular Cayley graphs that are cospectral and pairwise non-isomorphic. To the best of our knowledge, it is unknown whether the same result holds for $k = 3, 4, 5$.

\begin{prob}
For $k = 3, 4, 5$, are there arbitrarily large families of $k$-regular Cayley graphs that are cospectral and pairwise non-isomorphic?
\end{prob}

\begin{prob}
\label{pb:non-isom-cosp}
Study the existence of non-isomorphic cospectral Cayley graphs on a given group and determine all pairs of such graphs if they exist.
\end{prob}

Obviously, this problem is quite general and so a universal solution to it is unlikely to exist. One may thus focus on some special families of groups such as the dihedral group $D_{2n}$, the dicyclic group $\mathrm{Dic}_n$ of order $4n$, the above-mentioned group $U_{6n}$, etc. So we have the following subproblem of Problem \ref{pb:non-isom-cosp}.

\begin{prob}
Determine all pairs of non-isomorphic cospectral Cayley graphs on $\Ga$, where $\Ga$ is $D_{2n}$,  $\mathrm{Dic}_n$ or $U_{6n}$.
\end{prob}

This problem in turn contains the next problem as a subproblem which extends a question asked by Huang, Huang and Lu in \cite{HuangHL17}. Note that when restricted to cubic Cayley graphs this question has been answered by Huang, Huang and Lu in the same paper as seen in Theorem \ref{thm:CayDSD2p1}.

\begin{prob}
Determine which Cayley graphs on the dihedral group $D_{2n}$ are Cay-DS.
\end{prob}

\begin{prob}
Prove or disprove Conjecture \ref{SoConj} (\cite[Conjecture 7.3]{So06}): All integral circulant graphs are Cay-DS (or, quivalently, two integral circulant graphs are isomorphic if and only if they are cospectral).
\end{prob}

A necessary and sufficient condition for $D_{2p}$ to be Cay-DS, $p$ being a prime, was given in Theorem \ref{thm:CayDSD2p}. A similar result for a general dihedral group is unknown as far as we know.

\begin{prob}
\label{prob:CayDS1}
Determine all positive integers $n$ for which $D_{2n}$ is Cay-DS.
\end{prob}

In general, we have the following two very broad problems, each of which may be challenging even when restricted to special cases.

\begin{prob}
\label{prob:CayDS}
Determine which Cayley graphs on a given group are Cay-DS.
\end{prob}

\begin{prob}
\label{prob:CayDS2}
Determine which groups are Cay-DS.
\end{prob}

Recall that a group $\Ga$ is said to be determined by spectrum (DS) if any graph cospectral with a Cayley graph on $\Ga$ must be isomorphic to that Cayley graph. Theorem \ref{thm:CayDSSol} says that every DS group must be solvable. It is unknown to what extent or under what conditions the converse of this statement is true. So we have the following general problem.

\begin{prob}
\label{prob:CayDS3}
Under what conditions is a solvable group determined by spectrum?
\end{prob}

In view of Theorems \ref{Spec1mod4} and \ref{Spec3mod4}, if $R$ is a finite commutative ring as in Assumption \ref{as:1} such that $|R_i|/m_i \equiv 1~(\mod~4)$ for all but at most one $i$, and $|R_i|/m_i \equiv 3~(\mod~4)$ for this exceptional $i$, then the eigenvalues of the quadratic unitary Cayley graph $\mathcal{G}_R$ of $R$ are known. In all other cases the eigenvalues of $\mathcal{G}_R$ are unknown. So we pose the following problem.

\begin{prob}
Determine the eigenvalues of the quadratic unitary Cayley graph of any finite commutative ring.
\end{prob}

As seen in the first three subsections of Section \ref{sec:EnCay}, there is an extensive body of work on the energies of circulant graphs, especially integral circulant graphs. Nevertheless, we still do not know the energies of all integral circulant graphs, let alone all circulant graphs. So we have the next three problems of which the last one contains the other two as subproblems.

\begin{prob}
\label{pb:int-cir}
Determine the energies of all integral circulant graphs.
\end{prob}

\begin{prob}
\label{pb:int-cir-1}
Determine the minimum and maximum energies of an integral circulant graph of order $n$.
\end{prob}

\begin{prob}
\label{pb:en-cir}
Determine the energies of all circulant graphs.
\end{prob}

The next three problems are about Conjectures 6.1, 6.2 and 6.3 in \cite{Le12-1}, respectively.

\begin{prob}
Prove or disprove Conjecture \ref{conj:61}: For any integer $n \ge 3$, the minimum energy of an integral circulant graph of oder $n$ is equal to $2n(p-1)/p$, where $p$ is the smallest prime factor of $n$.
\end{prob}

\begin{prob}
Prove or disprove Conjecture \ref{conj:62}: For any integer $n \ge 3$ and any integral circulant graph $\ICG(n,D)$ of oder $n$ achieving the minimum energy, $D$ must be a multiplicative divisor set.
\end{prob}

\begin{prob}
Prove or disprove Conjecture \ref{conj:63}: For any integer $n \ge 3$ and any two multiplicative sets $D_1, D_2\subseteq D(n)\setminus\{n\}$ such that $\ICG(n,D_1)$ and $\ICG(n,D_2)$ are cospectral, we must have $D_1 = D_2$.
\end{prob}

The next problem is about the gcd graphs $\ICG(R/(c), D)$ as defined in Section \ref{subsec:eng-ring-gcd}, where $R$ is a unique factorization domain, $c$ is a nonzero nonunit element of $R$ such that $R/(c)$ is finite, and $D$ is a set of proper divisors of $c$.

\begin{prob}
\label{pb:RcD}
Determine the energy of the gcd graph $\ICG(R/(c), D)$.
\end{prob}

Theorem \ref{XnCor1} gives a necessary and sufficient condition for the unitary Cayley graph of $\ZZZ_n$ to be Ramanujan. We do not have a similar result for general integral circulant graphs, though partial results exist in the literature as seen in Section \ref{sec:RamCayAbelian}. Moreover, almost nothing is known about when a non-integral circulant graph is Ramanujan. So we have the following two problems.

\begin{prob}
Characterize Ramanujan integral circulant graphs.
\end{prob}

\begin{prob}
Characterize Ramanujan non-integral circulant graphs.
\end{prob}

Recall that Cayley graphs on dihedral groups are called dihedrants. Theorem \ref{RamnujanDH1} identifies all Ramanujan dihedrants in two special cases, but we do not have a similar result for general dihedrants. Similarly, we do not have a complete characterization of Ramanujan quadratic unitary Cayley graphs of finite commutative rings.

\begin{prob}
Characterize Ramanujan dihedrants.
\end{prob}

\begin{prob}
Characterize Ramanujan quadratic unitary Cayley graphs of finite commutative rings.
\end{prob}

The next two problems are about two conjectures in \cite{Allen98}.

\begin{prob}
Prove or disprove Conjecture \ref{AllenConj1}.
\end{prob}

\begin{prob}
Prove or disprove Conjecture \ref{AllenConj2}.
\end{prob}

The following problem is concerned with the signed reversal graph $SR_n$ (see Section \ref{subsec:cox} for definition).

\begin{prob}
Prove or disprove Conjecture 6.1 in \cite{CioabaRT20}: For any $n \ge 2$, the second largest eigenvalue of $SR_n$ is ${\choose{n}{2}}$.
\end{prob}

Recall from Section \ref{subsec:cox} that $T(n,I) = \{\sigma\in S_n: |\mathrm{supp}(\sigma)|\in I\}$ for $\emptyset \ne I \subseteq \{2,3,\ldots,n-1, n\}$ and that a Cayley graph on $S_n$ is said to have the Aldous property if its largest eigenvalue strictly smaller than its degree is attained by the standard representation of $S_n$. The next two problems posed in \cite{LiXZ21} arise from Theorems \ref{thm:LiXZ21a} and \ref{thm:LiXZ21} naturally.

\begin{prob}
\label{prob:LiXZa}
Give a necessary and sufficient condition for $\Cay(S_n,T(n,I))$ with $\{n-1\} \subset I \subset \{2,3,\ldots,n-2,n-1\}$ to have the Aldous property for sufficiently large $n$.
\end{prob}

\begin{prob}
\label{prob:LiXZb}
Give a necessary and sufficient condition for $\Cay(S_n,T(n,I))$ with $\{n\} \subset I \subseteq \{2,3,\ldots,n-2,n\}$ to have the Aldous property for sufficiently large $n$.
\end{prob}

The next problem is concerned with the Weyl group $W_n = W(B_n)$ and a specific permutation representation $\mathbf{P}_n$ of $W_n$ (see Section \ref{subsec:cox} for related definitions). For $A \subset \{1, 2, \ldots, n\}$, let $s_A$ be the element of $W_n$ which in the defining representation is given
by the diagonal matrix $\diag(x_i)_{i=1}^n$, where $x_i = -1$ if $i \in A$ and $x_i = 1$ if $i \notin A$.

\begin{prob}
\label{pb:Weyl}
Prove or disprove Conjecture 5.2 in \cite{Cesi20}: Let $w \in \mathbb{C}W_n$ be given by
$$
w = \sum_{A} a_{A} s_{A} + \sum_{(i, j)} b_{ij} (i, j),\;\, a_{A} \ge 0, b_{ij} \ge 0,
$$
where the first sum is running over all subsets $A$ of $\{1, 2, \ldots, n\}$ with odd cardinality, and the second sum is running over all transpositions $(i, j)$ in $S_n$. Then
$$
\mu_{W_n}(w) = \mu_{W_n}(w, \mathbf{P}_n).
$$
\end{prob}

The first part of the following question is Question 1 in \cite{CaoCL20}.

\begin{prob}
Is there any circulant admitting perfect state transfer which is isomorphic to a dihedrant? Determine all pairs of such graphs if they exist.
\end{prob}

\begin{prob}
Give a characterization of the gcd graphs of abelian groups admitting perfect state transfer.
\end{prob}

A few results about the existence of perfect state transfer in $\ICG(R/(c), D)$ were given in Section \ref{subsec:pst-rings}, where $R$, $c$ and $D$ are as in the paragraph immediately before Problem \ref{pb:RcD}. However, the picture is incomplete for this family of Cayley graphs. So we have the following general problem.

\begin{prob}
Study the existence of perfect state transfer in the gcd graphs $\ICG(R/(c), D)$.
\end{prob}

As seen in Section \ref{subsec:PSTdih}, we know exactly when a Cayley graph on a dihedral group admits perfect state transfer. Since dihedral groups are nilpotent and metacyclic, it is natural to consider the following two problems. (A metacyclic group is an extension of a cyclic group by a cyclic group.) Note that both problems are broad.

\begin{prob}
Study the existence of perfect state transfer in Cayley graphs on nilpotent groups.
\end{prob}

\begin{prob}
Study the existence of perfect state transfer in Cayley graphs on metacyclic groups.
\end{prob}

As seen in Section \ref{sec:pst}, if we want to construct regular graphs admitting perfect state transfer, integral Cayley graphs are good candidates. However, not every Cayley graph admitting perfect state transfer is integral, though we do not have many counterexamples so far.

\begin{prob}
Construct more families of non-integral Cayley graphs admitting perfect state transfer.
\end{prob}

Note that such a Cayley graph must be on a non-abelian group.

\begin{prob}
Characterize circulant graphs admitting pretty good state transfer. (Equivalently, characterize non-integral circulant graphs admitting pretty good state transfer.)
\end{prob}

\begin{prob}
In the case when $\Ga$ is abelian, a formula for the eigenvalues of the bi-Cayley graph $\Cay(\Ga; R, S, T)$ in terms of that of $\Cay(\Ga, R)$, $\Cay(\Ga, S)$ and $\Cay(\Ga, T)$ was given in Theorem \ref{thm:bi-Cay}. Explore whether similar results can be achieved for some other families of groups.
\end{prob}

Recall that bi-Cayley graphs and tri-Cayley graphs on cyclic groups are called bicirculants and tricirculants, respectively. The next two problems are related to the results in Section \ref{subsec:srnCay}.

\begin{prob}
Classify all strongly regular bicirculants.
\end{prob}

\begin{prob}
Classify all strongly regular tricirculants.
\end{prob}

The Petersen graph is an integral bicirculant with spectrum $((-2)^4, 1^5, 3)$. This observation prompts us to raise the next problem and its counterpart for tricirculants.

\begin{prob}
\label{pb:int-bicir}
Classify all integral bicirculants.
\end{prob}

\begin{prob}
\label{pb:int-tricir}
Classify all integral tricirculants.
\end{prob}

A characterization of integral Cayley graphs on abelian groups was given by Alperin and Peterson as seen in Theorem \ref{thm:Alperin12}. The next two problems call for counterpart results for bi-Cayley graphs and tri-Cayley graphs, respectively.

\begin{prob}
Characterize all integral bi-Cayley graphs on abelian groups.
\end{prob}

\begin{prob}
Characterize all integral tri-Cayley graphs on abelian groups.
\end{prob}

The next two problems are related to Problem \ref{pb:int-bicir}.

\begin{prob}
\label{pb:pst-int-bicir}
Classify all integral bicirculants which admit perfect state transfer.
\end{prob}

\begin{prob}
\label{pb:pgst-int-bicir}
Classify all integral bicirculants which admit pretty good state transfer.
\end{prob}

The next two problems are related to Problem \ref{pb:int-tricir}.

\begin{prob}
\label{pb:pst-int-tricir}
Classify all integral tricirculants which admit perfect state transfer.
\end{prob}

\begin{prob}
\label{pb:pgst-int-tricir}
Classify all integral tricirculants which admit pretty good state transfer.
\end{prob}

As mentioned in Section \ref{subsec:pgst}, a regular integral graph admits pretty good state transfer if and only if it admits perfect state transfer. So Problem \ref{pb:pgst-int-bicir} is identical to Problem \ref{pb:pst-int-bicir} for regular bicirculants, but the two problems are different for bicirculants which are not regular. Similarly, Problem \ref{pb:pgst-int-tricir} is different from Problem \ref{pb:pst-int-tricir} only for tricirculants which are not regular.

Finally, as seen in Section \ref{subsec:dis-reg-cay}, all distance-regular graphs in the following families of Cayley graphs have been classified: Cayley graphs on cyclic groups; Cayley graphs on abelian groups (with respect to minimal inverse-closed generating sets); Cayley graphs on dihedral groups. In this line of research one may attempt to classify distance-regular Cayley graphs (or, more restrictively, strongly regular Cayley graphs) for other families of groups such as metacyclic groups, generalized dihedral groups, etc.



\subsection*{Acknowledgements}

We are extremely grateful to the anonymous referee for many helpful comments which led to a number of improvements to this paper. In particular, Section \ref{sec:open} would not exist without the referee's suggestion. We would like to thank Cristina Dalf\'{o}, Mojtaba Jazaeri, Kazufumi Kimoto, Yuxuan Li, Daria Lytkin, Bogdan Nica, Paul Renteln, Cid Reyes-Bustos, Shohei Satake, Sergey V. Savchenko and Binzhou Xia for suggesting papers, clarifying points in their papers and/or informing us omissions and errors in earlier versions of this article.


\end{document}